\documentclass[12pt,a4paper]{amsart} 
\usepackage{amssymb, amsfonts, amscd}
\usepackage[all]{xy} 
\usepackage{amsmath}
\usepackage{amsthm}
\usepackage{epsfig}
\usepackage{url}

\usepackage{hyperref} 
\graphicspath{{figures-pdf/}}

\newcommand{\id}{{\mathrm {id}}}
\newcommand{\R}{{\mathbb R}}   
\newcommand{\Z}{{\mathbb Z}}   

\newcommand\Cat{\mbox{Cat}}

\newcommand{\inv}{{-1}}
\newcommand{\uline}[1]{{\underline{#1}}} 
\newcommand{\mac}[1]{{\mathcal {#1}}}

\newcommand{\Glax}{weak } 

\theoremstyle{definition} 
\newtheorem{theorem}{Theorem}[section]  
\newtheorem{lemma}[theorem]{Lemma}  
\newtheorem{conjecture}[theorem]{Conjecture} 
\newtheorem{definition}[theorem]{Definition}
\newtheorem{corollary}[theorem]{Corollary}
\newtheorem{example}[theorem]{Example}
\newtheorem{remark}[theorem]{Remark}

\newtheorem*{theorem*}{Theorem}
\newtheorem*{theoremc*}{Theorem/Conjecture}

\theoremstyle{remark}

\begin{document}

\title {Gray categories with duals and their diagrams}

\author{John W. Barrett}
\address{School of Mathematical Sciences, University
  Park, Nottingham NG7 2RD, UK} 
\email{john.barrett@nottingham.ac.uk}

\author {Catherine Meusburger}
\address {Department Mathematik,
  Friedrich-Alexander-Universit\"at Erlangen-N\"urnberg,
  Cauerstrasse 11,
  D-91058 Erlangen,
  Germany}
\email{catherine.meusburger@math.uni-erlangen.de}

\author {Gregor Schaumann}
\address {Universit\"at W\"urzburg, Institut f\"ur Mathematik, 
Lehrstuhl X,
Emil-Fischer-Strasse 40, 
D-97074 W\"urzburg, 
Germany}
\email{gregor.schaumann@uni-wuerzburg.de}

\date{16 May 2024} 
\thanks{ Gray categories with duals and their diagrams © 2024 by John W. Barrett, Catherine Meusburger and Gregor Schaumann is licensed under CC BY-NC-ND 4.0. To view a copy of this license, visit https://creativecommons.org/licenses/by-nc-nd/4.0/  \  Also published at https://doi.org/10.1016/j.aim.2024.109740\ }

\begin{abstract}
The geometric and algebraic properties of Gray categories with duals are investigated by means of a diagrammatic calculus. The diagrams are three-dimensional stratifications of a cube, with regions, surfaces, lines and vertices labelled by Gray category data. These can be viewed as a generalisation of ribbon diagrams. The Gray categories present two types of duals, which are extended to  functors of 2-strict tricategories with natural isomorphisms, and correspond directly to symmetries of the diagrams.  It is shown that these functors can be strictified so that the symmetries of a cube are realised exactly. A new condition on Gray categories with duals called the spatial condition is defined. A class of diagrams for which the evaluation for spatial Gray categories is invariant under homeomorphisms is exhibited. This relation between the geometry of the diagrams and structures in the Gray categories proves useful in computations and has  applications in topological quantum field theory. 
\end{abstract}

\maketitle

\section{Introduction}

The aim of this paper is to develop the theory of duals for Gray categories. 
The principal tool is a diagrammatic calculus introduced in this paper.   This can be viewed as  a higher-categorical, three-dimensional  analogue of the diagrams used for computations in  pivotal categories. Many of the algebraic results on Gray categories with duals can be understood in terms of the geometry of the corresponding diagrams.

Our main motivation is the applications of higher categories in (extended) topological quantum field theory, quantum geometry and conformal field theory  \cite{KapTop, Lurie}.  For instance, one would like  to construct  topological quantum field theories with `defects'.  These are theories in which certain embedded submanifolds are labelled with geometric data.  In three dimensions, it is natural that the data on these labelled submanifolds should arise from a tricategory \cite{CMS}.
An  example of this is the work on the relation between Reshetikhin-Turaev and Turaev-Viro invariants \cite{balsam_kirillov,TuraevNew}. In this case, the relevant higher category is the centre of a spherical category, which is a tricategory with a single object and a single 1-morphism.  
Another example arises from  bimodule categories, which plays an important role in conformal field theory, see e.g. \cite{KapSau, schweigert}.
The notion of duals is required to incorporate orientation. If  data from a  subcategory is used to label distinguished submanifolds, orientation reversal of these submanifolds  must be reflected in a corresponding structure in the tricategory, namely the duals.

Another important motivation is the benefit of diagrams for computations in higher categories that arises from a direct relation between  geometry and structures in the category. That the diagrams have a non-trivial geometric content is familiar from the example of
ribbon categories and knots, or more generally, ribbon graphs embedded in three-dimensional space. The Reshetikhin-Turaev  invariants  \cite{RT} define a functor which takes ribbon graphs in three-dimensional space labelled with data from a ribbon category and evaluates them in the category itself. 
The relations in the category imply invariance of the evaluation under homeomorphisms of three-dimensional space. In this way the homeomorphism invariance is the geometric expression of the relations in the category. 

This article does not consider general tricategories, but restricts attention to Gray categories. As every tricategory is triequivalent to a Gray category   \cite{GPS, gurski} and there is no stricter version of a tricategory with this property, Gray categories can be viewed as maximally strict tricategories.
The practical reason for using Gray categories  is  the wish to avoid a degree of complexity that makes algebraic manipulations nearly impossible. 
The deeper reason is that the coherence data for Gray categories is precisely that part of the coherence data for tricategories that can be given a diagrammatic meaning. 
An analogous situation arises for pivotal tensor categories. There is a weak and a strict notion of monoidal structure and of duality and the diagrams for pivotal tensor categories 
reflect precisely  the coherence data for a strict pivotal category with strict monoidal structure, while  the rest of its coherence data is not given a diagrammatic representation.

Section 2 of the paper introduces diagrams for Gray categories without duals. This is a generalisation of the diagrammatic calculus for braided monoidal categories introduced by Joyal and Street \cite{JS}. The diagrams for Gray categories are located in a cube with the three coordinate axes corresponding to the three compositions in a Gray category. These diagrams consist of 3-, 2-, 1-, and 0-dimensional strata which are labelled, respectively, with objects and 1-, 2-, 3-morphisms in a Gray category. The categorical axioms are introduced in an `unpacked' manner. We see this concreteness as an advantage, both in view of  possible applications in state sum models and  because this yields a direct link between categories and geometry.

 The familiar diagrams for  braided monoidal categories have  an evaluation to morphisms in the category and this evaluation is invariant under isotopies of diagrams.
  A braided monoidal category can be viewed as a Gray category with a single object and 1-morphism, see e.g. Lemma \ref{braidten}, and we generalize the evaluation of diagrams for braided monoidal categories to the evaluation of Gray category diagrams.  A  Gray category diagram is thereby evaluated to 
 a 3-morphism between  2-morphisms associated with the top and bottom of the cube. This evaluation of Gray category diagrams  is invariant under a similar set of moves as the evaluations of diagrams for braided monoidal categories and 
as a consequence, it is invariant under the following isotopies.

\begin{theorem*}(Theorem \ref{th_progressive_invariant})
  Let $D$, $D'$ be generic Gray category diagrams that are isotopic by a one-parameter family of isomorphisms of progressive diagrams. Then the evaluations of $D$ and $D'$ are equal.
\end{theorem*}

Section 3 introduces Gray categories with duals using the definition of Baez and Langford \cite{BL} but with some minor modifications. The diagrammatic representation of the data for these duals is explained in this section.
The Gray categories possess two types of duals, ${*}$ and $\#$,  which  correspond, in a sense explained in this paper, to $180$ degree rotations 
around two different coordinate axes. The $*$-duals are familiar from pivotal or ribbon categories, but the $\#$-duals are a feature of Gray categories that does not appear in the pivotal or ribbon cases. The coherence data for $\#$ is such that it matches the appearance of folds and cusps in projections of surfaces.   As in the case of  Gray categories, the  duals considered in this paper are not the most general ones, and their axioms could be weakened.  Again, the strictness of the axioms ensures that all their  coherence data has a diagrammatic representation.

Our first central result  in Section 4  concerns the algebraic structure of these duality operations. 
\begin{theorem*}[Theorem \ref{graycatduals}, Lemma \ref{lemma:doublestar-is-one}, Theorem \ref{graydualnat}]
  The duals extend in a canonical way to (partially contravariant)  functors of 2-strict tricategories $*,\#\colon\mac G\to\mac G$ with $**=1$ and define natural isomorphisms $\Gamma\colon *\#{*}\#\to 1$, $\Theta\colon\#\# \to 1$.  
\end{theorem*}

The structure maps for these natural isomorphisms are interpreted geometrically in terms of diagrams. By investigating a  closely-related natural isomorphism $\Delta\colon \#\to*\# *$ one obtains two  diagrams that are homeomorphic,  but whose evaluations 
are not necessarily equal.   
This motivates the 
definition of a spatial Gray category as a Gray category with duals in which
such an equality holds. This condition is a generalisation of the ribbon condition for a ribbon category.

Section 5 contains the second  important  result: a strictification  theorem for the  duals given in Theorem \ref{thm:strictification} and \ref{lem:graystrict}. While the  functor of 2-strict tricategories  $*\colon\mac G\to\mac G$ satisfies the identity $**=1$, 
which corresponds to its geometrical interpretation, the functor $\#\colon\mac G\to\mac G$  satisfies such an identity only up to a natural isomorphism. Similarly,  the functor $*\#{*}\#$ whose geometrical counterpart is the identity rotation of $\mathbb R^3$, is not equal to the identity functor but only isomorphic to it.  However, a spatial Gray category can be strictified (in the sense of \cite{gurski}) to one in which these duality functors do indeed satisfy the relations for 180 degree rotations around different coordinate axes: 

\begin{theorem*} [Theorems \ref{thm:strictification} and \ref{lem:graystrict}] Every spatial Gray category  with duals can be strictified to a Gray category  whose duals $*,\#\colon\mac G\to\mac G$  satisfy $**=1$, $\#\#=1$ and $*\#{*}\#=1$.
\end{theorem*}
Thus the geometrical interpretation of the action of $\#$ as a rotation is restored for higher morphisms, 
which justifies  the original set of duality axioms. The proof for this result is conceptually clear and  may be of independent interest. 

Section 6 explores in more depth the relation between Gray categories with duals and their diagrams.  The diagrams in this paper have no framing. This is adequate to express all of the axioms for the category and also the structure maps for the duality functors. However it does restrict the generality of the invariance results. Diagrams are labelled with category data using a generalisation of `blackboard framing' familiar from knot theory. These diagrams are called standard. 

The general invariance result in this section holds for a large class of Gray category diagrams, called surface diagrams,  whose 0-,  1- and 2-strata form a surface with boundary. 
It states that the evaluations of standard surface diagrams are invariant under a set of moves that are the PL counterparts of the moves induced by projecting an isotopy  in the smooth setting (Theorem \ref{invariance_surface}). Under the conjecture (Conjecture \ref{conjecture}) that these are also all the moves arising from projecting PL isotopies, it implies that  oriented isomorphisms of standard surface diagrams leave their evaluations invariant. While it appears to be a reasonable conjecture that the moves on folds and cusps in the PL setting are the same as in the smooth setting, we do not know of any previous work on this problem.

\begin{theoremc*} [Theorem \ref{invariance_surface}, Conjecture \ref{conjecture} ]
  Let $D$ and $D'$ be standard surface diagrams that are labelled with a spatial Gray category. Let $f\colon D\to D'$ be an oriented isomorphism of standard surface diagrams and the labels of $D'$ induced from $D$ by $f$. Then the evaluations of $D$ and $D'$ are equal.   
\end{theoremc*}

For surface diagrams, this indicates that there are no further conditions other than the spatial condition on a Gray category with duals that are needed to prove invariance under homeomorphisms. Essentially it arises because surface diagrams have a uniquely determined notion of framing. Extending this to all Gray category diagrams would require a general definition of framing; whilst this is an interesting problem we leave it for future work.

\tableofcontents

\section{Category, 2-category and Gray category diagrams}

The aim of this section is to develop a diagrammatic calculus for Gray categories  and to show that the evaluation of diagrams labelled with  Gray category data is invariant under certain mappings of diagrams.  These diagrams can be viewed as  higher-dimensional analogues of  spin network diagrams in the physics literature, string diagrams in the mathematics literature and tangle diagrams in knot theory. They are dual to the more common pasting diagrams of the category theory literature. 

The diagrams for an $n$-category are located in a geometrical `cube' $[0,1]^n$.  It should be possible to define diagrams analogous to the ones considered here for  arbitrary dimension  $n$.  However,  it is most practical to be guided by known examples rather than abstract formalism. Hence this work only considers the cases  up to $n=3$. 
The evaluation of an $n$-dimensional diagram labelled with $n$-category data will be defined inductively, in terms of a projection to  an $(n-1)$-dimensional diagram labelled with $(n-1)$-category data. 

For this reason, the definition of Gray category diagrams and their mappings, which corresponds to $n=3$,  requires a careful discussion of their lower-dimensional counterparts.  
We start with a discussion of one-dimensional diagrams, then the two-dimensional case before introducing diagrams for Gray categories.   At each stage, we discuss the $n$-dimensional diagrams, their mappings as well as their labelling with data from an $n$-category and their evaluation. 

\subsubsection{Piecewise-linear topology}\label{subsubsec:plbackground}
Throughout the article,  the piecewise-linear (PL) framework is used. The basic definitions of PL topology are taken from Rourke and Sanderson \cite{RS}, so that all spaces in this article are polyhedra and the mappings are piecewise-linear. A {\bf polyhedron} is a subset of $\R^n$, for some $n$, that is locally a cone over a compact subset of $\R^{n}$. As a consequence, each polyhedron is a locally finite union of simplexes in $\R^n$. If compact, it is the union of a finite set of simplexes in $\R^n$ that form a simplicial complex and conversely, every such simplicial complex determines a compact polyhedron. 

 A {\bf piecewise-linear map} $f\colon P\to Q$ is a map that is locally conical.
It follows that there is a locally-finite decomposition of $P$ into simplexes such that 
$f$ is (affine) linear on each simplex. If $P$ and $Q$ are the polyhedra of a simplicial complex, then a piecewise-linear map $P\to Q$ is linear on each simplex of some subdivision of $P$, and is in fact a simplicial map to some subdivision of $Q$.

Note that in the category of topological spaces, embeddings of polyhedra are classified into tame or wild  according to whether or not they are equivalent to PL embeddings \cite{Rushing}. Thus PL embeddings are by definition tame. 

An {\bf isotopy} of  $X$ is a PL isomorphism $\phi\colon X\times[0,1]\to X\times [0,1]$ such that $\phi(x,t)=(\phi_t(x),t)$ and $\phi(x,0)=(x,0)$.
This can be thought of as a PL isomorphism $\phi_1$ that is continuously connected to the identity map.

A polyhedron is a {\bf $k$-manifold (with boundary)} if it is locally PL isomorphic to $\R^k$ or $\R^k_+$; this is a property of the polyhedron and does not require any additional structure. 

\subsubsection{Comparison with the smooth case}
It would also be possible to work with smooth manifolds and smooth maps. In a general dimension this is significantly different, as the classic work of Kervaire and Milnor shows that in sufficiently high dimension there are PL manifolds that can have many smooth structures or none; however this does not occur in dimension three. Whether this implies that smooth diagrams for Gray categories are essentially the same as the PL ones is not entirely clear to us. There are some technical details that are particular to the PL case: the definition of a stratification in Section \ref{sec:graycatd} is simpler in the PL case as the `Whitney frontier conditions' that are standard in the smooth case are not necessary. Also, the `Alexander trick'  is used in theorem \ref{invariance_surface}; the Alexander trick is true
 for all dimensions in the PL case but fails for sufficiently high dimension in the smooth case. Note that the use of the Alexander trick could be avoided by simply replacing `oriented isomorphism' by `isotopy' in the statement of the theorem.

There are significant advantages in working in the PL setting. Since any compact polyhedron is the union of a finite number of simplexes one sees immediately that the constructions are combinatorial in nature and could be easily automated, e.g. by using simplicial complexes with rational coordinates for the vertices. The composition of diagrams is much simpler, as it is clear that 
the union of two 
 PL manifolds along a boundary component is  again a PL manifold; however this is not directly true for smooth manifolds due to the need to address the delicate technical issue of `smoothing corners'. 

One advantage of smooth diagrams is that they are often easier to draw and read. Thus in some places the diagrams are drawn smoothly for better legibility, but these should be understood as a representation of a finely triangulated PL diagram to a limited resolution on the page.

\subsection{Categories and diagrams}

We start by considering the one-dimensional case, which is required for the definitions in higher dimensions,  but also  can be viewed as a toy model  that exhibits the general features of the construction.

\begin{definition}
  A {\bf one-dimensional diagram} is a finite set of points, called {\bf vertices}, in the interior of the unit interval $[0,1]$. The complement of the vertices is a disjoint union of its connected components,  which are called {\bf regions} of the diagram. 
\end{definition}

The one-dimensional diagrams are a purely topological construction. They become \emph{category diagrams} once labelled with data from a category $\mathcal C$. A category diagram with a single vertex is called an elementary diagram. A general category diagram can then be defined in terms of the elementary diagrams.
\pagebreak[4] 

\begin{definition}
Let $\mac C$ be a category.
  \begin{enumerate}
  \item 
    An {\bf elementary category diagram} for  $\mathcal C$ is a one-dimen\-sio\-nal diagram with one vertex together with a morphism $f\colon A\to B$ in $\mac C$. The object $A$ is associated with the region containing 0, $B$ with the region containing 1, and $f$ with the vertex.
  \item 
    A {\bf category diagram} for  $\mac C$ is a one-dimensional diagram together with a labelling of each region with an object in $\mac C$ and a labelling of each vertex  with a morphism in $\mac C$. 
For every vertex it is required that there exists a neighbourhood $[a,b]$ such that the unique affine map $[a,b] \to [0,1]$, sending $a$ to $0$ and $b$ to $1$, maps $[a,b]$ to an elementary category diagram. 
  \end{enumerate}
\end{definition} 

The benefit of category diagrams and their higher-dimensional counterparts is that they 
allow one to visualise a  calculation in the category  $\mathcal C$. This calculation is called the evaluation of the diagram.

\begin{figure}
  \centering
  \includegraphics[scale=0.3]{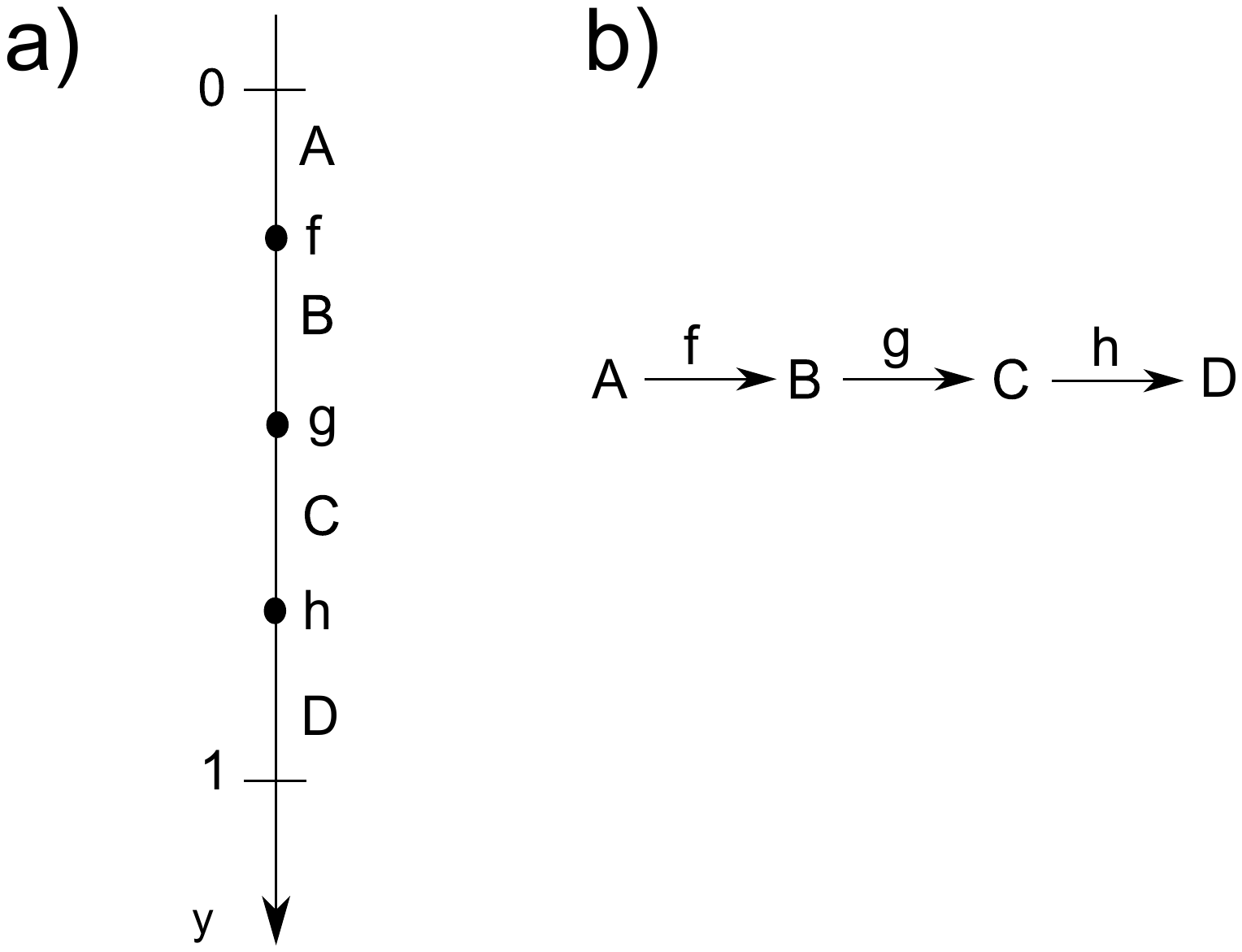}
  \caption{Category diagram (a) and its evaluation (b).}
  \label{1catdiagram}
\end{figure}

\begin{definition}
\label{1d_eval}
  The {\bf evaluation of a category diagram} is the product of the morphisms at the vertices in the order of increasing values of $v$. It maps the object for the region containing $0$ to the object for the region containing $1$ (see Figure \ref{1catdiagram}). A diagram without vertices  is called an {\bf identity diagram} and its evaluation is the identity on the single object. 
\end{definition}
The evaluation is well-defined by associativity. 
A category diagram and its evaluation are shown in Figure \ref{1catdiagram}.
The usefulness of diagrams for visualising calculations in  (higher) categories  is due to the fact  that their evaluation  is invariant under certain manipulations of diagrams such as homeomorphisms of diagrams and subdivisions.  A precise formulation of this idea requires the notion of a mapping of diagrams.
\pagebreak[3] 

\begin{definition}
$\quad$\label{1dmap}
  \begin{enumerate}
  \item A {\bf mapping of one-dimensional diagrams} $D\to D'$ is a PL embedding $m \colon [0,1]\to [0,1]$ such that  $m(v)$ is a vertex of $D'$ for each vertex $v$ of $D$.

  \item If the mapping has the property that $v$ is a vertex if and only if $m(v)$ is, then it is called a {\bf subdiagram}.

  \item  A mapping $m\colon D\to D'$ of one-dimensional diagrams is called a  {\bf homomorphism of diagrams} if  $m(0)=0$ and $m(1)=1$ and an {\bf isomorphism} if  $m^{-1}$ is also a homomorphism. 

  \item  If the mapping  $m\colon[0,1]\to[0,1]$ is the identity map, then the mapping is called a {\bf subdivision} of $D$. 
  \end{enumerate}
\end{definition}

\begin{figure}
  \centering
  \includegraphics[scale=0.4]{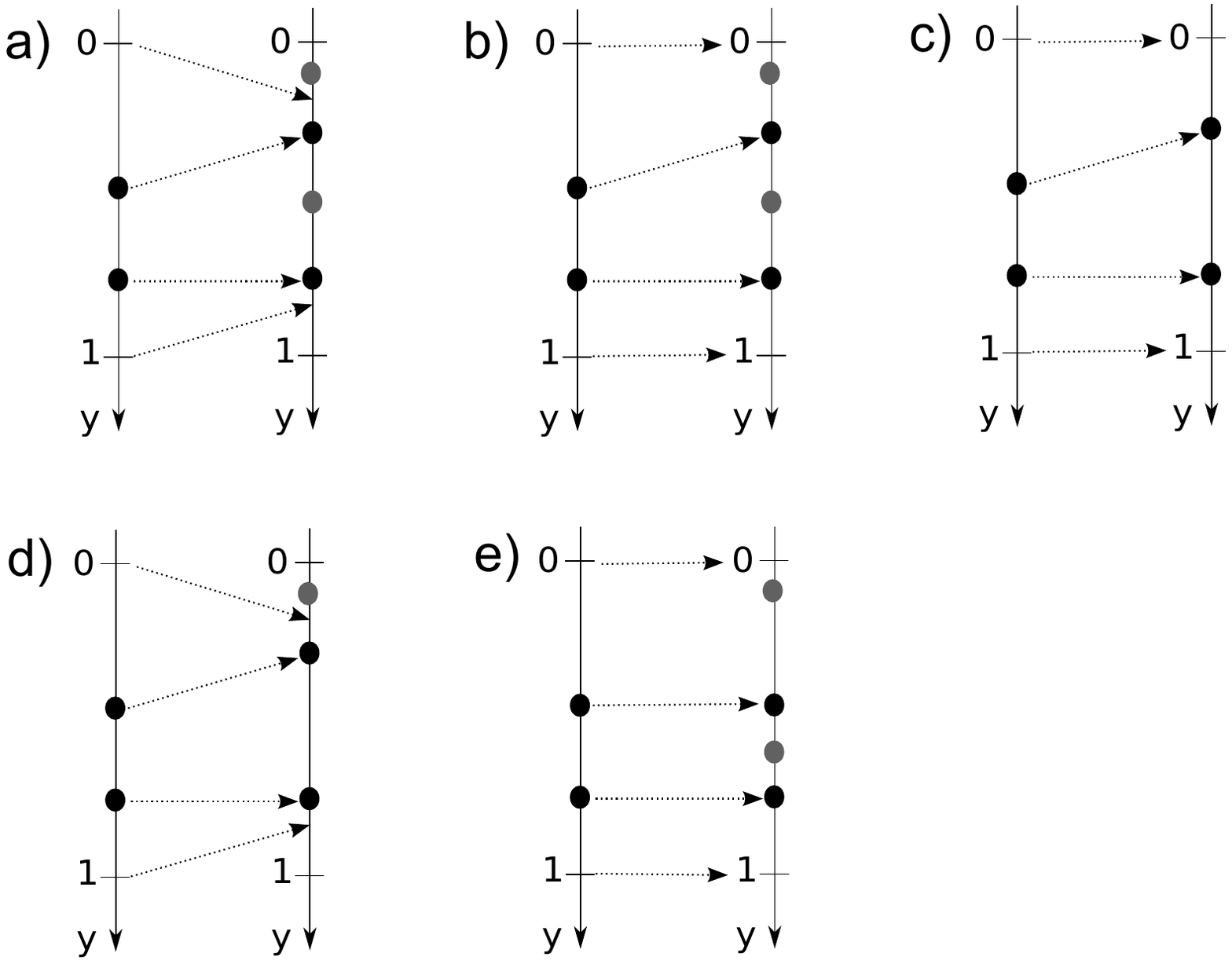}
  \caption{a), b), c), d), e) mappings of  diagrams,\newline
    b), c), e) homomorphisms of diagrams, \newline
    c) isomorphism of diagrams, d) subdiagram, \newline
    e) subdivision.
  }
  \label{mappings}
\end{figure}

Mappings, homomorphisms, isomorphisms, subdiagrams and subdivisions of one-dimensional diagrams are shown in Figure \ref{mappings}.
Note that if $m\colon D\to D'$ is a mapping,  a homomorphism, a subdiagram or a subdivision, the image $D'$ can have more vertices than $D$. 
In the case of an isomorphism, there is a  bijection between the vertices of $D$ and of $D'$.  As the map $m\colon[0,1]\to[0,1]$ is an embedding, every  homomorphism of diagrams is a homeomorphism and can be decomposed into an isomorphism of diagrams and a subdivision.  

The concept of a mapping can be extended to mappings of diagrams that are labelled with data from a category. This amounts to imposing certain relations between the label of a  point in the diagram and the label of its image.

\begin{definition}
Let $\mac C$ be a category.
  A {\bf mapping of category diagrams}  $D,D'$  for $\mac C$ is a mapping of one-dimensional diagrams $m\colon D\to D'$ that is orientation-preserving ($m(0)<m(1)$) and preserves the labelling:
  \begin{enumerate} 
  \item If $x\in D$ and $m(x)\in D'$ are both vertices or both points in a region, then their labels are equal.
  \item  If $x\in D$ is a point in a region  of $D$ labelled by $A$ and $m(x)$ is a vertex of $D'$, then $m(x)$ is labelled with the morphism $1_A\colon A\to A$. 
  \end{enumerate}
  A mapping of category diagrams is called a {\bf homomorphism}, {\bf isomorphism}, {\bf subdivision} or {\bf subdiagram} if the underlying mapping  $m\colon D\to D'$ of one-dimensional diagrams  is.
\end{definition}

It is clear that if $D$ and $D'$ are category diagrams and $m\colon D\to D'$ a homomorphism of category diagrams, then the evaluations of $D$ and $D'$ are equal.

\subsection{2-category diagrams}\label{2catdiagsection}

This section extends the notions of  labelled diagrams, their evaluation and their mappings  to two dimensions. In this case, the  diagrams are labelled with data from a 2-category.
Compactly, a  2-category can be defined as a category enriched in the category $\Cat$ of categories and functors \cite{Kelly-B, KS}. Unpacked, this amounts to the following definition.

\begin{definition}[\cite{StreetH}]
\label{2catdef}
  A {\bf 2-category} $\mac C$  consists of a set of objects, and for each pair of objects $A,B$ 
  a category $\mathcal C(A,B)$, whose objects are called {\bf 1-morphisms} and whose morphisms are called {\bf 2-morphisms}.
  The composition of morphisms in $\mac C(A,B)$ is denoted  $\cdot$ and called {\bf vertical composition}.
  It is required that for each triple of objects $A,B,C$ there is an assigned 
  a {\bf horizontal composition functor} $\circ\colon\mathcal C(B,C)\times \mathcal C(A,B)\to \mathcal C(A,C)$, which is strictly associative and unital.
  The unit 1-morphism for an object $A$ is denoted  $1_A$ and the unit $2$-morphism for a 1-morphism $f$ is denoted $1_f$.
\end{definition}
 It follows that the objects and 1-morphisms of a 2-category $\mac C$ form a category $\mathcal C_1$. 
It is possible to define a 2-category without using the horizontal composition of two 2-morphisms; all that is required is the horizontal composition of a 2-morphism with a 1-morphism.
The notion of a 2-category can then be generalised  by regarding the 
horizontal composition of two 2-morphisms as undefined and dropping the interchange law.
This will be called a pre-2-category (it is also called a sesqui-category in \cite{StreetH,Crans}), and will be useful for the discussion of Gray categories below. The definition of a pre-2-category is given next, followed by the interchange law. Taken together, these give an explicit definition of a 2-category.

\begin{definition}

  A {\bf pre-2-category} consists  of  a set of objects, and for each pair of objects $A$, $B$ a category $\mathcal C(A,B)$.
 Furthermore, there  is the datum of  a horizontal composition of 1-morphisms $\circ\colon\mathcal C_1(B,C)\times \mathcal C_1(A,B)\to \mathcal C_1(A,C)$ with units that makes the objects and 1-morphisms into a category $\mathcal C_1$ and extends to a left action of $\mathcal C_1(B,C)$ on $\mathcal C(A,B)$ and a right action of $\mathcal C_1(A,B)$ on $\mathcal C(B,C)$
  by functors. This means there are functors
  \begin{equation*}
f \circ - \colon      \mathcal C(A,B)\to \mathcal C(A,C) \quad \text{and} \quad - \circ g\colon \mathcal C(B,C)\to \mathcal C(A,C)
  \end{equation*}
for 1-morphisms $f \in C_{1}(B,C)$ and $g \in C_{1}(A,B)$.
  These actions are required to be unital and associative and to commute with each other.  
\end{definition}

The notation is the same as before: if $f\colon A\to B$ is a 1-morphism, $\Psi\in \mathcal C(B,C)$ a 2-morphism and $n\colon C\to D$ a 1-morphism, then $\Psi\circ f\in \mathcal C(A,C)$ and $n \circ \Psi \in \mathcal C(B,D)$  denote the horizontal composites.
In a pre-2-category there are two  possible  definitions for the horizontal composite of two 2-morphisms. For  1-morphisms $f,g\colon A\rightarrow B$, $h,k\colon B\rightarrow C$ and 2-morphisms $\Phi\in\mathcal C(f,g)$, $\Psi\in\mathcal C(h,k)$ these are 
$$(k\circ\Phi)\cdot (\Psi\circ f)\quad\quad\text{and}\quad\quad (\Psi\circ g)\cdot (h\circ \Phi).$$
\begin{equation*}
   \xymatrix{
  A
  \ar@/^2pc/[rr]_{\quad}^{f}="1"
  \ar@/_2pc/[rr]_{g}="2"
  &&B
  \ar@{}"1";"2"|(.2){\,}="7"
  \ar@{}"1";"2"|(.8){\,}="8"
  \ar@{=>}"7" ;"8"^{\Phi} 
  \ar@/^2pc/[rr]_{\quad}^{h}="9"
  \ar@/_2pc/[rr]_{k}="10"
  &&C
  \ar@{}"9";"10"|(.2){\,}="11"
  \ar@{}"9";"10"|(.8){\,}="12"
  \ar@{=>}"11" ;"12"^{\Psi} 
}
\end{equation*}
The interchange law in a 2-category states that these  2-morphisms are equal and define $\Psi\circ\Phi$. Thus a 2-category can be viewed as a pre-2-category with an interchange law.

The notion of a pre-2-category will be useful later: Definition \ref{graydef} of a Gray category can be regarded as a categorification of a 2-category, where the interchange law is replaced by a family of isomorphisms. 
\bigbreak

Two-dimensional diagrams are a direct generalisation of one-dimen\-sional diagrams. The unit interval $[0,1]$ is replaced by the unit square $[0,1]^2$, and the  only  additional condition  is  that  lines  meet the boundary of the square only at its top and bottom edge. 

\begin{definition}
$\quad$ \label{genprog2d}
  \begin{enumerate}
  \item
    A {\bf two-dimensional diagram}  is a set of closed  subspaces $\emptyset=X^{-1}\subset X^0\subset X^1\subset X^2=[0,1]^2$ called $k$-skeleta,  such that each  $X^k\setminus X^{k-1}$ is a PL manifold of dimension $k$, for $k=0,1,2$, $X^0$ lies in the interior of the square, and  all intersection points of $X^1$ with the boundary of the square are contained in $(0,1)\times\{0,1\}$.
    The connected components of  $X^2\setminus X^{1}$ 
    are called {\bf regions}, the connected components of $X^1\setminus X^0$ {\bf lines}, and the elements of $X^0$ {\bf vertices}.  

  \item A two-dimensional diagram is called {\bf progressive} if the projection of the diagram $p_1\colon (x,y)\mapsto y$ is regular, i.e., the mapping of each line is a PL isomorphism to its image in $\R$. 

  \item A progressive two-dimensional diagram is called {\bf generic} if the $y$-coordi\-na\-tes of any two different vertices are different.
  \end{enumerate}
\end{definition}

Note that because all spaces $X^k$ are compact polyhedra, the set  of vertices of a two-dimensional diagram is always finite, see the preliminaries in Subsection \ref{subsubsec:plbackground}.

Mappings between two-dimensional diagrams are defined in direct analogy to the one-dimensional case; they are  PL embeddings that  map $k$-skeleta to $k$-skeleta.

\begin{definition} \ 
\label{2dmappings} 
  \begin{enumerate}
  \item A {\bf mapping of two-dimensional diagrams} $D\to D'$ is  a PL embedding $m\colon [0,1]^2\to [0,1]^2$
    that preserves the $k$-skeleta, i.e., $m(X^k)\subset X'^k$ for $k=0,1$. 
  \item If a mapping has the property that $v$ is a vertex if and only if $m(v)$ is, then it is called a {\bf subdiagram}.
  \item If $m$ is a PL homeomorphism that is the identity map on $\partial [0,1]^2$, the mapping is called a {\bf homomorphism of diagrams} and an {\bf isomorphism} if $m^{-1}$ is also a homomorphism. 
 \item  If  $m$ is the identity map, then the mapping is called a {\bf subdivision} of $D$. 
  \end{enumerate}
\end{definition}

Let  $m\colon [0,1]^2\to [0,1]^2$ be an arbitrary  PL homeomorphism  that is the identity on the boundary. 
Applying $m$ to any diagram $D$ defines a 
diagram $D'$ by transporting the skeleta of $D$ along $m$. This gives $m$ the structure of  a homomorphism of diagrams. Further examples can be constructed by subdividing $D'$, i.e., adding additional lines or vertices. 

The next step is to label two-dimensional diagrams with data from general 2-categories. To 
this end we consider the class of \emph{progressive diagrams}.  
They were first studied in the context of monoidal categories by Joyal and Street \cite{JS} and
have an important local structure:  any vertex  $v=(x,y)$ in a progressive diagram has a rectangular neighbourhood $[x-\epsilon_1,x+\epsilon_1]\times[y-\epsilon_2,y+\epsilon_2]$  which is, after a rescaling of the coordinates, a subdiagram that is an `elementary' two-dimensional diagram with one vertex.  For this, one chooses
a sufficiently small 
neighbourhood of $v$ in which all line segments are linear and  $\epsilon_2$ sufficiently small  so that all line segments exit the rectangle through its top or bottom edge.  

In analogy to the one-dimensional case,  two-dimensional diagrams can be labelled with data from a 2-category. For a generic diagram, the definition of the labelling does not use the interchange law of a 2-category and hence can be stated in the more general context of a pre-2-category. Hence generic diagrams are the ones that are appropriate for pre-2-categories. This will turn out to be important subsequently in the definition of  Gray category diagrams since in that context the interchange law no longer holds as an equation but is `categorified' to a 3-morphism that interchanges lower morphisms.
 
Where (pre)-2-category appears, both cases are considered at once, the 2-category case being the one with `pre-' deleted everywhere. As in the one-dimensional case, the labelling of a two-dimensional diagram is defined in  terms of elementary diagrams. 

\begin{definition}
$\quad$ Let $\mac C$ be a (pre-)2-category.
  \begin{enumerate}
  \item An {\bf elementary (pre-)2-cate\-gory diagram} $D$ for $\mac C$ is a progressive two-dimensional diagram with exactly one vertex, which meets every line in the diagram, together with:
    \begin{itemize}
    \item  a labelling of each region with an object in $\mathcal C$, 
    \item a labelling of each line with a 1-morphism, 
    \item  a labelling of the vertex with a 2-morphism.
    \end{itemize}
    The top edge  $[0,1]\times \{0\}$  and the bottom edge $[0,1]\times\{1\}$  of the diagram are required to be category diagrams for $\mathcal C_1$ with labels induced by the labelling of regions and lines of $D$. They   evaluate to  1-morphisms which are, respectively,   the source and the target for the  2-morphism at  the vertex.

  \item A {\bf (pre-)2-category diagram} for a (pre-)2-category $\mathcal C$ is  a  progressive two-dimensional diagram
 together with a labelling of each region with an object in $\mathcal C$, 
a labelling of each line with a 1-morphism  and a labelling of each vertex with a 2-morphism, such that
\begin{itemize}
\item  the top and bottom edges are  category diagrams for $\mathcal C_1$,
\item  each vertex $v$ has a neighbourhood that is  an elementary (pre-)2-category diagram after affine linear rescaling. 
\end{itemize}
  \end{enumerate}
\end{definition}

\begin{figure}
  \centering
  \includegraphics[scale=0.3]{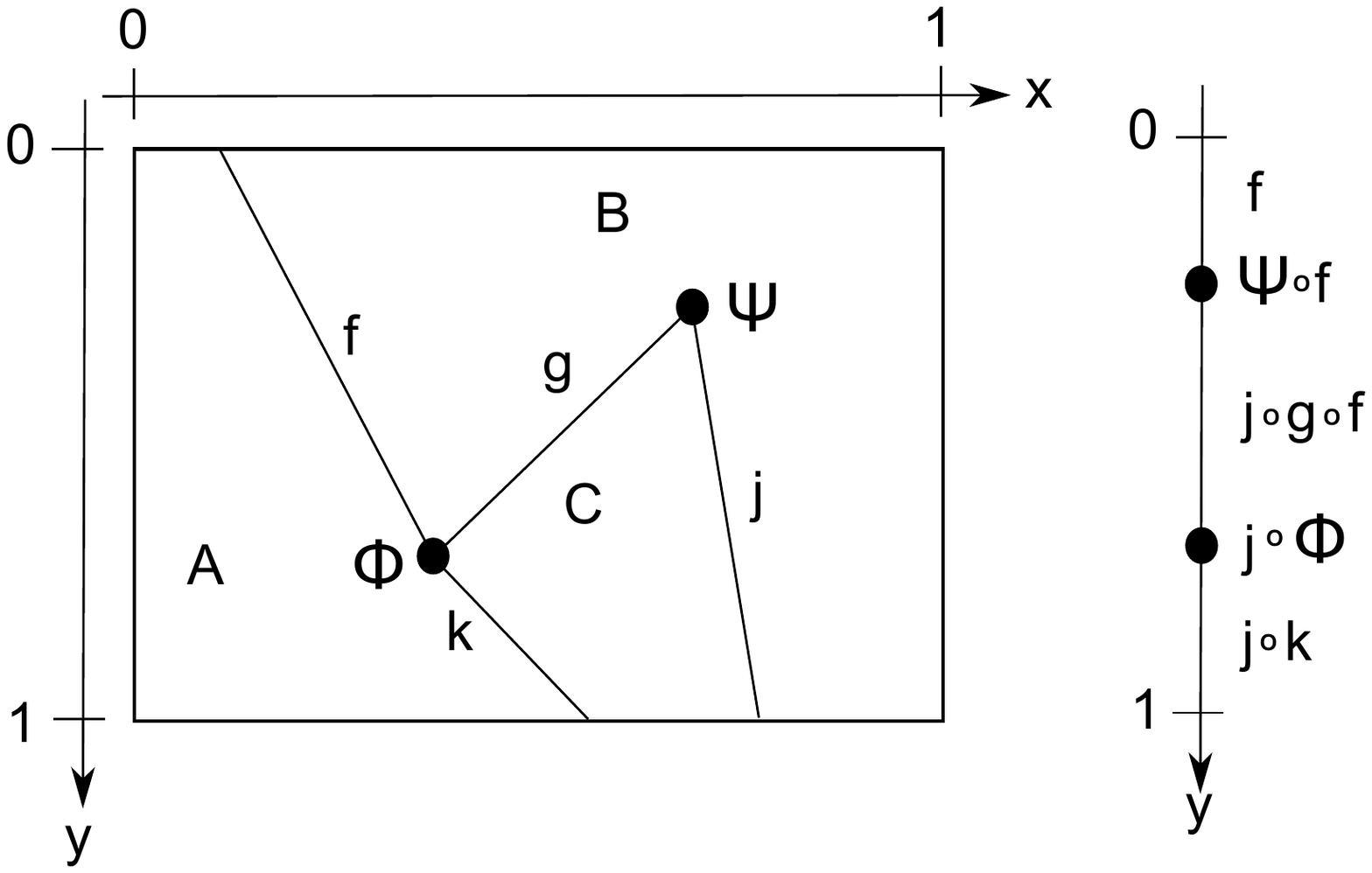}
  \caption{A 2-Category diagram together with its projection onto the $y$-axis.}
  \label{2catdiagram}
\end{figure}

An example of a (pre-)2-category diagram is shown in Figure \ref{2catdiagram}. The requirement that vertices are locally isomorphic to an elementary vertex enforces the condition that the source and targets of 1- and 2-morphisms match.  Important examples are the \emph{identity diagrams} which have a number of vertical lines and no vertices. More precisely, an identity diagram is a diagram of the form $1_D=D\times [0,1]$, where $D$ is a category diagram for $\mathcal C_1$. The regions and lines of $1_D$ correspond to the regions and vertices of $D$. 

\bigskip
Next we define the evaluation of a (pre-)2-category diagram. For a pre-2-category, only the evaluation of  \emph{generic} diagrams (see Definition \ref{genprog2d}) is defined and 
the 2-category diagrams are treated by perturbing them to generic diagrams. So we first treat evaluations of generic diagrams. 
 The evaluation of a generic (pre-)2-category diagram consists of  two steps. 
The first is to  project the (pre-)2-category diagram  to a category diagram   via the 
projection map  $p_1\colon (x,y)\mapsto y$. The second step is the  evaluation of the resulting category diagram. 

The category diagram $p_1D$ is obtained as follows.
Consider a generic (pre-)2-category diagram $D$ as in Figure \ref{2catdiagram},   
where the region containing the left-hand edge $\{0\}\times[0,1]$ is labelled with an object $A$ and the region containing the right-hand edge $\{1\}\times[0,1]$ is labelled with an object $B$ in a (pre-)2-category $\mac C$. Then the projection $p_1D$ is labelled  with data from the category $\mathcal C(A,B)$. The labelling of a point $y\in [0,1]$ depends on whether $p_{1}^\inv(y)$ contains a vertex: 
\begin{enumerate}
\item 
  If    $p_{1}^\inv(y)$  does not contain a vertex of $D$, then $y$ lies in a region of the category diagram  $p_1D$  and this  region is labelled with the horizontal composite  of the 1-morphisms in $p_1^{-1}(y)$, 
  composed as shown in Figure \ref{2catdiagram}.

\item  If  $p_{1}^\inv(y)$ contains a vertex, then $y$ is  a vertex of the category diagram $p_1D$ and is labelled with the horizontal composite  of the 1-morphisms and the single 2-morphism in $p_1^{-1}(y)$ as shown in Figure \ref{2catdiagram}. 
\end{enumerate}
The evaluation of the category diagram  $p_1D$ according to Definition \ref{1d_eval} is a morphism in $\mac C(A,B)$ and hence a 2-morphism in $\mac C$.  This defines the evaluation of a generic (pre-)2-category diagram.

\begin{definition}
\label{generic2diagram}
  The  evaluation of a  {\bf generic (pre-)2-category diagram} $D$ labelled with data from a (pre-)2-category $\mac C$   is the 2-morphism in $\mathcal C$ defined by the evaluation of the category diagram $p_1D$.
\end{definition}

As in the case of tangle diagrams, there are two products for diagrams labelled with data from a (pre-)2-category $\mathcal C$. Vertical composition is defined if the bottom edge of a diagram $D$, with its labelling, matches the top edge of another diagram $D'$  and consists of drawing one diagram above the other (along the $y$-axis). The evaluation of the composite diagram is then given by the vertical composite  of the evaluations of $D$ and $D'$. 

Horizontal composition consists of juxtaposing two diagrams along the $x$-axis and is defined only if the object on the left-hand side of one diagram matches the object on the right-hand side of the other. In the case of a pre-2-category, this product is defined only in the cases where one of the two diagrams is an identity diagram. 

By analogy to the category of tangles, it seems plausible to expect 
that by taking a suitable quotient by isotopies it would be possible to make the diagrams into a (pre-)2-category. However we do not develop this idea here.

\subsection{Invariance of the evaluation of 2-category diagrams}

For a 2-category $\mathcal C$, the 2-category diagrams are dual to the usual pasting diagrams considered in category theory. However the former contain  more information than the latter, namely the values of the $y$-coordinate. It is  therefore important to note that the evaluation is in fact independent of these values, which is a consequence of  the interchange law.  
Formulating this precisely requires an appropriate notion of mappings between generic progressive  diagrams  and a proof that the evaluation of  the diagrams is invariant under these mappings. 

The homomorphisms  in Definition \ref{2dmappings} are too general to give a meaningful notion of mappings between generic progressive diagrams and to preserve their evaluation. For instance, there are examples of isomorphisms that change the order of the lines incident at a vertex.
The appropriate notion of  mappings for generic progressive diagrams was determined by  Joyal and Street \cite{JS} for the case of monoidal categories, which are 2-categories with a single object. The relevant mappings are the ones that are 
determined by a PL isotopy from the identity mapping of the diagram.
As the action of an isotopy  on a 2-category diagram preserves all labels, this result has a direct generalisation to the context of 2-category diagrams.

\begin{theorem} \label{generic2invariance} The evaluation of a generic 2-category diagram is invariant under a piecewise-linear isotopy that starts at the identity mapping and is a one-parameter family of isomorphisms of 2-category diagrams.
\end{theorem}

\begin{figure}
  \centering
  \includegraphics[scale=0.35]{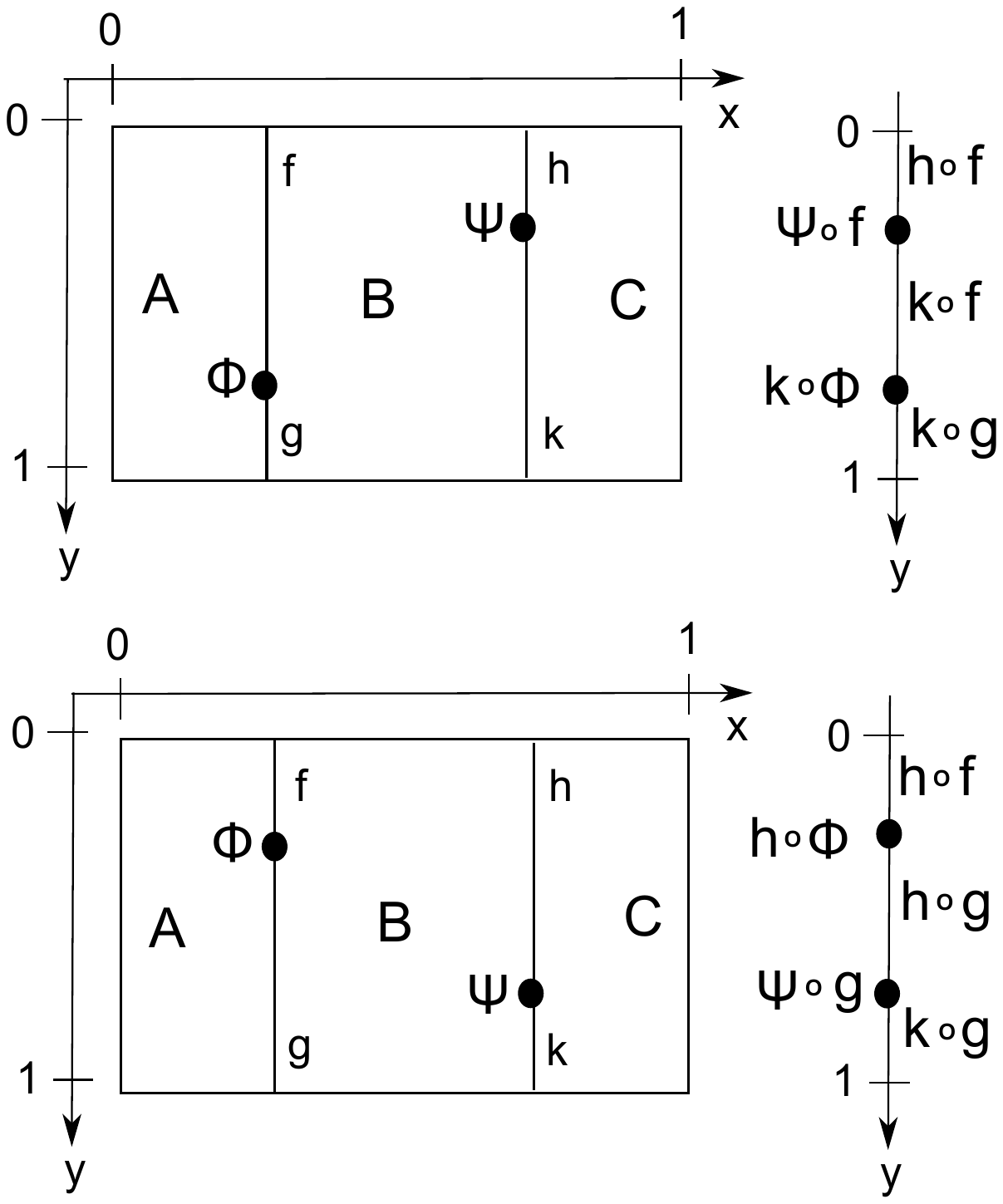}
  \caption{The interchange law: associated 2-category diagrams with associated projections.}
  \label{exchange}
\end{figure}

\begin{proof} This proof is essentially the one given in  \cite[Theorem 1.2]{JS}, where there is more detail. The important point in the statement of the theorem is that at every stage in the isotopy the diagram is progressive. Therefore the category diagram obtained by projection with $p_1$ changes only by an isotopy of $[0,1]$ if the order of the $y$-coordinates of the vertices does not change. In this case it follows that the evaluation is invariant under the isotopy. 

  If the isotopy does change the order of the $y$-coordinates of the vertices, then the isotopy can be perturbed slightly so that they change one at a time. This can be done by composing the isotopy with suitable isotopies of square neighbourhoods of each vertex. The invariance of the evaluation 
  under an isotopy that changes the order of two neighbouring vertices then follows from the
  interchange law (see Figure \ref{exchange}).\end{proof}

By means of  Theorem \ref{generic2invariance}, it is possible to extend the definition of 2-category diagrams to non-generic diagrams by  dropping the requirement that distinct vertices have different $y$-coordinates. This cannot be done for pre-2-category diagrams.

\begin{definition}
 \label{progdiag} The {\bf evaluation of a 2-category diagram} is defined by perturbing it by an isotopy to a generic 2-category diagram and evaluating the resulting diagram.
\end{definition}

  This is well-defined, since the result of the evaluation is independent of the choice of isotopy by Theorem \ref{generic2invariance}.
In fact it is easy to see that the definition could also be extended by allowing the product of more than one 2-morphism in the projection in Definition \ref{generic2diagram} by evaluating directly in the 2-category, and this would give the same result. 

For pre-2-category diagrams, there is a result similar to Theorem  \ref{generic2invariance} but where the isotopy is through generic diagrams, i.e., it is required to preserve the ordering of the vertices by the $y$-coordinate. This ensures that the interchange law is not required.

\subsection{Gray categories}
The  three-dimensional categories considered in this article are Gray categories, the principal example being 2Cat \cite{Gray}, consisting of  strict 2-categories, strict 2-functors,  pseudonatural transformations and modifications (see Appendix \ref{2functors}). Although tricategories present a more general notion of a three-dimensional category, Gray categories have the advantage that their coherence data is stricter than that of a general tricategory and, consequently,  the constructions are less involved. 
Note that the focus on Gray categories  is only a minor restriction on the generality of the constructions, since  every tricategory  $T$ is triequivalent to a Gray category $\mac G$ \cite{GPS, gurski}.

The standard definition of a Gray category \cite{GPS} is a 
category enriched over the monoidal category {\bf Gray}, which is constructed  using the Gray tensor product.  This one-sentence definition, which is not given precisely here, 
is equivalent to  `strict cubical tricategory', as shown by Gordon, Power and Street \cite{GPS}, see also \cite{gurski}. The cubical tricategory point of view turns out to be easier to work with and so forms the starting point for the definitions in this paper. In actual fact, our conventions turn out to be that of a strict opcubical tricategory, the difference between cubical and opcubical being purely presentational; one can convert one to the other by interchanging left and right in the definitions.

The main problem in presenting the definition is that the definition of a tricategory is long and complicated and is much more general than the strict version that is needed here. Accordingly, the definition of a tricategory in \cite{GPS} is specialised to the case of a strict opcubical tricategory in Appendix \ref{2functors}, Definition \ref{strictopcubical}. The definitions of functors and transformations for strict tricategories are also summarised in Appendix \ref{2functors}.  It is worth noting that in general there are some differences between the definition of tricategory in \cite{GPS} and `algebraic tricategory' in \cite{gurski}. However, for strict (op)cubical tricategories, and also their functors and transformations of functors, these definitions coincide. Thus the main definition is

\begin{definition}\label{tricatGray}
 A {\bf Gray category} is a strict opcubical  tricategory.
\end{definition}

This conceptual definition takes some work to unpack; this is done in \cite{Crans} and is summarised here, with some change in notation. The 
unpacked data and relations are called `Gray category data'. This is described explicitly so that this paper is self-contained in the sense that one does not need to know the definition of a tricategory. Lemma \ref{cubegray} shows that this approach is equivalent to the definition as a strict opcubical tricategory. The difference between the two is that in the tricategory approach the Gray product $\Box$ is defined on all 2- and 3-morphisms. This means one can discuss the Gray categories in familiar terms as a collection of 2-categories.

\begin{definition}
\label{graydef}
 {\bf Gray category data} $\mathcal G$ consists of a set of objects, and for any pair of objects $\mathcal C$, $\mathcal D$, of a 2-category $\mathcal G(\mathcal C,\mathcal D)$ of 1-, 2- and 3-morphisms.
  In this 2-category, the notation is as defined previously: $\circ$ for the horizontal composition and $\cdot$ for the vertical composition.

  The additional data is the {\bf Gray product} $\Box$ and the {\bf tensorator}. The Gray product defines a product  $G\Box F$ of 1-morphisms $F\colon \mathcal C\to \mathcal D$ and $G\colon \mathcal D\to \mathcal E$,  which extends to a product $ \Phi\Box F$ of a 1-morphism with a 2- or 3-morphism $\Phi\in\mathcal G(\mathcal D,\mathcal E)$ and to a product $G\Box \Psi$ of a 2- or 3-morphism $\Psi\in\mathcal G(\mathcal C,\mathcal D)$ with a 1-morphism $G$. These products  are required to determine strict 2-functors, $-\Box F$ and $G\Box -$, and the $\Box$ product is required to be strictly unital and associative.  The former means that each object $\mathcal C$ has a unit 1-morphism $1_\mathcal C$ and the 2-functors $-\Box 1_{\mathcal C}$ and $1_{\mathcal C}\Box -$ are the identity 2-functors.  The associativity condition requires that  all  $\Box$-composable morphisms $P,Q,R$, two of which are 1-morphisms and the third a 1- 2- or 3-morphism, satisfy
  $$(P\Box Q)\Box R=P\Box (Q\Box R).$$
  The {\bf tensorator} or {\bf braiding} is the categorification of the interchange law, as noted in Section \ref{2catdiagsection}. It consists of invertible 3-morphisms 
  \begin{align}\label{tens}\sigma_{\mu,\nu}\colon (\mu\Box F_2)\circ(G_1\Box\nu)\Rrightarrow (G_2\Box\nu)\circ(\mu\Box F_1),\end{align} for all composable 2-morphisms $\nu\colon F_1\Rightarrow F_2\in\mac G(\mac C, \mac D)$ and $\mu\colon G_1\Rightarrow G_2\in\mac G(\mac D,\mac E)$. It is required to be an identity 3-morphism if either $\mu$ or $\nu$ is an identity 2-morphism: 
  \begin{align}\label{idtens}
    \sigma_{\mu, 1_{F_1}}=1_{\mu\Box F_1}\qquad\sigma_{1_{G_1},\nu}=1_{G_1\Box \nu}
  \end{align}
  and is required to be natural in both arguments:  
  \begin{align}\label{nattens}
    &\sigma_{\mu,\nu'}\cdot\bigl((\mu\Box F_2)\circ(G_1\Box\Phi)\bigr)
    =\bigl((G_2\Box\Phi)\circ(\mu\Box F_1)\bigr)\cdot\sigma_{\mu,\nu}\\
    &\sigma_{\mu',\nu}\cdot\bigl((\Psi\Box F_2)\circ(G_1\Box\nu)\bigr)\nonumber
    =\bigl((G_2\Box\nu)\circ(\Psi\Box F_1)\bigr)\cdot\sigma_{\mu,\nu}\end{align}
  for all 3-morphisms  $\Phi\colon\nu\Rrightarrow\nu'$, $\Psi\colon\mu\Rrightarrow\mu'$.
  It  is also required to be compatible with the horizontal composition $\circ$ of 2-morphisms: 
  \begin{align}\label{tenshor}
    &\sigma_{\mu,\bar\nu\circ\nu}=
    \bigl((G_2\Box\bar\nu)\circ\sigma_{\mu,\nu}\bigr)\cdot\bigl(\sigma_{\mu,\bar\nu}\circ(G_1\Box\nu)\bigr),\\
    &\sigma_{\bar\mu\circ\mu,\nu}=
    \bigl(\sigma_{\bar\mu,\nu}\circ(\mu\Box F_1)\bigr)\cdot\bigl((\bar\mu\Box F_2)\circ\sigma_{\mu,\nu}\bigr)\nonumber\end{align}
  for all 2-morphisms $\bar\nu\colon F_2\Rightarrow F_3$, $\bar\mu\colon G_2\Rightarrow G_3$.  
  In addition, for all 2-morphisms $\mu$ and $\nu$, and 1-morphisms $F$, the following equations  are required to hold whenever the $\Box$ compositions are defined
  \begin{align}\label{tensbox}\sigma_{\mu\Box F,\nu}=\sigma_{\mu,F\Box\nu}\quad
    \sigma_{F\Box\mu,\nu}=F\Box\sigma_{\mu,\nu}\quad
    \sigma_{\mu,\nu\Box F}=\sigma_{\mu,\nu}\Box F.\end{align}
\end{definition}

\bigskip
Note that  Definition \ref{graydef}  implies the relations
\begin{align}\label{onebox}1_{\nu\Box G}=1_\nu\Box G\qquad 1_{F\Box \nu}=F\Box 1_{\nu}\end{align}
for all  1-morphisms $F, G$ and 2-morphisms $\nu$ for which these expressions are defined. 
Using this definition, it can be checked that the 0- 1- and 2-morphisms of $\mathcal G$ form a pre-2-category, which is denoted $\mathcal G_2$. The 0- and 1-morphisms form a category denoted 
$\mathcal G_1$. Where it is not ambiguous, the symbol $\Box$ may be omitted, so that the product of $G$ and $F$ may be written as just $GF$.

It is now shown that the Gray category data   in Definition \ref{graydef} is equivalent to the Definition \ref{strictopcubical}  of a strict cubical or opcubical tricategory. 
\begin{lemma} \label{cubegray} 
The set of strict  opcubical tricategories is in canonical bijection with 
the set of  Gray category data. Likewise, the set of strict cubical  tricategories is in canonical bijection with 
the set of  Gray category data.
\end{lemma}
\begin{proof}
  Let $\mac G$ be a strict  (op)cubical tricategory according to Definition \ref{strictopcubical} 
  with composition $\Box\colon \mac G(\mac D,\mac E)\times \mac G(\mac C,\mac D)\to\mac G(\mac C,\mac E)$ and coherence 3-morphisms
  $$
  \Box_{\mu,\nu}\colon  (\mu_1\Box\mu_2)\circ (\nu_1\Box \nu_2)\to (\mu_1\circ \nu_1)\Box (\mu_2\circ\nu_2)   
  $$
  for all $\Box$-composable pairs of 2-morphisms $\mu=(\mu_1,\mu_2)\colon (H_1,H_2)\to (K_1,K_2)$, $\nu=(\nu_1,\nu_2)\colon  (G_1,G_2)\to (H_1,H_2)$. Then $\Box$ defines the Gray product of 1-morphisms with 1-, 2- and 3-morphisms, and the tensorator   is given by
  $$
  \sigma_{\nu_1,\mu_2}=\Box^\inv_{(1_{H_1}, \mu_2),(\nu_1, 1_{H_2})}
  $$
  in case $\mac G$ is opcubical and by
  $$
  \sigma_{\mu_1,\nu_2}=\Box_{(\mu_1, 1_{H_2}),(1_{H_1}, \nu_2)}
  $$
  in case $\mac G$ is cubical, the (op)cubical condition being used to show that the source and target are correct.
  A direct computation shows that the axioms of a strict (op)cubical tricategory in Definition \ref{strictopcubical} imply  that the conditions in Definition \ref{graydef} are satisfied.

  Conversely, if $\mac G$ is a set of Gray category data, then one obtains a strict opcubical tricategory by   promoting the left-hand-side of \eqref{tens} to the product of 2- and 3-morphisms
  \begin{equation}\label{eq:product}
    \Psi\Box\Phi=\left(\Psi\Box F_2\right)\circ\left(G_1\Box\Phi\right)
  \end{equation}
  for all 3- or 2-morphisms $\Phi\in\mac G(F_1,F_2),\Psi\in\mac G( G_1, G_2)$ and 1-morphisms $F_1,F_2\colon \mac C\to\mac D$, $G_1,G_2\colon\mac D\to\mac E$. 
  The coherence morphisms for $\Box$ are then given by the collection of natural isomorphisms
  $$
  \Box_{\mu,\nu}=1_{\mu_1\Box K_2 }\circ \sigma^\inv_{\nu_1,\mu_2}\circ 1_{G_1\Box \nu_2}\colon (\mu_1\Box\mu_2)\circ (\nu_1\Box\nu_2)\to (\mu_1\circ \nu_1)\Box (\mu_2\circ\nu_2)
  $$
  for all $\Box$-composable pairs of 2-morphisms $\mu=(\mu_1,\mu_2)\colon(H_1,H_2)\to (K_1,K_2)$, $\nu=(\nu_1,\nu_2)\colon (G_1,G_2)\to (H_1,H_2)$. That this  determines a collection of weak 2-functors 
  $\Box\colon \mathcal G(\mathcal D,\mathcal E)\times\mathcal G(\mathcal C,\mathcal D)\to\mathcal G(\mathcal C,\mathcal E)$
  with strict units (see  Definition \ref{laks2func})
  is a direct consequence of  the axioms of the Gray category data. 
  Consistency condition (1) of  Definition \ref{laks2func} follows from Definition \ref{graydef} \eqref{idtens}, and consistency condition (2) of  Definition \ref{laks2func} from Definition \ref{graydef} \eqref{tenshor}. That the functor $\Box$ is opcubical follows directly from the definition.  

  Analogously, one obtains a strict cubical tricategory by  promoting the right-hand-side of \eqref{tens} to the product of 2- and 3-morphisms
  \begin{equation}\label{eq:product2}
    \Psi\Box\Phi=\left(G_2\Box \Phi\right)\circ\left(\Psi\Box F_1\right)
  \end{equation}
  for all 3- or 2-morphisms $\Phi\in\mac G(F_1,F_2),\Psi\in\mac G( G_1, G_2)$ and 1-morphisms $F_1,F_2\colon \mac C\to\mac D$, $G_1,G_2\colon\mac D\to\mac E$. 
  The coherence morphisms for $\Box$ are then given by the collection of natural isomorphisms
  $$
  \Box_{\mu,\nu}=1_{K_1\Box \mu_2}\circ \sigma_{\mu_1,\nu_2}\circ 1_{\nu_1\Box G_2}\colon (\mu_1\Box\mu_2)\circ (\nu_1\Box\nu_2)\to (\mu_1\circ \nu_1)\Box (\mu_2\circ\nu_2)
  $$
  for all $\Box$-composable pairs of 2-morphisms $\mu=(\mu_1,\mu_2)\colon(H_1,H_2)\to (K_1,K_2)$, $\nu=(\nu_1,\nu_2)\colon (G_1,G_2)\to (H_1,H_2)$. The proof that this defines a strict cubical tricategory is analogous to the opcubical case.
\end{proof}

The passage between a set of Gray category data and the associated cubical and opcubical tricategories can be viewed as a special case of the operation called ``nudging'' in \cite{GPS}, which allows one to pass between cubical and opcubical tricategories and functors. 
The details depend on the level of strictness of the tricategory and the functors. Here, the relevant notions are  2-strict and  strict (op)cubical  tricategories, which are strict tricategories with additional conditions, see Definition \ref{strictopcubical}, and  functors of 2-strict tricategories, see Definition \ref{grayfunc}.
This gives rise to the following statement. 

\begin{corollary}
  \label{lemma:cubical-to-opcubical-tricat}
  For every strict cubical (opcubical) tricategory $\mathcal{G}$, there exists a canonical strict opcubical (cubical) tricategory $\widehat{\mathcal{G}}$ and functors of 2-strict tricategories $\Sigma \colon \mathcal{G} \rightarrow 
  \widehat{\mathcal{G}}$, $\Sigma^{-1} \colon \widehat{\mathcal{G}} \rightarrow \mathcal{G}$, that are the identity mappings on all objects and morphisms and satisfy $\Sigma \circ \Sigma^{-1}=1$, $\Sigma^{-1} \circ \Sigma=1$. 
\end{corollary}
\begin{proof}
  Let $\mathcal{G}$ be a strict opcubical tricategory. Then by Lemma \ref{cubegray}, 
$\mathcal{G}$ defines a set of Gray category data. 
  Define  $\widehat{\mac G}$ as  the strict cubical tricategory determined by this set of Gray category data according to Lemma \ref{cubegray} and 
  define the functor $\Sigma\colon\mathcal{G} \rightarrow \widehat{\mathcal{G}}$ of 2-strict tricategories
  by taking the identity mappings on the objects and the identity functors for each 2-functor $\Sigma_{\mac C,\mac D}\colon \mac G(\mac C,\mac D)\to\widehat{\mac G}(\mac C,\mac D)$. 
  The only nontrivial data of  $\Sigma$ are the natural isomorphisms $\kappa_{\mu,\nu}\colon  \mu \widehat{\Box} \nu \rightarrow \mu \Box \nu $ from Definition \ref{grayfunc}, where $\Box$ and $\widehat\Box$ denote, respectively, the products in the tricategories $\mac G$ and $\widehat{\mac G}$. 
  These are given by the tensorator:
  $$\kappa_{\mu,\nu}=\sigma_{\mu,\nu}^{-1}$$ for all $\Box$- composable 2-morphisms $\mu,\nu$. It follows directly from the  properties of  the tensorator in Definition \ref{graydef} that this defines a functor $\Sigma\colon\mac G\to\widehat{\mac G}$ of 2-strict tricategories.
  By taking again the identity mappings on the objects and the identity 2-functors together with the coherence isomorphisms
  ${\kappa}_{\mu,\nu}^{-1} = \sigma_{\mu,\nu}$, we obtain a functor   that is strictly inverse to $\Sigma$.
\end{proof}

This corollary implies in particular that the tricategories 
obtained  via \eqref{eq:product} and \eqref{eq:product2}
are equivalent and one of these constructions can be chosen arbitrarily. 
The definition \ref{tricatGray} of a Gray category uses the strict {\em opcubical} tricategory. 

The notation for the $\Box$ product in a Gray  category used in the following is as follows. If $F$ is a 1-morphism and $\phi$ is a 2- or 3-morphism, then the preferred notation for the product is using the Gray category data convention rather than the tricategory convention, i.e., $F\Box\phi$ rather than $1_F\Box\phi$ or $1_{1_F}\Box \phi$ (and similarly for $F$ on the right). Also, this is abbreviated to $F\phi$ in places. In contrast, the $\circ$ product will continue to be written in full.
\medbreak

A special case of a Gray category is a braided monoidal category. The following is proved in \cite[Prop 8.6,8.7]{GPS}. 

\begin{lemma}\label{braidten}
  If $\mac G$ is a Gray category, then for every object $\mac C$ the category $\mac G(1_{\mac C}, 1_{\mac C})$ is a braided strict monoidal  category. 
  Conversely, a braided strict monoidal category is a Gray category with a single object and a single 1-morphism.
\end{lemma}

\subsection{Example: 2Cat and 2Cat(1)} \label{sec:2cat}
An example of a Gray category investigated in this paper is  2Cat. The objects of 2Cat are 2-categories, the 1-morphisms are strict functors of 2-categories, the 2-morphisms pseudo-natural transformations and 3-morphisms modifications.  In the following, these definitions,  the compositions  and the tensorator are recalled. For a proof that this defines a Gray category, see \cite[\S I,4.5]{Gray}. 

To simplify the notation, and because it is the only case used in the following, it is assumed throughout that the 2-categories have only one object,   denoted $\bullet$, and can thus be described as monoidal categories \cite{ChengGurski}. Thus the Gray category defined here is 2Cat(1) in the notation of Cheng and Gurski. The definitions extend easily to the case of a general  2-category. 

The objects of 2Cat(1) are strict monoidal categories.  The objects of the monoidal category 
correspond to the 1-morphisms of the corresponding 2-category with one object and  its morphisms to 2-morphisms of the corresponding 2-category with one object. The horizontal composition  is given by the monoidal product, denoted $\circ$ in the following, and the unit 1-morphism
in the associated 2-category corresponds to the monoidal unit $e$. The vertical composition of morphisms is denoted $\cdot$, as before.
The 1-morphisms in 2Cat(1) are strict monoidal functors.
\begin{definition}
\label{tensfunc}
  A {\bf strict monoidal functor} $F\colon \mathcal C \to \mathcal D$ between strict monoidal categories $\mac C$, $\mac D$ is a functor $F\colon\mac C\to\mac D$ with $F(e_{\mac C})=e_{\mac D}$  and $F(x\circ y)=F(x)\circ F(y)$ for all objects $x$ and $y$ of $\mac C$.
\end{definition}

The 2-morphisms in 2Cat(1) are pseudo-natural transformations between strict monoidal functors. 
They can be viewed  as a generalisation of natural transformations and are obtained by restricting the general definition of natural 2-transformations in Definition \ref{nat2transformations} to 2-categories with a single object.

\begin{definition}
  \label{psnat}
  A {\bf pseudo-natural transformation} $\mu\colon F\Rightarrow G$ between strict monoidal functors $F,G\colon\mathcal C\rightarrow\mathcal D$ 
  consists of an object  $x$ of $\mathcal D$, together with a collection of isomorphisms  $\mu_y\colon x\circ F(y)\rightarrow G(y)\circ x$ for all objects $y$ of $C$ such that
  \begin{enumerate}
  \item  $\mu_e=1_x\colon x\to x$ is the identity morphism.
  \item $\mu$ is natural in $y$:     for  all  morphisms $\alpha\colon y\rightarrow z$ in $\mac C$, the following diagram commutes
    \begin{align}\label{psnatdiag}
      \xymatrix{ x\circ F(y) \ar[d]^{1_x\circ F(\alpha)} \ar[r]^{\mu_y} & G(y)\circ x\ar[d]^{G(\alpha)\circ 1_x}\\
        x\circ F(z) \ar[r]^{\mu_{z}} & G(z)\circ x.}
    \end{align}
  \item $\mu$ is compatible with the monoidal product:  for all objects $y,z$ 
    \begin{align}\label{psnattens}\mu_{y\circ z}=(1_{G(y)}\circ \mu_z)\cdot (\mu_y\circ 1_{F(z)}).
    \end{align}
  \end{enumerate}
\end{definition}

The 3-morphisms in 2Cat(1) are modifications between pseudo-natural transformations. Their definition is obtained by restricting  Definition \ref{modi_gen}  to 2-categories with a single object.

\begin{definition}
\label{modidef} Let  $\mu, \nu\colon F\Rightarrow G$ be pseudo-natural transformations with  component morphisms
  $\mu_z\colon x\circ F(z)\rightarrow G(z)\circ x$, $\nu_z\colon y\circ F(z)\rightarrow G(z)\circ y$ for all objects $z$ of $\mathcal C$.
  A {\bf modification}  $\Phi\colon \mu \Rrightarrow \nu$ 
  is a morphism $\Phi\colon x\rightarrow y$ such that the following diagram commutes for all objects $z$ of $\mathcal C$
  \begin{align}
    \label{modi}
    \xymatrix{ x\circ F(z) \ar[d]^{\Phi\circ 1_{F(z)}} \ar[r]^{\mu_z} &  G(z)\circ x\ar[d]^{1_{G(z)}\circ\Phi}\\
      y\circ F(z)\ar[r]^{\nu_z} & G(z)\circ y.}
  \end{align}
\end{definition}

\begin{figure}
  \centering
  \includegraphics[scale=0.4]{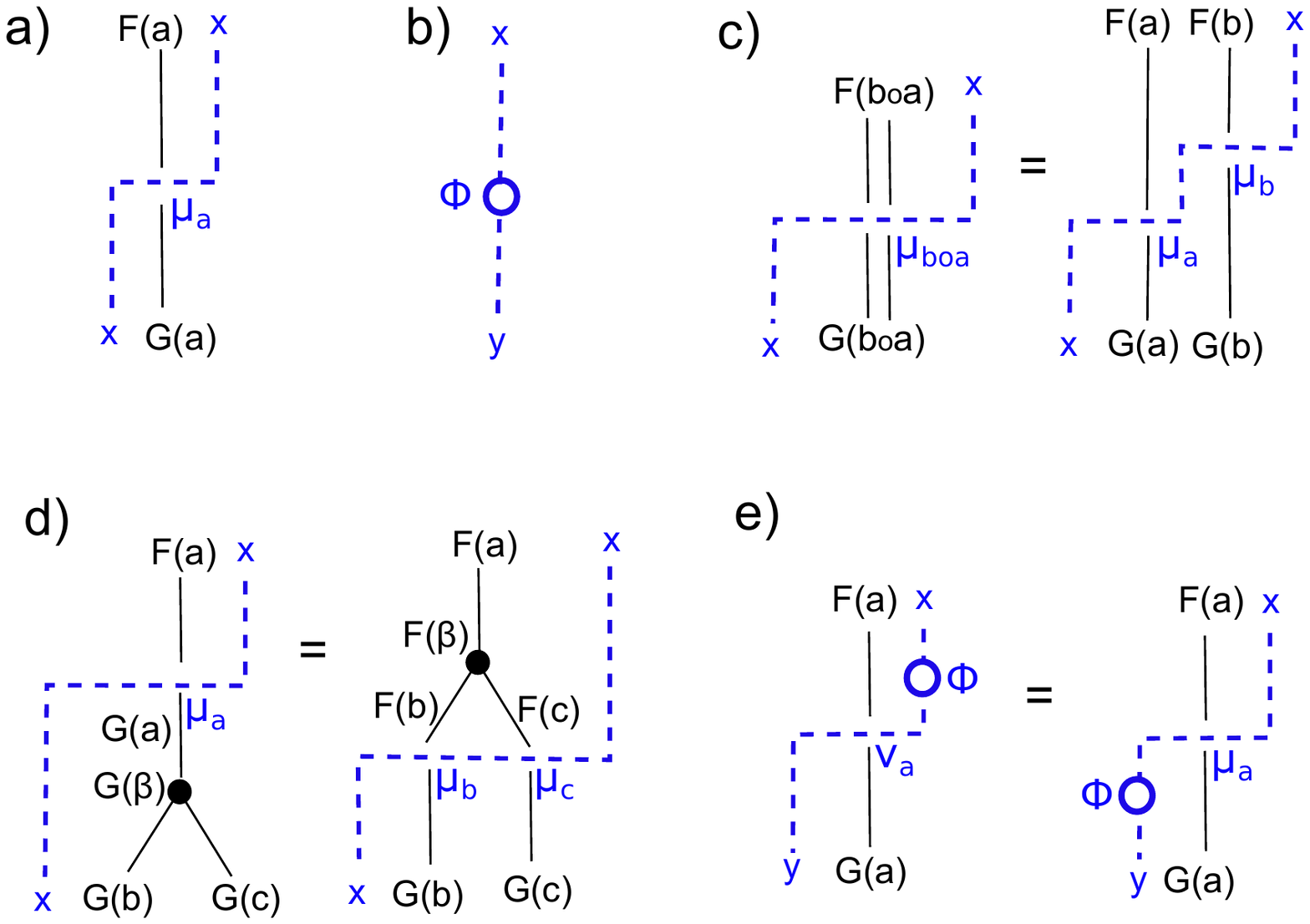}
  \caption{Diagrams for  pseudo-natural transformations and modifications.\newline
    a) Pseudonatural transformation $\mu\colon F\Rightarrow G$ with component morphisms
    $\mu_a\colon x\circ F(a)\to G(a)\circ x$.\newline
    b) Modification  $\Phi\colon \nu\Rrightarrow \mu$ between pseudo-natural transformations with component morphisms $\mu_a\colon x\circ F(a)\to G(a)\circ x$, $\nu_a\colon y\circ F(a)\to G(a)\circ y$.
    \newline c) Compatibility of pseudo-natural transformation with the monoidal product. 
    \newline d) Naturality property  of pseudo-natural transformations. \newline e) Defining property of modifications.
  }
  \label{psnatmodi}
\end{figure}

The defining properties of pseudo-natural transformations and modifications are depicted  in Figure  \ref{psnatmodi}.

\medskip
The product operations  and the tensorator of 2Cat(1)  and the tensorator are obtained by specialising the ones in 2Cat to the case of a single object and  are summarised  in the following definition.

\begin{definition}\label{moncat} The product operations and the tensorator of 2Cat(1) are as follows:
  \begin{enumerate}
  \item {\bf Gray product:} 

    \begin{itemize}
    \item The composition of functors $F\colon\mathcal B\rightarrow \mathcal C$, $G\colon\mathcal C\rightarrow\mathcal D$ defines the Gray product $G\Box F\colon\mathcal B\rightarrow\mathcal D$.

    \item The product $F\Box \mu$ of a functor $F$ with a pseudo-natural transformation $\mu\colon G\Rightarrow H$, $\mu_y\colon x\circ G(y)\rightarrow H(y)\circ x$,
      is the pseudo-natural transformation $F\mu\colon FG\Rightarrow FH$ with component morphisms $(F\Box \mu)_y=F(\mu_y)\colon F(x)\circ FG(y)\rightarrow FH(y)\circ F(x)$.
    \item The product $\mu\Box K$ of a functor $K$ with a pseudo-natural transformation $\mu\colon G\Rightarrow H$, $\mu_y\colon x\circ G(y)\rightarrow H(y)\circ x$,
      is the pseudo-natural transformation $\mu \Box K\colon GK\Rightarrow HK$ with component morphisms $( \mu\Box K)_y=\mu_{K(y)}\colon x\circ GK(y)\rightarrow HK(y)\circ x$.
    \item The product $F\Box \Phi$ of functor $F$ with a modification $\Phi\colon\mu\Rrightarrow \nu$, $\Phi\colon x\to y$  is  defined by the morphism $F(\Phi)\colon F(x)\rightarrow F(y)$.
    \item The product $\Phi\Box K$ of a functor $K$ with a modification $\Phi\colon\mu\Rrightarrow \nu$, $\Phi\colon x\to y$ is given by the morphism $\Phi\colon x\rightarrow y$.
    \end{itemize}

    \medskip
  \item {\bf horizontal composition:}

    \begin{itemize}
    \item The horizontal composite of pseudo-natural transformations $\mu\colon G\Rightarrow H$, $\mu_y\colon x\circ G(y)\rightarrow H(y)\circ x$ and $\nu\colon H\Rightarrow K$, $\nu_y\colon z\circ H(y)\rightarrow K(y)\circ z$ is the pseudo-natural transformation $\nu\circ\mu\colon G\Rightarrow K$ with component morphisms\newline $(\nu\circ\mu)_y=(\nu_y\circ 1_x)\cdot(1_z\circ\mu_y)\colon z\circ x\circ G(y)\to K(y)\circ z\circ x.$

    \item The horizontal composite $\Psi\circ \mu\colon \nu\circ \mu\Rrightarrow \rho\circ \mu$ of a pseudo-natural transformation 
      $\mu\colon F\Rightarrow G$, $\mu_y\colon x\circ F(y)\to G(y)\circ x$ and a modification $\Psi\colon \nu\Rrightarrow \rho$, $ w\to z$ 
      is given by the morphism $\Psi\circ 1_x\colon w\circ x\to z\circ x$.

    \item The horizontal composite $\mu\circ \Phi\colon \mu\circ \tau\Rrightarrow \mu\circ \rho$ of a pseudo-natural transformation 
      $\mu\colon H\Rightarrow K$, $\mu_y\colon x\circ H(y)\to K(y)\circ x$ and a modification $\Phi\colon \tau\Rrightarrow \rho$, $\Phi\colon w\to z$  
      is  given by the morphism $1_x\circ \Phi\colon x\circ w \to x\circ z$.
    \end{itemize}

    \medskip
  \item {\bf vertical composition:} The vertical composition  of two modifications $\Phi\colon \mu\Rrightarrow\nu$,  $\Phi\colon  w\to x$  and $\Psi\colon\nu\Rrightarrow\rho$, $\Psi\colon x\to z$  is given by the composition of the associated morphisms  $\Psi\cdot \Phi\colon w\xrightarrow{\Phi} x\xrightarrow{\Psi} z$.

  \item  {\bf tensorator:} The tensorator $\sigma_{\mu,\nu}$  of two pseudo-natural transformations $\nu\colon F_1\Rightarrow F_2$, $\mu\colon G_1\Rightarrow G_2$ with associated morphisms $\nu_z\colon x\circ F_1(z)\rightarrow F_2(z)\circ x$ and  $\mu_z\colon y\circ G_1(z)\rightarrow G_2(z)\circ y$
    is the modification $\sigma_{\mu,\nu}\colon (\mu\Box F_2)\circ (G_1\Box \nu)\Rrightarrow(G_2\Box \nu)\circ (\mu\Box F_1)$ given by the morphism $\mu_x\colon y\circ G_1(x)\rightarrow G_2(x)\circ y$.

  \end{enumerate}
\end{definition}

Particular examples of subcategories of 2Cat(1) that are  of interest are the subcategory  obtained by restricting attention to 2-functors on a single category $\mathcal C$, which is a monoidal 2-category, 
and the subcategory of this that consists of pseudo-natural transformations and modifications on a trivial functor. This is a braided monoidal category, called the centre of $\mac C$. Further examples arise when the monoidal categories  $\mathcal C$ have more structure.

Other important examples of Gray categories are Gray groupoids, which are obtained from 2-crossed  modules \cite{KP}  and, more generally, Gray categories obtained from the strictification of tricategories.

\subsection{Gray category diagrams}\label{sec:graycatd}

The definition of a diagrammatic calculus for Gray categories follows the pattern for categories and 2-categories. The diagrams are a three-dimensional generalisation of the two-dimensional diagrams defined above, and were previously studied informally by Trimble \cite{Tweb}.

Gray category diagrams   are located in the unit cube $[0,1]^3$ and consist of a number of points, lines, surfaces, etc.~in the cube. It seems that the clearest way to organise the definition is in terms of a PL stratification, which is just called a stratification here. 

\begin{definition}\label{def:stratification} A {\bf stratification} of an $n$-dimensional  
manifold $X$ is a set of closed subspaces $\emptyset=X^{-1}\subset X^0\subset X^1\subset X^2\subset\ldots\subset X^n=X$ called the $k$-skeleta, such that $X^k\setminus X^{k-1}$ is a manifold of dimension $k$, with boundary $\partial 
\bigl(X^k\setminus X^{k-1}\bigr)
\subset\partial X$. Each component of $X^k\setminus X^{k-1}$ is called a $k$-stratum. 
\end{definition} 

For each stratification of a compact manifold, there is a simplicial complex of which all $k$-skeleta of the stratification are  subcomplexes. \cite[Addendum 2.12]{RS}. This implies in particular that the set $X^0$ of vertices is finite.

A stratification is called homogeneous if it obeys an additional condition (adapted from \cite{Stratification},
where it is called CS, to allow boundaries), to ensure that each stratum has a cross-section 
that does not change along the stratum. This is defined formally in the following definition.

If $L$ is a compact polyhedron then the cone over $L$ is denoted $cL$ with vertex $v$. If $L$ is stratified, then $cL$ is stratified in an obvious way, with skeleta $(cL)_{k+1}=c(L_k)$ and $0$-skeleton $\{v\}$. Recall that $\R^n_+$ denotes the half-space $\{(x_1,\ldots,x_n)|x_n\ge0\}$. The product of $cL \times\R^n_+$ is also a stratified space, in which each stratum is a product of  $\R^n_+$ with a stratum of $cL$.

\begin{definition}\label{def:}
A $k$-stratum is called {\bf homogeneous} if there exists a stratification of an $(n-k-1)$-sphere $L$ and, for any point $x$ in the stratum, an isomorphism (of stratifications) of an open neighbourhood of $x$ in $X$ with
an open subset of $(cL)\times \R^k_+$. A {\bf homogeneous stratification} of a manifold $X$ is a stratification of $X$ such that each stratum is homogeneous.
\end{definition}

Note that the homogeneity condition is automatically satisfied for the $n$-strata of an $n$-manifold by the definition of a manifold, and also for $0$-strata by the definition of polyhedra. For $(n-1)$-strata, the homogenity condition is equivalent to local flatness, which  always holds for $n\le4$, by the Sch\"onflies theorem. Thus all stratifications are homogeneous for $n\le2$. For $n=3$, the condition can fail on the 1-strata.

\begin{figure}
  \centering
  \includegraphics[scale=0.3]{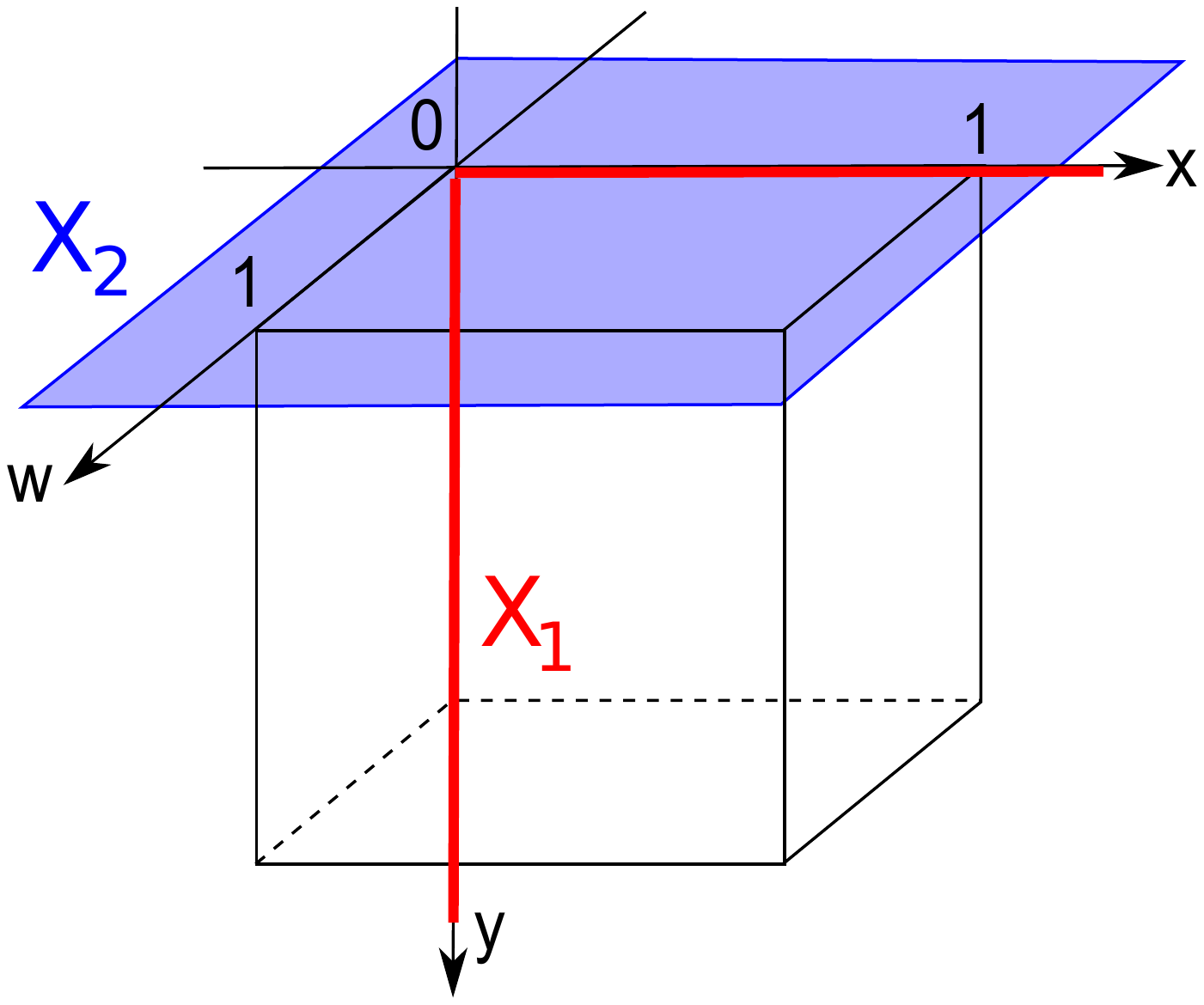}
  \caption{ A stratification of $\mathbb R^3$ that is not homogeneous.
  }
  \label{nothom}
\end{figure}

\begin{example} Let $X=\R^3$, $X_0=\emptyset$, $X_1=\{(0,0,y)|y\ge0\}\cup\{(0,x,0)|x\ge0\}$, $X_2=\{(w,x,0)|w,x\in\R\}$, as shown in Figure \ref{nothom}. This is a stratification of $\R^3$ that is not homogeneous, since the 1-stratum is not homogeneous.
\end{example}

Note that this example can be made homogeneous by putting $X_0=\{(0,0,0)\}$. In general, any stratification in three dimensions can be made homogeneous by adding extra vertices to the 1-skeleton.

\begin{definition}
\label{3ddiag} $\quad$
  \begin{enumerate}
  \item A {\bf three-dimensional diagram} is a homogeneous stratification of $[0,1]^3$ so that 
    \begin{enumerate}\item\label{3d1}
      each 
$k$-stratum 
$C^k$ satisfies
      $$ \partial [0,1]^3\cap C^k=C^k\cap\left((0,1)^{3-k}\times\partial[0,1]^k\right),$$
    \item\label{3d2} the side faces $[0,1]\times\{0\}\times[0,1] $ and $[0,1]\times\{1\}\times[0,1] $ are progressive two-dimensional diagrams.
    \end{enumerate} 
    The 0-, 1-, 2- and 3-strata are called {\bf vertices}, {\bf lines}, {\bf surfaces} and {\bf regions}, of the diagram, respectively.  The face $[0,1]^2\times\{0\}$ is called the {\bf source} and the face 
    $[0,1]^2\times\{1\}$ the {\bf target} of the diagram.

  \item A three-dimensional diagram is called {\bf progressive}, if the projection $p_2\colon (w,x,y)\mapsto (x,y)$ is a regular mapping of each surface, and the projection $p_1\circ p_2\colon (w,x,y)\mapsto y$ a regular mapping of each line.

  \item A progressive three-dimen\-sional diagram is called {\bf generic} if the following conditions on the image of the diagram under the projection $p_2\colon(w,x,y)\mapsto(x,y)$ hold:
    \begin{enumerate}
     \item The images of any two lines meet only at interior points of $[0,1]^2$,  and at every  point where they meet, they cross transversally.
     \item Crossing points and the projections of vertices do not coincide with the projections of any other vertices or points on lines.    
\end{enumerate}
  \end{enumerate}
\end{definition}

Just as the number of vertices, the number of crossing points in the projection of  a generic  three-dimensional diagram is finite. This follows because the image of a compact polyhedron under a PL map is a compact polyhedron \cite[Sec.~2.5]{RS}.  Applied  to the projection of the 1-skeleton to $\R^2$, this yields that projection is a compact polyhedron and hence can only have a finite number of crossing points.

Condition (1) (a) in this definition states that lines in the diagram can intersect the boundary of the unit cube only in its top face $[0,1]^2\times\{0\}$ and its bottom face $[0,1]^2\times\{1\}$ and that surfaces of the diagram cannot intersect its front face $\{1\}\times[0,1]^2$ or its back face $\{0\}\times [0,1]^2$. This is the three-dimensional analogue of the condition that lines in two-dimensional diagrams intersect only the top and bottom edge of the diagram.
Condition (1) (b) is a new feature that does not appear in lower dimensions; in dimension two the side edges are required to be empty. 
For a generic three-dimensional diagram, the source and target are generic two-dimensional diagrams.

\bigskip
Mappings of three-dimensional diagrams are defined in analogy to the one- and two-dimensional case.
\begin{definition}
 $\quad$
  \begin{enumerate}
  \item A {\bf mapping} $D\to D'$ {\bf of three-dimensional diagrams} is a PL embedding $m\colon [0,1]^3\to [0,1]^3$ that preserves the $k$-skeleta, i.e., $m(X^k)\subset X'^k$, for $k=0,1,\ldots 3$.
  \item A mapping of three-dimensional diagrams is called a {\bf homomorphism} if it is a PL homeomorphism and is the identity map on the boundary $\partial[0,1]^3$. An {\bf isomorphism} of three-dimensional diagrams  is a homomorphism that has an inverse.

  \item {\bf Subdivisions} and {\bf subdiagrams} are defined analogously to  the one- and two-dimensional cases (Definition \ref{1dmap}, Definition \ref{2dmappings}).
  \end{enumerate}\label{3dmappings}
\end{definition}

To define  diagrams labelled with data from a Gray category, it is necessary to restrict attention to progressive diagrams. As in the two-dimensional  case,  each vertex in a progressive three-dimensional diagram  has a neighbourhood that is isomorphic to a  diagram with a single vertex. The simplest type of one-vertex diagram has no crossings in its two-dimensional projection and is called an elementary diagram.

\begin{definition}\label{def:elementary}
 \label{graydiag} 
 Let $\mac G$ be a Gray category.
  \begin{enumerate}
  \item An {\bf elementary Gray category diagram} for $\mac G$ is a progressive three-dimensional diagram with one vertex  such that the images of its lines  under the projection $p_2\colon (w,x,y)\mapsto (x,y)$ do not intersect together with
    \begin{itemize}
    \item a labelling of each  region with an object in $\mathcal G$
    \item a labelling of  each surface with a 1-morphism in $\mac G$
    \item a labelling of  each line with a 2-morphism  in $\mac G$
    \item a labelling of the vertex with a 3-morphism in $\mac G$
    \end{itemize}
    The source and target of the diagram are required to be  pre-2-category diagrams for $\mathcal G_2$ and  evaluate to  2-morphisms which are  the source  and target for the vertex 3-morphism. 
  \item A {\bf Gray category diagram} for  $\mathcal G$ is  a progressive three-dimensional diagram 
together with a labelling of each region with an object in $\mathcal G$,
 each surface with a 1-morphism, each line with a 2-morphism and each vertex with a 3-morphism.
 The source and target are required to be (generic) pre-2-category diagrams for $\mathcal G_2$ and each vertex $v$ is required to have a rectangular neighbourhood that is an elementary Gray category subdiagram
after affine rescaling. 
  \end{enumerate}
\end{definition}

\begin{figure}
  \centering
  \includegraphics[scale=0.35]{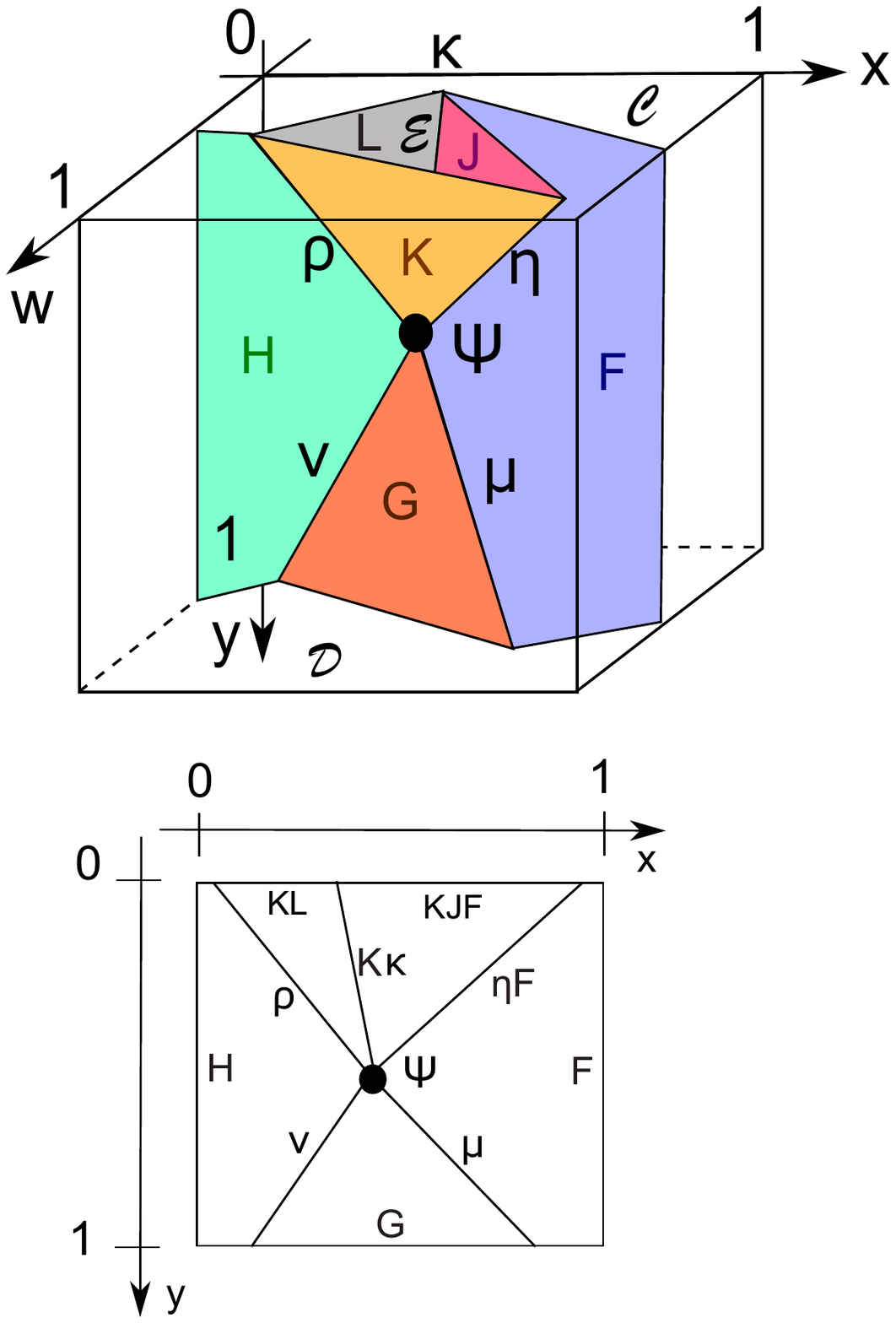}
  \caption{ Elementary Gray category diagram  and its projection to a 2-category diagram for 
    $\mathcal G(\mathcal C,\mathcal D)$: \newline
    $\bullet$ The regions are labelled with objects $\mathcal C,\mathcal D,\mathcal E$. \newline
    $\bullet$ The surfaces   with 1-morphisms $F,G,H\colon\mathcal C\rightarrow\mathcal D$, $J\colon\mathcal D\rightarrow\mathcal E$, $K\colon\mathcal E\rightarrow\mathcal D$ and $L\colon\mathcal C\rightarrow\mathcal E$.\newline
    $\bullet$ The lines with 2-morphisms $\rho\colon H\Rightarrow KL$, $\kappa\colon L\Rightarrow JF$, $\eta\colon KJ\Rightarrow 1_{\mathcal D}$, $\nu\colon H\Rightarrow G$, $\mu\colon G\Rightarrow F$.\newline
    $\bullet$ The vertex with a 3-morphism\newline  $\Psi\colon (\eta F)\circ (K\kappa)\circ \rho \Rrightarrow \mu\circ \nu$.
  }
  \label{graycatdiag}
\end{figure}

An elementary Gray category diagram and its projection are shown in Figure \ref{graycatdiag}. 
Note that the requirement on the vertex neighbourhoods is both a restriction on the topology at a vertex, so that the plane projections of lines are locally non-intersecting, and a restriction on the vertex label. Any two such rectangular neighbourhoods give isomorphic elementary vertex subdiagrams, so that the choice of rectangular neighbourhood does not matter.

Note also that a diagram with a plane projection that has a single crossing and no other vertices is a Gray category diagram but not an \emph{elementary} Gray category diagram. The requirement that the source and target are pre-2-category diagrams implies that the $x$-coordinates of  the intersection points of lines with the source and target of the diagram are all different.

As in two dimensions, the evaluation of a Gray category diagram is  obtained by projecting it to a two-dimensional 2-category diagram and then evaluating the resulting 2-category diagram  according to Definition \ref{generic2diagram}. The construction is described first for the case of a generic diagram. 

The image of the diagram under the  projection $p_2\colon(w,x,y)\mapsto(x,y)$ defines a two-dimensional diagram. This two-dimensional diagram has vertices given by the image of vertices or crossing points, and lines given by segments of images of lines between either vertices or crossing points.

\begin{figure}
  \centering
  \includegraphics[scale=0.45]{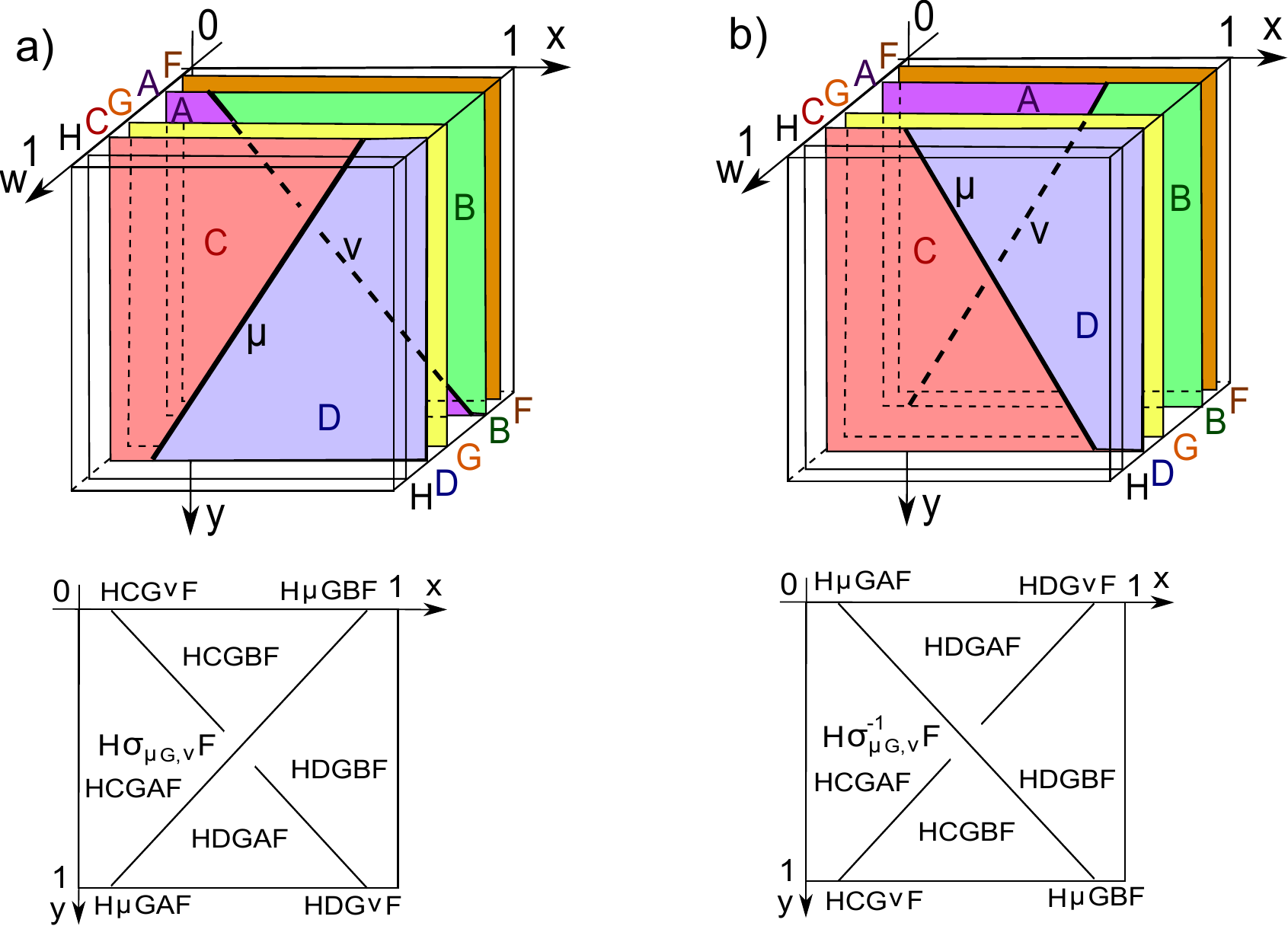}
  \caption{Gray category diagrams and projections for a crossing: {a)} the tensorator $H\sigma_{\mu G,\nu}F$ and {b)} its inverse.  The labelling of regions by objects is omitted. The shortened notation omitting $\Box$ is used.}
  \label{tensorator}
\end{figure}

Let $D$ be a generic Gray category diagram 
with initial region (containing the face $\{0\}\times[0,1]^2$) labelled by an object $\mathcal C$ in $\mac G$ and final region (containing the face $\{1\}\times[0,1]^2$) labelled with an object $\mathcal D$ of $\mac G$.  Then its projection $p_2D$ is  
a 2-category diagram
for the 2-category $\mathcal G(\mathcal C,\mathcal D)$ as follows: The label at a point 
$(x,y)\in[0,1]^2$  of $p_2D$ depends on whether the set  $p_2^{-1}(x,y)$ contains a vertex, an interior  point of a line, interior  points of two different lines, or none of these:
\begin{enumerate} 
\item If  $p_2^{-1}(x,y)$ contains a vertex,
  then the points in $p_2^\inv(x,y)$ define  a sequence  $F_1,F_2,\ldots,F_j,\Phi,G_1,G_2,\ldots G_k$, where $F_n$, $G_m$ are 1-morphisms labelling surfaces and $\Phi$ is 3-morphism labelling the vertex,  in order of increasing $w$ coordinate. The point $(x,y)$ is a vertex of $p_2D$ and is labelled with the 3-morphism
  $$G_k\Box\ldots G_2\Box G_1\Box\Phi\Box F_j\ldots F_2\Box F_1.$$

\item  If $p_2^{-1}(x,y)$ contains an interior point of a  line,  the labelling is analogous, but the 3-morphism  $\Phi$ is replaced by the 2-morphism $\nu$ labelling the line.  The point $(x,y)$ lies  on a line  of $p_2D$ and is labelled with a 2-morphism. 

\item  If $p_2^{-1}(x,y)$ contains no vertex and no interior points of lines, the labelling is as in (1) but with the 3-morphism $\Phi$ removed.  The point $(x,y)$ lies in a region  of $p_2D$ and is labelled with the corresponding 1-morphism.

\item If $p_2^{-1}(x,y)$ contains interior points of two different lines,  then the point $(x,y)$ is a crossing in $p_2D$.  In this case, the sequence associated with $p_2^\inv(x,y)$ is of the form  $$F_1,F_2,\ldots,F_j,\nu,G_1,G_2,\ldots G_k, \mu, H_1,H_2,\ldots H_l,$$ where $F_l, G_m, H_n$ are 1-morphisms in $\mac G$ that label the surfaces  and  $\nu\colon A\to B$,  $\mu\colon C\to D$ are the 2-morphisms labelling the two lines in the preimage of the crossing. 
  In this case, there are two possible diagrams, whose labellings are given in Figure \ref{tensorator} a) and b). The vertex of the 2-category diagram is labelled, respectively,  with $H\Box \sigma_{\mu\Box G,\nu}\Box F$ and with $H\Box \sigma_{\mu\Box G,\nu}^{-1}\Box F$ where $\sigma_{\mu,\nu}$  stands for the tensorator (see equation \eqref{tens}) and we abbreviate $F=F_j\Box\ldots F_2\Box F_1$, $G=G_k\Box\ldots G_2\Box G_1$ and $H=H_l\Box\ldots H_2\Box H_1$.    
\end{enumerate}

The 2-category diagram $p_2 D$  obtained in these four cases  defines the evaluation of a generic Gray category diagram: 
\begin{definition}
  Let $D$ be a generic Gray category diagram whose initial region is labelled with an object $\mathcal C$  and whose final region is labelled with an object  $\mathcal D$ in $\mac G$. 
The {\bf evaluation} of $D$ is  the 3-morphism in  $\mac G(\mac C,\mac D)$  that is the  evaluation of the 2-category diagram $p_2D$. 
\end{definition}

\begin{figure}
  \centering
  \includegraphics[scale=0.65]{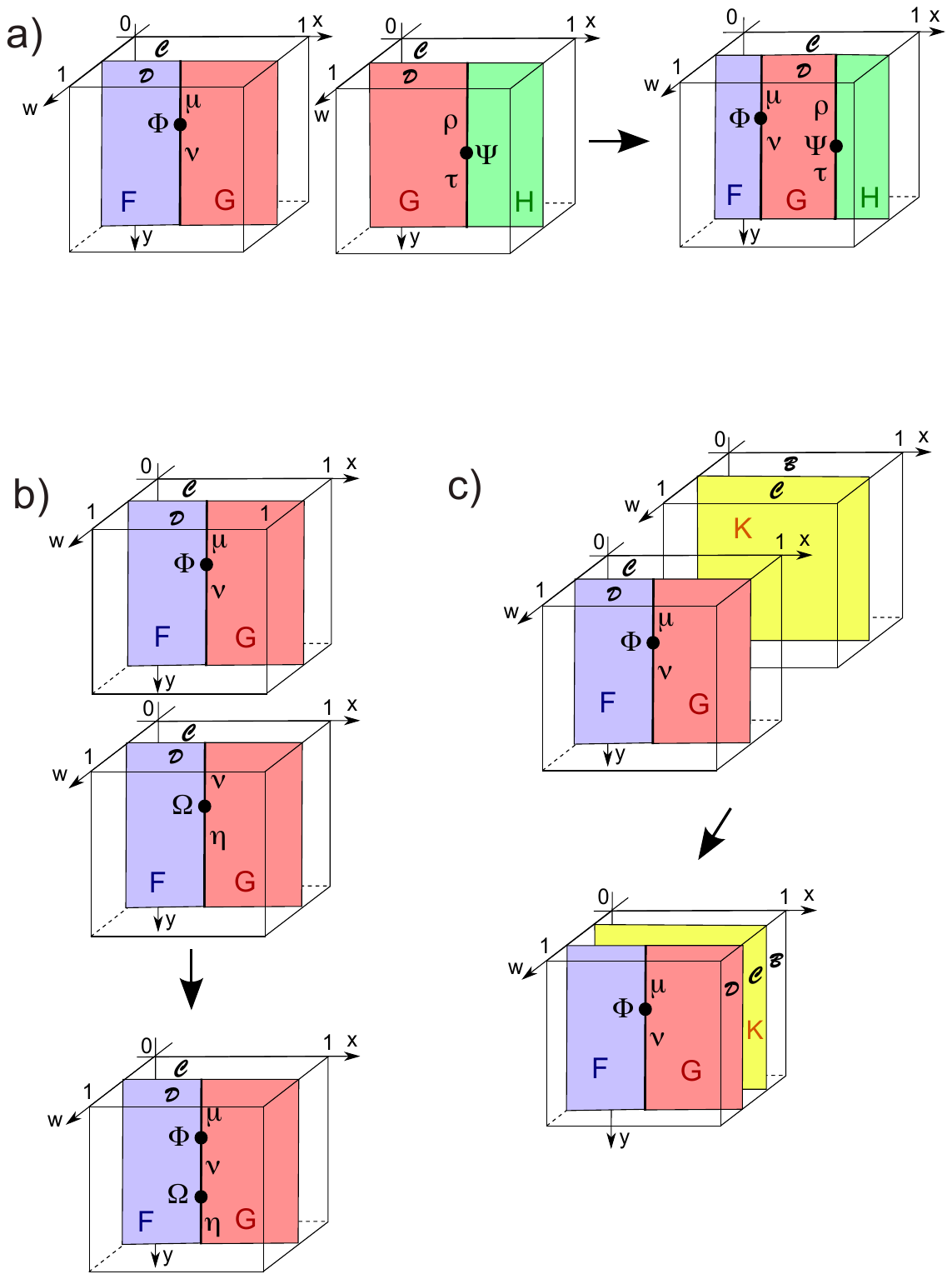}
  \caption{Composition of Gray category diagrams: \newline
    a) Horizontal composition $\circ$,\newline
    b) Vertical composition $\cdot$,\newline
    c) Gray product $\Box$.}
  \label{graycatcomp}
\end{figure}

As in the one- and two-dimensional case, there  is a relation between the composition of Gray category diagrams and the three compositions in a Gray category $\mac G$. 
Gray category diagrams  can be composed in the $w$-, $x$- and $y$-direction as depicted in Figure \ref{graycatcomp}.

The composition in the direction of the $w$-axis is defined if the object in  $\mac G$ labelling the initial region at the face $\{0\}\times[0,1]^2$ of one diagram agrees with the object labelling the final region at the face $\{1\}\times[0,1]^2$ of the other diagram, as shown in Figure \ref{graycatcomp} c).  This composition corresponds to the Gray product of two 3-morphisms in the Gray category $\mac G$. If $D, D'$ are generic progressive Gray category diagrams that can be composed in this way such that their composite diagram $\tilde D$ is again a generic progressive diagram, then the evaluation of $\tilde D$ is the Gray product  of the evaluation of $D$ and $D'$.

The composition  in the direction of the $x$-axis is defined if the labelled progressive two-dimensional diagram at the face $[0,1]\times \{0\}\times[0,1]$ of one of the diagrams matches the labelled diagram at the face $[0,1]\times\{1\}\times[0,1]$ of the other, as shown in Figure \ref{graycatcomp} a). It corresponds to the horizontal composition in $\mac G$. If two generic progressive
Gray category diagrams $D,D'$ are composable in this sense and the resulting diagram $\tilde D$ is again generic and progressive, then the evaluation of $\tilde D$ is the horizontal composite of the evaluation of $D$ and of $D'$.

The composition in the direction of the $y$-axis is defined if the labelled two-dimensional diagram at the face $[0,1]^2\times\{0\}$ of one of the two diagrams matches the  labelled two-dimensional at the face $[0,1]^2\times\{1\}$ of the other, as shown in Figure \ref{graycatcomp} b).  This composition of diagrams corresponds to the vertical composition in $\mac G$.
If $D,D'$ are generic progressive Gray category diagrams which can be composed such that the resulting composite diagram $\tilde D$ is again a progressive generic diagram, then the evaluation of $\tilde D$ is  the vertical composite of the evaluations of $D$ and $D'$.

It seems plausible that by considering Gray category diagrams up to suitable  isotopies, one  should obtain a Gray category  of Gray category diagrams, which generalises the well-known example of the tangle category. In this framework, the evaluation should define a functor from the Gray category  of  diagrams labelled with $\mac G$ to the Gray category $\mac G$.  The above relations between the composition of diagrams and the composition of their evaluations would then correspond to the axioms of a functor of 2-strict tricategories.  However, this aspect is not developed further in the paper.

\begin{figure}
  \centering
  \includegraphics[scale=0.6]{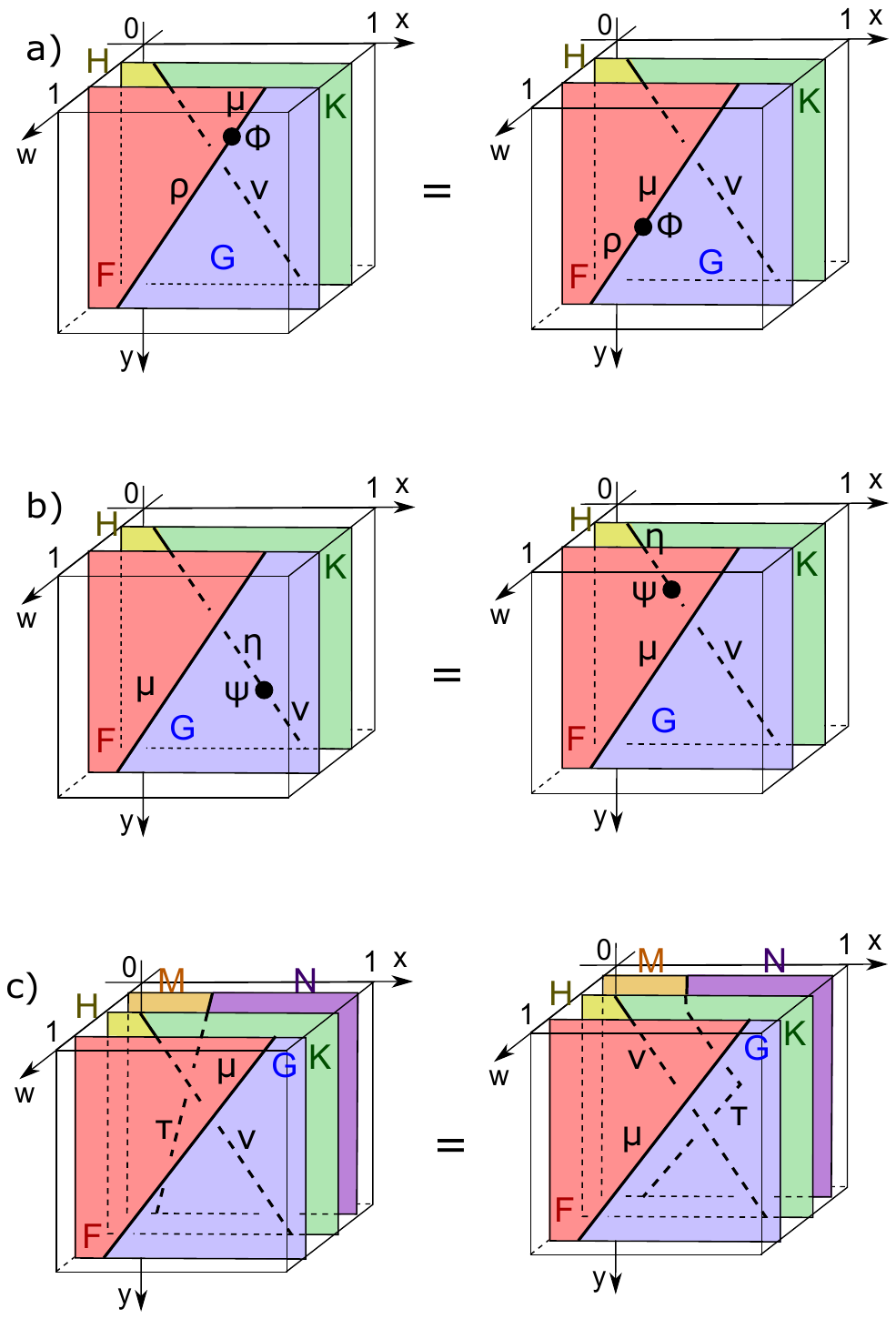}
  \caption{Properties of the tensorator $\sigma_{\mu,\nu}$:  \newline
    a)  naturality in the first argument,\newline
    b), c) naturality in the second argument.
  }
  \label{braiding}
\end{figure}

\medskip
In analogy to the lower-dimensional cases, each generic Gray category diagram  can be viewed as a calculation in a Gray category $\mac G$, which is given by the evaluation of the diagram.  The benefit of such a diagrammatic calculus is that calculations in $\mac G$ can be easily visualised. This requires a statement about the invariance of the evaluation under certain mappings of diagrams. 

Let $D$, $D'$ be Gray category diagrams and $m\colon D\to D'$ be an isomorphism of the underlying progressive three-dimensional diagrams according to Definition \ref{3dmappings}. Suppose $v$ is a vertex of $D$. Then $m(v)$ is a vertex of $D'$ and has an elementary diagram neighbourhood $C'$. It is possible to choose an elementary diagram neighbourhood $C$ of $v$ sufficiently small so that $m(C)\subset C'$. The map $m$ provides a bijection between the set of source lines of $C$ and the set of source lines of $C'$, and also a bijection between the set of target lines of $C$ and the set of target lines of $C'$.

The problem is that there are examples in which $m$ does not respect the orderings of the lines at $v$. Therefore the 3-morphisms at vertices $v$ and $m(v)$ cannot always be equal, and there is no obvious rule to relate them. In the cases where $m$ does respect the orderings at all the vertices, it is possible to set the labels on strata of $D$ and their images in $D'$ to be equal. However there are examples where the evaluations of $D$ and $D'$ are not equal.

A solution to these problems is to restrict to isomorphisms obtained from an isotopy of Gray category diagrams. A condition is added so that the orderings do not change with the isotopy parameter.

\begin{definition}\label{gcdisotopy} An isotopy of Gray category diagrams is an isotopy $\phi\colon [0,1]^3\times[0,1]\to [0,1]^3\times[0,1]$, such that each $\phi_t$ is a mapping of three-dimensional diagrams $D_0$ to $D_t$. This isotopy is required to satisfy a condition for each vertex of $D_0$. The condition on a vertex $v$ is that there is a rectangular neighbourhood $C$ of $v$ that is an elementary Gray category diagram, and for all $t$, the set $\phi_t(C)$ is contained in an elementary diagram neighbourhood of the vertex $\phi_t(v)$. If $s$ is a stratum of $D_0$ then the label of $\phi_t(s)$ in $D_t$ is required to be equal to the label of $s$ in $D_0$.
\end{definition}

This definition only makes sense if for all $t$, the source and targets for $\phi_t(v)$ are equal to those of $v$. This is guaranteed by the following lemma.
\begin{lemma} For each $t$, the map $\phi_t$ preserves the ordering of the source lines at $v$ and also the ordering of the target lines at $v$.  \label{lem:gcdordering}
\end{lemma}
\begin{proof} Recall that $p_2$ is the projection map $(w,x,y)\mapsto (x,y)$. The composite $p_2\circ \phi_t$ maps lines in $C$ to lines in $p_2(\phi_t(C))$. Due to the condition on $C$, these lines in the plane do not cross. Since these lines vary continuously with $t$, the orderings of the lines do not change with $t$. Since $\phi_0$ is the identity map, it follows that $\phi_t$ preserves the orderings of the lines.
\end{proof}

Examples of isotopies that satisfy Definition \ref{gcdisotopy} are the commonly-studied isotopies of diagrams with coupons, as introduced by Reshetikhin and Turaev \cite{RT}, with $C$ corresponding to the coupon. 
The invariance result follows.

\begin{theorem}\label{th_progressive_invariant}Let $D$, $D'$ be generic Gray category diagrams that are related by an isotopy of Gray category diagrams. Then the evaluations of $D$ and $D'$ are equal.
\end{theorem}

\begin{proof} The proof is a generalisation of the proof of Reidemeister's theorem for PL knots and the extension by Yetter to knotted graphs \cite{Y}.
 
 Consider a set of strata $\{C_i |i=1,\ldots n\}$ of the three-dimensional diagram whose projections $p_2(C_i)$ have a common intersection in $[0,1]^2$. The intersection is called transverse if it is a subset of $[0,1]^2$ of dimension $d$, with 
$2-d=\sum_i (2-\dim C_i)$, i.e., the codimensions are additive. 

A generic diagram $D$ is a diagram for which every intersection of strata of $p_2D$ is transverse. Any diagram can be perturbed to a generic diagram. The only possible transverse intersections satisfying the dimension constraint are
 \begin{enumerate}\item Two 1-strata  $(d=0)$
 \item A 0-stratum and a 2-stratum $(d=0)$ 
 \item A 1-stratum and a 2-stratum $(d=1)$ 
 \end{enumerate}
 plus, in each case, any number of additional 2-strata (these have codimension 0).
 
 Under an isotopy, the stratum $C_i$ becomes $C_i\times[0,1]\subset[0,1]^3\times[0,1]$, and the projections  $p_2(C_i)\times[0,1]$ intersect in a subset of dimension $d'$. The isotopy can be perturbed so that these intersections are transverse.
Then $3-d'=\sum_i (3-(1+\dim C_i))=\sum_i (2-\dim C_i)$.  
Therefore the intersections listed in (1)-(3) now have dimension $d'=d+1$ in $[0,1]^2\times[0,1]$, and so they persist in a one-dimensional interval of the `time' isotopy parameter $t$. 
  
 There are new intersections, in addition to (1)-(3), of dimension $d'=0$, namely of the $p_2(C_i)\times [0,1]$ with the $C_i$ given by
 \begin{enumerate}
 \setcounter{enumi}3
 \item A 0-stratum and a 1-stratum
  \item Three 1-strata 
 \end{enumerate}
 plus, as before, any number of 2-strata.
 Since these intersections have  $d'=0$, they can appear only at discrete values of $t$, i.e., are moves between generic  diagrams. The resulting set of moves for Gray category diagrams is depicted in Figure  \ref{braiding}, with a) and b) corresponding to (4), and c) corresponding to (5). The invariance of the evaluation under these moves then follows from the properties of the tensorator in Definition \ref{graydef}. 
  
 Now consider an intersection of strata $\{C_i\}$ with $d'\ge 1$. Suppose these intersect in a time interval $(a,b)$ but not for $t=b$. Then at $t=b$ there is an intersection involving the boundary of the closure of at least one of the $\{C_i\}$, 
 which is a stratum of lower dimension. This intersection in the diagram at $t=b$ is of type (4) or (5), and is not transverse. (A similar conclusion holds for $t=a$.) 
 
Consider a time interval  $[t_1, t_2]$ such that all of intersections have $d'\ge 1$. Since the neighbourhood of each vertex is preserved according to Definition \ref{gcdisotopy}, vertices and crossing points cannot merge. This means that each intersection persists for all values of $t$ in that interval. 
Then there is an isotopy of the projection square $[0,1]^2$ relating the diagrams at $t=t_1$ and $t=t_2$, due to the following argument. Firstly, there is an isotopy of $p_2(X^0)$ in $[0,1]^2$. Then by a result of Hamstrom \cite[Theorem 1]{H}  there is an isotopy in $[0,1]^2$ of the crossing points and the arc segments of $p_2(X^1)$ between vertices and/or crossing points.

 This isotopy  of the 1-complex formed by $p_2(X^0)$ and $p_2(X^1)$
 can be extended to an ambient isotopy of $[0,1]^2$ due to the fact that arc segments in two dimensions are unknotted. Invariance of the evaluation under three-dimensional isotopies that induce isotopies of the two-dimensional diagrams obtained via the projection follows directly from Theorem \ref{generic2invariance}. 
\end{proof}

 Joyal and Street \cite{JS} prove an analogue of this theorem for the special case of a braided monoidal category, using smooth diagrams. Here it is explained briefly how their smooth formalism corresponds to the PL case of this paper. 

The smoothness comes into their proof in only two places. The first is the use of the  definition \cite[Definition 3.1 (iii)]{JS} which requires that the plane projection of the tangent vectors of the curves incident at a vertex differ from each other. 
 In this paper, the analogue is the condition on an elementary diagram that the images of its lines under the plane projection $p_2$ do not intersect (Definition \ref{def:elementary}).

The second place where the smoothness is used is in \cite[Definition 3.2 (iii)]{JS} which is a condition on the derivative of a line in a deformation of diagrams. This ensures that the ordering of the edges arriving at a vertex does not change in the deformation parameter $t$. The PL analogue of this condition is the condition on each vertex in an isotopy in Definition \ref{gcdisotopy}. 

By means of Theorem \ref{th_progressive_invariant}, it is possible to extend the definition of the evaluation to progressive  Gray category diagrams which are not generic.  The idea is the same as in the two-dimensional case (see Theorem \ref{generic2invariance} and Definition \ref{progdiag}), namely to perturb a  non-generic progressive diagram into a generic progressive diagram by means of an isotopy.

\begin{definition}
\label{progdiaggray} 
The {\bf evaluation of a Gray category diagram} is defined by perturbing it, by an isotopy that fixes the boundary, to a generic progressive Gray category diagram and evaluating the resulting diagram. 
\end{definition}
This is independent of the choice of isotopy and hence well-defined by Theorem \ref{th_progressive_invariant}.

\section{Gray categories with duals}
\label{graycatdual}

This section introduces 2-categories and Gray categories with duals.  The main aim is to investigate the structure and the diagrammatic representation of Gray categories with duals.  As the associated diagrams and their evaluation are defined in terms of 2-category diagrams, this requires 
a careful investigation of diagrams for 2-categories with duals. Most of the material on 2-categories and the associated diagrams  is standard \cite{JS2},  
but the discussion of mappings of 2-dimensional diagrams contains the detail that is required in Section \ref{non-progressive}.

\subsection{2-categories with duals}
\label{2catdual}

As the duality operations introduced in this section reverse certain products in the Gray categories, we require two notions of opposites for 2-categories. We start by introducing the relevant notation. For a category $\mac C$ we denote by $\mac C^{op}$ the opposite category and for a morphism $f\in C(A,B)$ by $f^{op}\in (\mac C^{op})(B,A)$ the corresponding morphism with source and target reversed.  We denote by $\mathbin{\widetilde\cdot}$ the composition of morphisms in $\mac C^{op}$, e.g.,
$$f^{op}\mathbin{\widetilde\cdot} g^{op}= (g\cdot f)^{op}.$$ 
Similarly,  we denote for a functor $F\colon \mathcal C\to \mathcal D$  by  $F^{op}\colon \mathcal C^{op}\to\mathcal D^{op}$  the opposite functor with $F^{op}(f^{op})=(F(f))^{op}$ and  for a natural transformation $\nu\colon F\to G$ by
$\nu^{op}\colon G^{op}\to F^{op}$  the opposite natural transformation defined by $\nu^{op}(A)=\nu(A)^{op}$. In particular, if  $\nu$ is a natural isomorphism, then $(\nu^{op})^{-1}$ is a natural isomorphism from $F^{op}$ to $G^{op}$.

\begin{definition}
\label{op2catdef}
  Let $\mac C$ be a 2-category. Then  $\mac C^{op}$ denotes  the corresponding 2-category with both products reversed
  \begin{itemize} \item $(\mathcal C^{op})(A,B)=\mathcal C(B,A)^{op}$ for objects  $A$, $B$ 
  \item $\alpha^{op}\mathbin{\widetilde\circ}\beta^{op}=(\beta\circ\alpha)^{op}$ for composable 1- or 2-morphisms $\alpha$, $\beta$,
  \end{itemize}
  and $\mac C_{op}$  the 2-category with the same vertical but opposite horizontal product
  \begin{itemize}
  \item $(\mathcal C_{op})(A,B)=\mathcal C(B,A)$ for objects  $A$, $B$ 
  \item $\alpha_{op}\mathbin{\widetilde\circ}\beta_{op}=(\beta\circ\alpha)_{op}$ for composable 1- or 2-morphisms $\alpha$, $\beta$.
  \end{itemize}
 
\end{definition}

 Note that in the 2-category literature a more usual notation for $\mac C^{op}$ is $\mac C^{coop}$ and the more usual notation for $\mac C_{op}$ is $\mac C^{op}$. However, these notations would be inconvenient for the extension to Gray categories.  In the sequel we will abuse notation and simply denote a morphism $f^{op}$ by $f$ whenever it is clear from the context to which category $f$ belongs.

The appropriate notion of a 2-category with duals that  will be  used  later in the definition of a Gray category with duals is that of a planar 2-category.  
A planar 2-category is a direct generalisation of a  pivotal  category, which in turn can be regarded  as a planar 2-category with one object.

The definition uses a 2-category notion of contravariant functor. By definition, a contravariant functor $G$ between 2-categories $\mathcal B$ and $\mathcal C$ is a strict 2-functor $\overline G\colon\mathcal B \rightarrow \mathcal C^{op}$.
The product $FG$ of a contravariant functor $G$ followed by functor $F\colon\mathcal C \rightarrow \mathcal D$ is defined as the contravariant functor associated to ${F}^{op}\,\overline G$, and $GF$ is the contravariant functor associated to $\overline G F$. Similarly, if $F$ is 
instead a contravariant functor, then the product $FG=\overline{F}^{op}\,\overline G$ is a (covariant) functor.

In the following, it will be convenient to use the same symbol for $\overline G$ and $G$. As can be seen from the products of functors just described, this amounts to using the covariant functors and omitting the superscript $^{op}$ on functors where it occurs. These superscripts can be unambiguously re-inserted by examining the source category of each functor.

\begin{definition}
\label{plandef}
  A {\bf planar 2-category} is a 2-category $\mac C$  together with a contravariant functor $*\colon\mathcal C \rightarrow \mathcal C$, associated to a strict 2-functor  $\mathcal C \rightarrow \mathcal C^{op}$, 
that is the identity on objects. Additionally, it is required that there is  a collection of 2-morphisms $\epsilon_a\colon 1_{A'}\rightarrow a\circ a^*$ for all 1-morphisms $a\colon A\to A'$ of $\mathcal C$ such that:
  \begin{enumerate}
  \item $**=1_{\mathcal C}$ is  
    the identity functor 
  \item 
    for all 1-morphisms $a,b,c$  and  2-morphisms $\alpha\colon a\rightarrow b$ for which these expressions are defined:
    \begin{align}\label{pivcond}
      &(\alpha\circ 1_{a^*})\cdot \epsilon_a=(1_b\circ \alpha^*)\cdot \epsilon_b\qquad (1_{a}\circ \epsilon_{a^*}^*)\cdot(\epsilon_{a}\circ 1_{a})=1_{a}\\
      &(1_a\circ \epsilon_c\circ 1_{a^*})\cdot \epsilon_a=\epsilon_{a\circ c}.\nonumber
    \end{align}
  \end{enumerate}
A {\bf pivotal category} is a planar 2-category with one object, regarded as a monoidal category with additional structure.
\end{definition}

Note that the pivotal category in Definition \ref{plandef} is often called {\em strict pivotal category} in the literature. 
As we do not make use of weak pivotal categories,  the  word 'pivotal' always stands for  {\em strict pivotal} throughout the article.  

Note that the 2-functor $*$ and the collection of morphisms $\epsilon_a$ in a planar 2-category are not independent. For a given 1-morphism $a$, the 2-morphisms $\epsilon_a$ and $\epsilon_{a^*}^*$ determine an adjunction between $a$ and $a^*$ in the 2-category $\mac C$ \cite{KS}. The following lemma shows that this choice of data for every 1-morphism determines the action of the  functor $*$ on the 2-morphisms uniquely.

\begin{lemma}\cite{kelly_laplaza, BWSC} \label{dualmodi}
  For any 2-morphism $\alpha\colon a\rightarrow b$ in a planar 2-category, the dual $\alpha^*\colon b^*\rightarrow a^*$ is given by 
  \begin{align}\label{dual2morph}
    \alpha^*=&(\epsilon_{b^*}^*\circ 1_{a^*})\cdot(1_{b^*}\circ \alpha\circ 1_{a^*})\cdot(1_{b^*}\circ \epsilon_a)\\
    =&(1_{a^*}\circ \epsilon_{b}^*)\cdot(1_{a^*}\circ \alpha\circ 1_{b^*})\cdot(\epsilon_{a^*}\circ 1_{b^*}),\nonumber
  \end{align}
  depicted in Figure \ref{plancat}  f), 
  and the 2-morphism $\alpha$ satisfies the pivotal condition in Figure \ref{plancat} g):
  \begin{align}\label{pivotalc}
    \alpha=\; &
    (1_b\circ \epsilon_{a^*}^*)\cdot (1_b\circ 1_{a^*}\circ \epsilon_b^*\circ 1 _a)\cdot (1_b\circ 1_{a^*}\circ \alpha\circ 1 _{b^*}\circ 1 _{a})\\
    \cdot &(1_b\circ \epsilon_{a^*}\circ 1 _{b^*}\circ 1 _a)\cdot (\epsilon_b\circ 1 _a).\nonumber
  \end{align}
\end{lemma}

\begin{proof} 
  The proof is a direct generalisation of the corresponding proof for pivotal categories, see \cite{kelly_laplaza}, \cite[Lemma 2.2]{BWSC}.
  The identities in \eqref{dual2morph} follow from the first and second identity in \eqref{pivcond} together with the exchange law. The pivotal condition \eqref{pivotalc} is then obtained by applying \eqref{dual2morph} twice and using the identity  $**=1_{\mathcal C}$.
\end{proof}

\subsection{Diagrammatic representation of  the 2-category duals}\label{2catdualssection}

The $*$-duals in a planar 2-category are the extra data required to define 2-category diagrams that are not progressive. In this setting, the condition that a 2-category diagram $D$ is progressive can be relaxed to the condition 
that it is `piecewise progressive',  i.e., that there is a subdivision $m\colon D\to P$ such that
$P$ is progressive.  These exist if the diagram is generic in the following sense.

\begin{definition}
  A two-dimensional diagram is called {\bf generic} if the only singularities of the projection $p_1$ on lines are maxima and minima, and all vertices, maxima and minima have different  $y$-coordinates.
\end{definition}

By promoting the maxima and minima of a generic 2-dimensional diagram $D$ to vertices, one obtains a subdivision $m\colon D\to P$ called the  {\bf minimal progressive subdivision} in the following. $P$ is a `generic progressive diagram' in the sense of Definition \ref{genprog2d}. 

In the rest of this section, it is assumed that two-dimensional diagrams are generic where appropriate. Note, however, that an isotopy between generic diagrams  may fail to be generic at intermediate values of the isotopy parameter.

Subdividing generic two-dimensional diagrams to obtain  generic progressive diagrams
introduces a complication with the labelling. A line in $D$ that is labelled with an object $x$ and
zig-zags upwards and downwards with respect to the $y$-coordinate
corresponds to a collection of progressive line segments in $P$, whose labels vary
between $x$ and $x^*$, depending on their orientation. Keeping track of these labels and of the analogous labelling problems caused by the rotation of vertices motivated the introduction of a new structure into a diagram, namely a framing. This was introduced for knots by Kauffman \cite{K-RI} and for monoidal category diagrams by Reshetikhin and Turaev \cite{RT}.

Nevertheless,  it is possible to reformulate this theory in terms of unframed diagrams. 
 Our experience is that this is much simpler in the case of three-dimensional diagrams, and we will work with unframed diagrams throughout the paper.  However the issues that motivated the introduction of framing in previous works then appear in the action of mappings on diagrams. 

\begin{figure}
  \centering
  \includegraphics[scale=0.3]{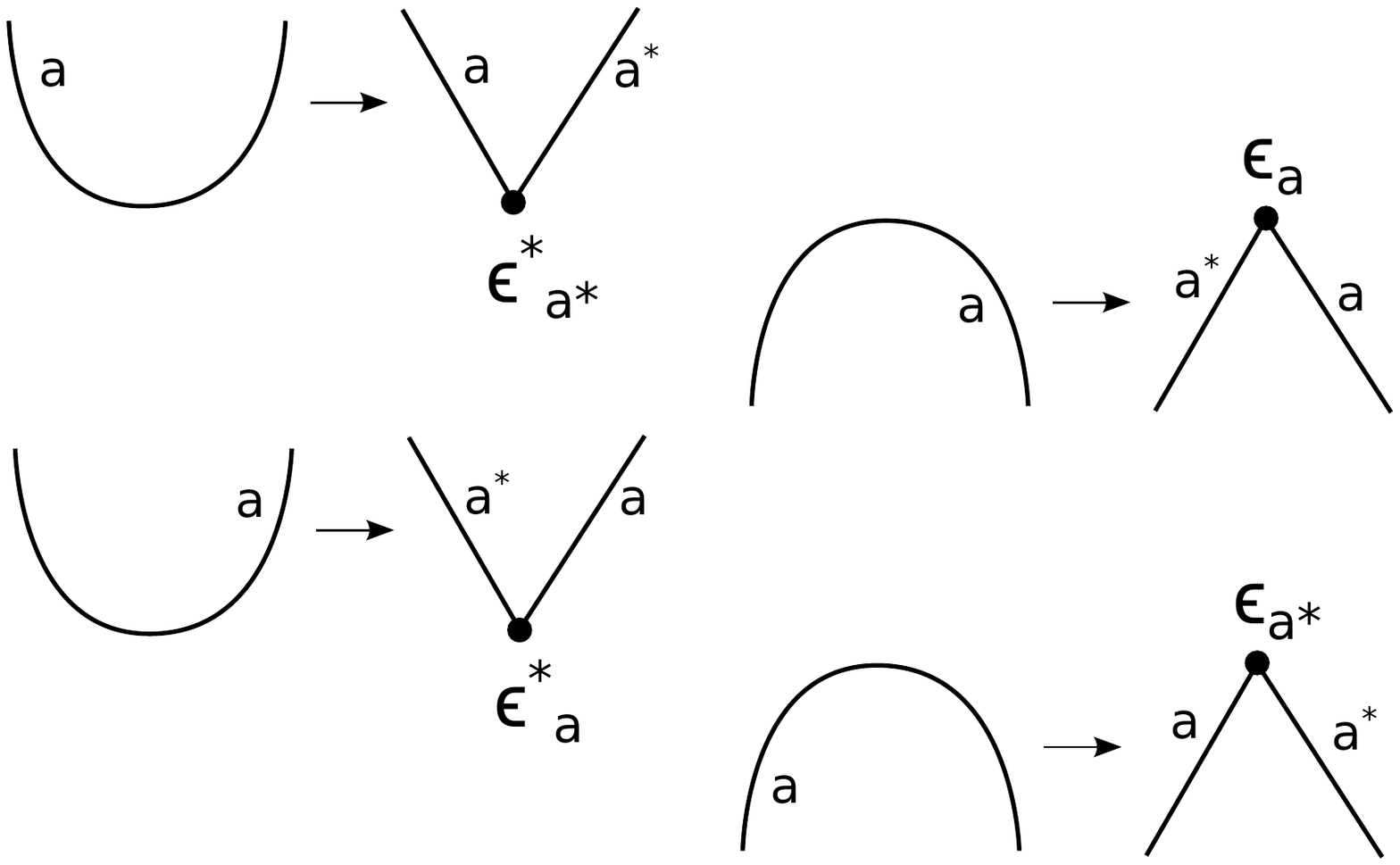}
  \caption{Labelling minima and maxima of a planar 2-category diagram.}
  \label{minmaxlabels}
\end{figure}

\begin{definition} \label{plan_sub}
Let $\mathcal C$ be a planar 2-category. 
  A {\bf planar 2-category diagram} for $\mac C$ is a generic two-dimensional diagram $D$ together with a
  labelling of the image of  its minimal progressive subdivision $D\to P$ with elements of $\mathcal C$ such that $P$  is a 2-category diagram. It is required that the additional vertices are labelled with the canonical 2-morphisms $\epsilon_a$ or  $\epsilon_a^*$  as shown in Figure \ref{minmaxlabels}. The {\bf evaluation} of $D$ is defined as the evaluation of $P$.
\end{definition}

\begin{figure}
  \centering
  \includegraphics[scale=0.6]{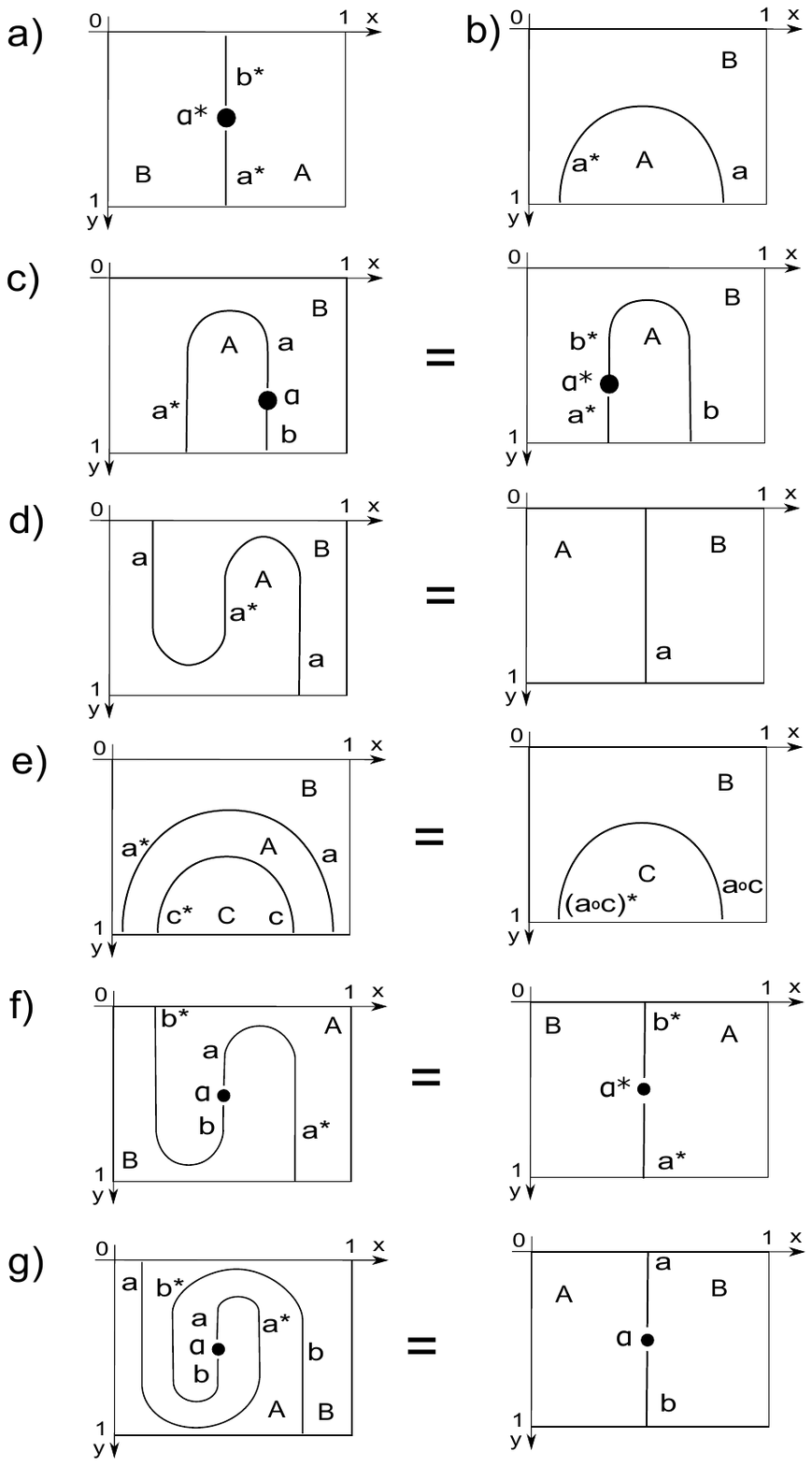}
  \caption{Planar 2-category diagrams. \newline
    a) 2-morphism $\alpha^*\colon b^*\rightarrow a^*$ for a 2-morphism $\alpha\colon a\rightarrow b$.\newline
    b) The 2-morphism $\epsilon_a\colon 1_B\rightarrow a\circ a^*$.\newline
    c) The identity $(\alpha\circ 1)\cdot \epsilon_a=(1\circ \alpha^*)\cdot \epsilon_b$ from \eqref{pivcond}.\newline
    d) The identity $(1_a\circ \epsilon_{a^*}^*)\cdot (\epsilon_a\circ 1_a)=1_a$ from \eqref{pivcond}.\newline
    e) The identity $(1_a\circ \epsilon_c\circ 1_{a^*})\cdot \epsilon_a=\epsilon_{a\circ c}$ from \eqref{pivcond}.\newline
    f) The identity $\alpha^*=(1_{a^*}\circ \epsilon_b^*)\cdot(1_{a^*}\circ \alpha\circ 1_{b^*})\cdot( \epsilon_{a^*}\circ 1 _{b^*})$.\newline
     g)  The condition \eqref{pivotalc} from Lemma \ref{dualmodi}.
  }
  \label{plancat}
\end{figure}

Examples of planar 2-category diagrams are given in Figure \ref{plancat}. \smallbreak

The aim now is to prove that the evaluation of planar 2-category diagrams  is invariant under a much more general class of mappings of diagrams  than those in Theorem \ref{generic2invariance}, namely under any homomorphism of diagrams. 
Unlike in the progressive case, 
this requires that the homomorphisms  act non-trivially on the labels of the minimal progressive subdivisions.
Hence, the first step is to define the relation between the labels of the source and target diagrams of a homomorphism of diagrams.

This can be reduced to specifying it for  isomorphisms of progressive diagrams by the following considerations.
If $f\colon D\to D'$ is a homomorphism of diagrams, then there are progressive diagrams $P$, $P'$ and subdivisions $s\colon D\to P$, $s'\colon D'\to P'$ such that $f$  induces  an {\em isomorphism} of diagrams $ f\colon P\to P'$. These subdivisions are obtained
 by first subdividing $D$ such that the homomorphism becomes an isomorphism, then adding vertices at the maxima and minima of both $D$ and $D'$ to make the diagrams progressive and finally then adding further vertices such  that $v$ is a vertex of the resulting diagram $P$ if and only if $f(v)$ is a vertex of the resulting diagram  $P'$. The additional vertices that are maxima or minima are labelled according to Definition \ref{plan_sub}.
 The additional vertices that are not maxima or minima are labelled with identity 2-morphisms, and the additional lines in $P$ are labelled with identity 1-morphisms. 

It is therefore sufficient to consider the case where $f\colon P\to P'$ is an isomorphism of  
diagrams. One can compare the configuration of lines and vertices with respect to the vertical projection $p_1\colon (x,y)\to y$. This projection induces an orientation on each line $l$ of $P$ and the corresponding line $f(l)$ of $P'$. The isomorphism $f$ either preserves or reverses the orientation. If the orientation is preserved, then the 1-morphism labels on $l$ and $f(l)$ are taken to be equal, and if the orientation is reversed, then the labels are related by $*$.

As explained in Section \ref{2catdiagsection}, the lines incident at a vertex $v$ of $P$ fall into two ordered sets, the  {\em  input lines}, and the {\em output lines}. The isomorphism $f$ exchanges lines between the two sets, preserving the cyclic order of the lines around the vertex. The combinatorial information in this exchange of lines is used to determine the relation of the 2-morphisms at the vertex.

\begin{figure}
  \centering
  \includegraphics[scale=0.5]{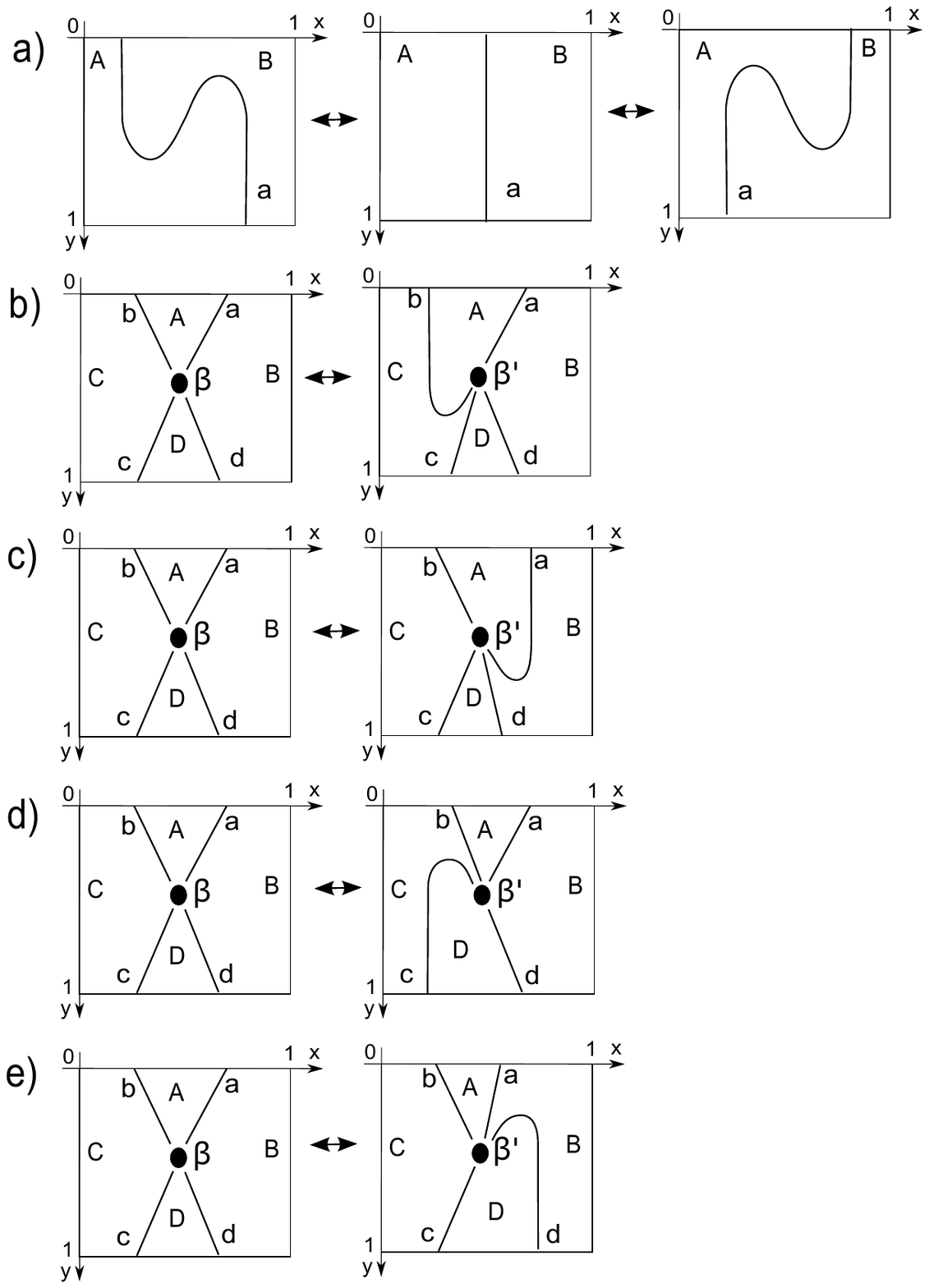}
  \caption{Diagrams for the moves from \cite[Proposition 1.8]{Y}.
  }
  \label{yetter}
\end{figure}

Some mappings of planar diagrams are shown in Figure \ref{yetter} b)-e). The diagrams on the left-hand side in Figure \ref{yetter} b)-e) show  vertices with  two input lines and two output lines, but there are obvious analogues of these moves for any number of lines. For each choice of the vertex 2-morphism $\beta'$ in the diagrams on the right-hand sides, the 2-morphism $\beta$ in the corresponding diagram on the left-hand side is then determined uniquely by the condition that the evaluations of the two diagrams are equal.  
This defines the action of an isomorphism on the vertex 2-morphisms.
Figure \ref{yetter} b) determines a map
\begin{equation}\phi_L\colon \mathcal C(a,d\circ c\circ b^*)\to\mathcal C(a\circ b,d\circ c),
\end{equation}
which is obtained by inserting the morphism $\beta'$ in the right-hand diagram and  evaluating the diagram. Similarly, 
Figure \ref{yetter} c) determines a map
\begin{equation}\mathcal\phi_R\colon \mathcal C( b,a^*\circ d\circ c)  \to \mathcal C(a\circ b,d\circ c).
\end{equation}
The other two moves, shown in  Figure \ref{yetter} d) and e), give the inverse maps. Again, there are obvious analogues of these maps with other numbers of input and output lines. For example, one can replace the line labelled by $a$ in Figure \ref{yetter} b) by any number of labelled lines. 

Consider now a vertex $v$ with $p$ input lines and $q$ output lines. Then the corresponding maps $\phi_L$ and $\phi_R$ are defined if 
$q>0$
and  their inverses $\phi_L^{-1}$ and $\phi_R^{-1}$ are defined if $p>0$. Any two of these four maps commute when their products in both orders are defined.
 Iterating these maps leads to the following definition of canonical maps on vertex 2-morphisms.
\begin{definition} Define the {\bf rotation} $r$ by
$$r=\begin{cases}\phi^{-1}_R\phi_L \quad \text{ if }&q>0\\ \phi_L\phi^{-1}_R &p>0\\1&p=q=0.\end{cases} $$
The {\bf canonical maps between sets of 2-morphisms} are defined for $n\in\Z$ and $-q\le k\le p$ by
$$c(k,n)=\phi_L^k\, r^n.$$
\end{definition}

\begin{lemma}\label{lem:vertexmove} Any iterate of $\phi_L$, $\phi_R$ and their inverses acting on a vertex with $p$ incoming and $q$ outgoing lines can be transformed to exactly one of the canonical forms $c(k,n)$ for $0\le n< p+q$ using the axioms of a planar 2-category.
\end{lemma}

\begin{proof}Consider a word of length $l$ in the generators $\phi_L$, $\phi_R$, $\phi_L^{-1}$, $\phi_R^{-1}$ and proceed by induction on $l$. For $l=0$ one has $1=c(0,0)$. The induction step is given by the relations $\phi_L \,c(k,n)=c(k+1,n)$, $\phi_L^{-1}\,c(k,n)=c(k-1,n)$, $\phi_R\, c(k,n)=c(k+1,n-1)$, $\phi_R^{-1}\,c(k,n)=c(k-1,n+1)$. The pivotal condition that is shown in Figure \ref{plancat} g)  implies that $c(k,n)=c(k,n+p+q)$.
Finally, the uniqueness is due to the fact that $k$ and $n$ can be read off from the powers of $\phi_L$ and $\phi_R$ in any word and the relations in a planar 2-category do not change $k$ and only can only change $n$ by the addition of a multiple of $p+q$.
\end{proof}

These definitions and results are summarised in the following formal definition.

\begin{definition} \label{rotatedvertex} A {\bf homomorphism of planar 2-category diagrams  $f\colon D\to D'$} is an orientation-preserving 
homomorphism of two-dimensional diagrams such that the associated isomorphism on subdivisions $f\colon P\to P'$ preserves the 2-category labels in the following sense: 
\begin{itemize} \item If $R\subset P$ is a region, then the objects labelling $R$ and $f(R)$ are equal.
\item If $l\subset P$ is a line, then the 1-morphisms labelling $l$ and $f(l)$ are equal if $f$ preserves the orientation and related by $*$ if the orientations are opposite.
\item If $v$ is a vertex of $P$ labelled with $\beta$, and the vertex $f(v)$ is labelled with $\beta'$, then $\beta=c(k,n)\beta'$, where $k,n$ are determined by the permutation of incoming and outgoing edges.
  \end{itemize}
\end{definition}

The main invariance result can now be stated.

\begin{theorem}\label{equality} Let $f\colon D\to D'$ be a homomorphism of planar 2-category diagrams.  Then the evaluations of the two diagrams are equal. 
\end{theorem} 

\begin{proof} The diagrams $D$ and $D'$ are replaced by progressive diagrams $P$ and $P'$ as described above.
The evaluations of $D$ and $P$ are equal, since they differ only by the insertion of identity morphisms. The same holds for the evaluations of $D'$ and $P'$. 

Thus it remains to consider the isomorphism $f\colon P\to P'$.  According to the result of Yetter \cite[Proposition 1.8]{Y}, there is a list of Reidemeister-type moves for a mapping of planar diagrams. These are the moves reproduced in Figure \ref{yetter}, plus the isotopies of progressive diagrams from Theorem \ref{generic2invariance}. As any mapping $f$ is a product of moves from the list, it is sufficient to check that the evaluation is invariant under these moves.  The invariance of the evaluation under the snake moves in Figure \ref{yetter} a) follows from the second condition in  \eqref{pivcond}.  The invariance of the evaluation under the vertex moves in Figure \ref{yetter} b)-e) follows from Lemma \ref{lem:vertexmove} and Definition \ref{rotatedvertex}, and the invariance of the evaluation and  under the isotopies of progressive diagrams follows from Theorem \ref{generic2invariance}.
\end{proof}

\subsection{Gray categories with duals}

The following definition is derived from the axioms that were first given by Baez and Langford \cite{BL} and Mackaay \cite{mackaay}. The main difference is that here only two independent duals are considered, whereas the previous authors defined three. The axioms are adapted from \cite{BL}. By referring to a planar 2-category these axioms can be cast into a more concise form. The definition has a conceptual justification
in Section \ref{grayfuncs} (Theorems \ref{graycatduals} and \ref{graydualnat}), where it is shown  
  that the data  for the duals determine functors from  that exhibit the symmetries of the cube. 
\begin{definition}
\label{dualgray}

  A {\bf Gray category with a duality structure} is a Gray category $\mathcal G$ with  the following additional structure:

  \begin{enumerate}
  \item \label{itm:GCstarduals}  
  $\mathcal{G}$ is a {\bf Gray category with *-duals}: 
for all objects $\mathcal C,\mathcal D$ of $\mathcal G$, the 2-category $\mathcal G(\mathcal C, \mathcal D)$ is  planar, and its dual $*$ is compatible with the Gray product:
    $$ (K\Box \mu\Box H)^*=K\Box \mu^*\Box H, \quad \text{and} \quad K\Box \epsilon_\mu\Box H=\epsilon_{K\Box \mu\Box H},$$ 
    for all 1-morphisms $H,K$ and  2-morphisms $\mu$ for which these expressions are defined.

  \item \label{itm:GCstarduals2}   For every 1-morphism $F\colon\mathcal C\rightarrow\mathcal D$, there is  a dual 1-morphism $F^\#\colon \mathcal D\rightarrow\mathcal C$, a 2-morphism $\eta_F\colon 1_{\mathcal D}\Rightarrow F\Box F^\#$, called {\bf fold},  and an invertible 3-morphism $T_F\colon(\eta_F^*\Box F)\circ(F\Box\eta_{F^\#})\Rrightarrow 1_{F}$, called the {\bf triangulator}, such that the following conditions are satisfied:
    \begin{enumerate}
    \item  $F^{\#\#}=F$ for all 1-morphisms $F\colon\mathcal C\rightarrow\mathcal D$,
      \smallskip
    \item $1_\mathcal C^\#=1_{\mathcal C}$, $\eta_{1_{\mathcal C}}=1_{1_\mathcal C}$,  $T_{1_\mathcal C}=1_{1_{1_\mathcal C}}$ for all objects $\mathcal C$,
      \smallskip
    \item  
      $(F\Box G)^\#=G^\#\Box F^\#$, $\eta_{F\Box G}=(F\Box\eta_G\Box F^\#)\circ \eta_F$,\\
      $T_{F\Box G}=(T_F\Box G\circ F\Box T_G)\cdot (1_{\eta^*_F\Box F\Box G}\circ \sigma_{F\Box \eta_G^*,\eta_{F^\#}\Box G}\circ 1_{F\Box G\Box \eta_{G^\#}})$
      for all composable 1-morphisms $F\colon\mathcal C\rightarrow\mathcal D$, $G\colon\mathcal B\rightarrow \mathcal C$, 
      \smallskip
    \item\label{eq:GCDaxiom} 
      $(1_{\eta_F^*} \circ T_F\Box F^\#)\cdot(\sigma_{\eta_F^*,\eta_F^*}\circ 1_{F\Box \eta_{F^\#}\Box F^\#})\cdot(1_{\eta_F^*}\circ F\Box T_{F^\#}^*)=1_{\eta_F^*}
      $.
      \bigskip
    \end{enumerate}
  \end{enumerate}
\end{definition}

Note that the duality structure for a Gray category in Definition \ref{dualgray} is more than just a higher version of duals. The requirement that the 2-categories $\mathcal G(\mathcal C,\mathcal D)$  are planar generalises a strict version of pivotality  for a rigid monoidal category, and condition (2) can be viewed as a strict higher analogue of pivotality. Additionally, we do not just impose that  fold 2-morphisms and cusp 3-morphisms as in Definition \ref{dualgray} exist, but a specification of such 2-morphisms and 3-morphisms is part of the structure. Nevertheless,  we use the shorter term
{\bf Gray category with duals} for  {\bf Gray category with duality structure}   in the following.

The relation with the notation of \cite{BL} is as follows. Firstly, the work \cite{BL} considers monoidal 2-categories, so what is termed here an $n$-morphism is the present paper is called an $(n-1)$-morphism there.
The duality on 2-morphisms (which they call 1-morphisms) is denoted $*$ in both works, and is extended here to 3-morphisms by the operation which is called an adjoint in \cite{BL}. The duality on 3-morphisms  in \cite{BL} (which they call 2-morphisms) does not appear in this paper, and hence their constraints that specify that the tensorator and triangulator are unitary are relaxed here to the conditions that these morphisms are invertible.

The duality on 1-morphisms is denoted $\#$ here and corresponds to the $*$-dual  on 1-morphisms  in \cite{BL} (which they call objects). There are two axioms in Definition \ref{dualgray} that have no analogue in \cite{BL}, namely the conditions
$1_\mathcal C^\#=1_{\mathcal C}$ and $\eta_{1_{\mathcal C}}=1_{1_\mathcal C}$. 
\medskip

Just as the $*$-duality can be regarded as a collection of adjunctions in a 2-category, the $\#$-duality can be regarded as a collection of biadjunctions in a Gray category.
For the case where the Gray category is 2Cat, these are particular types of the weak quasi-adjunctions defined in \cite[\S I.7]{Gray}, and, more specifically, in \cite{MacDonaldStone}. In the latter, the data determined by the folds and triangulators for the pair of 1-morphisms $F$ and $\# F$ are called there a `unit-counit diagram' and equation (\ref{eq:GCDaxiom}), 
sometimes also called `swallowtail equation',
is called there a `3-cell unit-counit equation'.

  \begin{remark}
    \label{remark:compareTricatDuals}
    We comment on weaker versions of higher categories with duals and their relation with the notions in this paper.
    \begin{enumerate}
    \item A 2-category $\mathcal C$ is said to {\bf have right duals}, if for every 1-morphism $a\colon A \rightarrow B$ there exists a morphism $a^{*} \colon B \rightarrow A$ together with evaluation and coevaluation morphisms
      $a \circ a^{*} \rightarrow 1_{B}$ and $1_{A} \rightarrow a^{*} \circ a$ that satisfy the  snake identities
      (see e.g.  \cite[Def I,6.1]{Gray}). Analogously one considers the case where $\mathcal C$ has left duals and $\mathcal C$ is said to have duals, if it has left and right duals. In the case where $\mathcal C$ has right duals,  $*$ extends
      naturally to a (weak) 2-functor $\mathcal C \rightarrow \mathcal C^{op}$, and by a
      {\bf weak pivotal structure} we mean a 2-isomorphism $** \rightarrow 1_{\mathcal C}$. 
  \item For a tricategory $\mathcal T$ the homotopy 2-category $h\mathcal T$ is defined by considering the objects and 1-morphisms of $\mathcal T$ together with the isomorphism classes of 2-morphisms of $\mathcal T$.
    We say that $\mathcal T$ has {\bf weak (left/right) $\#$-duals} if $h\mathcal T$ has (left/right) duals.
    Note  that it is only required that triangulators (isomorphisms replacing the snake identities) exist without 
    demanding those to satisfy identities analogous to equation (\ref{eq:GCDaxiom})  in Definition \ref{dualgray}. Thus in this formulation, the duals are just a property and no structure is specified.
    We conjecture that every tricategory which has $\#$-duals is equivalent to a tricategory with specified duality structures that satisfies the additional identities of Definition \ref{dualgray}.
   A proof, which would require adopting the methods of the proof of  \cite[Thm 7.2.1]{Schau}, is beyond the scope of this paper.  
  \item For a  tricategory $\mathcal T$  one can now combine the two duals and arrives at the following weak notions: 
    \begin{itemize}
    \item The tricategory $\mathcal T$ has $*$- and $\#$-duals if all 1- and 2-morphisms have corresponding duals. Here, no compatibilities are specified. Note that if one requires only $*$-duals and right $\#$-duals,
      a right $\#$-dual $F^{\#}$ of a 1-morphism $F$ is automatically also a left $\#$-dual of $F$: The $*$-dual of an evaluation 2-morphism for the right $\#$-dual is an coevaluation 2-morphism for $F^{\#}$, exhibiting it  as a left $\#$-dual.   
    \item One can consider a tricategory $\mathcal T$ that has  $*$- and $\#$-duals and a pivotal structure on all bicategories $\mathcal T(A,B)$
      for objects $A,B$, such that the compositions $a \Box -$ and $- \Box a$ for all 1-morphisms $a$ are pivotal functors, see \cite[Def 5.2.2]{Schau}. In \cite[Thm 7.2.1]{Schau}, it is shown that every such tricategory is
      equivalent to a Gray category with a duality structure. 
    \end{itemize}
  \end{enumerate}

 The main focus of this article is on three-dimensional diagrams and this requires a Gray category with duality structure as defined here.  As will be demonstrated in the following sections, the structures in a Gray category with duals reflect the singularities of the projections of 3d diagrams and are needed for a  consistent labelling of the corresponding diagrams. 

\end{remark}

\medskip
Gray categories with duals can be viewed as a generalisation of braided pivotal categories. In the literature, braided pivotal categories are also called balanced autonomous categories.  As a direct consequence of the axioms in Definition \ref{dualgray}, one obtains 
\begin{lemma}\label{pivten}
  If $\mac G$ is a Gray category with duals, then for every object $\mac C$ the category $\mac G(1_{\mac C}, 1_{\mac C})$ is a braided 
  pivotal  category. 
  Conversely, a braided   
  pivotal  category is a Gray category with duals with a single object and a single 1-morphism.
\end{lemma}
\begin{proof}
For every Gray category $\mathcal G$,   by Lemma   \ref{braidten} $\mac G(1_{\mac C}, 1_{\mac C})$ is a braided category. If $\mathcal G$ is a Gray category with duals, then $\mac G(1_{\mac C}, 1_{\mac C})$ is planar by Definition \ref{dualgray} (1) and hence a  pivotal category.  Conversely,   by Lemma  \ref{braidten} every braided  pivotal category is  a Gray category with $*$-duals with  a single object and a single morphism.  In this case, the $\#$-duals are determined uniquely 
   by the identities in
 Definition \ref{dualgray} (2)(b), and the conditions in Definition \ref{dualgray} (1), (2)(a), (2)(c) and (2)(d) are satisfied trivially.
\end{proof}

\medskip
As indicated by this lemma, some identities which are familiar from braided monoidal categories have  a direct analogue in Gray categories with duals. These similarities are also  apparent  in the diagrammatic calculus for Gray categories with duals  introduced in the next subsection. 
A specific example is the following lemma with the associated diagrams in Figure \ref{tensorator_identities}.
\begin{lemma} \label{sigma_ids}
  Let $\mac G$ be a Gray category with duals. Then for all 2-morphisms $\mu,\mu'\colon F\Rightarrow G$, $\nu\colon H\Rightarrow K$ and all 3-morphisms $\Phi\colon \mu\Rrightarrow\mu'$ for which these expressions are defined, one has
  \begin{align*}
    &(K\Box \Phi\Box H)^*=K\Box \Phi^*\Box H,\quad  \sigma^*_{\mu,\nu}=\sigma_{\mu^*,\nu^*},\\
    &(1_{(\mu\Box K)\circ (F\Box\nu)}\circ\epsilon^*_{\mu^*\Box H})\cdot (1_{\mu\Box K}\circ \sigma_{\mu^*,\nu}\circ 1_{\mu\Box H})\cdot(\epsilon_{\mu\Box K}\circ 1_{(G\Box\nu)\circ (\mu\Box H)})=\sigma_{\mu,\nu}^\inv,\\
    &(1_{(G\Box\nu)\circ(\mu\Box H)}\circ \epsilon^*_{F\Box\nu^*})\cdot (1_{G\Box\nu}\circ \sigma^\inv_{\mu,\nu^*}\circ 1_{F\Box\nu})\cdot(\epsilon_{G\Box\nu}\circ 1_{(\mu\Box K)\circ(F\Box\nu)})=\sigma_{\mu,\nu}.
  \end{align*}
\end{lemma}

\begin{proof}
  The first identity follows directly from the definition of the dual 3-morphisms in terms of the 3-morphism $\epsilon_\mu$  in equation \eqref{dual2morph} and from condition (1) in Definition \ref{dualgray}. The second identity follows from the
  third and the fourth. These two identities are  direct consequences of the properties of the tensorator together with condition (1) in Definition \ref{dualgray}.
\end{proof}

\subsection{Diagrammatic representation of the Gray category duals}
\label{dualdiags}

\begin{figure}
  \centering
  \includegraphics[scale=0.45]{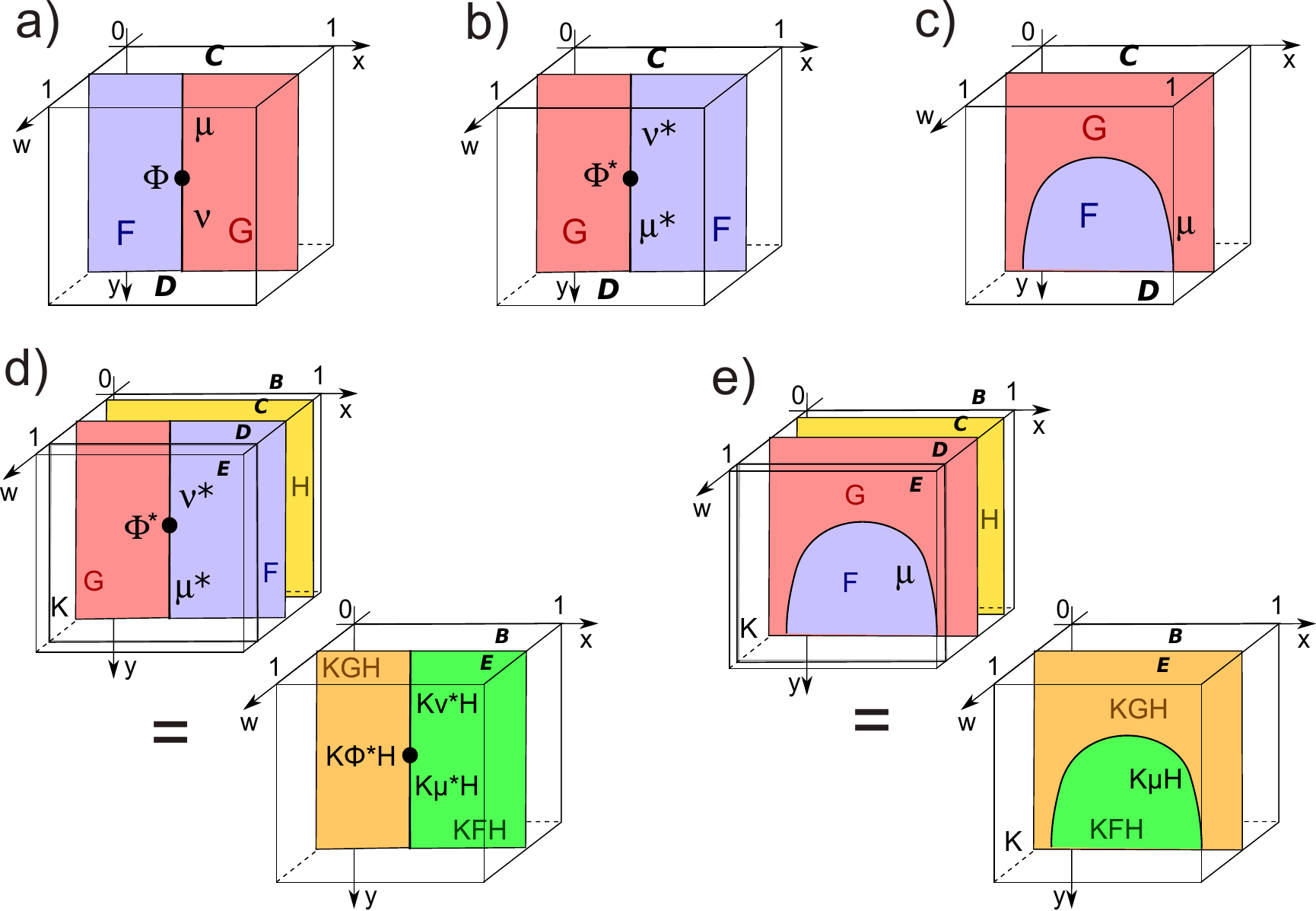}
  \caption{Gray category diagrams for the  $*$-duals:  \newline
    a) 3-morphism $\Phi\colon\mu\Rrightarrow \nu$,
    b) dual  $\Phi^*\colon\nu^*\Rrightarrow \mu^*$, \newline
    c) 3-morphism $\epsilon_\mu\colon 1_G\Rrightarrow \mu\circ \mu^*$, \newline 
    d) identity $(K\Box \Phi\Box H)^*=K\Box \Phi^*\Box H$,\newline
    e) identity  $\epsilon_{K\Box \mu\Box H}=K\Box \epsilon_\mu \Box H$.\newline
  }
  \label{graycatduals3}
\end{figure}

\begin{figure}
  \centering
  \includegraphics[scale=0.6]{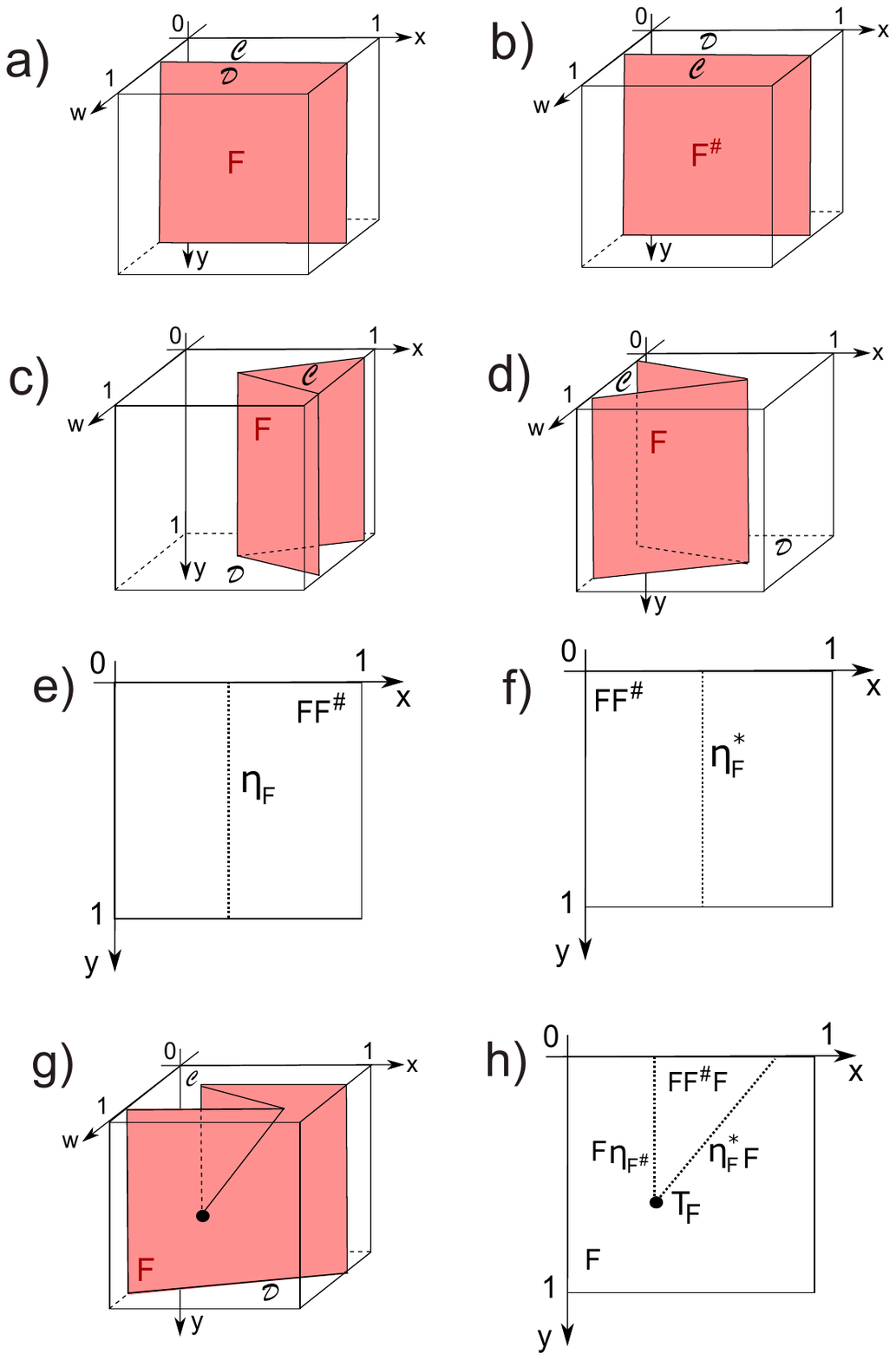}
  \caption{Diagrams for $\#$:\newline
    a) 1-morphism $F\colon \mathcal C\rightarrow\mathcal D$, b) its dual $F^\#\colon\mathcal D\rightarrow\mathcal C$,\newline
    c) Fold $\eta_F\colon 1_{\mathcal D}\Rightarrow F\Box F^\#$, d) its dual $\eta_F^*\colon F\Box F^\#\Rightarrow 1_{\mathcal D}$.\newline
    e) Projection of c), f) projection of d).\newline
    g) Triangulator $T_F\colon (\eta_F^*\Box F)\circ(F\Box\eta_{F^\#})\Rrightarrow 1_{F}$, h) its projection.}
  \label{baez}
\end{figure}

\begin{figure}
  \centering
  \includegraphics[scale=0.6]{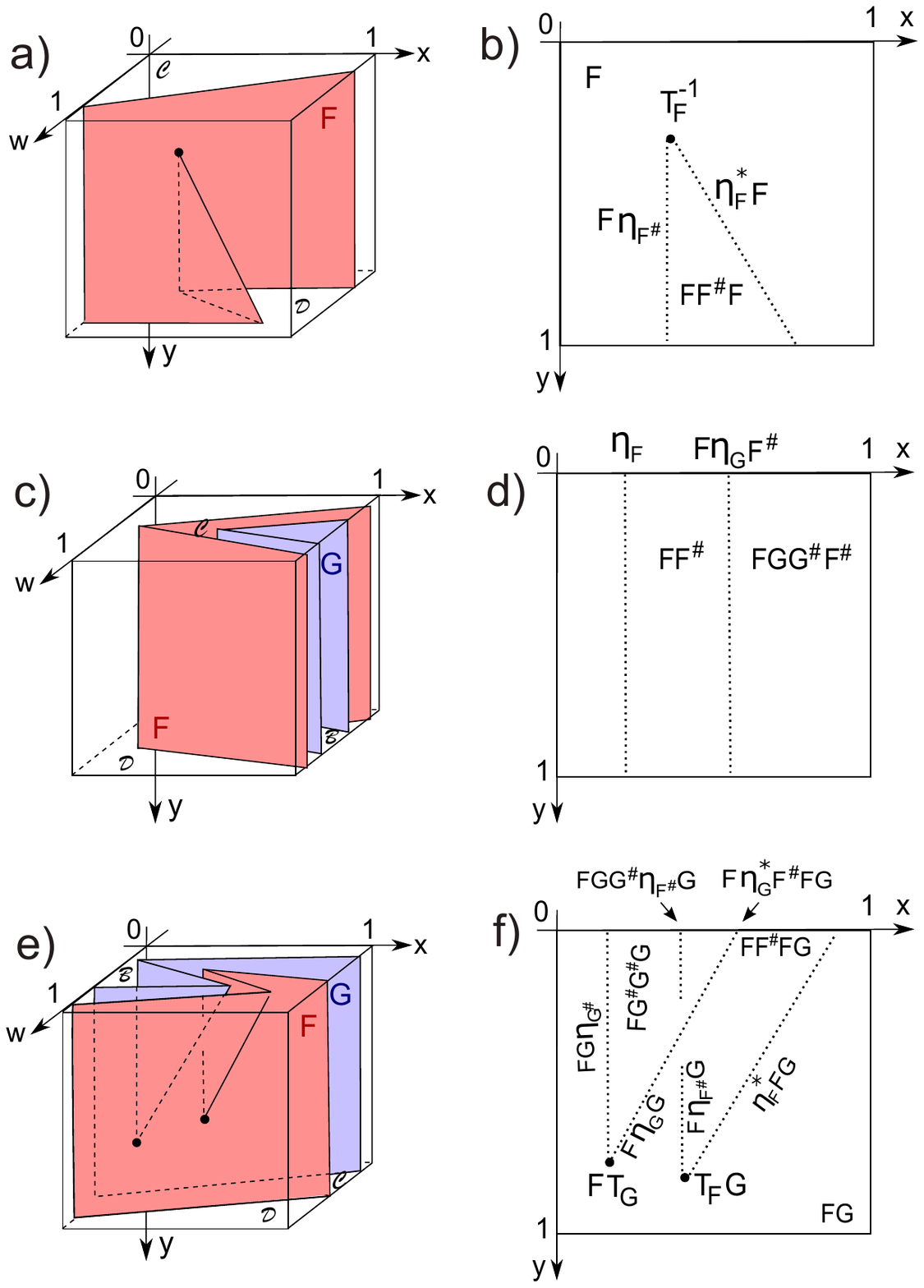}
  \caption{Diagrams for $\#$:\newline
    a) Inverse of the triangulator and b) its projection.\newline
    c) Fold $\eta_{F\Box G}$ and d) its projection.\newline
    e) Triangulator $T_{F\Box G}$ and  f) its projection.}
  \label{BL_consisteny}
\end{figure}

\begin{figure}
  \centering
  \includegraphics[scale=0.6]{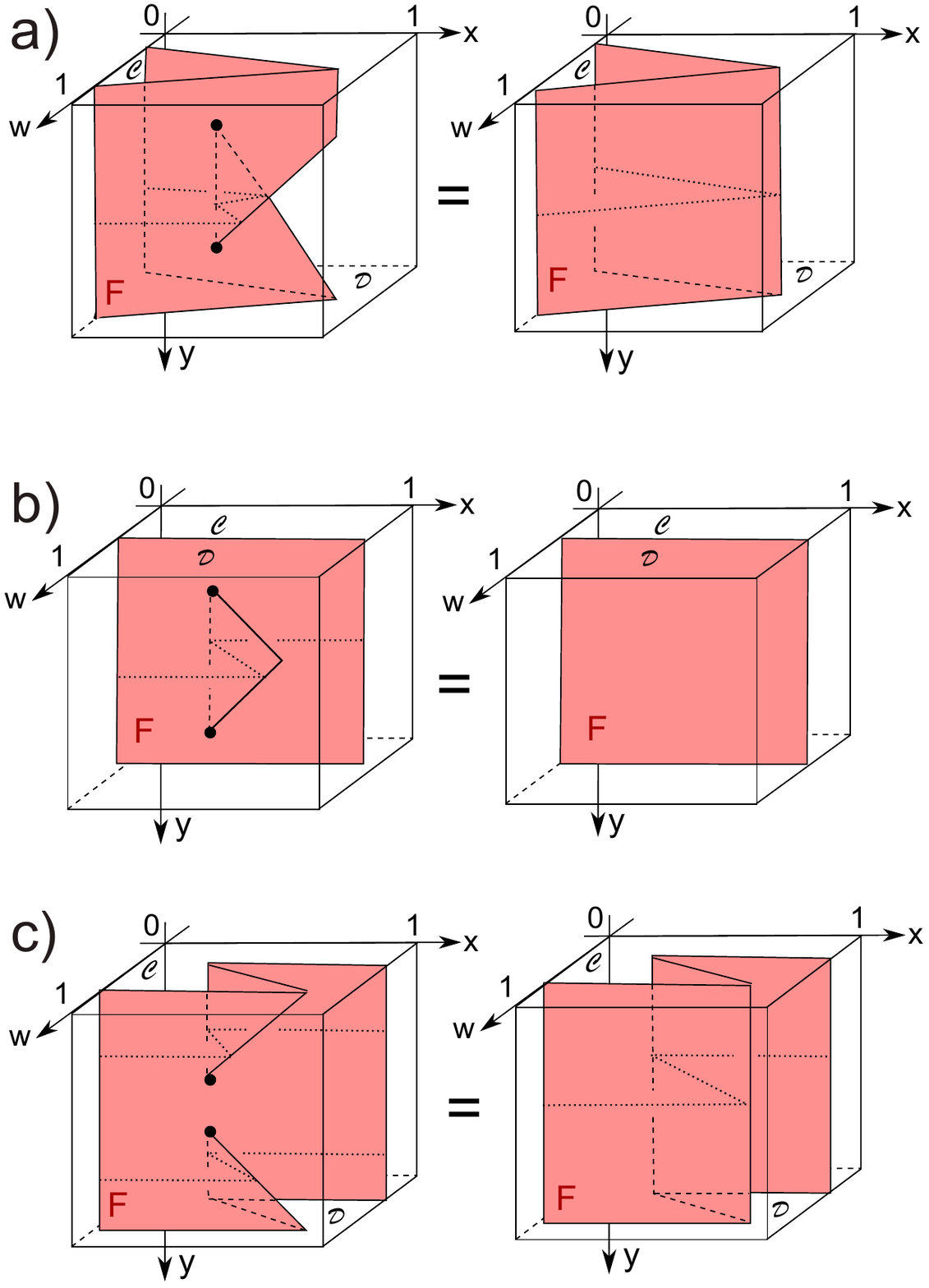}
  \caption{Consistency conditions for triangulators.\newline
    a) Diagrams for the identity  $(1_{\eta_F^*} \circ T_F\Box F^\#)\cdot(\sigma_{\eta_F^*,\eta_F^*}\circ F\Box 1_{\eta_{F^\#}}\Box F^\#)\cdot(1_{\eta_F^*}\circ F\Box T_{F^\#}^*)=1_{\eta_F^*}
    .$\newline
    b) Diagrams for the identity $T_F\cdot T_F^{-1}=1_{1_F}$.\newline
    c)  Diagrams for the identity $T_F^{-1}\cdot T_F=1_{\eta_F^*\Box F\circ F\Box \eta_{F^\#}}$.}
  \label{whitney}
\end{figure}

\begin{figure}
  \centering
  \includegraphics[scale=0.6]{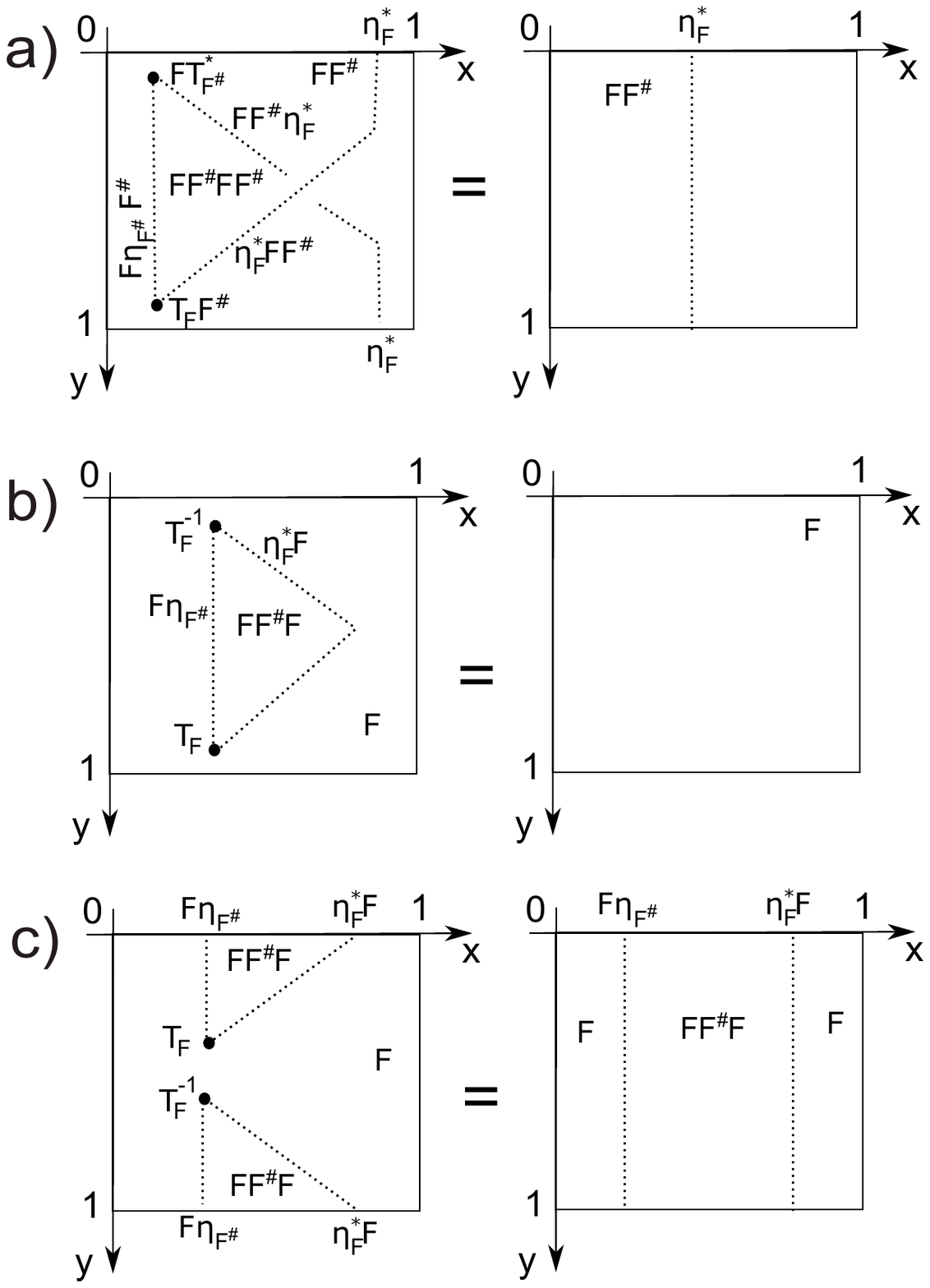}
  \caption{Consistency conditions for triangulators.\newline
    a) Projection of the diagrams in Figure \ref{whitney} a).\newline
    b) Projection of the diagrams in Figure \ref{whitney} b).\newline
    c) Projection of the diagrams in Figure \ref{whitney} c)}
  \label{whitney_proj}
\end{figure}

The set of progressive Gray category diagrams is sufficient to express the axioms of a Gray category with duals in diagrammatic form. 
 Non-progressive diagrams will be discussed in 
Section \ref{sec:diagrams_duals}.
Each of the canonical 2- and 3-morphisms in Definition \ref{dualgray} for a Gray category with duals determines a canonical diagrammatic element as follows:
\begin{itemize}

\item As for planar 2-categories, the 3-morphisms $\epsilon_\mu\colon 1_G\Rrightarrow \mu\circ \mu^*$ 
  are canonical vertices that correspond to  maxima and minima of  the lines. The Gray category diagram for  $\epsilon_\mu$ is  obtained by drawing the corresponding diagram for a planar  category  in Figure \ref{plancat} b) on a plane labelled with two 1-morphisms $F,G$ as shown in Figure \ref{graycatduals3} c).  As shown, the vertex is not labelled by any morphism. By convention this means that the morphism at this vertex is $\epsilon_\mu$. One of the lines meeting this vertex is labelled with $\mu$ and the other with $\mu^*$. It is therefore only necessary to show the label of one of these lines as the other label is then uniquely determined. This convention will be used in the following.

  The  compatibility condition (1) in Definition \ref{dualgray} which involves the Gray product of 1-morphisms with the 3-morphisms  $\epsilon_\mu\colon 1_G\Rrightarrow \mu\circ \mu^*$ is shown in Figure \ref{graycatduals3} e).

\item 
  The diagram for the fold
  2-morphism $\eta_F\colon 1_\mathcal D\Rightarrow F\Box F^\#$  is a line with two planes attached to the right, as shown in Figure \ref{baez} c).
  The convention is that the label for this line is not shown, and  only one of the two surfaces incident to the line is labelled. This diagram is in fact an identity diagram $1_D$, where $D$ is the elementary pre-2-category diagram for  $\eta_F$ in $\mathcal G_2$.
  The $*$-dual $\eta_F^*$ is represented by a diagram with two planes on the left, as shown in Figure \ref{baez} d).  The condition $\eta_{1_{\mac C}}=1_{1_{\mac C}}$ states that the 2-morphism $\eta_{1_{\mac C}}$ corresponds to an empty diagram.
  For better legibility of the diagrams, the 2-morphisms $\eta_F$ and their $*$-duals will also sometimes be drawn as rounded lines in the following. This does not affect any of the results.

\item The invertible 3-morphism $T_F\colon (\eta_F^*\Box F)\circ (F\Box \eta_{F^\#})\Rrightarrow 1_{F}$ for each 1-morphism $F\colon\mathcal C\rightarrow\mathcal D$
  corresponds to  the  Gray category diagram in Figure \ref{baez} g), and
  its inverse is depicted in Figure \ref{BL_consisteny} a).  As in the case of the 2-morphisms $\eta_F$, it is not necessary to label the vertex in this diagram or the lines incident at the vertex, and only one of the incident surfaces is labelled.

  The compatibility of the 2-morphisms $\eta_F$ and the 3-morphism $T_F$ with the Gray product $\Box$ (condition c) in Definition \ref{dualgray}), relates the 2-morphisms $\eta_{F\Box G}$ and the 3-morphisms $T_{F\Box G}$ to the corresponding 2- and 3-morphisms $\eta_F,\eta_G$ and $T_F,T_G$. 
  It states that the two diagrams  in Figure \ref{BL_consisteny} c), e) have a well-defined evaluation.  
  Condition (d) in Definition \ref{dualgray} and  the invertibility of the 3-morphisms $T_F$  are depicted, respectively,  in Figure \ref{whitney} a), b) and c) and their projections in Figure \ref{whitney_proj}. 
  
\end{itemize}

\begin{figure}
  \centering
  \includegraphics[scale=0.6]{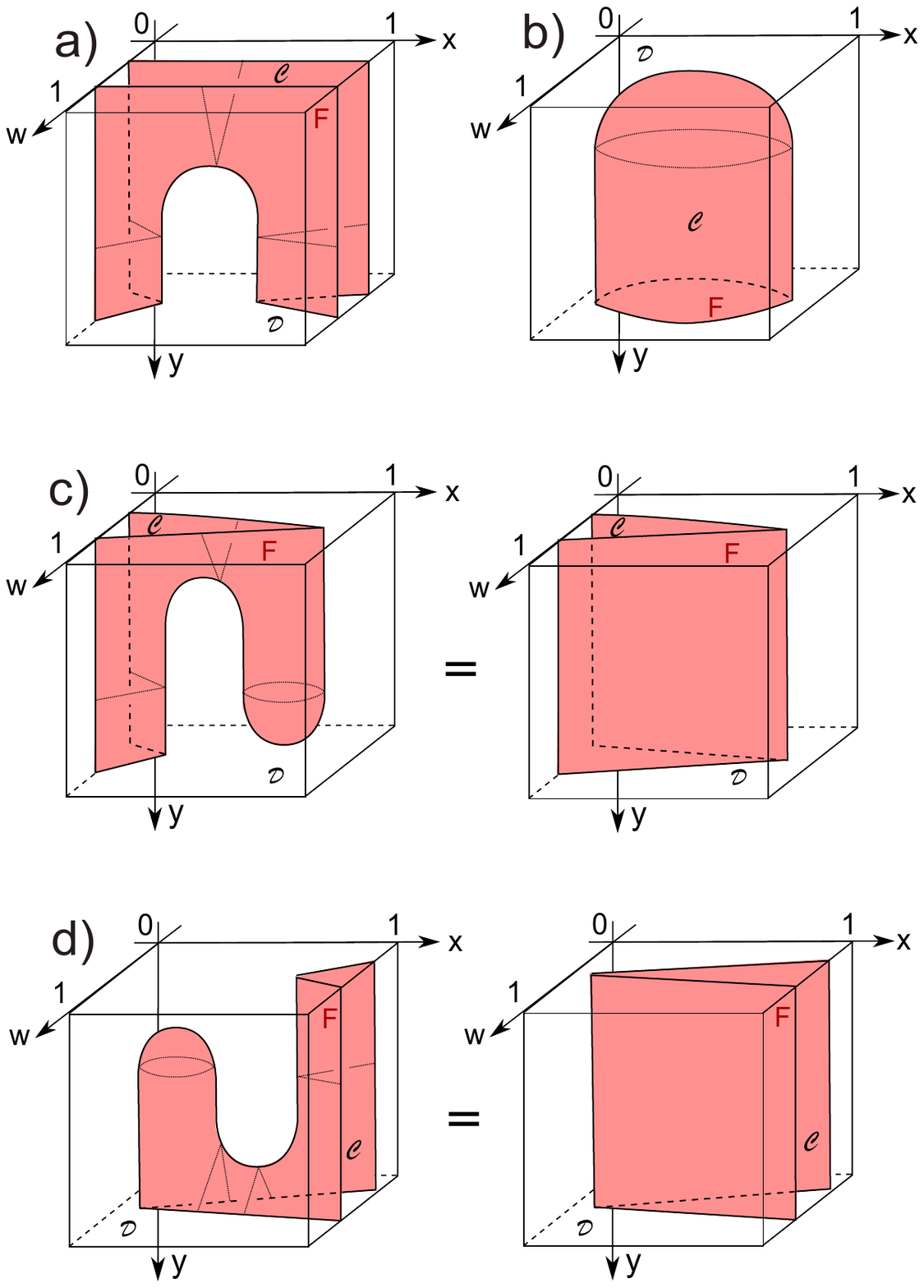}
  \caption{Gray category diagrams:\newline
    a) 3-morphism $\epsilon_{\eta_F}\colon 1_{F\Box F^\#}\Rrightarrow \eta_F\circ
    \eta_F^*$,\newline
    b) 3-morphism $\epsilon_{\eta_F^*}\colon 1_{1_D}\Rrightarrow \eta_F^*\circ
    \eta_F$,\newline
    c) identity $(\epsilon^*_{\eta_F^*}\circ 1_{\eta_F^*})\cdot
    (1_{\eta_F^*}\circ \epsilon_{\eta_F})=1_{\eta_F^*}$,\newline
    d) identity $(\epsilon^*_{\eta_F}\circ 1_{\eta_F})\cdot(1_{\eta_F}\circ
    \epsilon_{\eta_F^*})=1_{\eta_F}$.
  }
  \label{saddle}
\end{figure}

\begin{figure}
  \centering
  \includegraphics[scale=0.3]{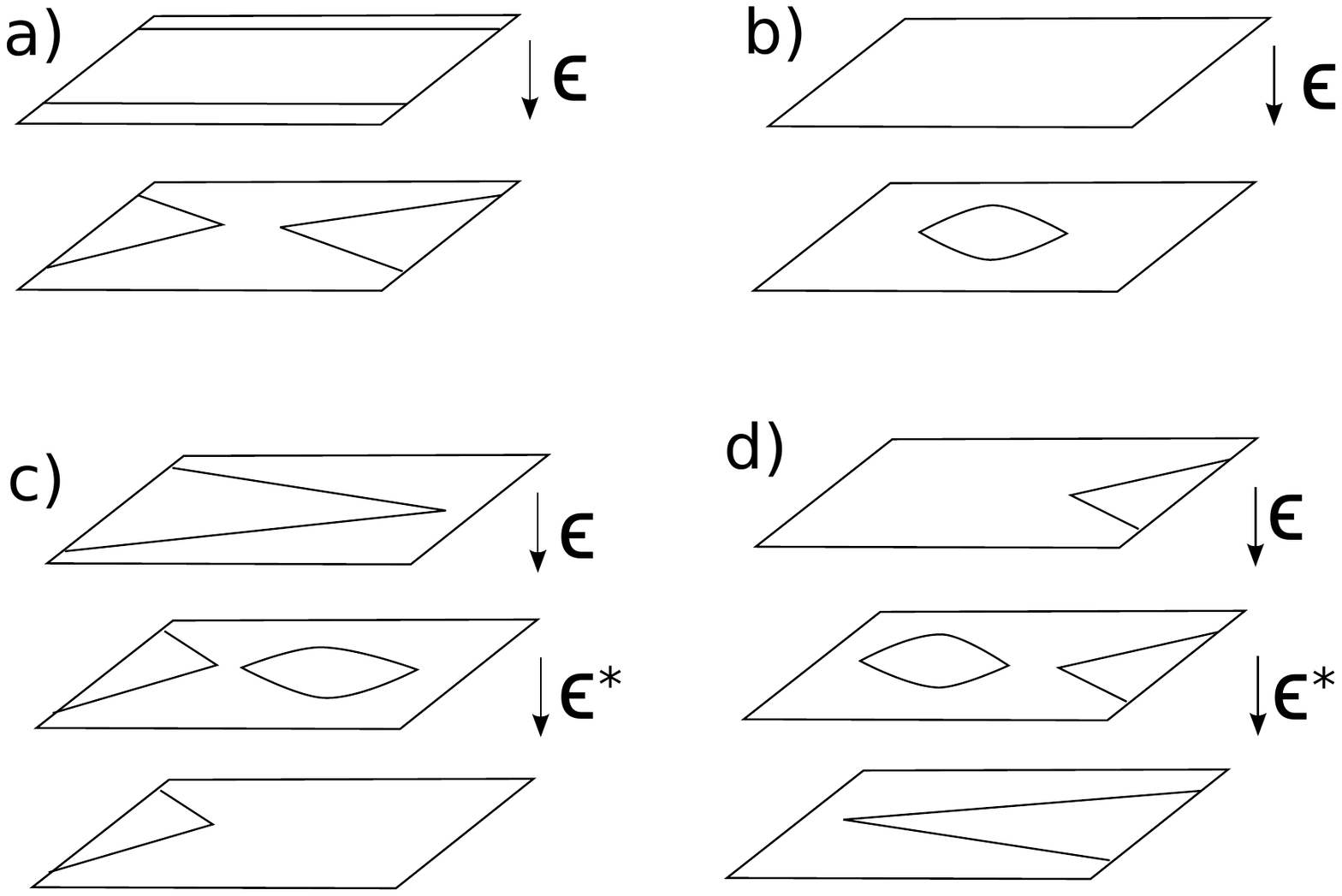}
  \caption{Gray category diagrams from Figure \ref{saddle} in their movie representation obtained by taking constant height slices of the diagrams in Figure \ref{saddle}:\newline
    a) 3-morphism $\epsilon_{\eta_F}\colon 1_{F\Box F^\#}\Rrightarrow \eta_F\circ
    \eta_F^*$ (Figure \ref{saddle} a),\newline
    b) 3-morphism $\epsilon_{\eta_F^*}\colon 1_{1_D}\Rrightarrow \eta_F^*\circ
    \eta_F$ (Figure \ref{saddle} b),\newline
    c) 3-morphism $(\epsilon^*_{\eta_F^*}\circ 1_{\eta_F^*})\cdot
    (1_{\eta_F^*}\circ \epsilon_{\eta_F})$ (Figure \ref{saddle} c),\newline
    d) 3-morphism $(\epsilon^*_{\eta_F}\circ 1_{\eta_F})\cdot(1_{\eta_F}\circ
    \epsilon_{\eta_F^*})$ (Figure \ref{saddle} d).\newline
    Some labels are omitted for legibility.
  }
  \label{saddle_slice}
\end{figure}

By composing these diagrammatic elements, one obtains diagrams for all structural data and relations of  a Gray category with duals.  For instance, the diagrams for the canonical 3-morphisms $\epsilon_{\eta_F}\colon 1_{F\Box F^\#}\Rrightarrow \eta_F\circ \eta_F^*$ and $\epsilon_{\eta_F^*}\colon 1_{\mac D}\Rrightarrow \eta_F^*\circ \eta_F$ for a 1-morphism $F\colon\mac C\to\mac D$ are given in Figures \ref{saddle} and \ref{saddle_slice} a), b).  The identities  
\begin{align*}
  &(\epsilon^*_{\eta_F^*}\circ 1_{\eta_F^*})\cdot (1_{\eta_F^*}\circ \epsilon_{\eta_F})=1_{\eta_F^*},\quad
  (\epsilon^*_{\eta_F}\circ 1_{\eta_F})\cdot(1_{\eta_F}\circ
  \epsilon_{\eta_F^*})=1_{\eta_F}
\end{align*}
from condition \eqref{pivcond} in the definition of a planar 2-category correspond to the diagrams in Figure \ref{saddle}  and \ref{saddle_slice} c), d).

\begin{figure}
  \centering
  \includegraphics[scale=0.6]{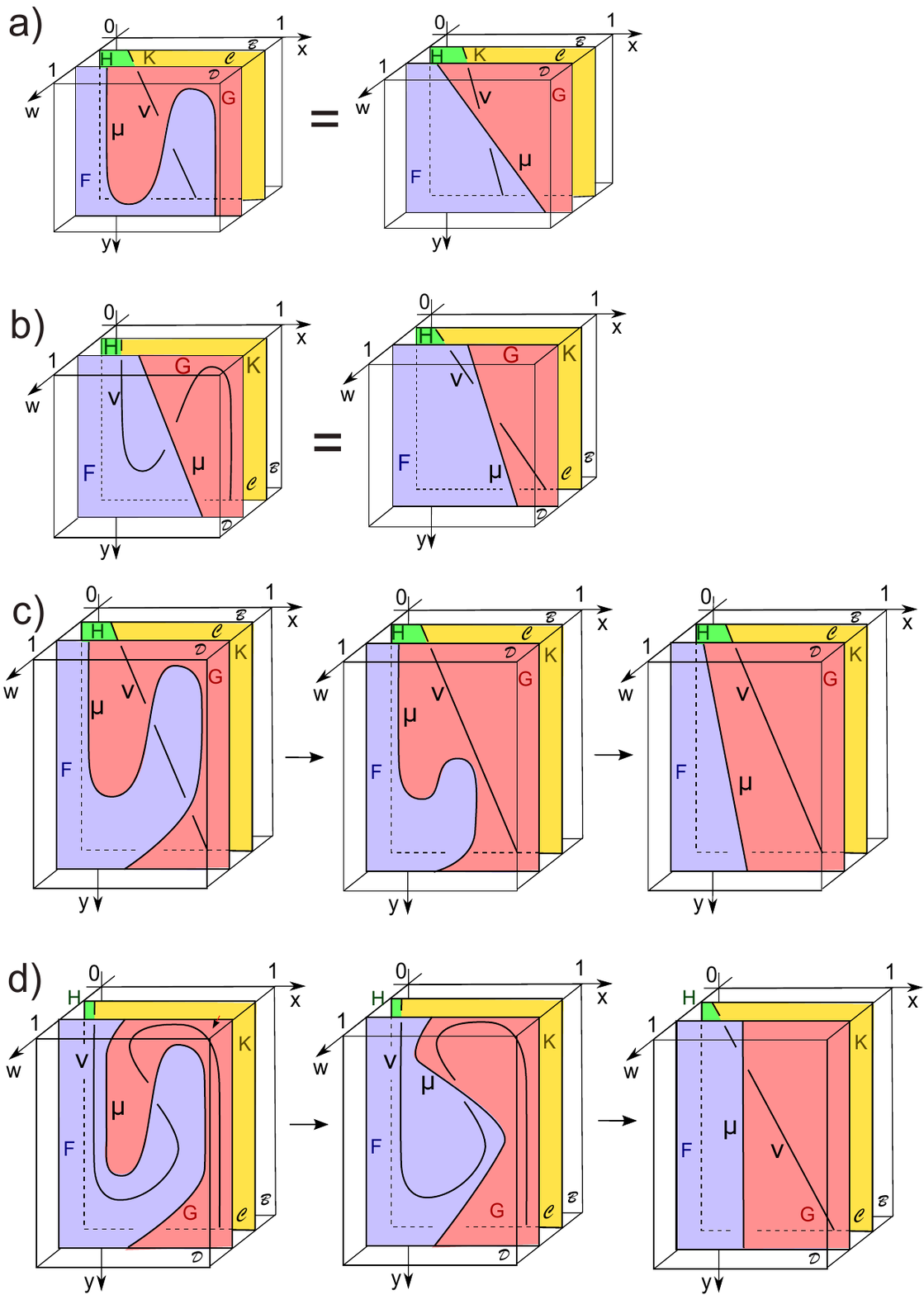}
  \caption{Gray category diagrams for Lemma \ref{sigma_ids}:  \newline
    a)  identity $(1_{(\mu\Box K)\circ (F\Box\nu)}\circ\epsilon^*_{\mu^*\Box
      H})\cdot (1_{\mu\Box K}\circ \sigma_{\mu^*,\nu}\circ 1_{\mu\Box
      H})\cdot(\epsilon_{\mu\Box K}\circ 1_{(G\Box\nu)\circ (\mu\Box
      H)})=\sigma_{\mu,\nu}^\inv$,\newline
    b) identity $(1_{(G\Box\nu)\circ(\mu\Box H)}\circ
    \epsilon^*_{F\Box\nu^*})\cdot (1_{G\Box\nu}\circ
    \sigma^\inv_{\mu,\nu^*}\circ 1_{F\Box\nu})\cdot(\epsilon_{G\Box\nu}\circ
    1_{(\mu\Box K)\circ(F\Box\nu)})=\sigma_{\mu,\nu}$,\newline
    c) diagrammatic proof of an identity equivalent to a),\newline
    d) diagrammatic proof of the identity
    $*\sigma_{\mu^*,\nu^*}=\sigma_{\mu,\nu}$.
  }
  \label{tensorator_identities}
\end{figure}

As in the case of  diagrams for monoidal categories,  Gray category diagrams prove  useful for computations in a Gray category with duals. An example is  the proof of  the last three identities in Lemma \ref{sigma_ids}, which is performed diagrammatically in Figure \ref{tensorator_identities}. 
This diagrammatic proof clearly exhibits the  structural similarities  between  Gray categories with duals and braided pivotal  categories.  

It also becomes apparent that diagrammatic calculations are much simpler than their symbolic counterparts. For these reasons, they will be used extensively in the following sections. Note, however, that at this stage the diagrammatic calculus is to be understood  as a  calculation in a Gray category with duals that is based on the evaluation of {\em progressive} Gray category  diagrams.

Although the diagrams for Gray categories with duals involve lines with maxima and minima, planes with folds (denoting the 2-morphisms $\eta_F$) and cusps (denoting the 3-morphisms $T_F$), these diagrams are considered as being progressive.  The maxima and minima, folds and cusps simply indicate a canonical labelling  of the lines and vertices in a progressive Gray category diagram.   

The diagrammatic  calculations use the fact that the evaluation of a generic progressive Gray category diagram is invariant under certain  isomorphisms of {\em progressive} diagrams (Theorem \ref{th_progressive_invariant}) and  the axioms of a Gray category with duals.
In particular, the diagrams in Figure \ref{whitney} represent relations between certain 3-morphisms in a Gray category with duals and are not to be interpreted as invariance of the  evaluation  under certain  homomorphisms of diagrams at this stage.

Non-progressive Gray category diagrams and the associated  homomorphisms of diagrams will be investigated in Section \ref{sec:diagrams_duals}. In particular, it will be shown there that the diagrams in Figure \ref{whitney} have an interpretation as Whitney moves relating surface projections.

\medskip

\begin{figure}
  \centering
  \includegraphics[scale=0.6]{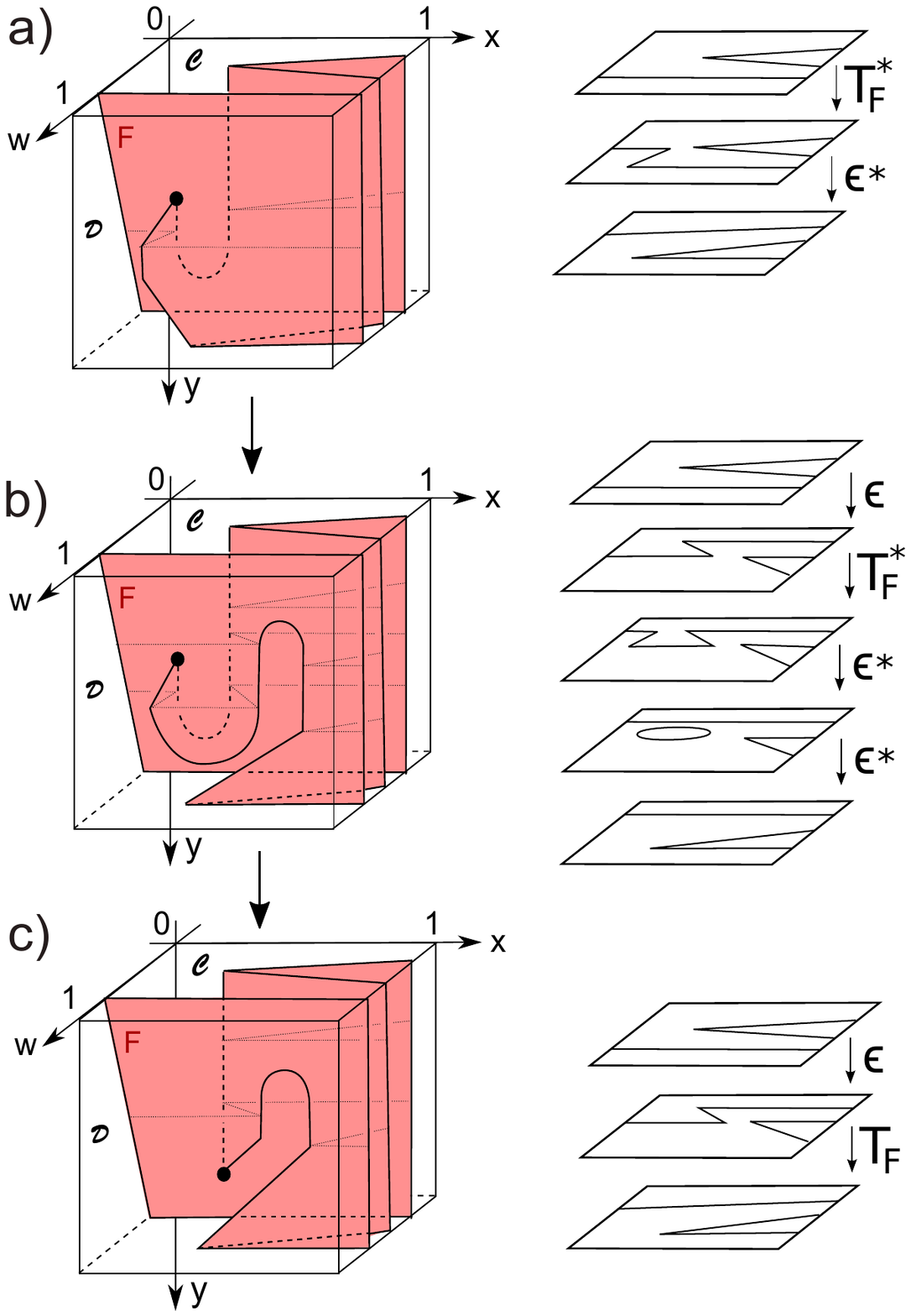}
  \caption{The horizontal cusp identity:\newline 
  a) left-hand side of the equation in Lemma \ref{lem:horcusp},\newline
  b) intermediate stage between  a) and c),  
\newline 
  c) right-hand side of the equation in Lemma \ref{lem:horcusp}.
  }
  \label{horcusp}
\end{figure}

This section is concluded with an application of the diagrammatic calculus. 
In the literature on higher categories sometimes an additional identity, the \emph{horizontal cusp identity}, is required as a condition for the interplay of two different dualities (see Figure \ref{horcusp}).
For instance, it appears in the approach to adjoint equivalences in the context of bicategories in \cite[Appendix C]{Bart} and is used in the discussion of the 
 sphere inversion in a higher categorical context in \cite[Figure 7.9]{ScottC}. 
In our framework it follows directly from the axioms.
\begin{lemma}\label{lem:horcusp}
  Every 1-morphism  $F\colon \mac C \rightarrow \mac D $  in a Gray category $\mac G$ with duals satisfies the \emph{horizontal cusp identity}
$$
\left(    (F \Box \epsilon_{\eta_{F^{\#}}}^{*}) \circ 1_{\eta_{F} \Box F } \right) 
\cdot (
1_{F\Box\eta_{F^{\#}}}   
\circ T_{F}^{*})
=(1_{\eta_{F}\Box F} \circ T_{F}) 
\cdot  \left( (\epsilon_{\eta_{F}} \Box F ) \circ 1_{F \Box \eta_{F^{\#}}}  \right)
$$
\end{lemma}
\begin{proof} A graphical proof is given in Figure \ref{horcusp}.
\end{proof}

\begin{remark}
  Very similar diagrams to the Gray category diagrams representing the duals in a Gray category with duals appeared in \cite[Figure 3.13]{Schom}. However, the interpretation of the diagrams there is different:
  They  present the oriented 2-dimensional extended bordism category $\mathrm{Bord}_{2}^{\mathrm or}$ as a symmetric monoidal computad \cite[Sec. 2.10]{Schom}, i.e.~they list
  objects (the two oriented points $(+)$ and $(-)$) and 1- and 2-morphisms (corresponding evaluation and coevaluation bordisms) as well as relations between the 2-morphisms that generate   $\mathrm{Bord}_{2}^{\mathrm or}$
  as a symmetric monoidal bicategory. Thus, the diagrams are \emph{not} embedded 3-dimensional diagrams as in our case, but abstract diagrams providing relations between the generators. As they represent generators of a  \emph{symmetric} monoidal bicategory, the diagrams are 3-dimensional, but the interpretation is independent of the embedding.

  Still, there is one notable connection between our diagrams and the results of \cite{Schom}. Suppose that $\mac G$ is a Gray category with duals, such that the corresponding braided monoidal category
 $\mac G(1_{\mac C}, 1_{\mac C})$ from Lemma \ref{braidten} is actually a symmetric monoidal category. Then each object $F$ of $\mac C$ provides an oriented fully extended 2-dimensional TFT, that is a symmetric monoidal functor $\mathrm{Bord}_{2}^{\mathrm or} \rightarrow \mac C$:
  By \cite[Thm. 3.50, Thm. 2.78 and Prop. 2.74]{Schom} such a functor is given by  providing objects, 1- and 2-morphisms of $\mac C$ that fulfill the relations of \cite[Figure 3.13]{Schom}. These data and relations are all present in a Gray category with duals: To the two objects $(+)$ and $(-)$ of $ \mathrm{Bord}_{2}^{\mathrm or}$
  one assigns the objects $F$ and $F^{\#}$ of $\mac C$ and to the evaluation and coevaluation 1- and 2-morphisms
the corresponding morphism $\eta_{F}$, $\eta_{F}^{*}$ and $\epsilon$-morphisms of $\eta$ from the data of a Gray category with duals. The relations between the morphisms in \cite{Schom} then correspond to the axioms of a Gray category with duals together with the horizontal cusp identity in Lemma \ref{lem:horcusp}. 
\end{remark}

\section{Duals as functors of 2-strict tricategories}
\label{grayfuncs}

In this section it is shown that  the duals $*$ and $\#$  in Definition \ref{dualgray}
define functors of 2-strict tricategories (see Definition \ref{grayfunc} in Appendix \ref{2functors}). As these functors will reverse different products in the Gray categories, this requires a notion of different opposites for Gray categories and functors thereof.  

\subsection{Higher opposites} The opposite 2-categories $\mathcal C^{op}$, $\mac C_{op}$ and 
contravariant functors associated to strict 2-functors $\mathcal C\to\mathcal D^{op}$ were introduced in Definition  \ref{op2catdef}. Here, the discussion is extended to  weak 2-functors. These constructions lead to the definition of two notions of the opposite of a Gray category. 

A weak 2-functor of 2-categories $F\colon\mathcal C\to\mathcal D$ (Definition  \ref{laks2func}) has the following data 
\begin{itemize}
\item A function $F_0\colon \text{Ob}(\mathcal C)\rightarrow \text{Ob}(\mathcal D)$.
\item For all objects $G,H$ of $\mathcal C$, a functor $F_{G,H}\colon \mathcal C_{G,H}\rightarrow \mathcal D_{F_0(G), F_0(H)}$.
\item For all 1-morphisms $\nu\colon G\rightarrow H$, $\mu\colon H\rightarrow K$, a 2-isomorphism $\Phi_{\mu,\nu}\colon F_{H,K}(\mu)\circ F_{G,H}(\nu)\rightarrow F_{G,K}(\mu\circ\nu)$. 
\item For all objects $G$, a 2-isomorphism $\Phi_G\colon 1_{F_0(G)}\to F_{G,G}(1_G)$.
\end{itemize}

The opposite $F^{op}\colon\mathcal C^{op}\to\mathcal D^{op}$ is determined by the following data
\begin{itemize}
\item $(F^{op})_0=F_0$
\item $(F^{op})_{H,G}=(F_{G,H})^{op}$ 
\item $(\Phi^{op})_{\nu,\mu}=(\Phi_{\mu,\nu}^{-1})^{op}$
\item $(\Phi^{op})_G=(\Phi_G^{-1})^{op}$,
\end{itemize}
where the right-hand involves the 1-categorical opposites defined in Section \ref{2catdual}.
The corresponding  opposite weak 2-functor $F_{op}\colon\mathcal C_{op}\to\mathcal D_{op}$ is determined by
\begin{itemize}
\item $(F_{op})_0=F_0$
\item $(F_{op})_{H,G}=F_{G,H}$ 
\item $(\Phi_{op})_{\nu,\mu}=\Phi_{\mu,\nu}$
\item $(\Phi_{op})_G=\Phi_G.$
\end{itemize}
Note that the choice of the coherence data for $F^{op}$ and $F_{op}$ is determined unambiguously  by  the coherence data of $F$ and the source and target of the coherence isomorphism. Hence it is justified to abuse notation and denote the functors
$F^{op}$ and $F_{op}$ again $F$. In the following,  the notation  $F^{op}$ and $F_{op}$ will be used only to emphasise their relations. 

These notions allow the definition of two types of opposite Gray category. The key point is the opposite of a Gray category is a tricategory but it is not necessarily a Gray category in the sense of Definition \ref{tricatGray} since, depending on which products are reversed, 
the resulting strict tricategory can be cubical instead 
of opcubical.

\begin{definition} Let $\mathcal G$ be a Gray category. 
  \begin{enumerate}
  \item The {\bf tricategory $\mathcal G^{op}$} has the same Gray product but opposite horizontal and vertical composition. Thus the 2-categories are $\mathcal G^{op}(\mathcal C,\mathcal D)=(\mathcal G(\mathcal C,\mathcal D))^{op}$ and the Gray product is the collection of opposite weak 2-functors $\Box^{op}$. Thus, in the abbreviated notation introduced above and in Definition \ref{op2catdef},  compositions of $\mathcal{G}^{op}$  are denoted $\Box$, $\widetilde{\circ}$ and $\widetilde\cdot$.

  \item The {\bf tricategory $\mathcal G_{op}$} has the same vertical composition but opposite Gray product and horizontal composition:  $\mathcal G_{op}(\mathcal C,\mathcal D)=(\mathcal G(\mathcal D,\mathcal C))_{op}$.  The Gray product for $\mathcal G_{op}$ is the collection of weak 2-functors $\widetilde{\Box}_{op}$ where $\Phi \widetilde{\Box}\Psi=\Psi\Box\Phi$ for all composable 2- and 3-morphisms $\Phi,\Psi$. Thus, in abbreviated notation, the compositions of $\mathcal{G}_{op}$ are $\widetilde{\Box}$, $\widetilde{\circ}$ and $\cdot$.
  \end{enumerate}
\end{definition}
There is the analogous definition of the opposites $\mathcal{T}^{op}$ and $\mathcal{T}_{op}$ of a strict cubical tricategory $\mathcal{T}$.
It is clear that $(\mathcal{T}^{op})^{op}=\mathcal{T}$ and $(\mathcal{T}_{op})_{op}=\mathcal{T}$.
In the following, the objects and morphisms of the opposite categories are not distinguished  in our notation but the compositions in $\mathcal{T}^{op}$ and $\mac T_{op}$ are denoted with the appropriate $op$-label.

The expressions for the coherence isomorphisms and  tensorators of the opposite strict  tricategories
are obtained by unpacking these definitions. 
Using the notation  introduced above, in Lemma \ref{cubegray} and Corollary \ref{lemma:cubical-to-opcubical-tricat}, one finds  that if $\mac T$ is a Gray category with
 structure isomorphisms $\Box_{\mu,\nu}$, then $\mathcal T^{op}$ has the structure isomorphism $(\Box^{op})_{\mu,\nu}=(\Box^{-1}_{\nu,\mu})^{op}$, and $\mathcal T_{op}$ has $\widetilde{\Box}_{op\,\mu,\nu}=\Box_{\nu',\mu'}$, with $\mu'=(\mu_2,\mu_1)$, $\nu'=(\nu_2,\nu_1)$. 
The following statement then follows directly from the definitions.
\begin{lemma}
  Let $\mathcal{T}$ be a strict cubical (opcubical) tricategory. Then $\mathcal{T}_{op}$ is again a strict cubical (opcubical) tricategory, while $\mathcal{T}^{op}$ is a strict  opcubical (cubical) tricategory.
\end{lemma}

Finally, there are notions of opposite functors for 2-strict tricategories. Let $F\colon T\to S$ be a functor of 2-strict tricategories. This has data (Definition \ref{grayfunc})
\begin{itemize}
\item a function $F_0\colon\text{Ob}(\mathcal T)\rightarrow \text{Ob} ({\mathcal S})$,
\item weak 2-functors $F_{\mathcal C,\mathcal D}\colon \mathcal T(\mathcal C,\mathcal D)\rightarrow \mathcal S(F_0(\mathcal C), F_0(\mathcal D))$ for all objects $\mathcal C,\mathcal D$ of $\mathcal T$
\item   For $\Box$-composable 2-morphisms $\mu$ and $\nu$, 3-isomorphisms \\
  $\kappa_{\mu,\nu}\colon F(\mu)\Box F(\nu)\to F(\mu\Box\nu)$.
\end{itemize}

The first opposite is the functor of 2-strict tricategories $F^{op}\colon T^{op}\to S^{op}$ with data
\begin{itemize}
\item $(F^{op})_0=F_0$
\item $(F^{op})_{\mathcal C,\mathcal D}=(F_{\mathcal C,\mathcal D})^{op}$, using the opposite of weak 2-functors
\item $(\kappa^{op})_{\nu,\mu}=(\kappa^{-1}_{\nu,\mu})^{op}.$
\end{itemize}

The second opposite is the functor of 2-strict tricategories $F_{op}\colon T_{op}\to S_{op}$ with data
\begin{itemize}
\item $(F_{op})_0=F_0$
\item $(F_{op})_{\mathcal C,\mathcal D}=(F_{\mathcal C,\mathcal D})_{op}$, using the opposite of weak 2-functors
\item $(\kappa_{op})_{\nu,\mu}=\kappa_{\mu,\nu}.$
\end{itemize}
Again, all coherence data is unambiguous and we denote the functors $F^{op}$ and $F_{op}$ by  $F$.

\subsection{Duals as functors}

To define a functor of 2-strict tricategories $*\colon \mac G\to\mac G^{op}$, the dual $*$
in the planar 2-categories $\mathcal G(\mathcal C,\mathcal D)$  is extended trivially to  objects and 1-morphisms of  $\mathcal G$. Similarly,  the dual $\#$ is  the identity on objects of $\mathcal G$. To extend it to 2- and 3-morphisms,
 define  for each 2-morphism $\nu\colon F\Rightarrow G$ and 3-morphism $\Psi\colon\nu\Rrightarrow\mu$  the associated 
$\#$-duals $\#\nu\colon G^\#\Rightarrow F^\#$, $\#\Psi\colon \#\nu\Rightarrow \#\mu$ by
\begin{align}\label{hashdualformula}
  &\#\nu=(F^\#\Box \eta_G^*)\circ (F^\#\Box\nu\Box G^\#)\circ (\eta_{F^\#}\Box G^\#)\\
  &\#\Psi=1_{F^\#\Box \eta_G^*}\circ (F^\#\Box \Psi\Box G^\#)\circ 1_{ \eta_{F^\#}\Box G^\#}.\nonumber
\end{align}
The diagrams for $\#\nu$ and $\#\Psi$ are given in Figure \ref{crossdiags} a). They correspond to folding the plane segment labelled by the 1-morphism $F$ to the front and the plane segment labelled by $G$ to the back of the cube.

The operations $*$ and $\#$ reverse some of the products and so extend to (partially) contravariant functors.   
The passage from a duality operation to a functor $F$ representing the duality is as follows. All the mappings in the definition of the functor $F$, i.e., $F_{0}$ and $F_{\mathcal{C},\mathcal{D}}$ in the notation of Definition 
\ref{grayfunc}, are given directly by the duality operation, with the result regarded  as an object or  morphism in the appropriate  opposite category. For example, for the functor $* \colon \mathcal{G} \rightarrow \mathcal{G}^{op}$, $\nu^{*}$ is regarded 
as $(\nu^{*})^{op}$ in $\mathcal{G}^{op}$ for a 2-morphism $\nu$ in $\mathcal{G}$.

\begin{theorem}[Duals as  functors of 2-strict tricategories] 
$\quad$\label{graycatduals}
  \begin{enumerate}
  \item The duality operation $*$ extends to a strict functor of 2-strict
    tricategories 
    $*\colon\mathcal G\rightarrow\mathcal G^{op}$ in the sense of Definition \ref{strictdef}.
  \item The duality operation $\#$ extends to a functor of 2-strict tricategories $\#\colon\mathcal G\rightarrow\mathcal G_{op}$ in the sense of Definition \ref{grayfunc}.
  \end{enumerate}
\end{theorem}

\begin{proof} $\quad$\newline\noindent
  1. The data for the functor $ *$ is
  \begin{itemize}\item The identity mapping on objects,
  \item The strict 2-functors $ *_{\mathcal C,\mathcal D}\colon\mathcal G(\mathcal C,\mathcal D)\to\mac G^{op}(\mathcal C, \mathcal D)$ defined by $\phi\mapsto (\phi^*)^{op}$ using the planar structure of $\mathcal G(\mathcal C,\mathcal D)$,
  \item natural isomorphisms  $\rho_{\mac B,\mac C,\mac D}\colon \Box^{op}\,( *_{\mac C,\mac D}\times  *_{\mac B,\mac C})\to  *_{\mac B,\mac D}\,\Box$ of weak 2-functors defined by components $\rho_{\mu,\nu}=(\sigma_{\mu^*,\nu^*})^{op}$.
  \end{itemize}
  To check that this is a strict 
  functor of 2-strict tricategories, we verify the axioms from Definition \ref{grayfunc} and Definition \ref{strictdef}. Since for composable 2-morphisms $\mu\colon F \Rightarrow G$ and $\nu\colon G \Rightarrow K$ in $\mathcal{G}$,
  $$ (*(\mu \Box \nu))^{op}=((F \Box \nu^{*}) \circ (\mu^{*} \Box K))^{op}=((\mu^{*})^{op} \widehat{\Box^{op}} (\nu^{*})^{op}),$$
  with $\widehat{\Box}$ defined as in  Corollary  \ref{lemma:cubical-to-opcubical-tricat},   this shows that $*$ defines a functor of 2-strict tricategories $\tilde{*}\colon \mathcal{G} \rightarrow \widehat{ \mathcal{G}^{op}}$. The coherence data of the functor $*$  is obtained by  composing  
  $\tilde{*}$ with  the functor $\Sigma\colon \widehat{\mathcal{G}^{op}} \rightarrow \mathcal{G}^{op}$ from Corollary 
  \ref{lemma:cubical-to-opcubical-tricat}: $ *= \Sigma \circ \tilde{*}$. This
  produces precisely the natural isomorphism  $\rho_{\mac B,\mac C,\mac D}$
  given above and shows that $*$  is indeed a  strict 
  functor of 2-strict tricategories in the sense of Definitions \ref{grayfunc} and \ref{strictdef}.

  \medskip\noindent
  2. As the duality $\#$ is the identity on the objects  of $\mathcal G$, it is sufficient to show that for each pair of objects $\mac C,\mac D$  
  the dual $\#$ defines weak 2-functors $\#_{\mac C,\mac D}\colon\mac G(\mac C,\mac D)\to\mac G_{op}(\mac C,\mac D)$ 
  and there are natural isomorphisms  $\kappa_{\mac B,\mac C,\mac D}\colon \widetilde{\Box}(\#_{\mac C,\mac D}\times \#_{\mac B,\mac C})\to \#_{\mac B,\mac D}\Box$  of weak 2-functors  
  satisfying the conditions specified in Definition \ref{grayfunc}. 

  To show that $\#$  defines a weak 2-functor $\#\colon\mathcal G(\mathcal C,\mathcal D)\rightarrow\mathcal G_{op}(\mathcal C,\mathcal D)$, we note that for all objects $F$ in $\mathcal G(\mathcal C,\mathcal D)$
  $$
  \#1_{F}= 
  (F^\#\Box \eta^*_F)\!\circ\! (\eta_{F^\#}\Box F^\#).
  $$
  This implies that the\ $*$-dual of the triangulator  defines an
  invertible 3-morphism  $\Phi_F=T^*_{F^\#}\colon 1_{F^\#}\Rrightarrow \#1_F$ with $\Phi_{1_\mathcal C}=1_{1_\mathcal C}$ in $\mathcal{G}_{op}$.   For each pair of composable 1-morphisms $\nu\colon F\Rightarrow G$, $\rho\colon E\Rightarrow F$ in $\mathcal G(\mathcal C,\mathcal D)$, one obtains a 2-morphism  $\Phi_{\rho,\nu}\colon \#\rho\circ\#\nu\Rrightarrow\#(\nu\circ\rho)$ in $\mathcal{G}_{op}(\mathcal{C},\mathcal{D})$ by
  \begin{align*}
    &\Phi_{\rho,\nu}=(1_{E^\#\eta_G^*\circ E^\#\nu G^\#}\circ E^\#T_{F}G^\#\circ 1 _{E^\#\rho G^\#\circ \eta_{E^\#}G^\#})\\
    &\cdot(\sigma_{E^\#\eta_F^*, \eta_G^*\circ\nu G^\#}\circ 1 _{E^\#F\eta_{F^\#}G^\#}\circ 1_{E^\#\rho G^\#\circ \eta_{E^\#}G^\#})\cdot(1_{E^\#\eta_F^*}\circ \sigma_{E^\#\rho\circ \eta_{E^\#} ,\#\nu}).
  \end{align*}
  The corresponding Gray diagram and its projection are given in Figure \ref{crossdiags} c), d) and in Figure \ref{crossdiags_slice} a). 
  It follows directly from the invertibility of the triangulator and the tensorator that  $\Phi_{\rho,\mu}$ is invertible, and the naturality of the tensorator implies that it is natural in both arguments.
  It remains to prove the identities
  \begin{align}\label{compat1}
    \Phi_{1_F,\nu}\cdot(\Phi_F\circ 1 _{\#\nu})=1_{\#\nu} = \Phi_{\nu, 1_G}\cdot (1_{\#\nu}\circ \Phi_G)
  \end{align}
  and the commutativity of the diagram
  \begin{align}\label{compat2}
    \xymatrix{\#\rho\circ\#\nu\circ \#\mu  \ar[r]^{1_{\#\rho}\circ \Phi_{\nu,\mu}} \ar[d]^{\Phi_{\rho,\nu}\circ 1_{\#\mu}} & \#\rho\circ \#(\mu\circ \nu)\ar[d]^{\Phi_{\rho,\mu\circ\nu}}\\
      \#(\nu\circ\rho)\circ \#\mu \ar[r]^{\Phi_{\nu\circ \rho, \mu}} & \#(\mu\circ\nu\circ\rho)  
    }
  \end{align}
  which correspond to the two consistency conditions in Definition \ref{laks2func}. For this,  note that for a 1-morphism $F\colon\mac C\to\mac D$,  the 3-morphism $\Phi_{1_F,\nu}$ is given by
  \begin{align*}
    \Phi_{1_F,\nu}=&(1_{F^\# \eta_{G}^*\circ F^\#\nu G^\#} \circ F^\#T_F G^\#\circ 1_{\eta_{F^\#} G^\#})\\
    &\cdot (\sigma_{F^\#\eta_F^*\,,\, \eta_G^*\circ \nu G^\#}\circ 1_{F^\#F\eta_{F^\#}G^\#}\circ 1_{\eta_{F^\#}G^\#})\cdot (1_{F^\#\eta_{F}^*}\circ \sigma_{F^\#\eta_{F^\#}\,,\, \#\nu}).
  \end{align*}
  Composing this expression with $T^*_{F^\#}\circ 1 _{\#\nu}$, one finds that
  the  conditions in \eqref{compat1} follow from the naturality of the tensorator $\sigma$, together with the invertibility of $T_F$  and  identity (d) in Definition \ref{dualgray} - see also the third diagram in Figure \ref{whitney}. A diagrammatic proof is given in Figure \ref{P_consist}. The commutativity of the diagram \eqref{compat2} follows from the naturality of the tensorator $\sigma$ together with the invertibility of $T_F$ and the exchange law for 2-categories. A diagrammatic proof is given in Figure \ref{prhomuconsist}. 
  This shows that for all objects $\mathcal C,\mathcal D$ of $\mathcal G$, the duality $\#$ defines a weak 2-functor  $\#_{\mathcal{C},\mathcal{D}}\colon\mathcal G(\mathcal C,\mathcal D)\rightarrow\mathcal G(\mathcal C,\mathcal D)_{op}$.

  \begin{figure}
    \centering
    \includegraphics[scale=0.6]{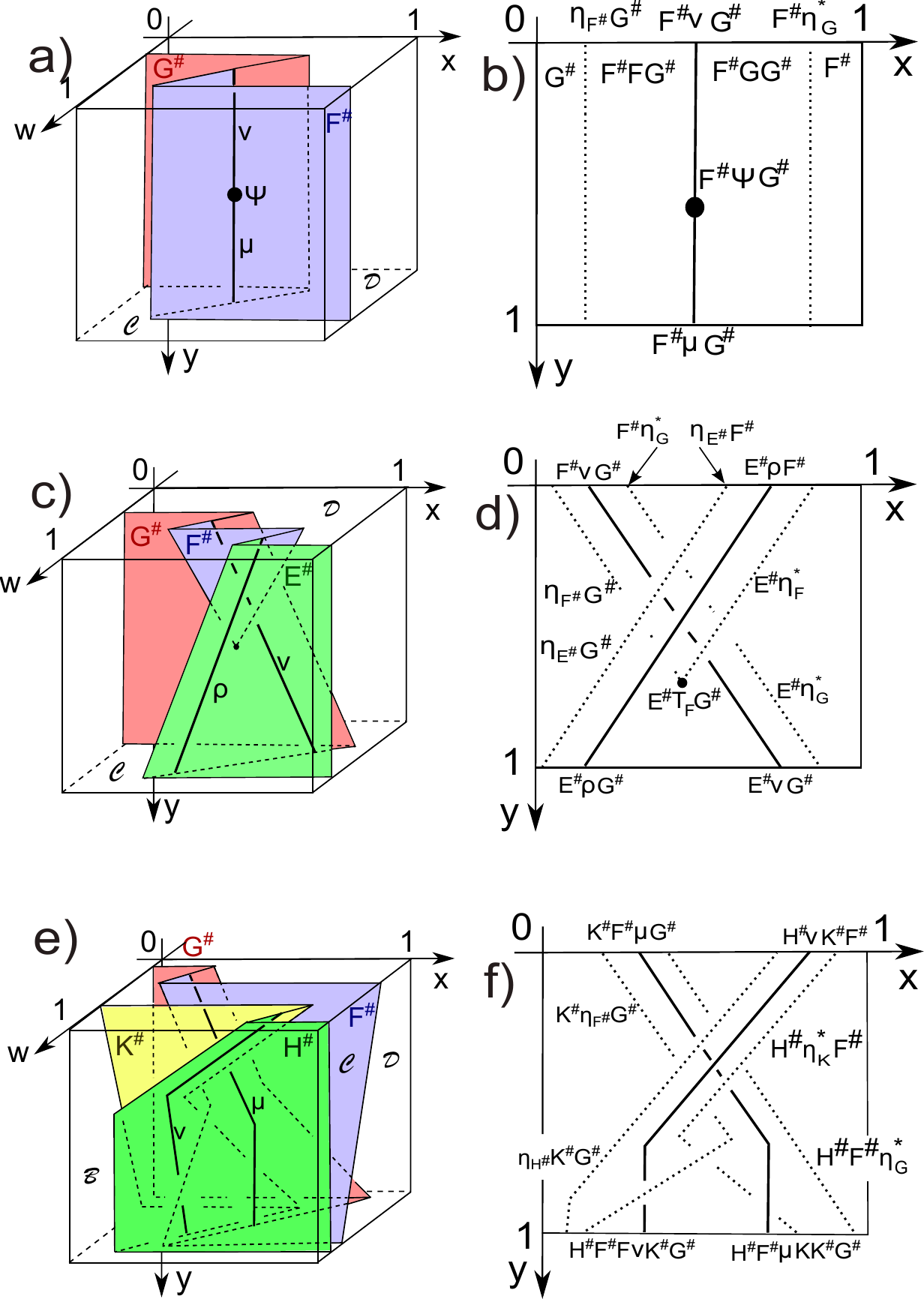}
    \caption{Diagrams for $\#$:\newline
      a),b) 3-morphism  $\#\Psi\colon \#\nu\Rrightarrow\#\mu$ and its projection.\newline
      c),d) 3-morphism  $\Phi_{\rho,\nu}\colon \#\rho\circ\!\#\nu\!\Rrightarrow\! \#(\nu\circ \rho)$, its projection.\newline
      e),f)  3-morphism  $\kappa_{\mu,\nu}\colon \#\nu\Box\#\mu\Rrightarrow \#(\mu\Box\nu)$ and its projection.
      Some labels are omitted for legibility.
    }
    \label{crossdiags}
  \end{figure}

  \begin{figure}
    \centering
    \includegraphics[scale=0.35]{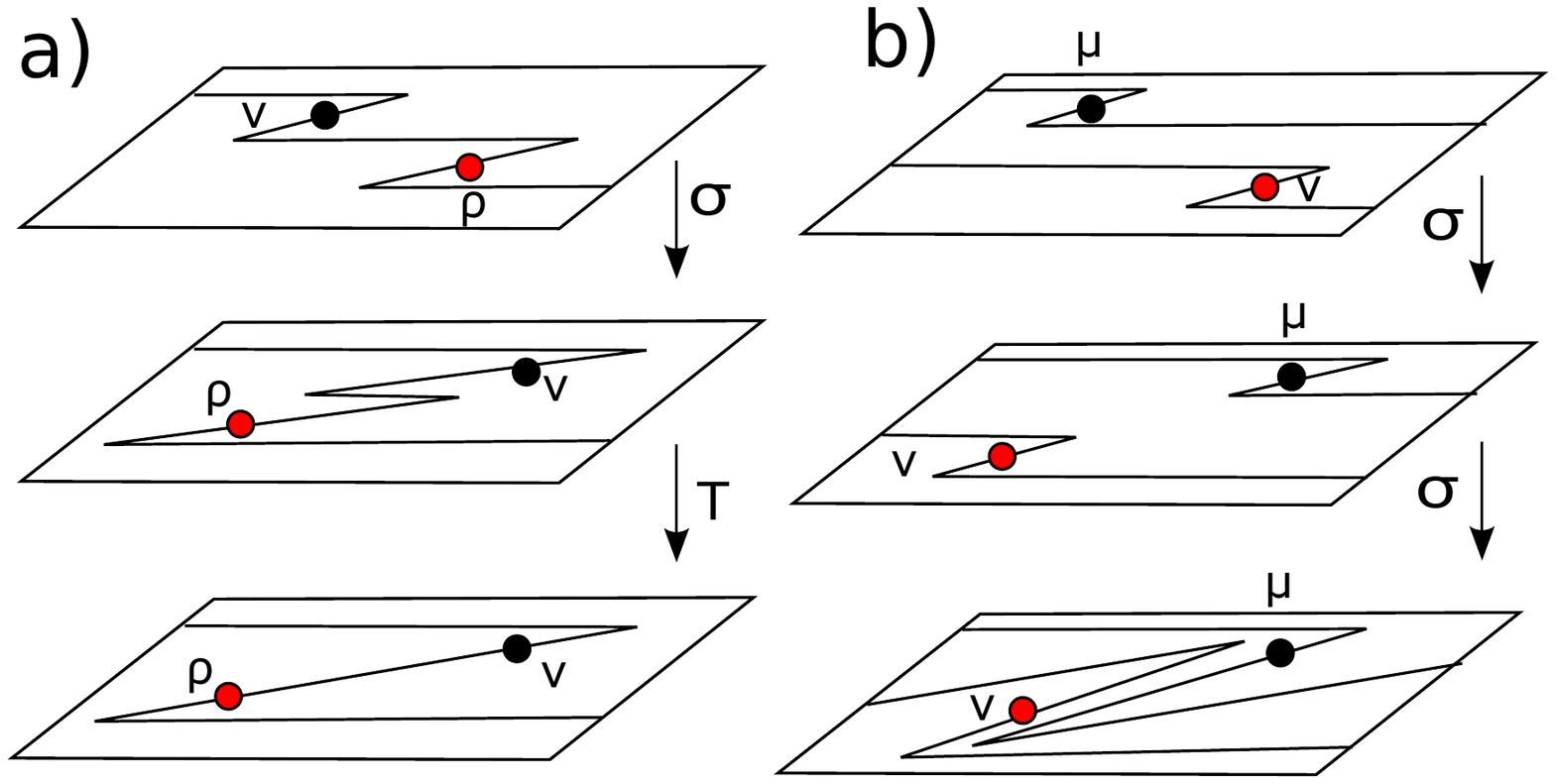}
    \caption{Diagrams in Figure \ref{crossdiags}  c),e) in their movie  representation obtained by taking constant height slices:\newline
      a) 3-morphism $\Phi_{\rho,\nu}\colon\#\rho\circ\!\#\nu\!\Rrightarrow\! \#(\nu\circ \rho)$  (Figure \ref{crossdiags} c).\newline
      b)  3-morphism  $\kappa_{\mu,\nu}\colon \#\nu\Box\#\mu\Rrightarrow \#(\mu\Box\nu)$   (Figure \ref{crossdiags} e).
      Some labels are omitted for legibility.
    }
    \label{crossdiags_slice}
  \end{figure}

  \begin{figure}
    \centering
    \includegraphics[scale=0.65]{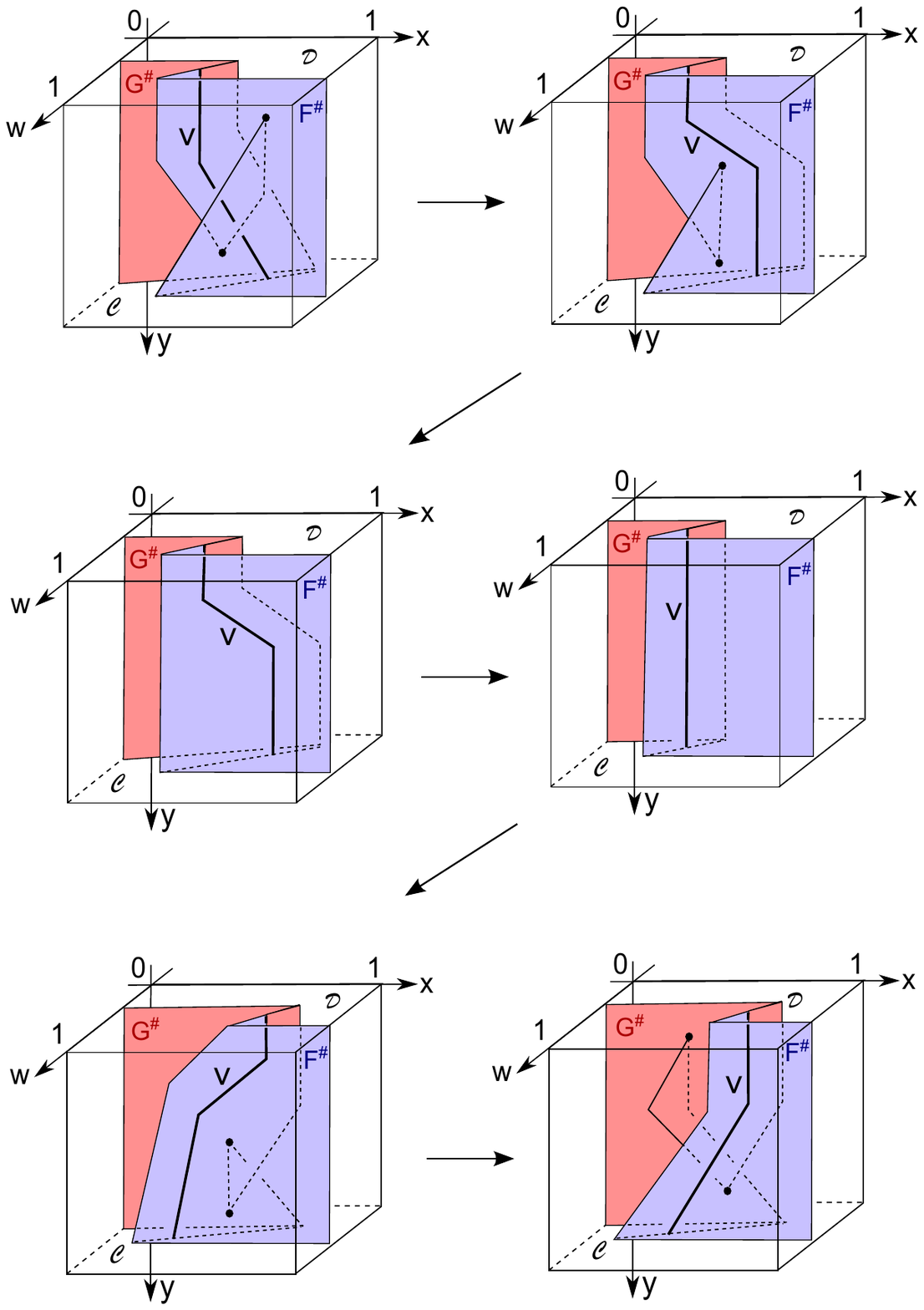}
    \caption{Diagrammatic proof of the identity \newline
      $\Phi_{1_F,\nu}\cdot(\Phi_F\circ 1 _{\#\nu})=1_{\#\nu}= \Phi_{\nu, 1_G}\cdot (1_{\#\nu}\circ \Phi_G)$.  
    }
    \label{P_consist}
  \end{figure}

  \begin{figure}
    \centering
    \includegraphics[scale=0.65]{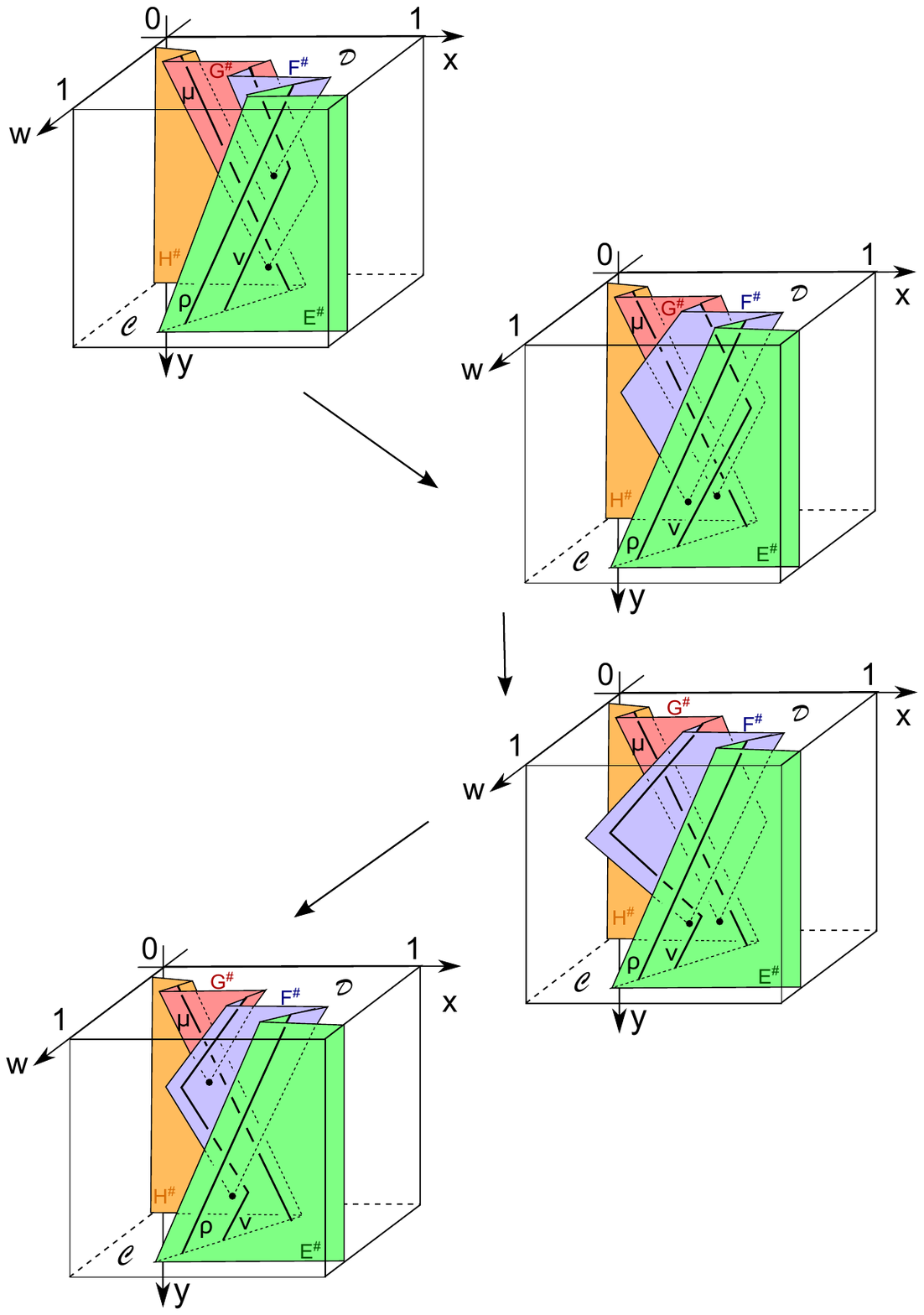}
    \caption{Diagrammatic proof of the identity\newline
      $\Phi_{\nu\circ \rho,\mu}\cdot(\Phi_{\rho,\nu}\circ1_{\#\mu})=\Phi_{\rho,\mu\circ \nu}\cdot (1_{\#\rho}\circ \Phi_{\nu,\mu})$. }
    \label{prhomuconsist}
  \end{figure}

  \begin{figure}
    \centering
    \includegraphics[scale=0.69]{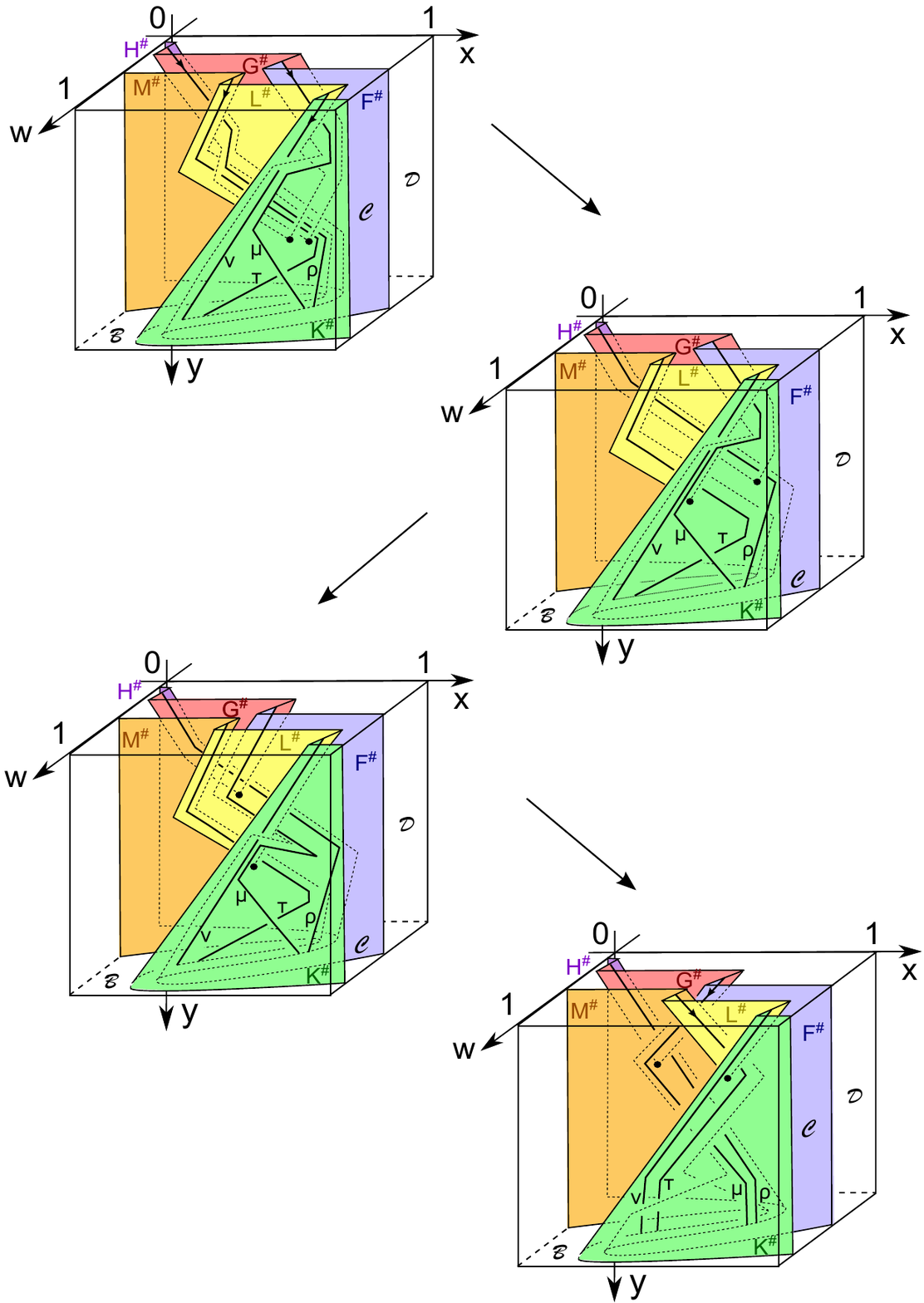}
    \caption{Diagrammatic proof of the identity \eqref{kappanatisom}: \newline
      $\#(1\Box \sigma^\inv_{\mu,\tau}\Box 1)\cdot \Phi_{\mu\Box\nu, \rho\Box\tau}\cdot (\kappa_{\mu,\nu}\circ\kappa_{\rho,\tau})=(\kappa_{\rho\circ\mu,\tau\circ\nu})\cdot (\Phi_{\nu,\tau}\Box\Phi_{\mu,\rho})\cdot \sigma^\inv_{\#\tau,\#\mu}$. }
    \label{kappa_consist0}
  \end{figure}

  \begin{figure}
    \centering
    \includegraphics[scale=0.65]{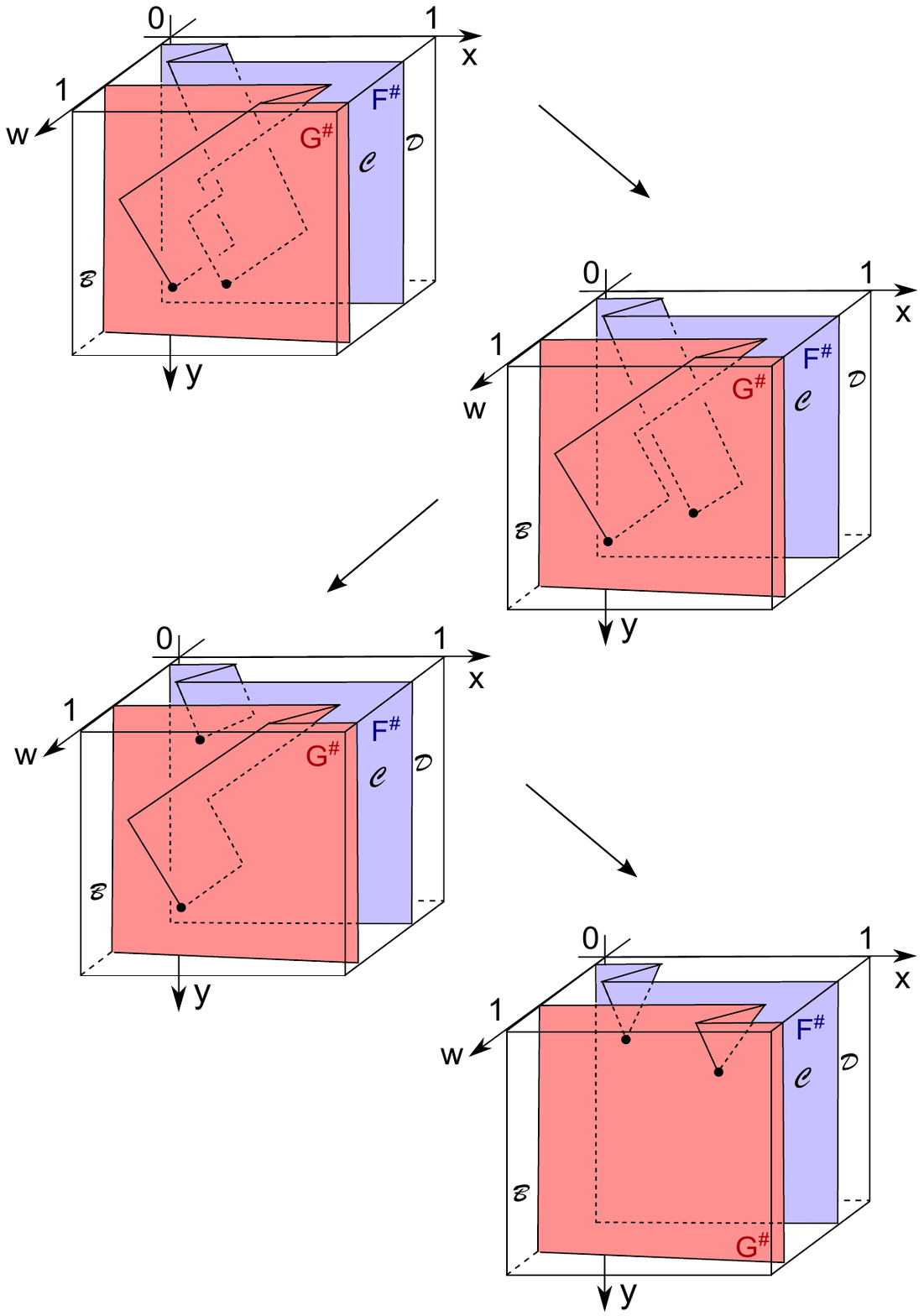}
    \caption{Diagrammatic proof of the identity \eqref{kappaunit}:\newline
      $\Phi_{F\Box G}^\inv\cdot \kappa_{1_G,1_{F}}=\Phi_F^\inv\Box\Phi_G^\inv$. }
    \label{kappa_consist1}
  \end{figure}

  \begin{figure}
    \centering
    \includegraphics[scale=0.68]{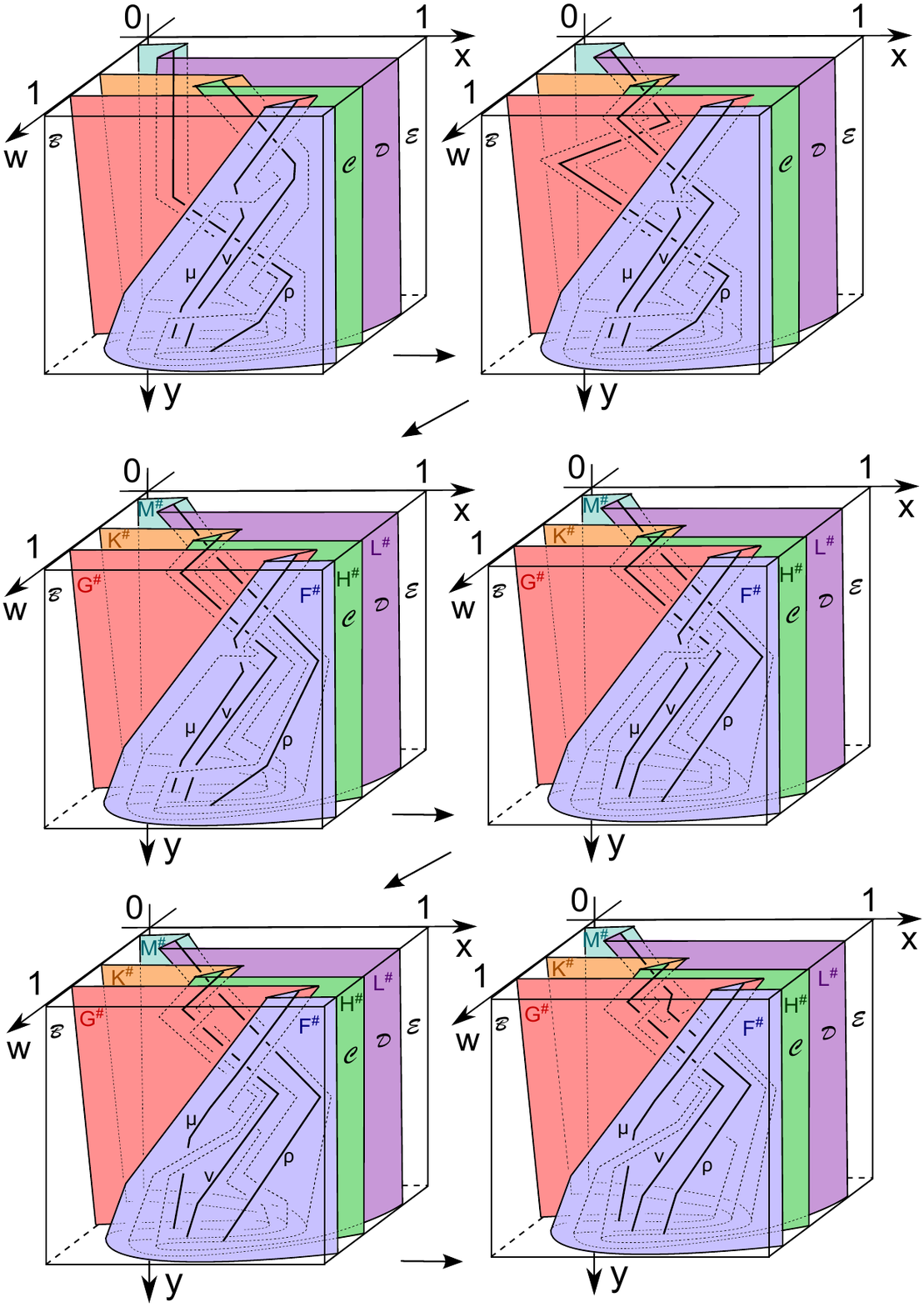}
    \caption{Diagrammatic proof of the identity \eqref{kappa_assoc}: \newline
      $\kappa_{\rho,\nu\Box\mu}\cdot (\kappa_{\nu,\mu}\Box 1)=\kappa_{\rho\Box\nu,\mu}\cdot (1\Box \kappa_{\rho,\nu})$. }
    \label{kappa_consist2}
  \end{figure}

  \medskip\noindent
  3. To show that the four consistency conditions in Definition \ref{grayfunc} are satisfied for the functor $\#$,  note that the operation $\#$ satisfies $1_{\mathcal C}^\#=1_{\mathcal C}$. The natural isomorphisms 
  $\kappa_{\mac B,\mac C,\mac D}\colon \widetilde{\Box}(\#_{\mac C,\mac D}\times \#_{\mac B,\mac C})\to \#_{\mac B,\mac D}\Box$ from Definition \ref{grayfunc}
  are thus specified by their component 3-morphisms 
  $$\kappa_{\mu,\nu}\colon  (\# \mu) \widetilde{\Box} (\# \nu)=(\# \nu)\Box(\# \mu)\Rrightarrow \# (\mu\Box\nu).$$
  These 3-morphisms define natural isomorphisms of weak 2-functors if and only if  they are natural in both arguments, invertible and the following two diagrams commute
  \begin{align}\label{kappanatisom}
    &\xymatrix{ (\#\nu\Box \#\mu)\circ (\#\tau\Box\#\rho) \ar[d]_{1\circ\sigma^\inv_{\#\tau,\#\mu}\circ 1} \ar[r]^{\quad\kappa_{\mu,\nu}\circ \kappa_{\rho,\tau}} & \#(\mu\Box\nu)\circ\#(\rho\Box\tau) \ar[d]^{\Phi_{\mu\Box\nu,\rho\Box\tau}}\\
      (\#\nu\circ\#\tau)\Box(\#\mu\circ\#\rho) \ar[d]_{\Phi_{\nu,\tau}\Box\Phi_{\mu,\rho}} & \#((\rho\Box\tau)\circ(\mu\Box\nu))\ar[d]^{\#(1\circ\sigma^\inv_{\mu,\tau}\circ 1)}\\
      \#(\tau\circ\nu)\Box\#(\rho\circ\mu) \ar[r]_{\kappa_{\rho\circ\mu, \tau\circ\nu}}& \#((\rho\circ\mu)\Box(\tau\circ\nu))
    }\\[+2ex]
    &\xymatrix{1_{G^\#}\Box 1_{F^\#}=1_{(F\Box G)^\#} \ar[d]_{\Phi_G\Box \Phi_F} \ar[rd]^{\Phi_{F\Box G}}\\
      \#1_G\Box \#1_F \ar[r]_{\kappa_{1_F,1_G}\qquad} & \#(1_F\Box 1_G)=\#(1_{F\Box G}).
    }\label{kappaunit},
  \end{align}
  where  the two vertical arrows labelled with tensorators arise from the definition of the Gray products in  Lemma \ref{cubegray} and Corollary   \ref{lemma:cubical-to-opcubical-tricat}.

  Condition (1) in Definition \ref{grayfunc} is equivalent to the commutativity of diagram
  \begin{align}\label{kappa_assoc}
    \xymatrix{
      \#\mu\Box\#\nu\Box\#\rho \ar[d]_{\kappa_{\nu,\mu}\Box 1}\ar[r]^{1\Box\kappa_{\rho,\nu}} & \#\mu\Box \#(\rho\Box\nu) \ar[d]^{\kappa_{\rho\Box\nu, \mu}}\\
      \#(\nu\Box \mu)\Box\#\rho  \ar[r]_{\kappa_{\rho,\nu\Box\mu}}  & \#(\rho\Box\nu\Box\mu),
    }
  \end{align}
  and conditions (2), (3) to the equations
  $
  \kappa_{\mu,1_{1_\mathcal C}}=\kappa_{1_{1_\mathcal D},\mu}=1_\mu
  $ for all 2-morphisms $\mu \in \mac G( \mac C, \mac D)$.

  \medskip
  For $\Box$-composable  2-morphisms $\mu\colon F\Rightarrow G$, $\nu\colon H\Rightarrow K$,  define a  3-morphism $\kappa_{\mu,\nu}\colon \#\nu\Box\#\mu\Rrightarrow\#(\mu\Box \nu)$  by
  \begin{align*}
    &\kappa_{\mu,\nu}=\\
    &(1_{H^\#F^\#\eta_G^*}\circ 1_{H^\#F^\#G\eta_K^*G^\#}\circ 1_{H^\#F^\#\mu KK^\#G^\#}\circ H^\#\sigma_{\eta_{F^\#},\nu}K^\#G^\#\circ 1_{\eta_{H^\#}K^\#G^\#})\nonumber\\
    &\cdot (1_{H^\#F^\#\eta_G^*}\circ H^\#\sigma_{F^\#\mu\circ \eta_{F^\#}, \eta_K^*}G^\#\circ 1_{H^\#\nu K^\#G^\#}\circ 1_{\eta_{H^\#}K^\#G^\#})\cdot \sigma_{\#\nu,\#\mu}\nonumber.
  \end{align*}
  This 3-morphism and its projection are shown in Figure \ref{crossdiags} e), f) and Figure \ref{crossdiags_slice} b). It follows directly from the definition of the tensorator that the 3-morphisms $\kappa_{\mu,\nu}$ are invertible, satisfy the conditions $
  \kappa_{\mu,1_{1_\mathcal C}}=\kappa_{ 1_{1_\mathcal D},\mu}=1_\mu
  $ and are natural in both arguments. 
  It is therefore sufficient to establish the commutativity of the diagrams in \eqref{kappanatisom}, \eqref{kappaunit} and \eqref{kappa_assoc}.
  A diagrammatic proof of these identities  is given, respectively, in Figures \ref{kappa_consist0}, \ref{kappa_consist1} and \ref{kappa_consist2}.  
\end{proof}

\medskip\noindent

\begin{lemma}
  \label{lemma:doublestar-is-one}
  The functor of 2-strict tricategories $*\colon\mathcal{G} \rightarrow \mathcal{G}^{op}$ satisfies $**=1$.
\end{lemma}
\begin{proof}
  Note that strictly speaking this identity should be written $*^{op} * =1_{\mathcal{G}}$,  identifying  $(\mathcal{G}^{op})^{op}=\mathcal{G}$.
  That the mappings of the functor $**$ are given by the identity follows directly from the fact that $*$ is trivial on the objects and 1-morphisms and the 2-categories $\mathcal G(\mathcal C,\mathcal D)$ are planar.
  It remains to show that the coherence morphisms of $**$ are the identities. Recall that the coherence data of $*$ is given by $\sigma$, i.e., the components of the natural transformation
  $\rho_{\mac B,\mac C,\mac D}\colon \Box^{op}\,( *_{\mac C,\mac D}\times  *_{\mac B,\mac C})\to  *_{\mac B,\mac D}\,\Box$ are $\rho_{\mu,\nu}=(\sigma_{\mu^*,\nu^*})^{op}$. According to the definition of $F^{op}$ for a functor $F$ of 2-strict tricategories,
  the components of the corresponding transformation for $*^{op}$ are given by $\sigma_{\mu^*,\nu^*}^{-1}$ and hence the composition $**$ is the identity functor.   
\end{proof}

Theorem \ref{graycatduals} gives a more conceptual understanding of the duals in terms of functors of 2-strict tricategories rather than the concrete axioms in Definition \ref{dualgray}. These functors of 2-strict tricategories are related to certain symmetries of the cube. The $*$-dual does not affect the 1-morphisms in $\mac G$ and corresponds to a $180$ degree rotation around the $w$-axis of the diagrams, and the  dual $\#$ corresponds to a $180$ degree  rotation around the $y$-axis.

Lemma \ref{lemma:doublestar-is-one} states that  the operation $*$ satisfies $**=1$, as expected for the composite of two $180$ degree rotations around the same axis. However  $\#\#$ and $*\#{*}\#$ are not equal to the identity, although the associated compositions of rotations are. A rationale for this is that in a higher category, one can in general only  expect such relations to hold up to higher morphisms. In the case at hand, one obtains natural isomorphisms of functors of 2-strict tricategories (Definition \ref{graynat}). 

\begin{figure}
  \centering
  \includegraphics[scale=0.65]{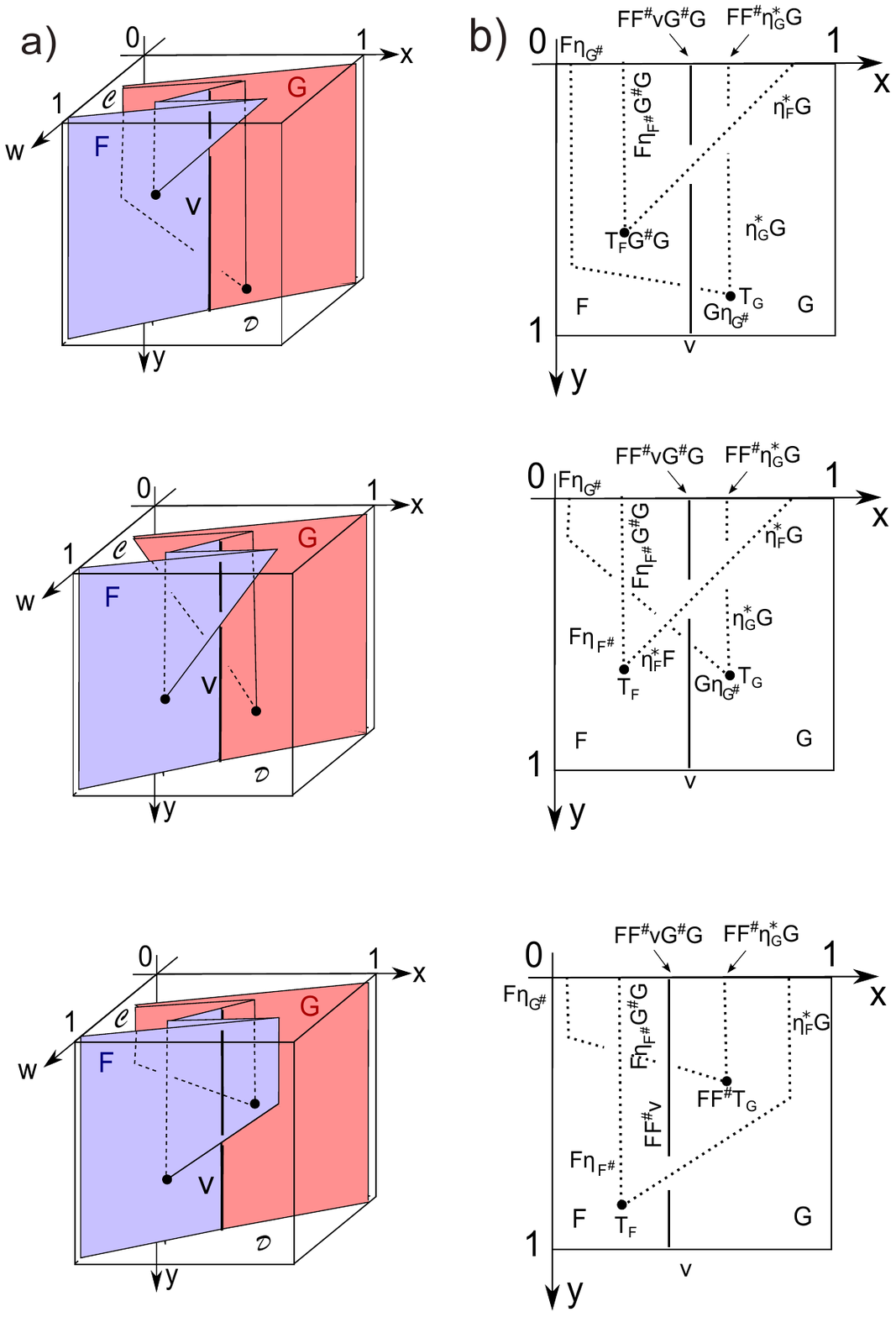}
  \caption{3-morphism $\Gamma_\nu\colon *\#{*}\#\nu\Rrightarrow \nu$ a) and b) its projection. The three Gray category diagrams have the same evaluation.}
  \label{gammanu}
\end{figure}

\begin{figure}
  \centering
  \includegraphics[scale=0.4]{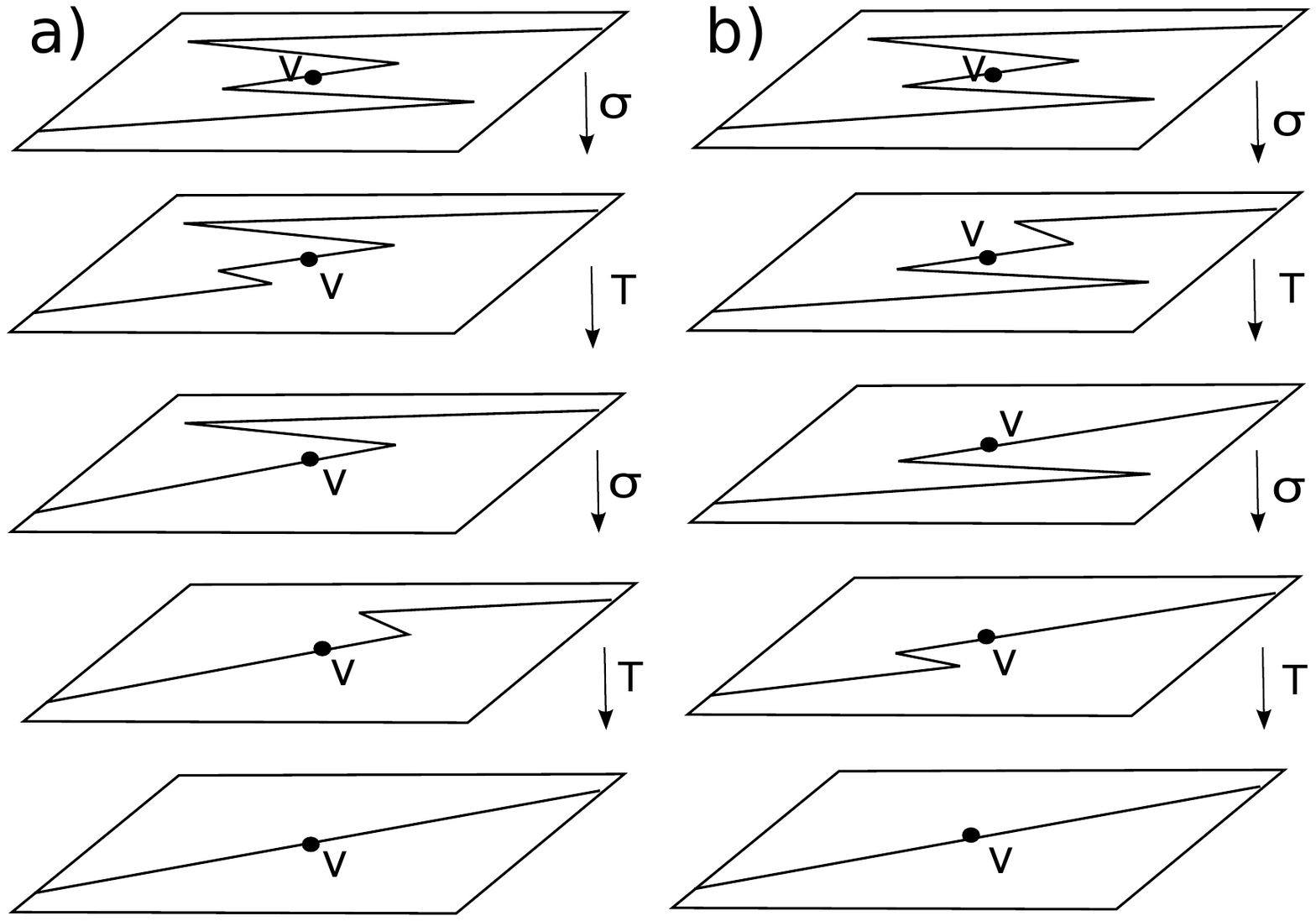}
  \caption{Gray category diagrams for the 3-morphism $\Gamma_\nu\colon *\#{*}\#\nu\Rrightarrow \nu$ in their movie representation, obtained from Figure \ref{gammanu} a) and c) by taking constant height slices. Some labels are omitted for legibility.}
  \label{gammanu_slice}
\end{figure}

\begin{figure}
  \centering
  \includegraphics[scale=0.5]{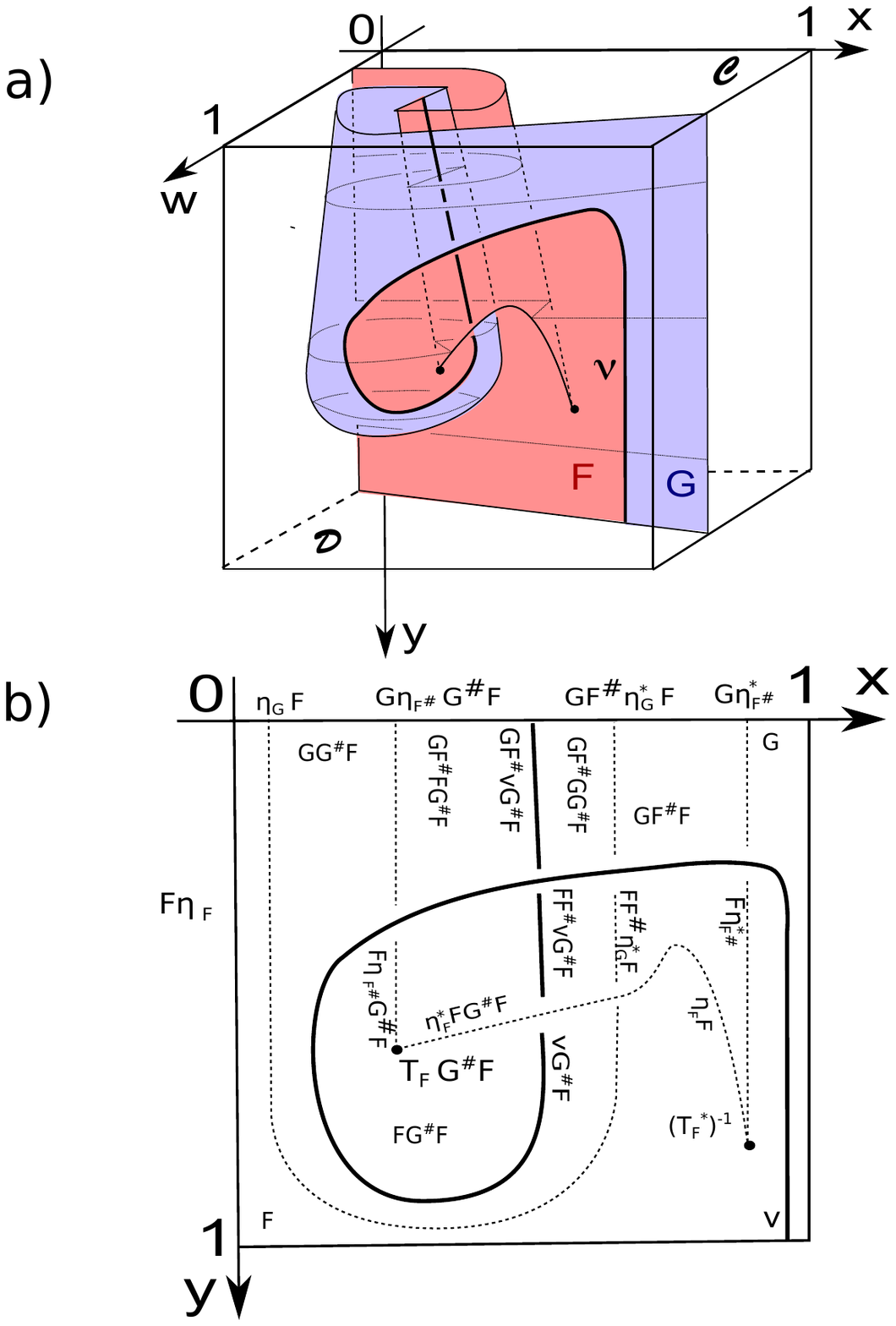}
  \caption{a) Gray category diagram for the 3-morphism $\Theta_\nu\colon \#\#\nu\Rrightarrow\nu$ and b) its projection.
  }
  \label{twist_simple}
\end{figure}

\begin{figure}
  \centering
  \includegraphics[scale=0.45]{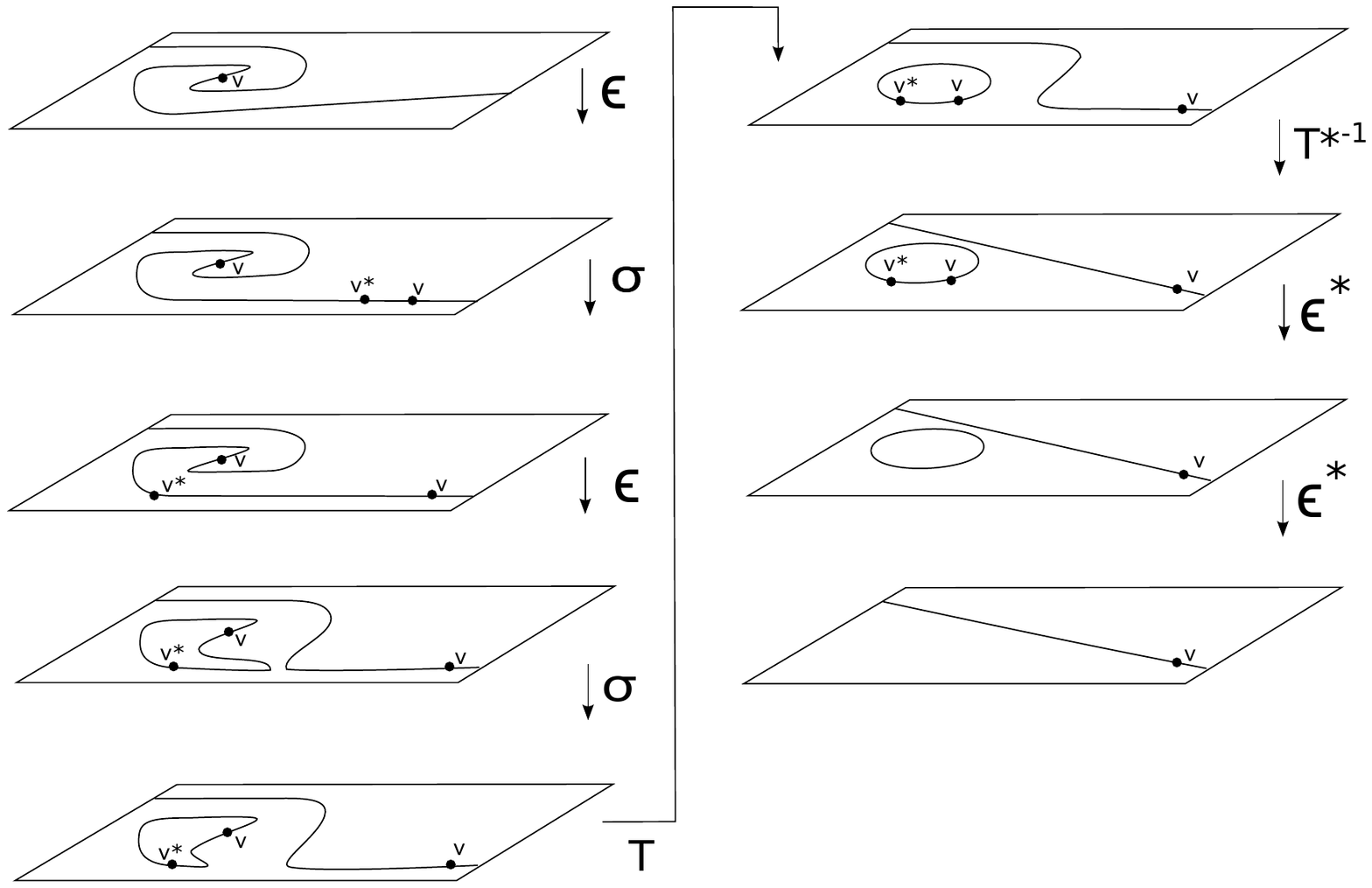}
  \caption{a) Gray category diagram for the 3-morphism $\Theta_\nu\colon \#\#\nu\Rrightarrow\nu$ in movie representation, obtained by taking constant height slices of the diagram in Figure \ref{twist_simple} a). Some labels are omitted for legibility.
  }
  \label{twist_slice}
\end{figure}

As the functors $*\#{*}\#*$ and  $\#\#$ of 2-strict tricategories  act trivially on the objects and 1-morphisms of $\mathcal G$,  natural isomorphisms $\Theta\colon \#\#\to 1$ and $\Gamma\colon *\#{*}\#\to 1$ are determined  by the collection of their component  3-morphisms $\Gamma_\nu\colon *\#{*}\#\nu\Rrightarrow \nu$ and $\Theta_\nu\colon \#\#\nu\Rrightarrow\nu$ for each 2-morphism $\nu$.  For a  2-morphism $\nu\colon F\Rightarrow G$, we define
\begin{align}\label{eq:gammanudef}
  \Gamma_{\nu}=&(T_G\circ 1 _\nu)\cdot (1_{\eta_G^* G}\circ \sigma_{\nu,\eta_{G^\#}})\cdot (1_{\eta_G^*G\circ \nu G^\#G}\circ T_FG^\#G\circ 1 _{F\eta_{G^\#}})\\
  \cdot&(\sigma_{\eta_{F}^*, \eta_G^*G\circ \nu G^\#G} \circ 1_{F\eta_{F^\#}G^\#G\circ F\eta_{G^\#}}),\nonumber
\end{align}
\begin{align}
  \label{eq:thetanudef}
  & \Theta_\nu=(1_\nu\circ \epsilon^*_{\eta_G^*F})\cdot (1_{\nu\circ \eta_G^*F}\circ \epsilon^*_{ \nu G^{\#}F}\circ 1_{\eta_GF})\\
  &\cdot (1_\nu\circ (T^*_F)^\inv\circ 1_{\eta_G^*F\circ\nu G^{\#}F\circ \nu^*G^{\#}F\circ \eta_G F})\nonumber\\
  &  \cdot (1_{\nu\circ F\eta_{F^\#}^*\circ\eta_FF\circ \eta^*_G F\circ \nu G^{\#}F} \circ T_FG^{\#}F\circ 1_{\nu^*G^{\#}F\circ\eta_GF})\nonumber\\
  &\cdot (1_{\nu\circ F\eta_{F^{\#}}^*\circ \eta_FF}\circ \sigma_{\eta_F^*F, FF^\#\eta_G^*F\circ FF^\#\nu G^\#F }\circ 1_{F\eta_{F^\#}G^\#F\circ\nu^*G^\#F\circ\eta_GF})\nonumber\\
  &\cdot (1_{\nu\circ F\eta^*_{F^\#}}\circ \epsilon_{\eta_FF}\circ 1_{F\#\nu F\circ\nu^*G^\#F\circ \eta_GF})\cdot (1_\nu\circ \sigma_{\nu^*, \eta_{F^\#}^*\circ\ \#\nu F}\circ 1_{\eta_GF})\nonumber\\
  &\cdot (\epsilon_\nu\circ 1_{\#\#\nu}).\nonumber
\end{align}
The Gray category diagram  for the 3-morphism $\Gamma_\nu\colon *\#{*}\#*\nu\Rrightarrow\nu$,  its projection and its movie representation are given in  Figures \ref{gammanu} and \ref{gammanu_slice},
the corresponding diagrams for  the 3-morphism $\Theta_\nu\colon\#\#\nu\Rrightarrow\nu$ in Figures \ref{twist_simple} and \ref{twist_slice}.

\begin{theorem} \label{graydualnat}The 3-morphisms $\Gamma_\nu\colon *\#{*}\#\nu\Rrightarrow \nu$  and $\Theta_\nu\colon \#\#\nu\Rrightarrow\nu$ in \eqref{eq:gammanudef} and \eqref{eq:thetanudef} define natural isomorphisms of functors of 2-strict tricategories $\Gamma\colon *\#{*}\#\rightarrow 1$, $\Theta\colon \#\#\rightarrow 1$.  
\end{theorem}

\begin{figure}
  \centering
  \includegraphics[scale=0.4]{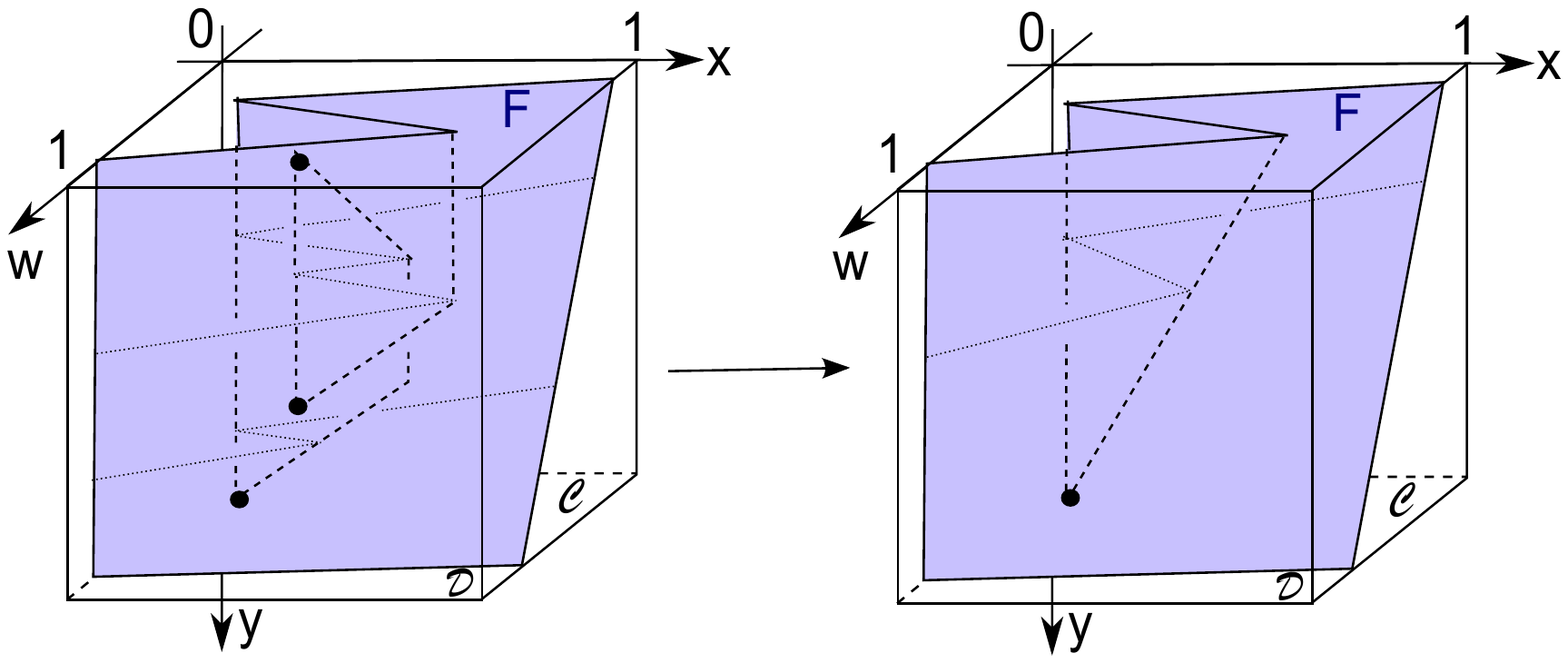}
  \caption{$\quad$\newline Diagrammatic proof of  identity \eqref{gamuncomp}: \newline $\Gamma_{1_F}\cdot (*\#*\Phi_F)=*\Phi_{F^\#}$.}
  \label{gammaid}
\end{figure}

\begin{figure}
  \centering
  \includegraphics[scale=0.67]{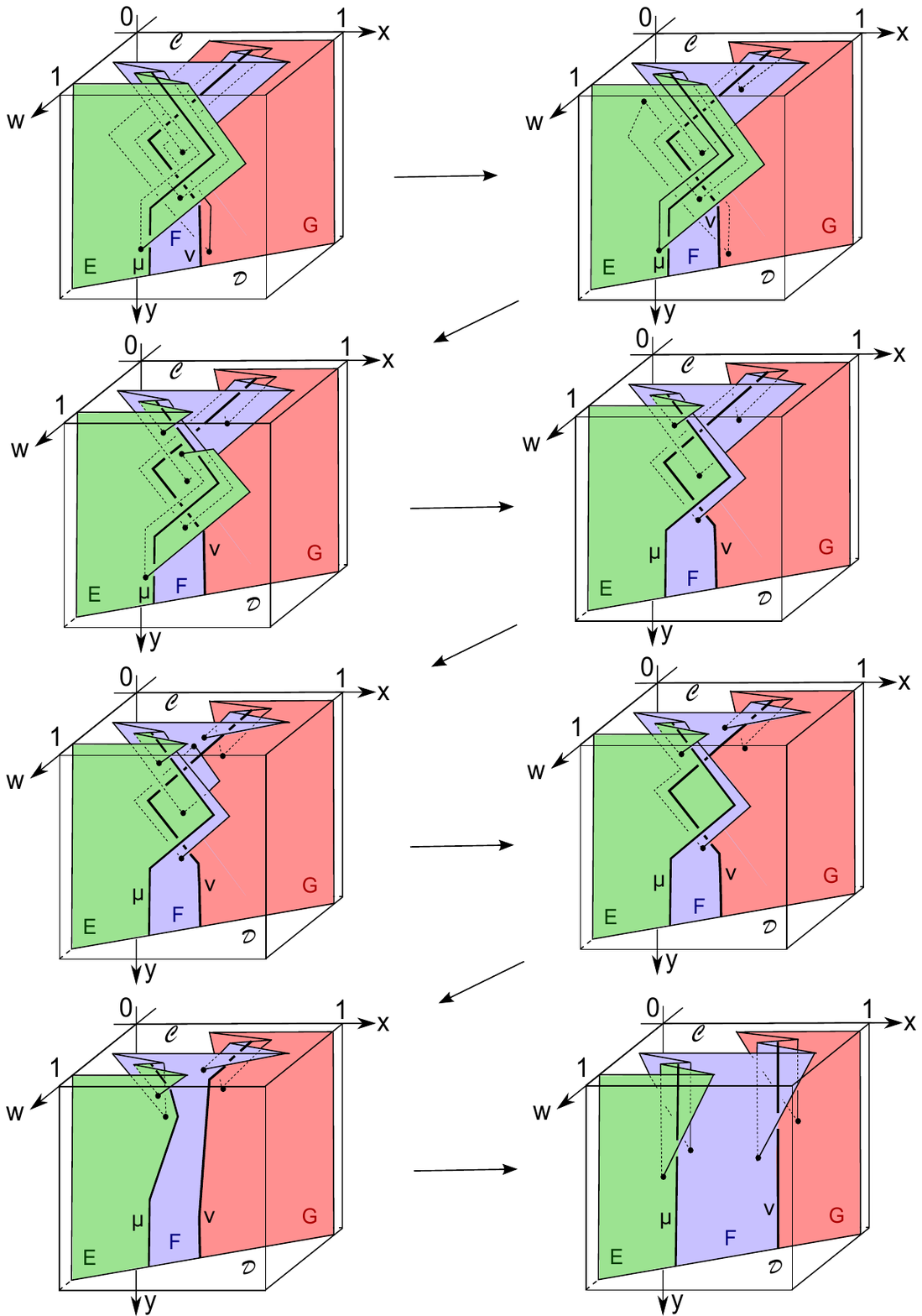}
  \caption{Diagrammatic proof of  identity \eqref{gamhorcompeq}: \newline $\Gamma_{\mu\circ\nu}\cdot (*\#*\Phi_{\mu,\nu})^{-1}\cdot (*\Phi_{*\#\nu,*\#\mu})^{-1}=\Gamma_\mu\circ\Gamma_\nu$.}
  \label{gammahorcomp}
\end{figure}

\begin{figure}
  \centering
  \includegraphics[scale=0.67]{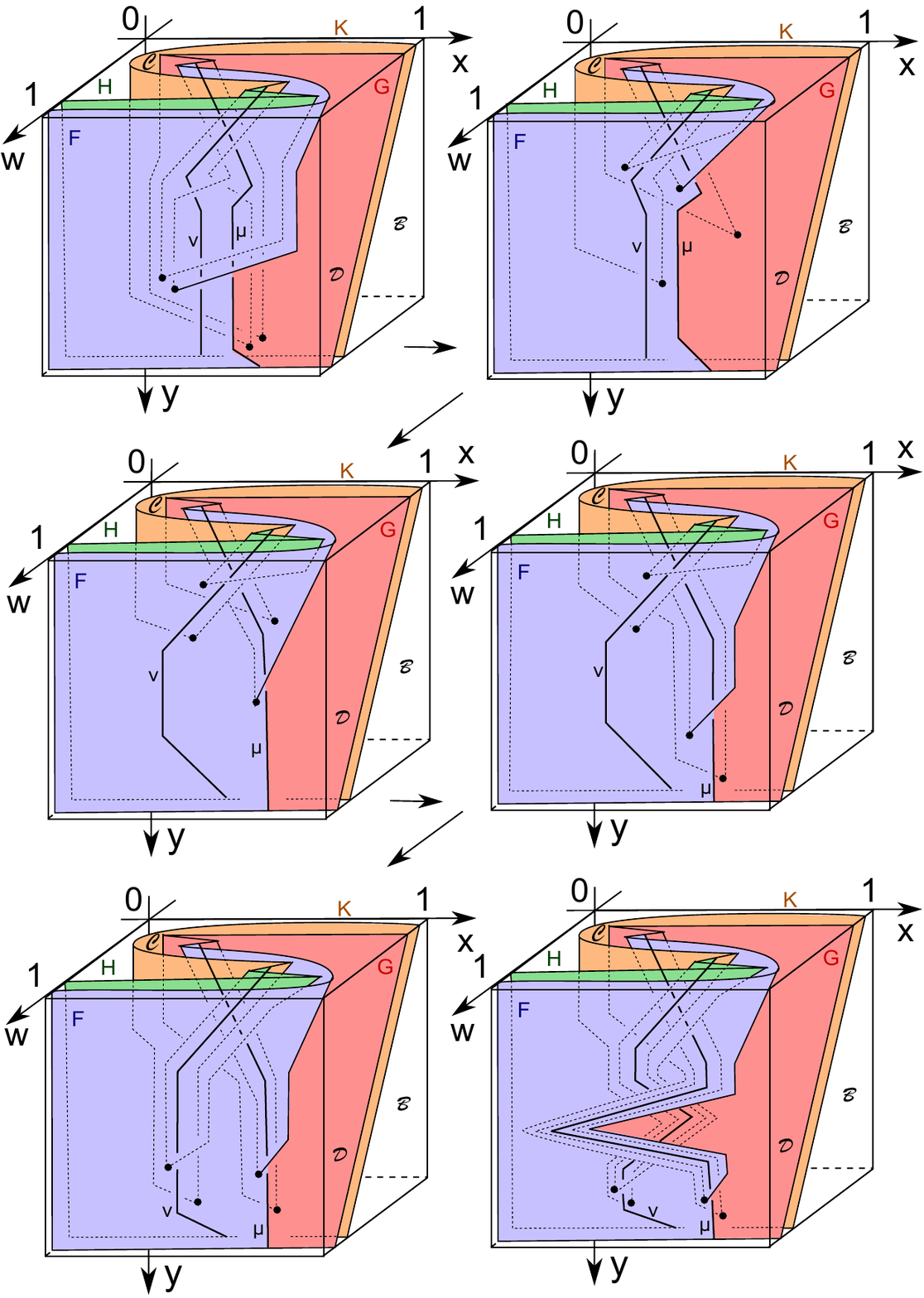}
  \caption{Diagrammatic proof of  identity \eqref{gamgraycomp}: \newline $\Gamma_{\mu\Box \nu}\cdot (*\#*\kappa_{\mu,\nu})=(\Gamma_\mu\Box\Gamma_\nu)\cdot \sigma^\inv_{*\#{*}\#\mu,*\#{*}\#\nu}\cdot *\kappa_{ *\#\nu,*\#\mu}\cdot *\#\sigma_{*\#\nu,*\#\mu}$.}
  \label{gammabox}
\end{figure}

\begin{proof} $\quad$ \newline
  1. As the functor of 2-strict tricategories $*\#{*}\#\colon\mathcal G\rightarrow\mathcal G$ acts trivially on objects and 1-morphisms,  it is sufficient to show that the 3-morphisms $\Gamma_\nu\colon *\#{*}\#\nu\Rrightarrow \nu$ are invertible for each 2-morphism $\nu$ and satisfy the following conditions:
  \begin{itemize}
  \item {\em naturality:} for each 3-morphism $\Psi\colon \mu\Rrightarrow\nu$ the following diagram commutes
    $$
    \xymatrix{ {*}\#{*}\#\mu \ar[r]^{\quad\Gamma_\mu} \ar[d]_{*\#{*}\#\Psi} & \mu \ar[d]^{\Psi}\\
      {*\#{*}\#} \nu \ar[r]_{\quad\Gamma_\nu} & \nu.
    }
    $$

  \item {\em compatibility  with the unit 2-morphisms:} for all 1-morphisms $F$, the following diagram commutes
    \begin{align}\label{gamuncomp}
      \xymatrix{ {*} {\#} * {\#} 
        1_F  \ar[rd]_{\Gamma_{1_F}} & \# 1_{F^\#} \ar[l]_{\quad\quad*\#* \Phi_F} \ar[d]^{* \Phi_{F^\#}}\\
        & 1_F. 
      }
    \end{align}

  \item {\em compatibility with the horizontal composition:} for all composable 2-morphisms $\mu$, $\nu$, the following diagram commutes
    \begin{align}\label{gamhorcompeq}
      \xymatrix{( {*}\#*\#\mu)\circ (* \#*\#\nu)  \ar[d]^{\Gamma_\mu\circ \Gamma_\nu} & \quad {*}\#*(\#\nu\circ \#\mu) \quad \ar[l]_{\quad  {*}\Phi_{{*}\#\nu,{*}\#\mu}} \ar[d]^{{*}\#*\Phi_{\mu,\nu}}\\
        \mu\circ \nu & {*}\#*\#(\mu\circ \nu) \ar[l]_{\Gamma_{\mu\circ \nu}}.
      }
    \end{align}

  \item
    {\em compatibility with the Gray product:} for all composable 2-morphisms \linebreak $\mu\colon F\Rightarrow G$, $\nu\colon H\Rightarrow K$, the following diagram commutes
    \begin{align}\label{gamgraycomp}
      \xymatrix{ {*}\#{*}\#\mu\Box {*}\#{*}\#\nu \ar[r]^{\qquad\Gamma_\mu\Box\Gamma_\nu} & \mu\Box\nu\\
        {*}(\#{*}\#\mu\Box \#{*}\#\nu) \ar[u]^{\sigma^\inv_{*\#{*}\#\mu,*\#{*}\#\nu}} & {*}\#{*}\#(\mu\Box\nu)\ar[u]_{\Gamma_{\mu\Box\nu}}\\
        {*\#(*\#\nu\Box *\#\mu)}\quad \ar[u]^{*\kappa_{*\#\nu,*\#\mu}}& {\quad *}\#{*}(\#\nu\Box\#\mu) \ar[l]^{*\#\sigma_{*\#\nu,*\#\mu}}  \ar[u]_{*\#*\kappa_{\mu,\nu}}.
      }
    \end{align}
  \end{itemize}
  Note that the four arrows labelled with the invertible 3-morphisms $*\#*\kappa_{\mu,\nu}$, $*\#\sigma_{*\#\nu, *\#\mu}$, $*\kappa_{*\#\nu, *\#\mu}$, $\sigma^\inv_{*\#{*}\#\mu, *\#{*}\#\nu}$ in this diagram compose to the coherence 3-morphism $\kappa^{*\#{*}\#}_{\mu,\nu}\colon \#*\#*\mu{\Box}*\#{*}\#\nu\Rrightarrow *\#{*}\#(\mu\Box\nu)$ of the functor $*\#{*}\#\colon\mac G\to\mac G$.

  That $\Gamma_\nu\colon *\#{*}\#\nu\Rrightarrow \nu$ is invertible for all 2-morphisms $\nu\colon F\Rightarrow G$ follows directly from its definition in equation \eqref{eq:gammanudef} and Figure \ref{gammanu}. The naturality of $\Gamma_\nu$ follows directly from the naturality of the tensorator.

  To show the  compatibility of $\Gamma_\nu$ with the unit 2-morphisms, recall that for each 1-morphism $F$, the tensorator $\sigma_{1_F, \eta_{F^\#}}$ is trivial. The associated 3-morphism $\Gamma_{1_F}$ therefore reduces to:
  \begin{align*}
    \Gamma_{1_F}=&T_F\cdot (1_{\eta_F^* F}\circ  T_FF^\#F\circ 1 _{F\eta_{F^\#}})
    \cdot(\sigma_{\eta_{F}^*, \eta_F^*F} \circ 1_{F\eta_{F^\#}F^\#F\circ F\eta_F^\#}).
  \end{align*}
  Composing this 3-morphism with the 3-morphism $*\#{*}\Phi_F=*\# T_{F^\#}$, one obtains the Gray category diagram in Figure \ref{gammaid}. The commutativity of the diagram 
  \eqref{gamuncomp} is then a direct consequence of identity (d) in Definition \ref{dualgray}. A diagrammatic proof is given in  Figure \ref{gammaid}.

  The compatibility condition  \eqref{gamhorcompeq} between  $\Gamma$ and  the horizontal composition follows from the definitions together with the invertibility of the triangulator, identity (d) in Definition \ref{dualgray}, the naturality of the tensorator and the exchange identity. 
  As the calculations are lengthy and technical, we give a diagrammatic proof  in Figure \ref{gammahorcomp}.

  The compatibility of $\Gamma$ with the Gray product  in equation \eqref{gamgraycomp} again follows from the definitions together with the naturality of the tensorator, the properties of the triangulator and the exchange identity. A diagrammatic proof is given in Figure \ref{gammabox}. This concludes the proof that the 3-morphisms $\Gamma_\nu\colon *\#{*}\#\nu\Rrightarrow\nu$ define a natural isomorphism of functors of 2-strict tricategories $\Gamma\colon *\#{*}\#\rightarrow 1$.

  \bigskip\bigskip\noindent
  2. As the functor of 2-strict tricategories $\#\#\colon\mathcal G\rightarrow\mathcal G$  acts trivially on the objects and 1-morphisms of $\mathcal G$, it is sufficient to show that the 3-morphisms $\Theta_\nu\colon \#\#\nu\Rrightarrow\nu$ are invertible for each 2-morphism $\nu\colon F\rightarrow G$ and  satisfy the following conditions:
  \begin{itemize}
  \item {\em naturality:} for each 3-morphism $\Psi\colon \mu\Rrightarrow\nu$ the following diagram commutes:
    $$
    \xymatrix{ \#\#\mu \ar[r]^{\Theta_\mu} \ar[d]_{\#\#\Psi} & \mu\ar[d]^\Psi\\
      \#\#\nu \ar[r]_{\Theta_\nu} &\nu.
    }
    $$
  \item  {\em compatibility with the unit 2-morphisms:} for all 1-morphisms $F$ the following diagram commutes:
    \begin{align}\label{twistun}
      \xymatrix{ \#\# 1_F   \ar[rd]_{\Theta_{1_F}} \ar[r]^{(\#\Phi_F)^\inv} & \# 1_{F^\#}  \ar[d]^{\Phi^\inv_{F^\#}} \\
        & 1_F
      }
    \end{align}

  \item {\em compatibility with the horizontal composition:} for all composable 2-morphisms $\mu$, $\nu$, the following diagram commutes:
    \begin{align}\label{twisthorcomp}
      \xymatrix{ \#\#\mu\circ \#\#\nu \ar[r]^{\Phi_{\#\mu,\#\nu}} \ar[d]^{\Theta_\mu\circ \Theta_\nu} & \#(\#\nu\circ \#\mu) \ar[d]^{\#\Phi_{\nu,\mu}}\\
        \mu\circ \nu & \#\#(\mu\circ \nu) \ar[l]_{\Theta_{\mu\circ \nu}}
      }
    \end{align}

  \item
    {\em compatibility with the Gray product:} for all  composable 2-morphisms $\mu\colon F\Rightarrow G$, $\nu\colon H\Rightarrow K$  the following diagram commutes:
    \begin{align}\label{twistbox}
      \xymatrix{\#\#(\mu\Box\nu) \ar[r]^{\Theta_{\mu\Box\nu}} & \mu\Box\nu\\
        \#(\#\nu\Box\#\mu) \ar[u]_{\#\kappa_{\mu,\nu}} & \#\#\mu\Box\#\#\nu \ar[u]_{\Theta_\mu\Box\Theta_\nu} \ar[l]_{\kappa_{\#\nu,\#\mu}}.
      }
    \end{align}
  \end{itemize}

\begin{figure}
  \centering
  \includegraphics[scale=0.3]{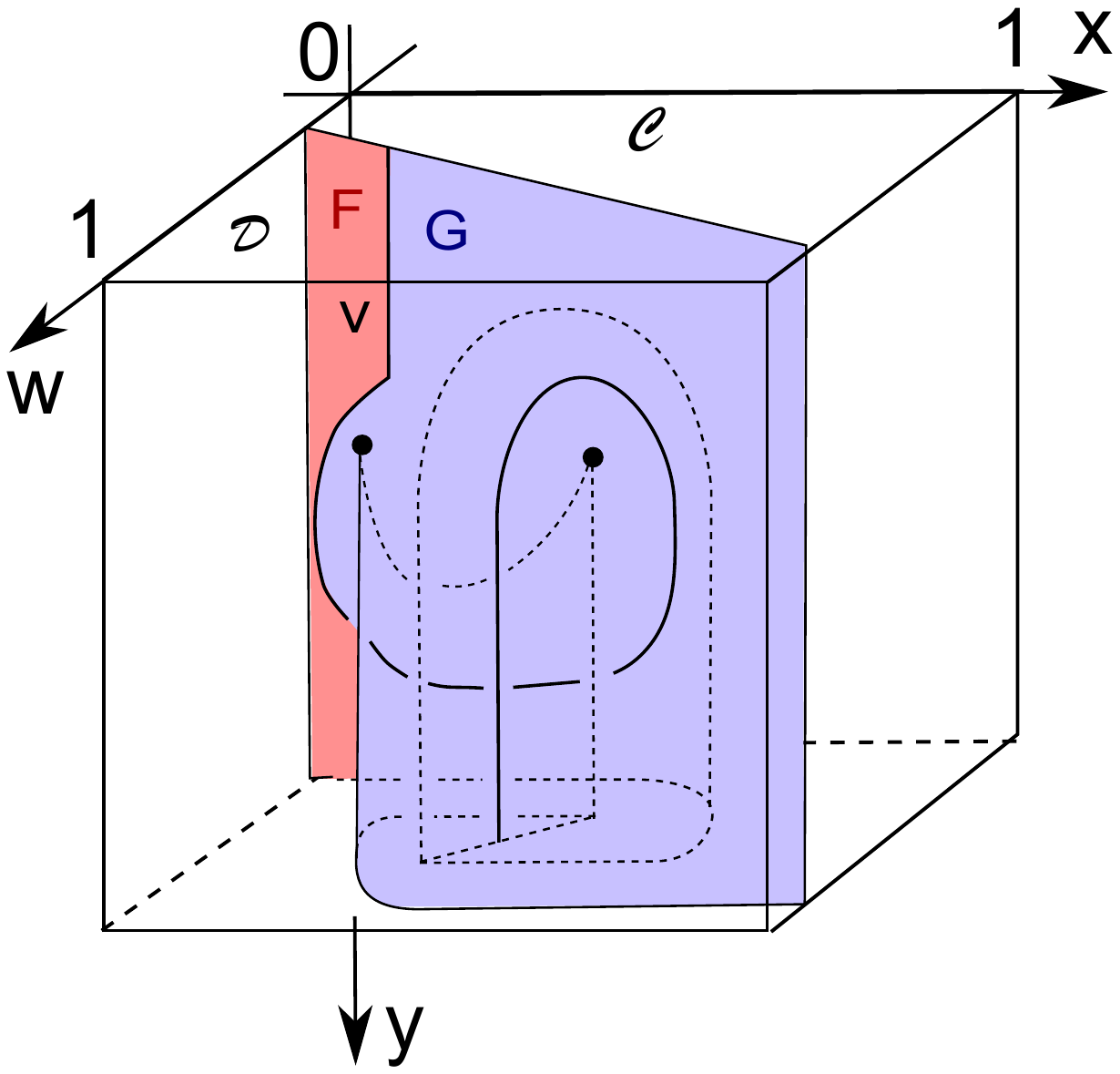}
  \caption{Inverse $\Theta_\nu^{-1}\colon \nu\Rrightarrow\#\#\nu$ of the 3-morphism $\Theta_\nu\colon\#\#\nu\Rrightarrow\nu$.}
  \label{twistinv}
\end{figure}

\begin{figure}
  \centering
  \includegraphics[scale=0.65]{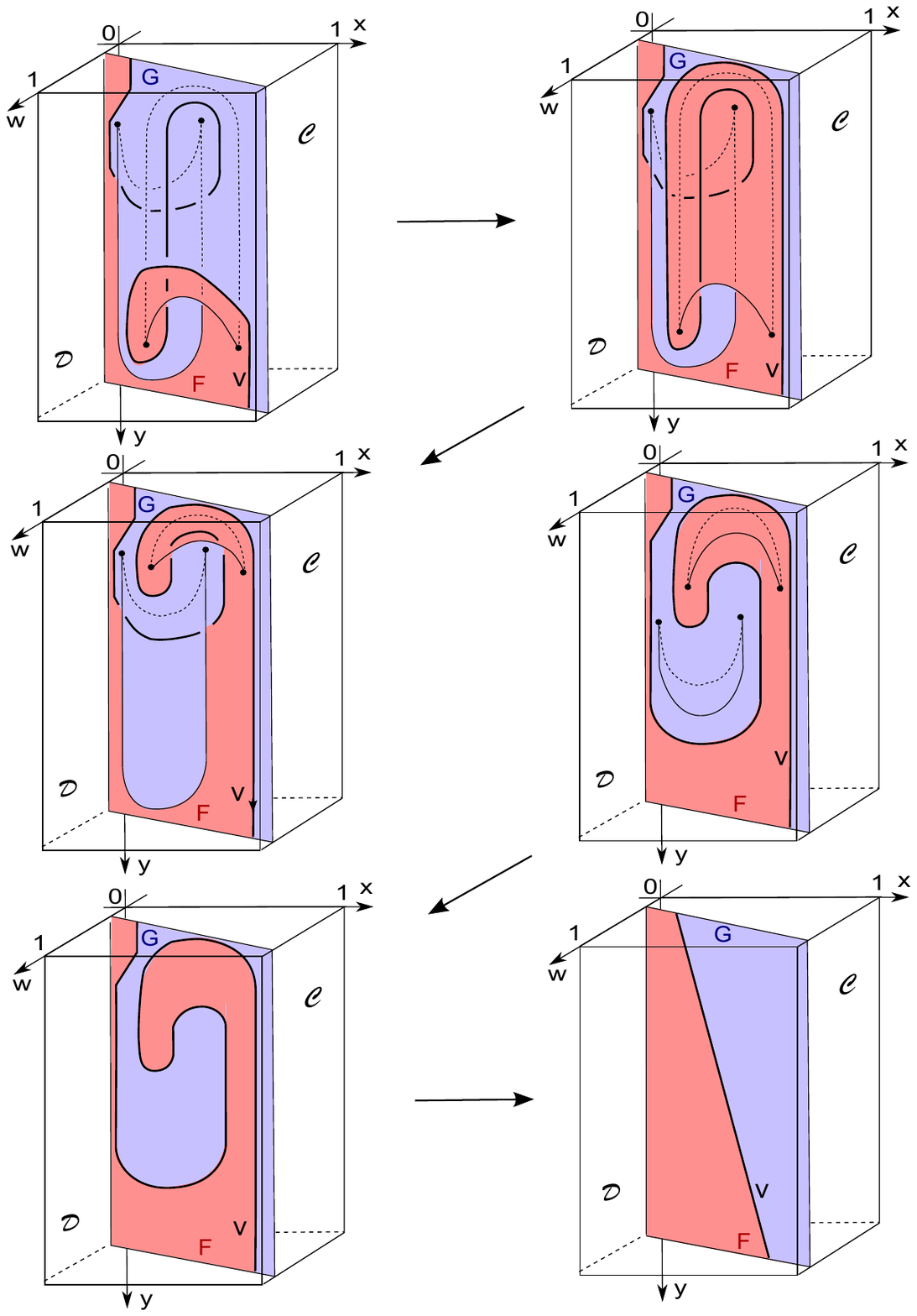}
  \caption{Diagrammatic proof of the relation \newline $\Theta_\nu\cdot \Theta_\nu^{-1}=1_\nu$.}
  \label{twistinvrel}
\end{figure}

  The naturality of  the 3-morphisms $\Theta_\nu$ is a direct consequence of the naturality  of the tensorator \eqref{nattens} together with the first condition in \eqref{pivcond} in the definition of a planar 2-category.  Also, it follows from  the invertibility of the triangulator, the invertibility of the  tensorator and the identities satisfied by the 3-morphisms $\epsilon_\nu$ that $\Theta_\nu$ is invertible 
  with inverse
  \begin{align*}
    &\Theta_\nu^{-1}=
    (1_{G\eta^*_{F^\#}}\circ \epsilon_{GF^\#\nu^*}^*\circ 1_{G\#\nu F\circ \eta_GF})
    \cdot(1_{G\eta_{F^\#}^*\circ GF^\#\nu^*}\circ \sigma_{G\#\nu \circ \eta_G\,,\, \nu})\\
    &\cdot(1_{G\eta_{F^\#}^*\circ GF^\#\nu^*\circ GF^\#\eta_G^*G\circ GF^\#\nu G^\#G\circ G\eta_{F^\#}G^\#G}\circ \epsilon_{G\eta_{G^\#}^*}^*\circ 1_{\eta_G G\circ \nu})\\
    & \cdot  (1_{G\eta_{F^\#}^*\circ GF^\#\nu\circ GF^\#\eta_G^*G}\circ \sigma^{-1}_{GF^\#\nu\circ G\eta_{F^\#}\,,\,\eta_{G^\#}^*}\circ 1_{G\eta_{G^\#}\circ \eta_GG\circ\nu})\\
    &\cdot(1_{G\eta^*_{F^\#}\circ GF^\#\nu^*}\circ GF^\#T_G\circ 1_{GF^\#\nu\circ G\eta_{F^\#}}\circ (T_G^*)^{-1}\circ 1_\nu)\\
    &\cdot (\epsilon_{G\eta_{F^\#}^*\circ GF^\#\nu^*}\circ 1_{\nu}).
  \end{align*}
  The 3-morphism $\Theta_\nu^{-1}\colon \nu\Rrightarrow\#\#\nu$ is  depicted in Figure \ref{twistinv} and a diagrammatic proof of the identity $\Theta_\nu\cdot \Theta_\nu^{-1}=1_\nu$  is given in  Figure \ref{twistinvrel}.

\begin{figure}
  \centering
  \includegraphics[scale=0.65]{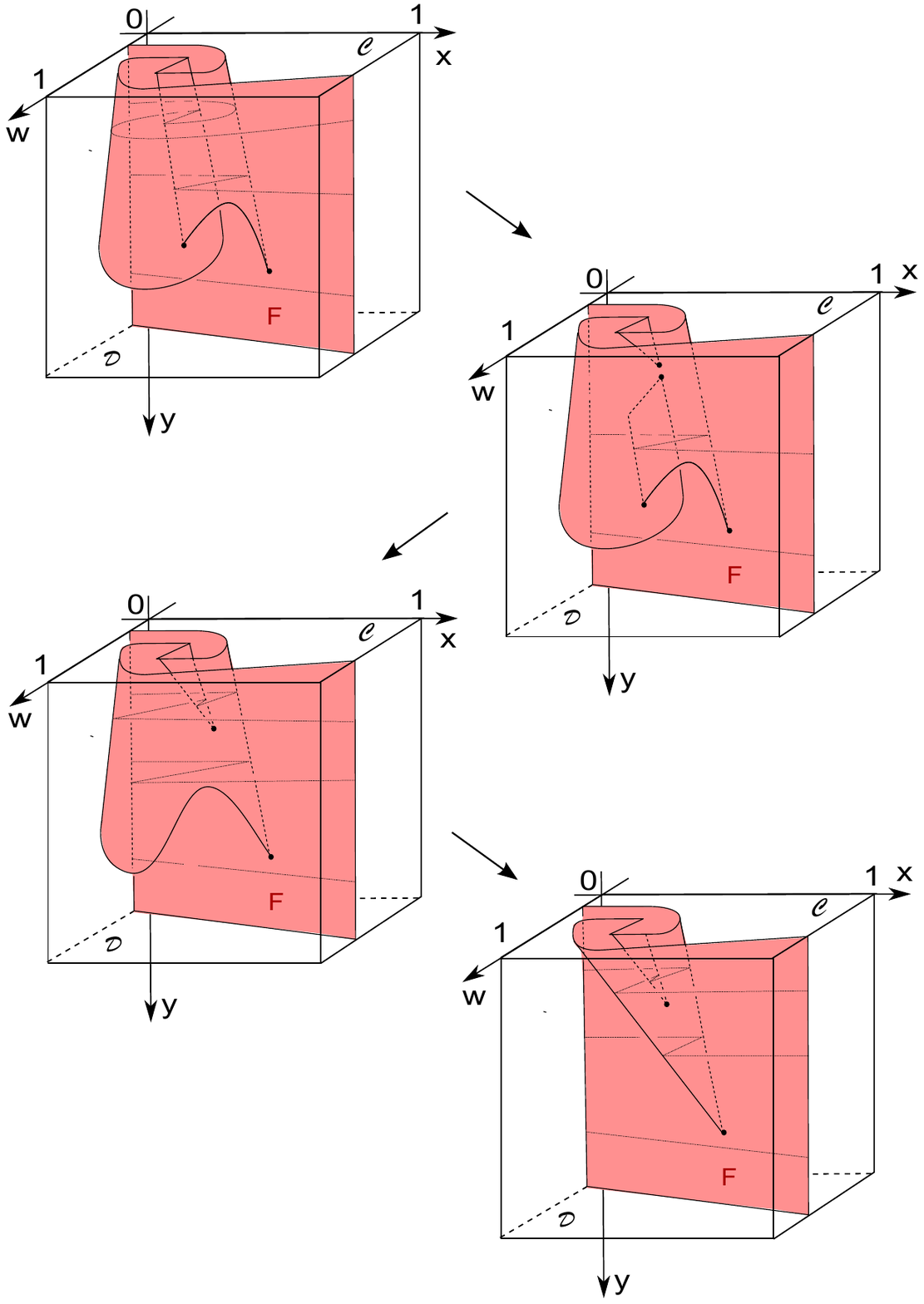}
  \caption{Diagrammatic proof of  identity \eqref{twistun}:\newline
    $\Theta_{1_F}=\Phi_{F^\#}^{-1}\cdot (\#\Phi_F)^{-1}$.
  }
  \label{twistunit}
\end{figure}

  To verify condition \eqref{twistun}  on the compatibility of $\Theta$ with the unit 2-morphisms $1_F\colon F\Rightarrow F$, we note that the 3-morphism $\Theta_{1_F}\colon \#\#1_F\Rrightarrow 1_F$ is given by
  \begin{align*}
    &\Theta_{1_F}= \epsilon^*_{\eta_F^*F}\cdot(1_{ \eta_F^*F}\circ T_FF^\#F\circ 1_{\eta_FF})
    \cdot ((T_F^*)^{-1}\circ\sigma_{\eta_F^*\,,\,\eta_F^*F}\circ 1_{F\eta_{F^\#}F^\#F\circ \eta_F F})\\
    &\quad\cdot(1_{ F\eta_{F^\#}^*}\circ \epsilon_{\eta_FF}\circ 1_{ \eta_FF}).
  \end{align*}
  By applying condition (d) in Definition \ref{dualgray} together with the invertibility of the triangulator   and the properties of the morphisms $\epsilon_{\nu}$ in a planar 2-category, one finds that the right hand side is equal to the 3-morphism $\Theta_{F^\#}^{-1}\cdot (\#\Theta_F)^{-1}$. A diagrammatic proof is given in Figure \ref{twistunit}.

\begin{figure}
  \centering
  \includegraphics[scale=0.65]{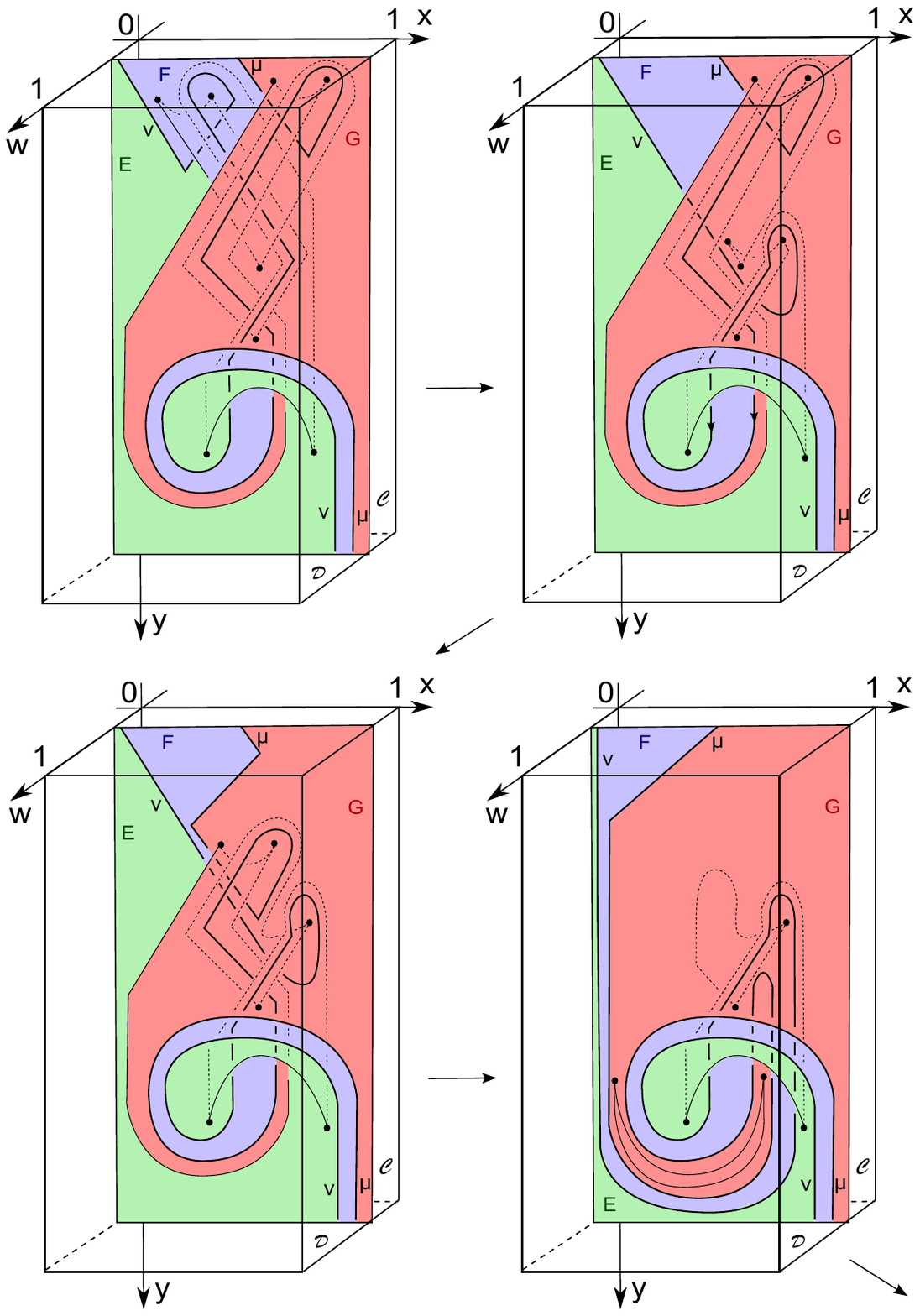}
  \caption{Diagrammatic proof of  identity \eqref{twisthorcomp}:\newline
    $\Theta_{\mu\circ\nu}\cdot \#\Phi_{\nu,\mu}\cdot \Phi_{\#\mu,\#\nu}\cdot(\Theta_\mu^\inv\circ \Theta_\nu^\inv)=1_{\mu\circ\nu}
    $ - continued in Figure \ref{thetacirc2}.
  }
  \label{thetacirc}
\end{figure}

\begin{figure}
  \centering
  \includegraphics[scale=0.65]{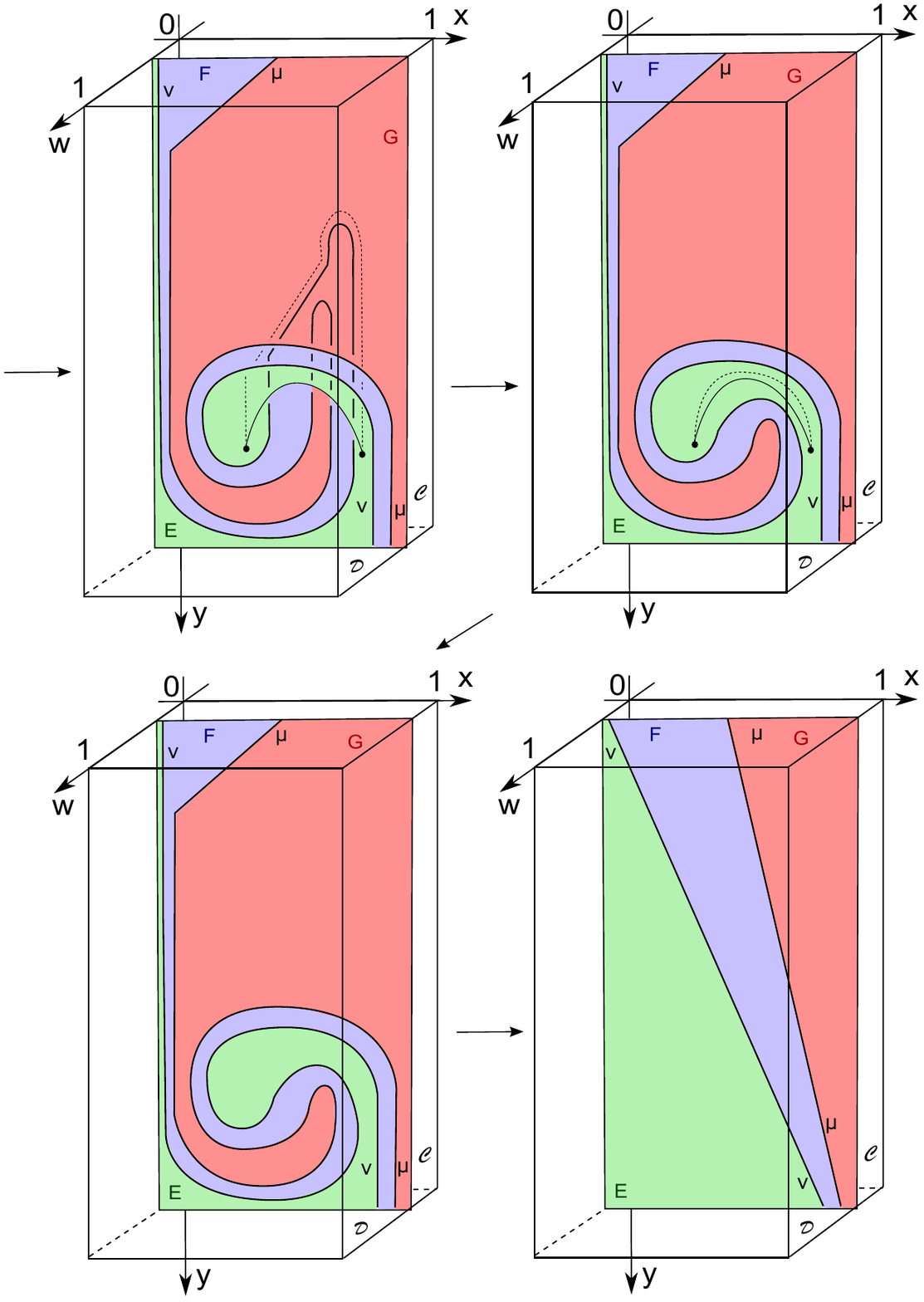}
  \caption{Diagrammatic proof of  identity \eqref{twisthorcomp}:\newline
    $\Theta_{\mu\circ\nu}\cdot \#\Phi_{\nu,\mu}\cdot \Phi_{\#\mu,\#\nu}\cdot(\Theta_\mu^\inv\circ \Theta_\nu^\inv)=1_{\mu\circ\nu}
    $ - continuation from  Figure \ref{thetacirc}.
  }
  \label{thetacirc2}
\end{figure}

  The condition \eqref{twisthorcomp}  on the compatibility of the 3-morphisms $\Theta_\nu$ with the horizontal composition of 3-morphisms is more involved. It is a consequence of the properties of the triangulator, the naturality of the tensorator and the properties of the  morphisms $\epsilon_{\nu}$ in a planar 2-category. A  diagrammatic proof  is given in Figure \ref{thetacirc} and \ref{thetacirc2}.

\begin{figure}
  \centering
  \includegraphics[scale=0.7]{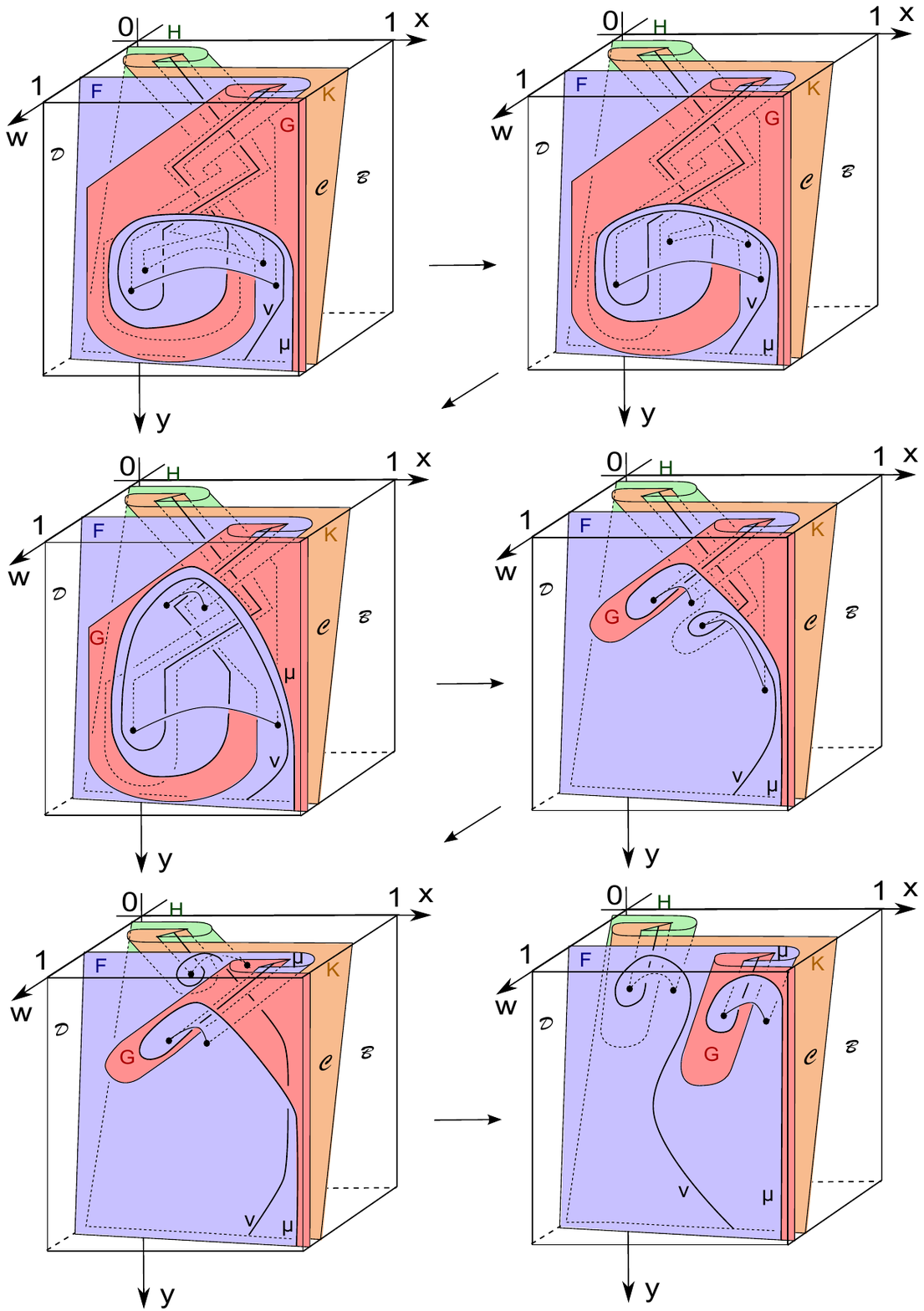}
  \caption{Diagrammatic proof of  identity \eqref{twistbox}:\newline
    $\Theta_{\mu\Box\nu}\cdot \#\kappa_{\mu,\nu}\cdot \kappa_{\#\nu,\#\mu}=\Theta_\mu\Box\Theta_\nu$.
  }
  \label{thetabox}
\end{figure}

  To prove identity \eqref{twistbox} which states the compatibility of $\Theta$ with the Gray product, consider the Gray category diagram for  $\Theta_{\mu\Box\nu}\cdot \#\kappa_{\mu,\nu}\cdot \kappa_{\#\nu,\#\mu,}$ in the upper left of Figure \ref{thetabox}. A diagrammatic proof that the associated 3-morphism  agrees with $\Theta_{\mu}\Box\Theta_{\nu}$ is given in Figure \ref{thetabox}.
\end{proof}

\subsection{Coherence properties of the duals}
\label{subsec}

This section investigates the interaction of the functors of 2-strict tricategories $*\colon \mac G\to\mac G^{op}$, $\#\colon\mac G\to\mac G_{op}$ with the natural transformations $\Theta\colon\#\#\to 1$, $\Gamma\colon *\#{*}\#\to 1$. The  results are needed in the strictification of these functors in Section \ref{sec:strictification}. 

The first result can be regarded as a coherence result for the functors of 2-strict tricategories $\#$ and $*$. By composing the natural transformations $\Theta\colon\#\#\to 1$ and $\Gamma\colon*\#{*}\#\to 1$ on the left and right with, respectively, the functors $\#$ and $*\#$, one obtains pairs of natural transformations
$\#\Theta, \Theta\#\colon \#\#\#\to\#$ and $\#*\Gamma^{\inv}, \Gamma *\#\colon \*\#{*}\#{*}\#\to *\#$. The following lemma shows that these natural transformations are equal.

\begin{lemma}\label{coherence} 
  The natural isomorphisms   $\Gamma\colon *\#{*}\#\to 1$  and  \linebreak $\Theta\colon \#\#\to1$ satisfy
  $$
  \#\Theta=\Theta\#,  \quad (*\#\Gamma)\cdot (\Gamma*\#)=1,
  $$
  and there is a natural isomorphism $\Delta\colon \#\to *\#*$ of functors of 2-strict tricategories  such that the following diagram 
  commutes
  \begin{align}\label{Deltagamma}
    \xymatrix{ 
      \#\#
      \ar[d]_{\Delta\#}  
      \ar[r]^{\#\Delta} 
      \ar[rd]^{\Theta} & \#*\#* \\
      {*\#{*}\#} \ar[r]_{\Gamma} & 1. \ar[u]_{*\Gamma*}
    }
  \end{align}
\end{lemma}

\begin{figure}
  \centering
  \includegraphics[scale=0.65]{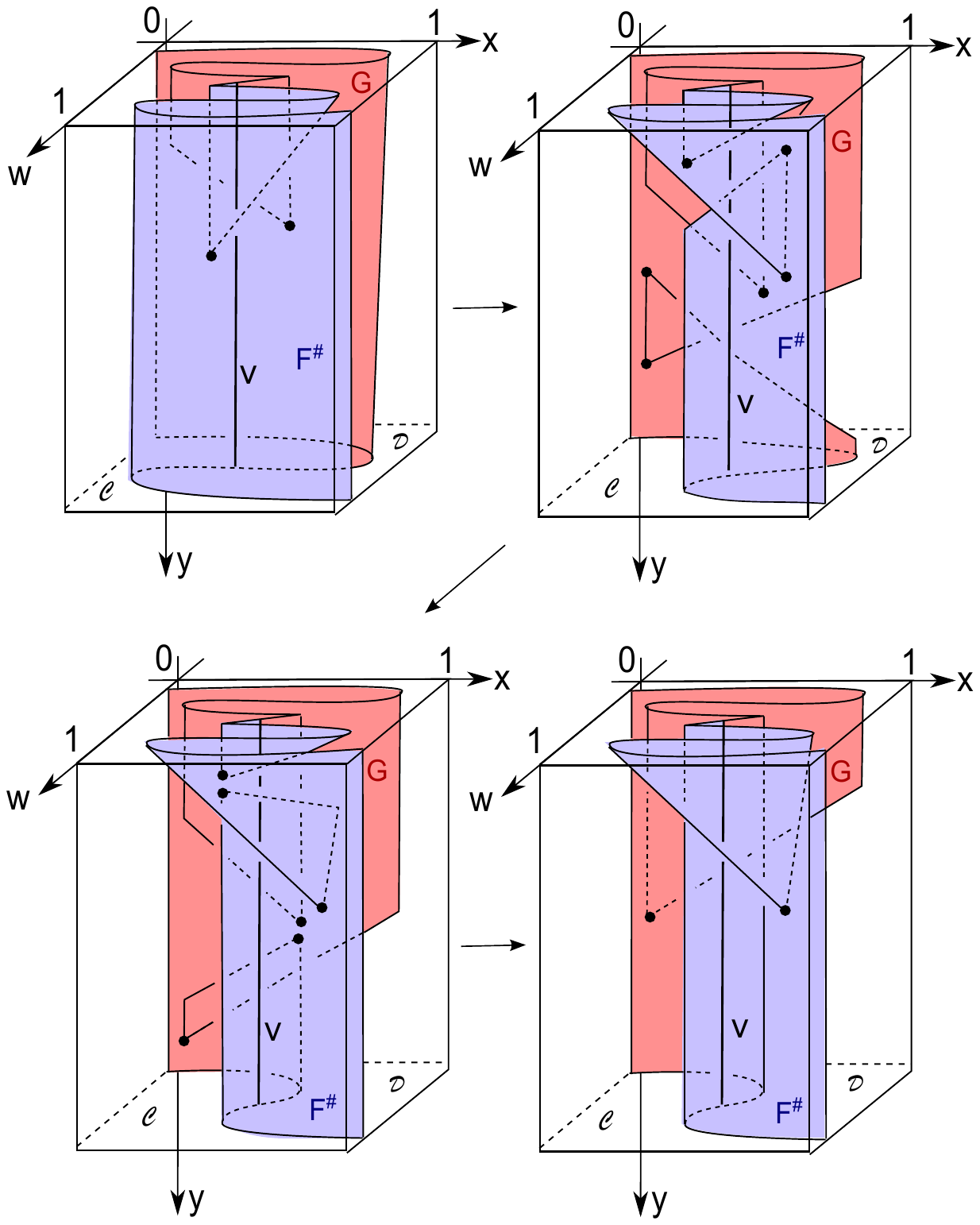}
  \caption{Diagrammatic proof of the  identity $\#\Gamma_\nu=(\Gamma_{*\#\nu}^*)^{-1}$.}
  \label{gammaassociativity}
\end{figure}

\begin{figure}
  \centering
  \includegraphics[scale=0.65]{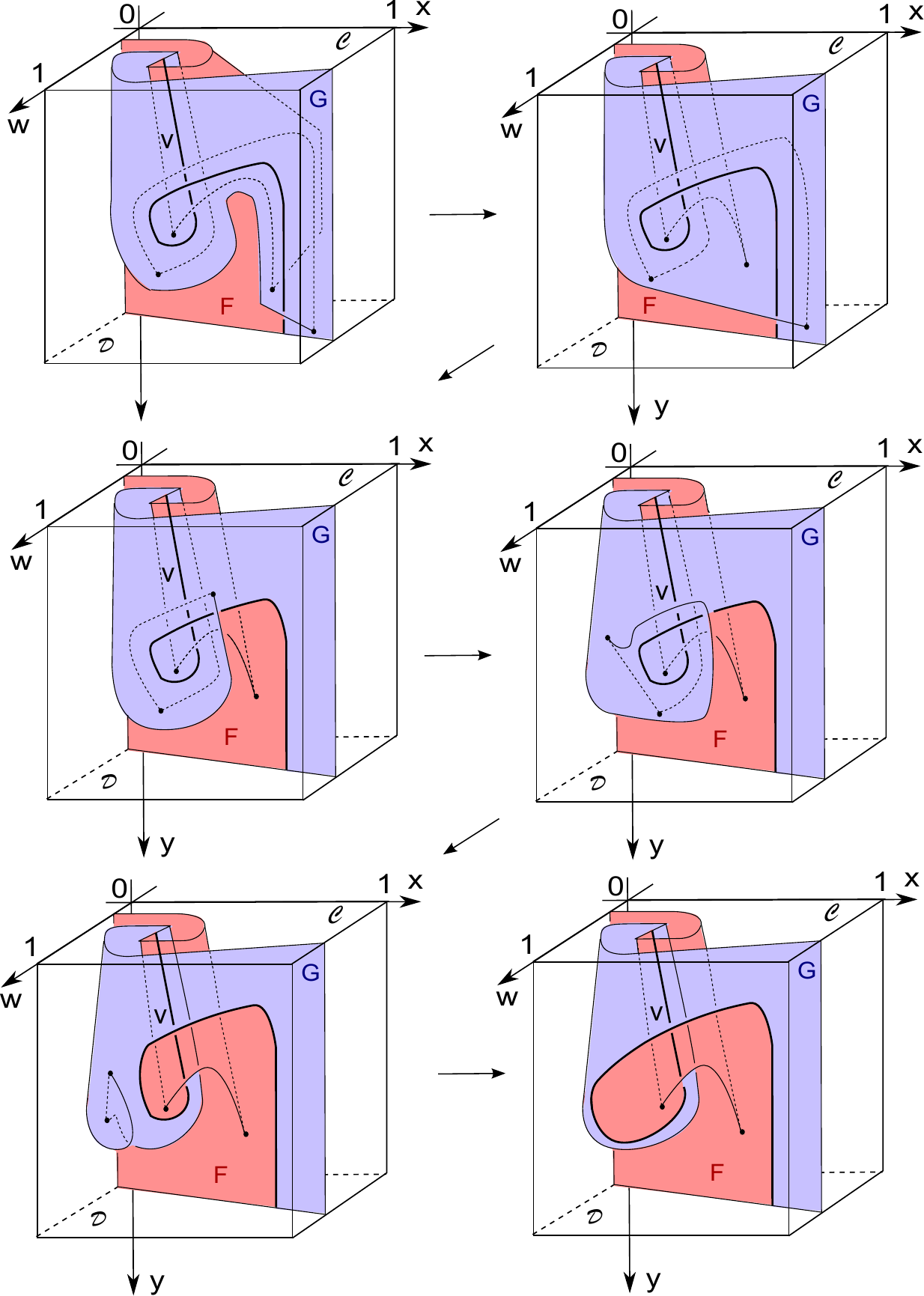}
  \caption{Diagrammatic proof of the identity \eqref{thetadeltarel}:\newline
    $\Theta_\nu=(\Gamma_{\nu^*}^*)^{-1}\cdot \#\Delta_{\nu}$.
  }
  \label{thetadelta}
\end{figure}

\begin{figure}
  \centering
  \includegraphics[scale=0.65]{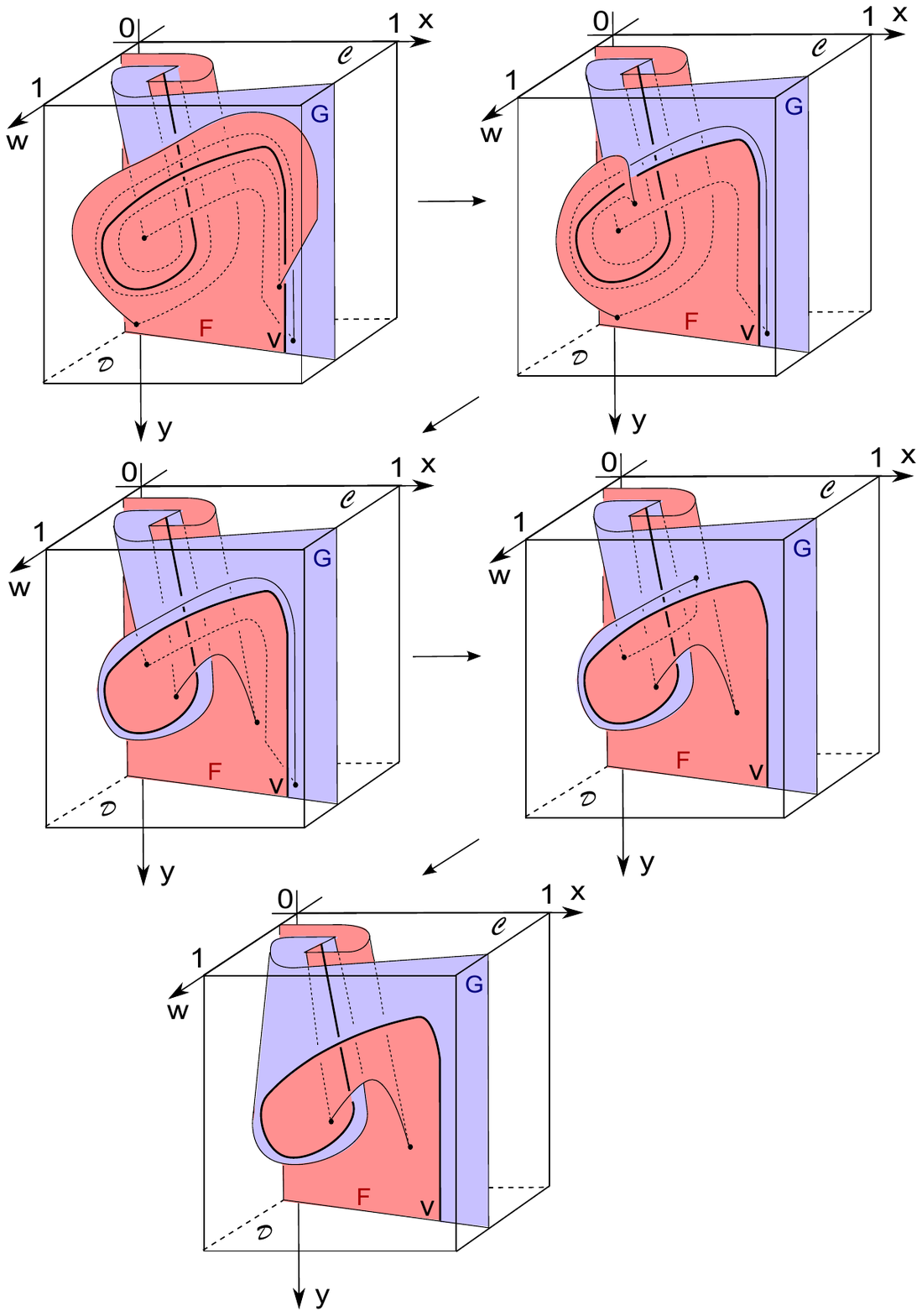}
  \caption{Diagrammatic proof of the identity \eqref{thetadeltarel}:
    $\Theta_\nu=\Gamma_\nu\cdot \Delta_{\#\nu}$.
  }
  \label{thetaassociativity}
\end{figure}

\begin{figure}
  \centering
  \includegraphics[scale=0.4]{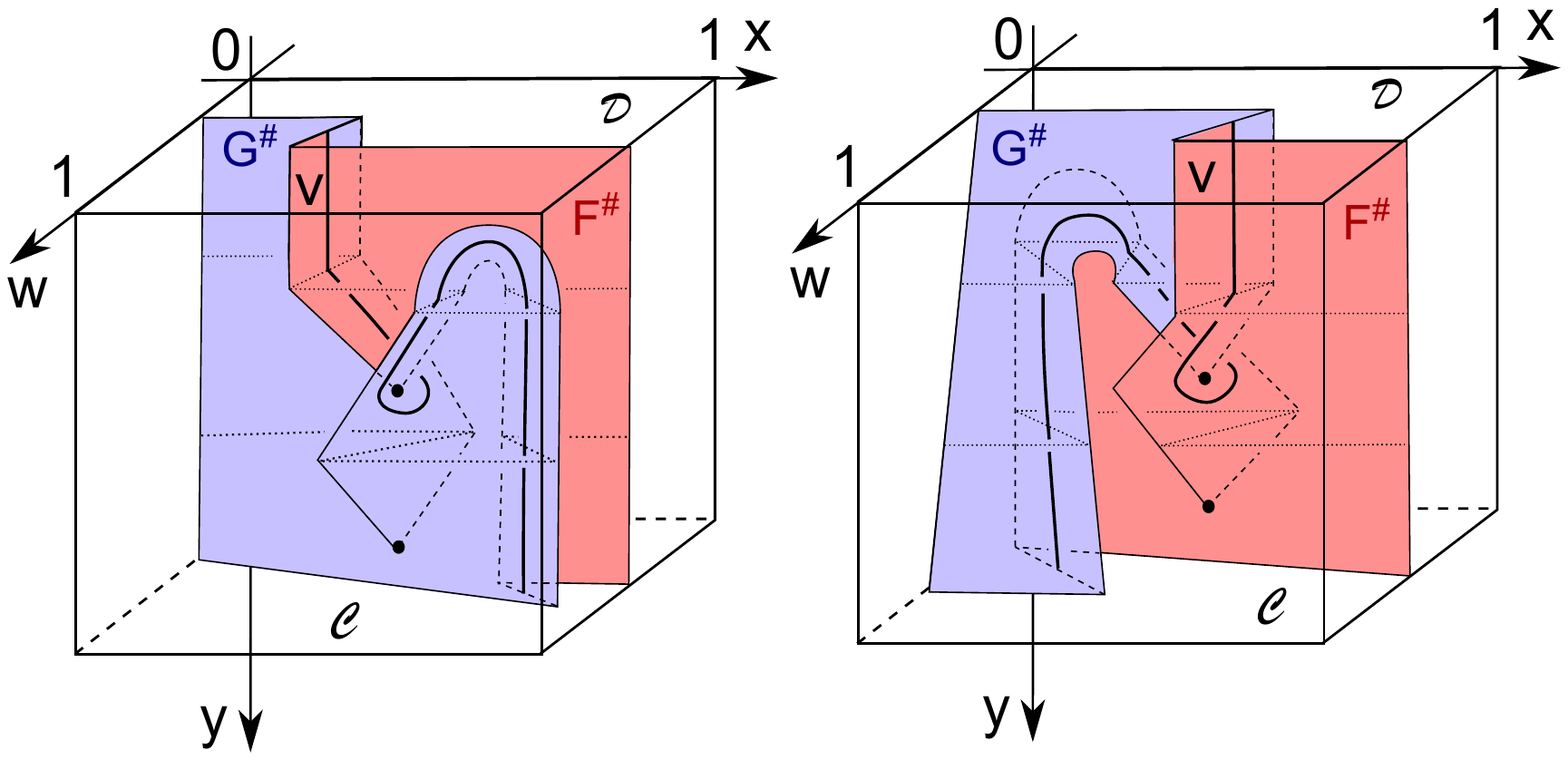}
  \caption{Gray category diagram for the 3-morphisms $\Delta_\nu\colon\#\nu\Rrightarrow *\#*\nu$ (on the left) and  $\Delta_{\nu^*}^*\colon \#\nu\Rrightarrow *\#*\nu$ (on the right).}
  \label{twistid}
\end{figure}

\begin{figure}
  \centering
  \includegraphics[scale=0.4]{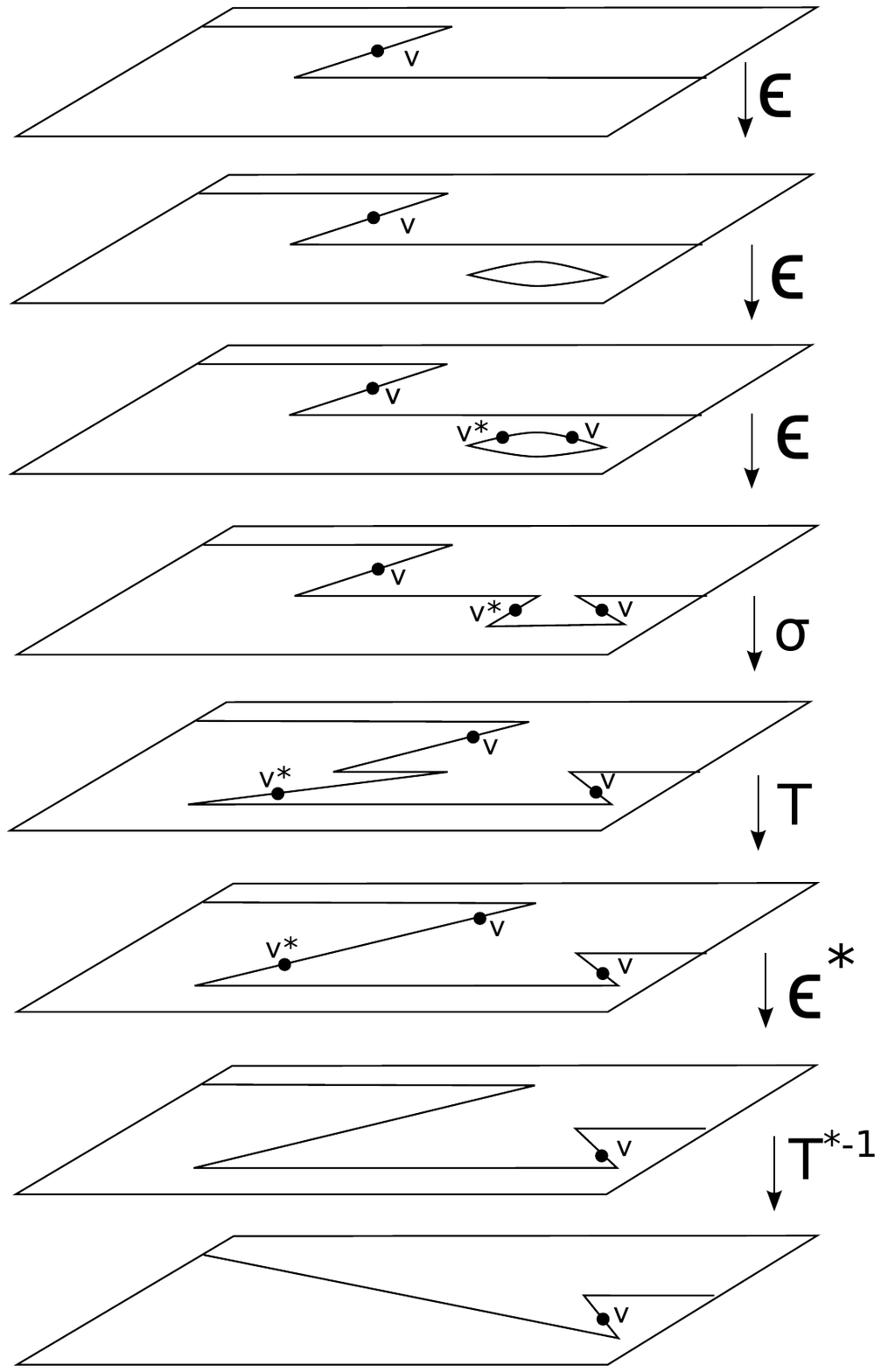}
  \caption{Gray category diagram for the 3-morphism $\Delta_\nu\colon\#\nu\Rrightarrow *\#*\nu$ in its movie representation obtained by taking constant height slices in Figure \ref{twistid} a). Some labels are omitted for legibility.}
  \label{twistid_slice}
\end{figure}

\begin{proof}

  In terms of the associated component 3-morphisms $\Theta_\nu\colon\#\#\nu\Rrightarrow \nu$ and $\Gamma_{\nu}\colon *\#{*}\#\nu\Rrightarrow \nu$, the first two relations in the lemma read
  $$
  \#\Theta_\nu=\Theta_{\#\nu},\quad *\#\Gamma_\nu=\Gamma_{*\#\nu}^\inv.
  $$
  A diagrammatic proof of the second  relation   is given in Figure \ref{gammaassociativity}. 

  To construct the natural transformation $\Delta\colon\#\to *\#*$, it is sufficient to specify its component 3-morphisms $\Delta_\nu\colon \#\nu\Rrightarrow *\#{*}\nu$ for each 2-morphism $\nu$ and to show that the following diagram commutes
  $$
  \xymatrix{ \#\#\nu \ar[d]_{\Delta_{\#_\nu}} \ar[r]^{\#\Delta_\nu} \ar[rd]_{\Theta_\nu} & \#*\#*\nu \ar[d]^{(\Gamma_{\nu^*}^*)^{-1}}\\
    {*}\#*\#\nu \ar[r]_{\Gamma_\nu} & \nu.
  }
  $$
  For a 2-morphism $\nu\colon F\Rightarrow G$, we define the  3-morphism $\Delta_\nu\colon \#\nu\Rrightarrow *\#{*}\nu$ as the composite
  \begin{align}\label{deltadef}
    &\Delta_\nu=(1_{*\#*\nu}\circ T^{*-1}_{G^\#} )
    \cdot(1_{*\#*\nu}\circ 1_{G^\#\eta_{G}^*}\circ G^\#\epsilon_{\nu}^*G^\#\circ 1 _{\eta_{G^\#}G^\#})\\
    &\cdot(1_{*\#*\nu}\circ 1_{G^\#\eta_G^*\circ G^\#\nu G^\#}\circ G^\#T_FG^\#\circ 1_{G^\#\nu^*G^\#\circ \eta_{G^\#}G^\#})\nonumber\\
    &\cdot (1_{*\#*\nu}\circ\sigma_{G^\#\eta_{F}^*\,,\,\eta_G^*\circ \nu G^\#} \circ 1_{G^\#F\eta_{F^\#}G^\#\circ G^\#\nu^*G^\#\circ \eta_{G^\#}G^\#})\nonumber\\
    &\cdot(1_{*\#*\nu}\circ 1_{G^\#\eta_F^*}\circ \sigma_{G^\#\nu^*\circ \eta_{G^\#}\,,\, \#\nu})\cdot ( \epsilon_{*\#*\nu}\circ 1_{\#\nu}).\nonumber
  \end{align}
  The Gray category diagram for the 3-morphism $\Delta_\nu$ is given by the left diagram in Figure \ref{twistid} and in Figure \ref{twistid_slice}.
  After some computations, which are performed diagrammatically in Figures \ref{thetadelta} and \ref{thetaassociativity},  one finds that the  3-morphism $\Theta_\nu\colon\#\#\nu\Rrightarrow\nu$ is given in terms of $\Delta_\nu$ by
  \begin{align}\label{thetadeltarel}
    \Theta_{\nu}=(\Gamma_{\nu^*}^*)^{-1}\cdot \#\Delta_\nu=\Gamma_\nu\cdot \Delta_{\#\nu}.
  \end{align}
  This implies the commutativity of the  diagram in \eqref{Deltagamma}.  It also follows directly that the 3-morphisms $\Delta_\mu$ define a natural isomorphism of functors of 2-strict tricategories.

  By combining diagram \eqref{Deltagamma} with the relation $*\#\Gamma_\nu=(\Gamma_{*\#\nu})^{-1}$, one obtains for all 2-morphisms $\nu$:
  $$
  \qquad\quad \Theta_{\#\nu}=(\Gamma_{*\#\nu}^*)^\inv\!\cdot \#\Delta_{\#\nu}=\#\Gamma_{\nu}\cdot \#\Delta_{\#\nu}=\#(\Gamma_\nu\cdot \Delta_{\#\nu})=\#\Theta_\nu,
  $$
  which proves the first identity in the lemma.
\end{proof}

Lemma \ref{coherence} has direct implications for the categories $\mathcal G(F,G)$ associated to 1-morphisms $F,G\colon\mathcal C\to\mathcal D$. 
The categories $\mathcal G(F,F)$ have a canonical structure as strict monoidal categories with the monoidal product given by the horizontal composition and the monoidal unit by the 2-morphism $1_F\colon F\Rightarrow F$. 
The functors of 2-strict tricategories $*$ and $\#$ 
induce functors 
$*\colon\mathcal G(F,G)\to\mathcal G(F,G)^{op}$ and $\#\colon\mathcal G(F,G)\to\mathcal G(G^\#,F^\#)$.

\begin{corollary} \label{cor_equiv_cat} For all 1-morphisms $F,G\colon\mathcal C\to\mathcal D$,
  the functors \linebreak $*\colon\mathcal G(F,G)\to\mathcal G(F,G)^{op}$, $\#\colon \mathcal G(F,G)\to\mathcal G(G^\#,F^\#)$ are equivalences of categories. When $\mathcal G(F,F)$ is equipped with its canonical monoidal structure, 
  then $*$ defines a  pivotal structure on $\mathcal G(F,F)$,
  $\#\colon\mathcal G(F,F)\to \mathcal G(F^\#,F^\#)_{op}$ is a strong monoidal  functor to the monoidal category $\mac G(F,F)_{op}$ with the opposite monoidal product,  and the 3-morphisms $\Delta_{\mu},\Delta^*_{\mu^*}\colon \#\mu\Rrightarrow *\#{*}\mu$ define natural isomorphisms $\#\to *\#*$.
\end{corollary}
\begin{proof} The functor $*\colon\mathcal G(F,G)\to\mathcal G(F,G)^{op}$ is an equivalence of categories since it is invertible: $**=1_{\mathcal G(F,G)}$. It  follows directly  from the axioms of a planar 2-category that $*$ equips each monoidal category $\mathcal G(F,F)$ with a  pivotal structure.

  To see that the functor  $\#\colon\mathcal G(F,G)\to\mathcal G(G^\#,F^\#)$ is essentially surjective, note that for each object $\mu$ of $\mathcal G(G^\#, F^\#)$,  the 3-morphism  $\Theta_\mu^\inv\colon \mu\Rrightarrow\#\#\mu$ defines an isomorphism in $\mac G(G^\#, F^\#)$ from $\mu$ to an object in the image of $\#$.
  That $\#\colon\mathcal G(F,G)\to\mathcal G(G^\#,F^\#)$ is fully faithful follows from the fact that $\Theta\colon\#\#\to 1$ defines a natural isomorphism $\#\#\cong 1_{\mac G(F,G)}$.

  To prove that $\#\colon\mathcal G(F,F)\to\mathcal G(F^\#,F^\#)_{op}$ is a strong monoidal functor,  consider the isomorphisms $\Phi_F\colon 1_{F^\#}\to \#1_{F}$ and the isomorphisms $\Phi_{\nu,\mu}\colon \#\nu\circ \#\mu\to \#(\mu\circ\nu)$ from the proof of Theorem \ref{graycatduals}. Identities \eqref{compat1} and \eqref{compat2} in the proof of Theorem \ref{graycatduals}  coincide with the axioms for a strong monoidal functor.
\end{proof}

The last  structural property of a Gray category with duals that will be required in the following is a relation between 
the  natural transformation $\Delta\colon\#\to *\#*$  from Lemma \ref{coherence} and its double $*$-dual
$*\Delta*\colon\#\to *\#*$.
The Gray category diagrams for their component morphisms are depicted in Figure \ref{twistid}.   If one restricts attention to  2-morphisms $\nu\colon 1_{\mathcal C}\Rightarrow 1_{\mathcal C}$ between trivial 1-morphisms, the folds and triangulators in Figure \ref{twistid} become trivial, and the diagrams reduce to two well-known diagrams from knot theory, which are required to be equal in a ribbon category.  It is thus natural to impose that the natural transformations $\Delta\colon\#\to *\#*$ and $*\Delta*\colon\#\to *\#*$
are equal, i.e., that the two Gray category diagrams in Figure \ref{twistid} have the same evaluation. This condition can be regarded as a horizontal categorification of the ribbon condition in a ribbon category.

\begin{definition}
\label{spatial} 
A Gray category $\mac G$ with duals is called {\bf spatial} if the natural transformations $\Delta\colon\#\to *\#*$ and $*\Delta*\colon\#\to *\#*$
  are equal.
\end{definition}

\begin{corollary}\label{cor:spatial} If $\mac G$ is a spatial Gray category, then for each object $\mac C$ of $\mac G$, the category $\mac G(1_{\mac C}, 1_{\mac C})$ is a ribbon category. Conversely, a ribbon category is a spatial Gray category with one object and one 1-morphism.
\end{corollary}
\begin{proof}
  If $\mac G$ is a Gray category with duals, then by  Lemma \ref{pivten} the category $\mac G(1_{\mac C}, 1_{\mac C})$ is a braided pivotal  category. As all Gray products with 1-morphism $1_\mac C$,  the 2-morphisms $\eta_{1_\mac C}$  and the 3-morphism $T_{1_{\mac C}}$ are trivial, it follows 
  that the 2-functor $\#\colon\mac G(1_{\mac C}, 1_{\mac C})\to\mac G_{op}(1_{\mac C}, 1_{\mac C})$ is trivial, and that the 3-morphisms from Theorem \ref{graydualnat} and Lemma \ref{coherence} satisfy $\Gamma_\mu=1_\mu$ and $\Theta_\mu=\Delta_\mu$ for all 2-morphisms $\mu$ in $\mac G(1_{\mac C}, 1_{\mac C})$.  For each object $\mu$, the 3-morphism  $\Theta_\mu=\Delta_\mu\colon \mu\Rrightarrow\mu$ reduces to the twist in a pivotal braided category. The condition that   
  $\mac G$ is spatial ensures that the twist satisfies the condition that makes $\mac G(1_{\mac C}, 1_{\mac C})$ into a ribbon category.
\end{proof}

\subsection{Geometrical interpretation of the duals}

\label{subsec:geom-interpret}

\begin{figure}
  \centering
  \includegraphics[scale=0.4]{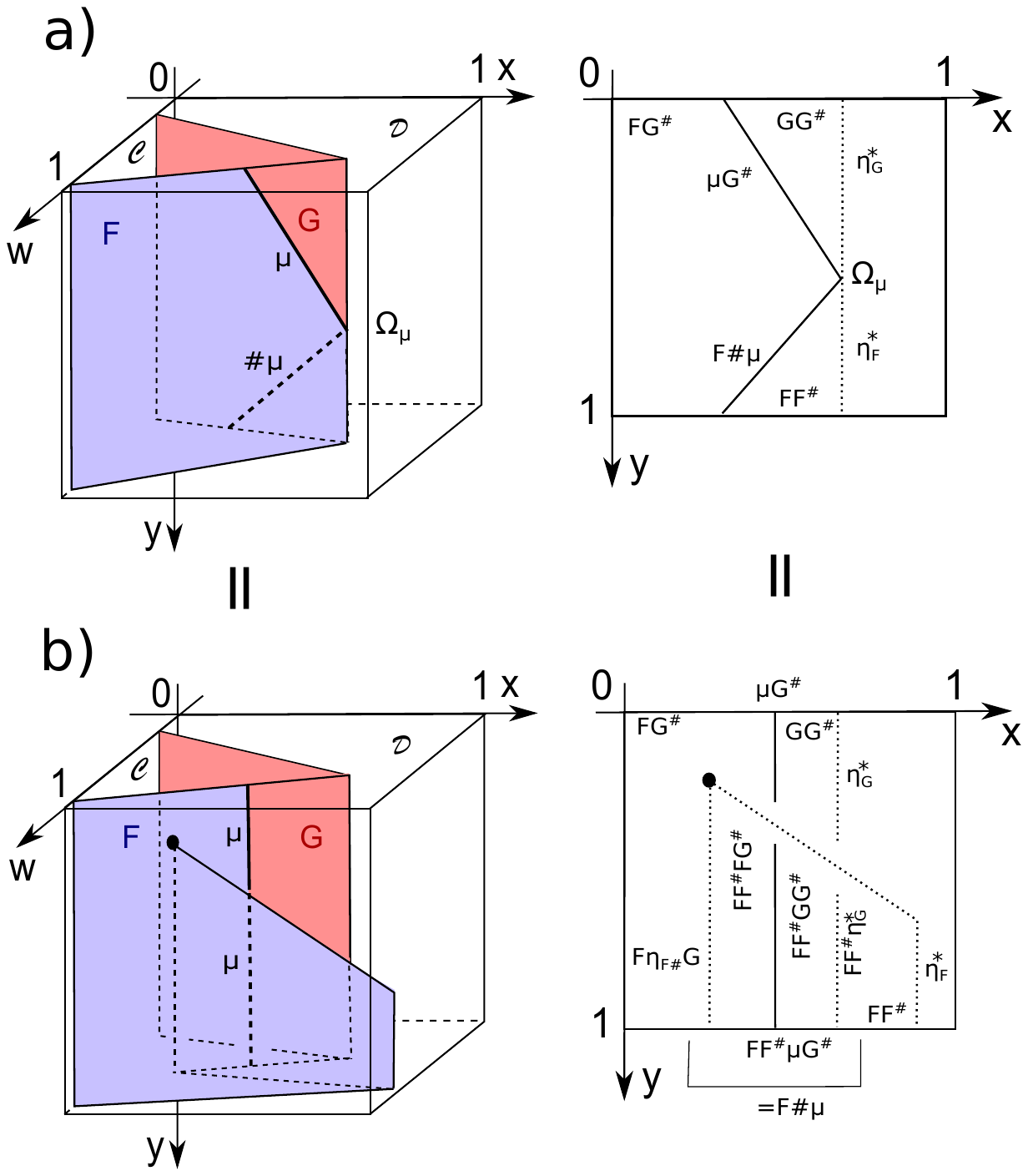}
  \caption{Gray category diagrams for the 3-morphism $\Omega_\mu\colon \eta_G^*\circ (\mu\Box G^\#)\Rrightarrow \eta_F^*\circ (F\Box \#\mu)$.}
  \label{foldcross3}
\end{figure}

\begin{figure}
  \centering
  \includegraphics[scale=0.62]{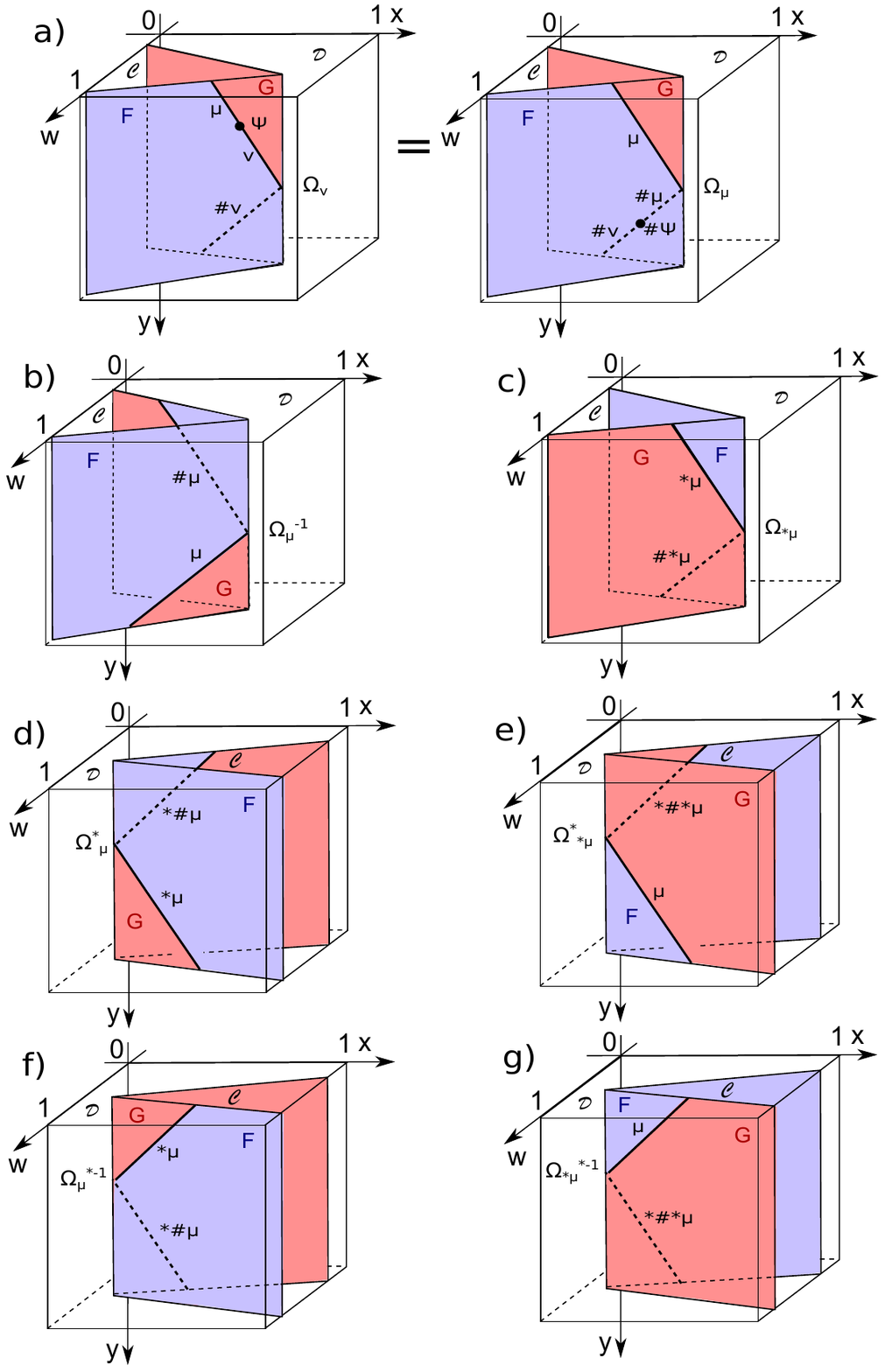}
  \caption{Gray category diagrams for  $\Omega_\mu$:\newline
    a) naturality condition,\newline 
b) inverse 3-morphism  $\Omega_\mu^\inv\colon \eta_F^*\circ (F\Box\#\mu)\Rrightarrow \eta_G^*\circ (\mu\Box G^\#)$\newline
    c) 3-morphism $\Omega_{\mu^*}\colon \eta_F^*\circ (\mu^*\Box F^\#)\Rrightarrow \eta_G^*\circ (G\Box \#*\mu)$,\newline
    d) dual 3-morphism $\Omega_\mu^*\colon (F\Box *\#\mu)\circ \eta_F\Rrightarrow (\mu^*\Box G^\#)\circ\eta_G$,\newline
    e) 3-morphism $\Omega^*_{\mu^*}$, f) 3-morphism $\Omega_{\mu}^{*-1}$, \newline g) 3-morphism $\Omega_{\mu^*}^{*-1}$.
  }
  \label{foldcross5}
\end{figure}

The functors of 2-strict tricategories $*\colon\mathcal G\to\mathcal G^{op}$, $\#\colon\mathcal G\to\mathcal G_{op}$  and the 
natural isomorphisms $\Gamma\colon *\#{*}\#\to 1$ and $\Theta\colon \#\#\to 1$ have a direct geometrical interpretation in terms of Gray category diagrams.  To see this, consider for each 2-morphism $\mu\colon F\Rightarrow G$ the 3-morphism 
$\Omega_\mu\colon \eta_G^*\circ (\mu\Box G^\#)\Rrightarrow \eta_{F}^*\circ (F\Box\#\mu)$ 
\begin{align}\label{omega_mu}
  &\Omega_\mu=(\sigma^\inv_{\eta_F^*, \eta_G^*\circ \mu G^\#}\circ 1_{F\eta_{F^\#}G^\#})\cdot (1_{\eta_G^*\circ \mu G^\#}\circ T_F^\inv G^\# ).
\end{align}
The  Gray category diagram for $\Omega_\mu$  and its projection are given in Figure \ref{foldcross3} b). By neglecting the expression for $\Omega_\mu$ in terms of the  data of a Gray category with duals and considering only its source and target 2-morphisms, one obtains the Gray category diagram in Figure \ref{foldcross3} a). The 3-morphism $\Omega_\mu$ thus allows one to let the lines labelled by 2-morphisms  cross folds. 
It follows directly from the properties of the tensorator that  $\Omega_\mu$ satisfies the naturality condition $$\Omega_\nu\cdot (1_{\eta_G^*}\circ (\Psi\Box G^\#))=(1_{\eta_{F}^*}\circ (F\Box \#\Psi))\cdot \Omega_\mu$$
for all 2-morphisms $\mu,\nu\colon F\Rightarrow G$ and all 3-morphisms $\Psi\colon \mu\Rrightarrow\nu$. This corresponds to sliding the dots labelled by 3-morphisms  over the folds as shown in Figure \ref{foldcross5} a). Moreover, the 3-morphism $\Omega_\mu$ is invertible, with inverse
$$
\Omega_\mu^\inv=(1_{\eta_G^*\circ \mu G^\#}\circ T_FG^\#)\cdot(\sigma_{\eta_F^*,\eta_G^*\circ \mu G^{\#}}\circ 1_{F\eta_{F^\#}G^\#}).
$$
The Gray category diagram for $\Omega_\mu^\inv$ is given in Figure \ref{foldcross5} b),  
and changing the orientation of the line leads to the Gray category diagram in Figure \ref{foldcross5} c). 
The Gray category diagrams  involving a fold that opens in the other direction are determined by the $*$-dual of $\Omega$ and are given in Figure \ref{foldcross5} d) to g).

\begin{figure}
  \centering
  \includegraphics[scale=0.65]{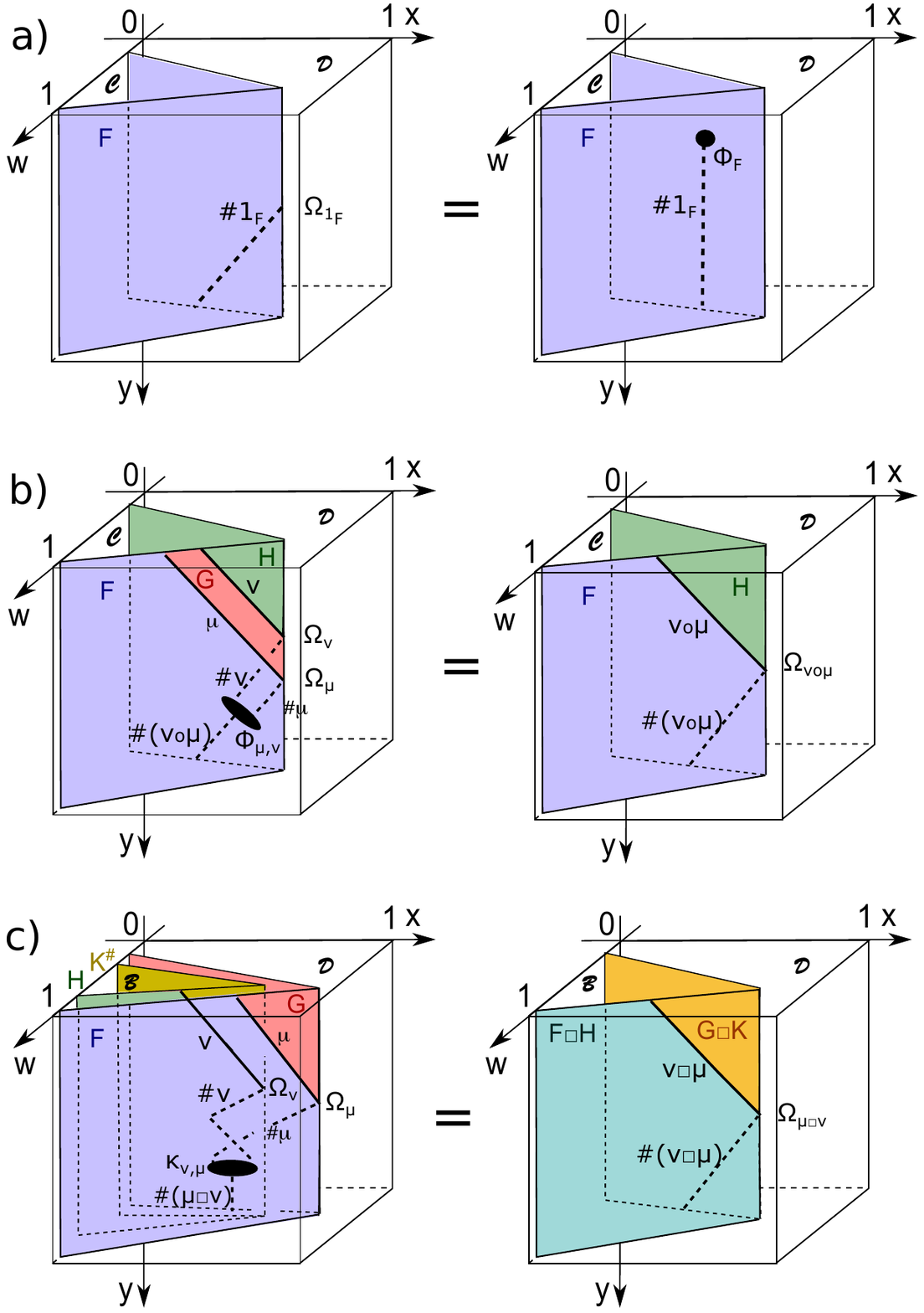}
  \caption{Gray category diagrams for the identities a) \eqref{foldnothing},  b) \eqref{foldcircle}  and c) \eqref{foldbox}. The 3-morphisms $\Phi_F$, $\Phi_{\mu,\nu}$ and $\kappa_{\mu,\nu}$ arise when lines corresponding to identity 2-morphisms or composite 2-morphisms cross a fold.}
  \label{foldcrossconsist3}
\end{figure}

\begin{figure}
  \centering
  \includegraphics[scale=0.65]{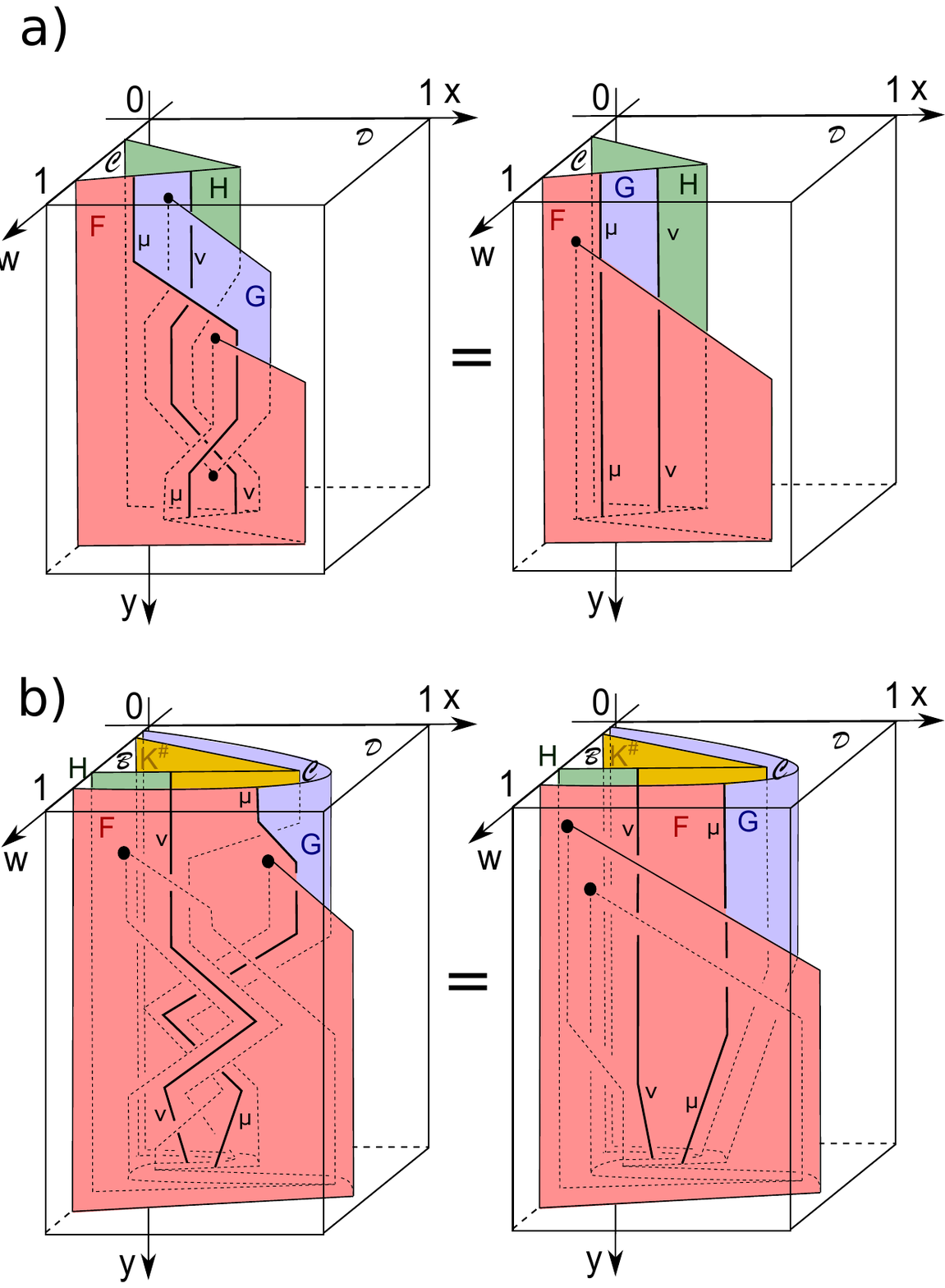}
  \caption{Diagrammatic proof of the identities a) \eqref{foldcircle}  and b) \eqref{foldbox}.}
  \label{foldcrossconsist}
\end{figure}

It remains to investigate the interaction of the 3-morphisms $\Omega_\mu$  with the unit 2-morphisms $1_F$, the horizontal composition $\circ$ and the Gray product $\Box$. For this, note that for every 1-morphism $F\colon\mathcal C\to\mathcal D$, the associated 3-morphism $\Omega_{1_F}\colon \eta_F^*\Rrightarrow \eta_F^*$ is given  in terms of the 3-morphism $\Phi_F\colon 1_{F^\#}\Rrightarrow \# 1_F$ from Theorem \ref{graycatduals}
\begin{align}\label{foldnothing}
  \Omega_{1_F}=\left(\sigma^\inv_{\eta_F^*,\eta_F^*}\circ 1_{F\eta_{F^{\#}}F^\#}\right)
\cdot \left(1_{\eta^*_F}\circ T_F^\inv F^\#\right)=1_{\eta_F^*}\circ F\Phi_F
\end{align}
which follows directly from the definition of $\Omega$ and the last property of the triangulator in Definition \ref{dualgray}. The corresponding Gray category diagram is given in Figure \ref{foldcrossconsist3} a).

The interaction of the 3-morphisms $\Omega_\mu$ with the horizontal composition is determined by the 3-morphisms $\Phi_{\mu,\nu}\colon \#\mu\circ \#\mu\Rrightarrow \#(\nu\circ\mu)$ from Theorem \ref{graycatduals}. A direct calculation shows that for all composable 2-morphisms
$\mu\colon F\Rightarrow G$, $\nu\colon H\Rightarrow K$, the following diagram commutes
\begin{align}\label{foldcircle}
  \xymatrix{\eta_H^*\circ \nu H^\#\circ \mu H^\# \ar[d]_{\Omega_\nu\circ 1_{\mu H^\#}} \ar[r]^{\Omega_{\mu\circ \nu}} & \eta_F^*\circ F\#(\nu\circ \mu)\\
    \eta_G^*\circ G\#\nu\circ \mu H^\# \ar[r]_{1\circ \sigma^\inv_{\mu,\#\nu}} & \eta_G^*\circ \mu G^\#\circ F\#\nu \ar[r]_{\Omega_\mu\circ 1_{F\#\nu}} & \eta_F^*\circ F\#\mu\circ F\#\nu \ar[lu]_{1_{\eta_F^*}\circ \Phi_{\mu,\nu}}.
  }
\end{align}
A diagrammatic proof of this identity is given in Figure \ref{foldcrossconsist} a). The corresponding Gray category diagrams are given in Figure \ref{foldcrossconsist3} b).

Similarly,  the interaction of the 3-morphism $\Omega_\mu$ with the Gray product  is governed by the 3-morphisms $\kappa_{\mu,\nu}\colon \#\nu\Box\#\mu\Rrightarrow\#(\mu\Box\nu)$ from Theorem \ref{graycatduals}. A direct computation which is performed diagrammatically in Figure \ref{foldcrossconsist} b) shows that  the following diagram commutes for all composable 2-morphisms $\mu\colon F\Rightarrow G$, $\nu\colon H\Rightarrow K$

\begin{align}\label{foldbox}
  \xymatrix{
    {\eta_{GK}^*\circ (\mu\Box\nu)K^\#G^\#}
    \ar[d]_{1_{\eta_G^*}\circ (\sigma^\inv_{\mu, \eta_K^*}G^\#)\circ 1 _{F\nu K^\#G^\#}} \ar[r]^{\Omega_{\mu\Box\nu}} & \eta_{FH}^*\circ  FH\#(\mu\Box\nu)\\
    \eta_G^*\circ \mu G^\#\circ F\eta_K^*G^\#\circ F\nu K^\#G^\# \ar[d]_{\Omega_\mu\circ  F\Omega_\nu G^\#} & \eta_{FH}^*\circ FH(\#\nu\Box\#\mu) \ar[u]_{1_{\eta^*_{FH}} \circ FH\kappa_{\mu,\nu}}\\
    \eta_F^*\circ F\#\mu\circ F\eta_H^*G^\# \circ FH\#\nu G^\# \ar[ru]_{\qquad \qquad 1_{\eta_F^*}\circ F\sigma^\inv_{\eta_{H}^*\circ H\#\nu, \#\mu}}.
  }
\end{align}
The corresponding Gray category diagrams are given in Figure \ref{foldcrossconsist3} c).

\begin{figure}
  \centering
  \includegraphics[scale=0.5]{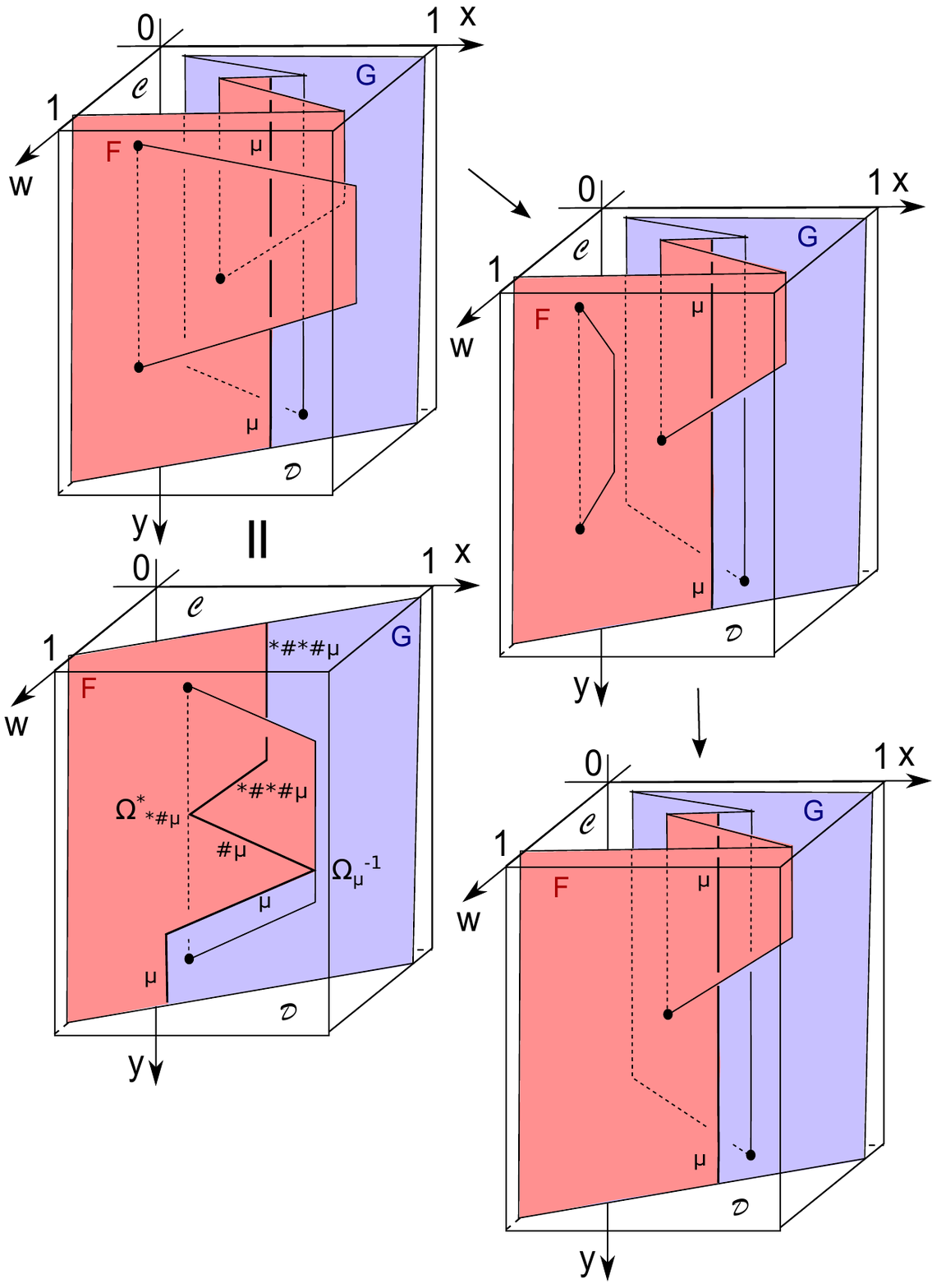}
  \caption{Diagrammatic proof of identity \eqref{gammaid_eq}.}
  \label{foldcrossgamma2}
\end{figure}

\begin{figure}
  \centering
  \includegraphics[scale=0.5]{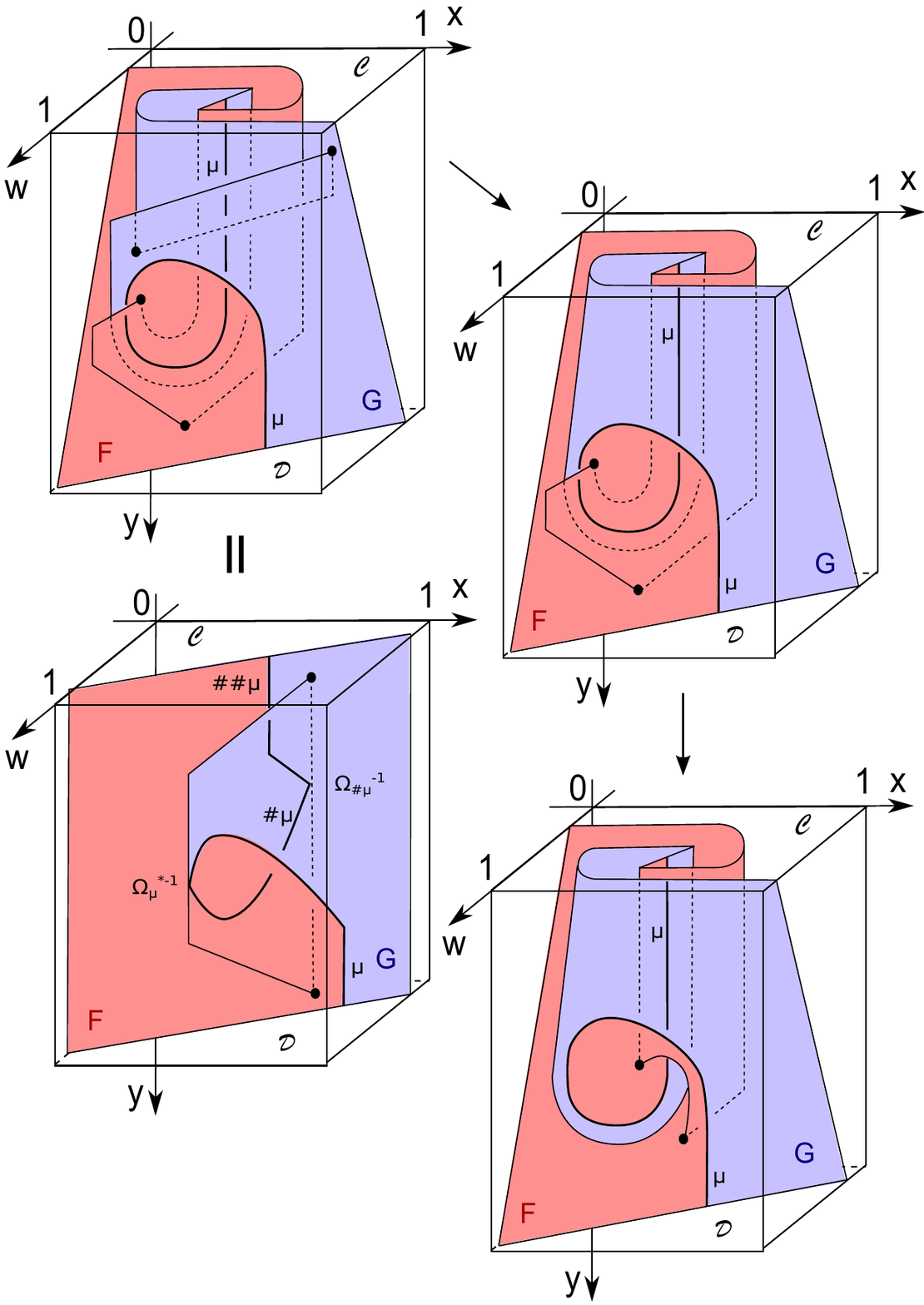}
  \caption{Diagrammatic proof of the identity \eqref{thetaid}}
  \label{foldcrosstheta}
\end{figure}

\medskip
To clarify the geometrical interpretation of the natural isomorphisms $\Gamma\colon *\#{*}\#\to 1$, $\Theta\colon \#\#\to 1$,  consider the Gray category diagrams given in the lower left of Figures \ref{foldcrossgamma2} and \ref{foldcrosstheta}. These diagrams contain two triangulators and two associated folds which together form a `pocket'   in the associated plane.  A line runs along the outside of this pocket, starting below and crossing its two folds to the top.  The difference between the diagrams in Figure \ref{foldcrossgamma2} and Figure \ref{foldcrosstheta} is that  it is a straight  line  in Figure \ref{foldcrossgamma2} while the diagram in   Figure \ref{foldcrosstheta} involves a line with maxima and minima.

Inserting the definition of the 3-morphism $\Omega_\mu$ (see Figure \ref{foldcross3}) in the lower left diagrams in Figures \ref{foldcrossgamma2} and \ref{foldcrosstheta} yields the diagrams in the upper left of these figures. It is shown in Figures \ref{foldcrossgamma2}, \ref{foldcrosstheta}   that the 3-morphisms obtained by projecting these diagrams are precisely the component 3-morphisms of the natural isomorphisms $\Gamma\colon *\#{*}\#\to 1$ and $\Theta\colon \#\#\to 1$. 
When expressed in terms of the 3-morphisms in the Gray category with duals, this corresponds to the  identities
\begin{multline}\label{gammaid_eq}
  \Gamma_\mu=(T_G\circ 1_\mu)\cdot (1_{\eta_G^*G}\circ \sigma_{\mu,\eta_{G^\#}})\cdot(\Omega_\mu^\inv G\circ 1_{F\eta_{G^{\#}}})\\
\cdot (1_{\eta_F^*G} \circ F\Omega^*_{*\#\mu})
  \cdot(\sigma^\inv_{\eta_F^*,*\#{*}\#\mu}\circ 1_{F\eta_{F^\#}})\cdot (1_{*\#{*}\#\mu}\circ T_F^\inv)
\end{multline}
\begin{multline} 
 \label{thetaid}
  \Theta_\mu=(1_\mu\circ T_F^{*-1})\cdot (\sigma^\inv_{\mu,\eta_{F^\#}^*}\circ 1_{\eta_F F}) \cdot (1_{G\eta^*_{F^\#}\circ \mu F^\#F}\circ F\epsilon_{\#\mu}^*F\circ 1_{\eta_FF})\\
  \cdot (1_{G\eta^*_{F^\#}\circ \mu F^\#F\circ F\#\mu F}\circ \Omega_\mu^{*-1}F)\cdot (1_{G\eta^*_{F^\#}}\circ \sigma^\inv_{\mu,\#\mu}F\circ 1_{\mu^*G^\#G\circ \eta_GF})\\
  \cdot (1_{G\eta^*_{F^\#}\circ G\#\mu F}\circ \epsilon_\mu G^\#F\circ 1 _{\eta_GF})\cdot(G\Omega^\inv_{\#\mu}\circ 1_{\eta_GF})\\
  \cdot  (1_{G\eta^*_{G^\#}}\circ \sigma_{\eta_G,\#\#\mu})\cdot(T_G^*\circ 1_{\#\#\mu}).
\end{multline}
The natural isomorphisms $\Gamma\colon*\#{*}\#\to 1$ and $\Theta\colon \#\#\to 1$ thus have a direct geometrical interpretation. Each of them relates  two Gray category diagrams that are obtained from each other by sliding a line over a cusp and relabelling one end of the line accordingly. 
The coherence data of the functors of 2-strict tricategories and the natural transformations $\Gamma$ and $\Theta$ thus arises naturally when one considers Gray category diagrams which involve lines crossing folds.

\section{Strictification for Gray categories with duals}
\label{sec:strictification}

In this  section,  it is shown that  for every spatial Gray category $\mac G$, the  functors of 2-strict tricategories $*:\mac G\to\mac G^{op}$
and $\#:\mac G\to\mac G_{op}$ can be strictified to strict functors of 2-strict tricategories   $\uline *:\uline{\mac G}\to \uline{\mac G}^{op}$ and $\uline\#: \uline{\mac G}\to \uline{\mac G}_{op}$
in the sense of Definitions \ref{grayfunc} and \ref{strictdef}, which satisfy $\uline\#\uline\#=1$, $\uline*\uline*=1$ and $\uline*\uline\#\uline*\uline\#=1$. The Gray category $\uline{\mac G}$ is  a Gray category with duals and equivalent to $\mac G$ as a   
2-strict
tricategory
(see Definition \ref{Equivalence-Gray-Cat}). The difference between $\mac G$ and $\uline{\mac G}$ 
is that the structures from  Definition \ref{dualgray} extend to strict functors of 2-strict tricategories.   This motivates the following definition.

\begin{definition}
\label{graystrictdual} 
A {\bf Gray category with strict duals} is a Gray category  with duals $\mac G$ such that the operations $\#$ and $*$ extend to
  strict functors of 2-strict tricategories $*:\mac G\to\mac G^{op}$, $\#:\mac G\to {\mac G}_{op}$ satisfying  
  $**=1$, $*\#{*}\#=1$, $\#\#=1$, and for all 2-morphisms $\mu$,
  $$\#\epsilon_\mu=\epsilon_{*\#\mu}.$$
\end{definition}

To strictify the spatial Gray with duals $\mac G$ to a Gray category $\uline{\mac G}$ with strict duals, we first construct a  Gray category $\uline{\mac G}$  with strict  functors of 2-strict tricategories 
$\uline *:\uline{\mac G}\to \uline{\mac G}^{op}$ and $\uline\#: \uline{\mac G}\to \uline{\mac G}_{op}$
satisfying $\uline\#\uline\#=1$, $\uline*\uline*=1$ and $\uline*\uline\#\uline*\uline\#=1$ that is equivalent to $\mac G$.
In a second step, we then show that the Gray category $\uline{\mac G}$ is a Gray category with duals in 
the sense of Definition \ref{dualgray} and satisfies the conditions in Definition \ref{graystrictdual}. 
This can be viewed as a higher analogue of Joyal and Street's result stating that every pivotal category is equivalent to a pivotal category with $**=\id$ \cite{JS2}. 
The proof is analogous, but in our  case, the non-trivial part is the identities $\uline\#\uline\#=1$, $\uline*\uline*=1$ and $\uline*\uline\#\uline*\uline\#=1$  while the identity  $**=\id$ already holds by assumption.

\begin{theorem}[Strictification] \label{thm:strictification}Let $\mathcal G$ be a spatial Gray category with duals. Then there exists a Gray category $\uline {\mathcal G}$ with the following properties:
  \begin{enumerate}

  \item $\uline{\mac G}$ is equipped with
    strict functors of 2-strict tricategories $\uline \#:\uline {\mathcal G}\to \uline{ \mathcal G}_{op}$, $\uline *: \uline {\mathcal G}\to\uline {\mathcal G}^{op}$ that satisfy $\uline\#\,\uline\#=1$, $\uline *\,\uline *=1$, $\uline*\underline\#\underline *\underline \#=1$.

  \item $\uline{\mac G}$ is equivalent to $\mac G$ as a Gray category: there  are \Glax functors of 2-strict tricategories $e: \uline{\mac G}\to{\mac G}$ and $f: \mac G\to\uline{\mac G}$ with $ef=1_{\mac G}$   and a natural isomorphism of \Glax functors of 2-strict tricategories $\eta:fe\to 1_{\uline{\mac G}}$. The natural isomorphism satisfies $e\eta=1_{e}$, $\eta fe=1_{fe}$,  and there is  an invertible modification $\Psi:\eta f \Rightarrow 1_f$ with $e\Psi=1_{{1_\mac G}}$.

  \item The \Glax functor of 2-strict tricategories $f:\mac G\to\uline {\mac G}$ satisfies $\uline * f=f*$, and there are
    natural isomorphisms of \Glax functors of 2-strict tricategories  $\xi: *e\to e\uline *$, $\chi: \#e\to e\uline\#$ and $\tilde\chi:\uline\#f\to f\#$.

  \end{enumerate}

\end{theorem}

\begin{proof} $\quad$\newline
  1.  We construct $\uline{\mac G}$ in analogy to the strictification proof for tricategories in \cite{gurski}. Note also that our proof can be viewed as a generalisation of the strictification proof for pivotal categories  in \cite{ng}.
  
  The objects of $\uline {\mac G}$
  are the objects of $\mac G$.  A  {basic 1-morphism} in $\uline{\mac G}$ from $\mac A$ to $\mac B$  is a pair $\uline F=(F,z)$ of a number $z\in
  \{1,-1\}$ and a 1-morphism
  $F: \mac A\to\mac B$ in $\mac G$ if $z=1$ or  a 1-morphism $F:\mac B\to\mac A$ in $\mac G$ if $z=-1$. The  1-morphisms of $\uline{\mac G}$ are composable strings $\uline F=(\uline {F_n},\ldots,\uline {F_1})$ of basic 1-morphisms, including the empty strings $\emptyset_{\mathcal A}:\mathcal A\to\mathcal A$. The evaluation of a 1-morphism $\uline F:\mathcal A\to\mathcal B$ in $\uline{\mac G}$ is the 1-morphism $e(\uline F):\mac A\to\mac B$ in $\mac G$ determined by
  \begin{align*}
    &e(F,1)=F,\quad e(F,-1)=\#F,\quad e(\emptyset_{\mac A})=1_{\mathcal A},\\
    &e(\uline {F_n},\ldots,\uline {F_1})=e(\uline {F_n})\Box \ldots \Box e(\uline {F_1}).
  \end{align*}
  A basic 2-morphism in $\uline{\mac G}$ is a pair $\uline\alpha=(\alpha, z): \uline F\Rightarrow \uline G$ consisting of a number $z\in\{1,-1\}$ and a 2-morphism $\alpha: e(\uline F)\Rightarrow e(\uline G)$ in $\mac G$  if $z=1$ or a 2-morphism $\alpha: \#e(\uline G)\Rightarrow \# e(\uline F)$ if $z=-1$. The 2-morphisms in $\uline{\mac G}$ are $\circ$-composable strings $\uline\alpha=(\uline{\alpha_n},\ldots,\uline{\alpha_1})$ of basic 2-morphisms in $\uline{\mac G}$, including the empty strings $\emptyset_{\uline F}: \uline F\Rightarrow \uline F$. The evaluation of a 2-morphism $\uline \alpha:\uline F\Rightarrow \uline G$ in $\uline{\mac G}$ is the 2-morphism $e(\uline\alpha): e(\uline F)\Rightarrow e(\uline G)$ in $\mathcal G$ determined  by
  \begin{align*}
    &e(\alpha,1)=\alpha,\quad e(\alpha,-1)=\#\alpha,\quad e(\emptyset_{\uline F})=1_{e(\uline F)},\\
    & e(\uline{\alpha_n},\ldots,\uline{\alpha_1})=e(\uline{\alpha_n})\circ\ldots\circ e(\uline{\alpha_1}),
  \end{align*}
  for all basic 2-morphisms $\uline{\alpha_{1}}, \ldots, \uline{\alpha_{n}}$.
  A 3-morphism $\uline\Gamma:\uline\alpha\Rrightarrow\uline \beta$ in $\uline{\mac G}$ is given by  a 3-morphism $e(\uline\Gamma): e(\uline\alpha)\Rrightarrow e(\uline\beta)$ in $\mac G$.

  The vertical composition 
  of 3-morphisms in $\uline{\mac G}$ is  the vertical composition in $\mac G$. The horizontal composition 
  of 2-morphisms in $\uline{\mac G}$ is  the concatenation of strings. This implies  that  the horizontal and vertical composition are strictly associative and the unit 2- and 3-morphisms are strict. As the evaluation  is also strictly compatible with the horizontal and vertical composition and the units,
  the horizontal composition of 3-morphisms in $\uline{\mac G}$ is given by the horizontal composition of 3-morphisms in $\mac G$. This shows that for each pair of objects $\mac A$, $\mac B$, $\uline{\mac G}(\mac A,\mac B)$ is a strict 2-category.

  The Gray product $\uline \Box$ in $\uline{\mac G}$ is defined as the concatenation of strings on 1-morphisms. 
  This implies  that  it is strictly associative and the unit 1-morphisms are strict:
  $$
  (\uline F\,\uline \Box \,\uline G)\uline\Box\,\uline H=\uline F\,\uline\Box(\uline G\,\uline \Box\,\uline H),\quad \uline F\,\uline \Box \emptyset_{\mac C}=\emptyset_{\mac D}\uline \Box\, \uline F=\uline F
  $$
  for all 1-morphisms $\uline F:\mac C\to\mac D$, $\uline G:\mac B\to\mac C$, $\uline H:\mac A\to\mac B$ in $\uline{\mac G}$. 
  It also follows that the Gray product of 1-morphisms is compatible with the evaluation. All composable 1-morphisms $\uline F,\uline G$ satisfy the relation
  $e(\uline F\,\uline \Box\, \uline G)=e(\uline F)\Box e(\uline G)$.

  For 1-morphisms $\uline F\in \uline{\mac G}(\mac C,\mac D)$, $\uline G\in \uline{\mac G}(\mac A,\mac B)$  and  a 2-morphism  $\uline\alpha \in  \uline{\mac G}(\mac B,\mac C)$ we define the Gray product by
  \begin{align*}
    &\uline F\,\uline \Box (\alpha,1)=(e(\uline F)\Box \alpha,1), &  &\uline F\,\uline \Box (\alpha,-1)=(\alpha\Box \#e(\uline F), -1),\\
    &(\alpha,1)\uline\Box\,\uline G=(\alpha\Box e(\uline G), 1), & &(\alpha,-1)\uline\Box\,\uline G=(\#e(\uline G)\Box \alpha,-1)\\
    &\uline F\,\uline\Box \emptyset_{\uline G}=\emptyset_{\uline F\,\uline \Box \,\uline G}, & &\uline F\,\uline\Box (\uline{\alpha_n},...,\uline{\alpha_1})=(\uline F\,\uline\Box\,\uline{\alpha_n},..., \uline F\,\uline\Box\,\uline{\alpha_1})\\
    &\emptyset_{\uline F}\uline\Box\,\uline G=\emptyset_{\uline F\,\uline\Box\,\uline G},
    & &(\uline{\alpha_n},...,\uline{\alpha_1})\uline\Box\,\uline G=(\uline{\alpha_n}\uline\Box\,\uline G,..., \uline{\alpha_1}\uline\Box\,\uline G).
  \end{align*}
  This  determines the Gray product of composable 2-morphisms
  $\uline\alpha:\uline F\Rightarrow\uline G$, $\uline\beta:\uline K\Rightarrow\uline H$ in $\uline{\mac G}$,  
  which is given by
  \begin{align}\label{gray2}
    \uline\alpha\,\uline\Box\,\uline\beta=(\uline\alpha\,\uline\Box\,\uline K)\uline \circ(\uline F\,\uline \Box\,\uline\beta).
  \end{align}
  As a direct consequence of these definitions, the Gray product  of 1- and 2-morphisms is again strictly associative, strictly compatible with the unit 2-morphisms and strictly compatible with the horizontal composition of 2-morphisms. 
  To define the Gray product of two 3-morphisms, we  compute the evaluation
  \begin{align*}
    &e(\uline F\,\uline\Box (\alpha,1))=e(\uline F)\Box e(\alpha,1), &  &e((\alpha,1)\uline\Box\, \uline G)=e(\alpha,1)\Box e(\uline G)\\
    &e(\uline F\,\uline\Box \emptyset_{\uline H})=e(\uline F)\Box e(\emptyset_{\uline H}), & &e(\emptyset_{\uline H}\uline\Box\,\uline G)=e(\emptyset_{\uline H})\Box e(\uline G),\\
    &e(\uline F\,\uline \Box (\alpha,-1))=\#(\alpha\Box \#e(\uline F)), & &e(\uline F)\Box e(\alpha,-1)=e(\uline F)\Box \#\alpha,\\
    &e((\alpha,-1)\uline\Box\, \uline G)=\#(\#e(\uline G)\Box\alpha), &  &e(\alpha,-1)\Box e(\uline G)=\#\alpha\Box e(\uline G).
  \end{align*}
  The Gray product of two 3-morphisms $\uline\Gamma:\uline\alpha\Rrightarrow\uline\alpha'$, $\uline\Psi:\uline\beta\Rrightarrow\uline\beta'$ between 2-morphisms $\uline\alpha, \uline\alpha': \uline F\Rightarrow \uline G$, $\uline\beta, \uline\beta': \uline H\Rightarrow \uline K$ is defined as
  $$
  e(\uline\Phi\,\uline\Box\,\uline\Psi)=( \tilde\iota_{\uline \alpha',\uline K}\circ \tilde \iota_{\uline F,\uline\beta'})
  \cdot (e(\uline\Phi)\Box e(\uline\Psi))\cdot (\tilde\iota_{\uline\alpha,\uline K}^\inv\circ \tilde \iota_{\uline F,\uline\beta}^\inv),
  $$
  where 
  $\tilde\iota_{\uline\alpha,\uline K}: e(\uline\alpha)\Box e(\uline K)\Rrightarrow e(\uline\alpha\,\uline\Box\,\uline K)$ and $\tilde\iota_{\uline L,\uline\alpha}: e(\uline L)\Box e(\uline\alpha)\Rrightarrow e(\uline L\,\uline\Box\,\uline \alpha)$  are the invertible 3-morphisms given by  
  \begin{align*}
    &\tilde \iota_{\emptyset_{\uline F}, \uline K}=1_{1_{e(\uline F)\Box e(\uline K)}}, \; \tilde \iota_{\uline L,\emptyset_{\uline F}}=1_{1_{e(\uline L)\Box e(\uline F)}}, &
    &\tilde \iota_{(\alpha,1), \uline K}\!=\!1_{\alpha\Box e(\uline K)}, \; \tilde \iota_{\uline L,(\alpha,1)}\!=\! 1_{e(\uline L)\Box \alpha},\\
    &\tilde \iota_{(\alpha,-1), \uline K}=\kappa_{ 1_{\#e(\uline K)},\alpha}\cdot (1_{\#\alpha}\Box \Phi_{\#e(\uline K)}), &
    &\tilde \iota_{\uline L,(\alpha,-1)}=\kappa_{\alpha, 1_{\#e(\uline L)}}\cdot (\Phi_{\#e(\uline L)}\Box 1_{\#\alpha}),
    \\
    &\tilde \iota_{(\uline\alpha_n,...,\uline\alpha_1), \uline K}=\tilde \iota_{\uline\alpha_n,\uline K}\circ\ldots\circ \tilde\iota_{\uline\alpha_1,\uline K}, &  &\tilde \iota_{\uline L, (\uline\alpha_n,...,\uline\alpha_1)}=\tilde\iota_{\uline L,\uline\alpha_n}\circ\ldots\circ \tilde\iota_{\uline L,\uline\alpha_1},
  \end{align*}
  with the 3-morphisms $\kappa_{\mu,\nu}: \#\nu\Box\#\mu\Rrightarrow \#(\mu\Box\nu)$ and $\Phi_F: 1_{\#F}\Rrightarrow \#1_{F}$ from  the proof of Theorem \ref{graycatduals}. In this, we used the identity $\#\#F=F$ in a Gray category with duals.
  It follows that the Gray product of 3-morphisms is strictly compatible with their vertical composition and with the unit 3-morphisms. 
  The  Gray product  is compatible with the horizontal composition of 3-morphisms 
  if and only if  the following two commutative diagrams and their counterparts with a 1-morphism on the left are equal
  \begin{align}\label{graycirc}
    &\xymatrix{
      e((\uline\alpha\,\uline\circ\,\uline\beta)\uline\Box\,\uline K) 
      \ar[rr]^{e((\uline \Psi\uline\circ \uline \Phi)\uline \Box\uline K)}
      \ar[d]_{\tilde\iota^\inv_{\uline\alpha\uline\circ\uline\beta,\uline K}}
      & & e((\uline\alpha'\uline\circ\uline\beta')\uline \Box K)  \\
      e(\uline\alpha\uline\circ\uline\beta)\Box e(\uline K) 
      \ar[rr]^{e(\uline\Psi\uline\circ\uline\Phi)\Box e(\uline K)}
      & & e(\uline\alpha'\uline\circ\uline\beta')\Box e(\uline K) 
      \ar[u]_{\tilde \iota_{\uline\alpha'\uline\circ\uline\beta',\uline K}}
    }\\
    &\xymatrix{
      e((\uline\alpha\,\uline\circ\,\uline\beta)\uline\Box\,\uline K) 
      \ar[rrr]^{e((\uline \Psi\uline\Box\uline K)\uline\circ (\uline \Phi\uline \Box\uline K))}
      \ar[d]_{\tilde\iota^\inv_{\uline\alpha,\uline K}\circ \tilde\iota^\inv_{\uline\beta,\uline K}}
      & & & e((\uline\alpha'\uline\circ\uline\beta')\uline \Box K)  \\
      e(\uline\alpha\uline\circ\uline\beta)\Box e(\uline K) 
      \ar[rrr]^{(e(\uline\Psi) \Box e(\uline K))\circ ( e(\uline\Phi)\Box e(\uline K))}
      & &  &e(\uline\alpha'\uline\circ\uline\beta')\Box e(\uline K) 
      \ar[u]_{\tilde \iota_{\uline\alpha',\uline K}\circ \tilde\iota_{\uline\beta', \uline K}}.
    }\nonumber
  \end{align}
  As the horizontal composition of 2-morphisms is strictly compatible with the evaluation, 
  the equality of these two diagrams follows directly  from the identities 
  \begin{align}\label{iotacomp}
    \tilde \iota_{\uline\alpha\,\uline\circ\,\uline\beta, \uline K}=\tilde\iota_{\uline\alpha, \uline K}\circ \tilde \iota_{\uline\beta, \uline K}, \quad \tilde \iota_{\uline L,\uline\alpha\,\uline\circ\,\uline\beta}=\tilde\iota_{\uline L,\uline\alpha}\circ \tilde \iota_{\uline L,\uline\beta}
  \end{align}
  which are satisfied by definition.  That the Gray product is associative amounts to the commutativity of the  following  diagrams
  \begin{align}
    &\xymatrix{e(\uline\alpha\, \uline\Box\, \uline F\,\uline\Box\,\uline G)  \ar@/ _7pc/[dd]_{\tilde\iota^\inv_{\uline\alpha,\uline F\uline \Box\uline G}} \ar[d]_{\tilde\iota^\inv_{\uline\alpha\uline\Box\uline F,\uline G}}\ar[r]^{e(\uline\Psi\uline\Box\uline F\uline\Box\uline G)} 
      & e(\uline\alpha'\uline\Box\uline F\uline\Box\uline G)
      \\
      e(\uline\alpha\uline\Box\uline F)\Box e(\uline G) 
      \ar[d]_{\tilde\iota^\inv_{\uline\alpha, \uline F}\Box e(\uline G)}  \ar[r]^{e(\uline\Psi \uline \Box \uline F)\Box e(\uline G)}& e(\uline\alpha'\uline\Box\uline F)\Box e(\uline G) \ar[u]_{\tilde\iota_{\uline\alpha'\uline\Box \uline F,\uline G}}
      \\
      e(\uline\alpha)\Box e(\uline F)\Box e(\uline G) \ar[r]^{e(\uline\Psi)\Box e(\uline F)\Box e(\uline G)}\;\;\; & \;\;\; e(\uline\alpha')\Box e(\uline F)\Box e(\uline G) \ar[u]_{\tilde\iota_{\uline\alpha',\uline F}\Box e(\uline G)} \ar@/ _7pc/[uu]_{\tilde \iota_{\uline\alpha',\uline F\uline \Box\uline G}}
    }\label{assoc_strict_1}\\
    \intertext{}
    &\xymatrix{ 
      e(\uline H\uline \Box\uline\alpha\uline\Box\uline F) 
      \ar@/ _5pc/[ddd]_{\tilde\iota^\inv_{\uline H,\uline\alpha\uline\Box\uline F}} 
      \ar[d]^{\tilde\iota^\inv_{\uline H\uline\Box\uline\alpha,\uline F}}
      \ar[r]^{e(\uline H\uline\Box\uline\Psi\uline\Box\uline F)}  
      & e(\uline H\uline \Box\uline\alpha'\uline\Box\uline F)  
      \\
      e(\uline H\uline\Box\uline\alpha)\Box e(\uline F) 
      \ar[d]^{\tilde\iota^\inv_{\uline H,\uline\alpha}\Box e(\uline F)} 
      \ar[r]^{e(\uline H\uline\Box \uline\Psi)\Box e(\uline F)} 
      &  e(\uline H\uline\Box\uline\alpha')\Box e(\uline F)  
      \ar[u]^{\tilde \iota_{\uline H\uline\Box \uline\alpha',\uline F}} \\
      e(\uline H)\Box e(\uline\alpha)\Box e(\uline F)\quad \ar[r]^{e(\uline H)\Box e(\uline\Psi)\Box e(\uline F)} &  \quad e(\uline H)\Box e(\uline\alpha')\Box e(\uline F) \ar[u]^{\tilde\iota_{\uline H,\uline\alpha'}\Box e(\uline F)} \ar[d]_{e(\uline H)\Box \tilde \iota_{\uline\alpha',\uline F}}\\
      e(\uline H)\Box e(\uline \alpha\uline\Box\uline F)  \ar[u]_{e(\uline H)\Box \tilde\iota^\inv_{\uline\alpha,\uline F}} \ar[r]^{{e(\uline H)}\Box e(\uline\Psi\uline\Box\uline F)}  &  e(\uline H)\Box e(\uline \alpha'\uline\Box\uline F)  \ar@/ _5pc/[uuu]_{\tilde\iota_{\uline H,\uline\alpha'\uline\Box\uline F}}
    } \label{assoc_strict_2}
  \end{align}
  and the analogue of  diagram  \eqref{assoc_strict_1} for  the composition with 1-morphisms from the left. In the  diagram \eqref{assoc_strict_1}, all squares commute by definition of the Gray product. 
  It remains  to prove the identities
  $$
  \tilde\iota_{\uline\alpha, \uline F\uline\Box\uline G}=\tilde\iota_{\uline\alpha\uline \Box\uline F,\uline G}\cdot (\tilde\iota_{\uline \alpha,\uline F}\Box 1_{e(\uline G)})\qquad  \tilde\iota_{\uline H\uline \Box\uline K,\uline\alpha}=\tilde\iota_{\uline H,\uline K\uline\Box\uline\alpha}\cdot(1_{e(\uline H)}\Box \tilde\iota_{\uline K,\uline\alpha})
  $$
  for all  2-morphisms  $\uline\alpha$ and 1-morphisms $\uline F,\uline G, \uline H,\uline K$ for which these expressions are defined. Due to the identities \eqref{iotacomp}, it is sufficient to prove this for basic 2-morphisms  and the empty string of 2-morphisms. For the latter and for basic 2-morphisms $\uline\alpha=(\alpha,1)$, this follows directly from the definition. For basic 2-morphisms 
  $\uline\alpha=(\alpha,-1)$,  inserting the definition of $\tilde\iota_{\alpha,\uline F}$  into these equations shows that this is the case for the first equation  if and only if the outer paths in the diagram
  $$
  \xymatrix{\#\alpha\Box  F\Box G 
    \ar[d]_{{\#\alpha}\Box\Phi_{\#F}\Box {G}} 
    \ar[r]^{{\#\alpha}\Box \Phi_{\# G\Box \# F}} 
    & \#\alpha\Box \#1_{\# G\Box \# F}  
    \ar[d]_{{\#\alpha}\Box \kappa^\inv_{ 1_{\# G},1_{\# F}}}  
    \ar[r]^{\kappa_{ 1_{\# G\Box \# F},\alpha}}
    & \#(1_{\# G\Box\# F}\Box \alpha)
    \\
    \#\alpha\Box \#1_{\# F}\Box  G \quad
    \ar[d]_{\kappa_{ 1_{\#F},\alpha}\Box G} 
    \ar[r]^{1_{\#\alpha}\Box \#1_{\# F}\Box \Phi_{\# G}}
    & \qquad \#\alpha\Box \#1_{\# F}\Box \#1_{\#G} 
    \ar[d]_{\kappa_{ 1_{\#F},\alpha}\Box \#1_{\#G}}
    \\
    \#(1_{\#F}\Box\alpha)\Box  G \quad
    \ar[r]^{{\#(1_{\#F}\Box\alpha)}\Box\Phi_{\#G}}
    & \qquad \#(1_{\# F}\Box\alpha)\Box \#1_{\# G} 
    \ar@/ _/[ruu]_{\kappa_{1_{\#G},1_{\#F}\Box\alpha}}
  }
  $$
  are equal.  
  The rectangle on the lower left commutes. The subdiagram on the right commutes due to  equation \eqref{kappa_assoc}, and the upper left rectangle due to the compatibility condition \eqref{kappaunit} between the 3-morphisms $\kappa$ and $\Phi_F$, $\Phi_G$. This shows that the diagram \eqref{assoc_strict_1} commutes. The proof for the commutativity of the corresponding diagram with the 1-morphisms on the left is analogous.

  In  diagram \eqref{assoc_strict_2}, the three rectangles in the middle commute by definition  of the Gray product. It is therefore sufficient to prove that the two subdiagrams with curved arrows on the left and right  commute.
  Using again the identities \eqref{iotacomp}, one finds that  it is sufficient to show that this is the case for the basic 2-morphisms   and the empty string of 2-morphisms. In the cases $\uline\alpha=\emptyset_{\uline F}$ and $\uline\alpha=(\alpha,1)$ the commutativity of the subdiagrams is obvious. For $\uline \alpha=(\alpha,-1)$, we
  insert the definition of $\tilde\iota_{\uline\alpha,\uline F}$ and $\tilde \iota_{\uline H,\uline \alpha}$  and obtain the following  diagram, in which we abbreviate $F=e(\uline F)$ and $H=e(\uline H)$
  
  \begin{align} \label{diagmove}
    \xymatrix{ & \#(\# H\Box \alpha\Box \#  F) \ar[ld]_{\kappa^\inv_{1_{ \#H\Box\alpha},\# F}} \ar[rd]^{\kappa^\inv_{1_{\# H},\alpha\Box \# F}}
      \\
      \#1_{\# F}\Box \#(\# H\Box \alpha) \ar[rd]^{\#1_{\# F}\Box \kappa^\inv_{ 1_{\# H},\alpha}} \ar[d]_{\Phi^\inv_{\# F}\Box {\#(\# H\Box\alpha)}}& &  \#(\alpha\Box \# F)\Box \#1_{\# H} \ar[ld]_{\kappa^\inv_{\alpha,1_{\# F}}\Box \#1_{\# H}} \ar[d]^{{\#(\alpha\Box \# F)}\Box \Phi^\inv_{\# H}} 
      \\
      { F}\Box \#(\#H\Box \alpha)  \ar[d]_{{F} \Box\kappa^\inv_{1_{\# H},\alpha} } &  \#1_{\# F}\Box \#\alpha\Box \#1_{\# H}   \ar[ld]^{\Phi^\inv_{\# F}\Box {\#\alpha\Box \#1_{\# H}}} \ar[rd]_{{\#1_{\#F}\Box {\#\alpha}}\Box \Phi^\inv_{\# H}} & \#(\alpha\Box \# F)\Box  H \ar[d]^{\kappa^\inv_{\alpha,1_{\#F}}\Box {H}} 
      \\
      { F}\Box \#\alpha\Box \#1_{\# H} \ar[rd]_{{ F\Box \#\alpha}\Box \Phi^\inv_{\# H}} & & \#1_{\# F}\Box \#\alpha\Box { H} \ar[ld]^{\Phi^\inv_{\# F}\Box {\#\alpha\Box  H}}\\
      &  F\Box \#\alpha\Box  H.
    }
  \end{align} 
  The outer paths in  diagram \eqref{diagmove} correspond to the 3-morphisms 
  $\tilde \iota_{\uline H, \uline\alpha\uline\Box\uline F}\cdot (1_{e(\uline H)}\Box \tilde\iota_{\uline\alpha,\uline F})$ and $\tilde\iota_{\uline H\uline\Box\uline\alpha,\uline F}\cdot(\tilde \iota_{\uline H,\uline\alpha}\Box 1_{e(\uline F)})$.
  It is directly apparent that the lower parallelogram in the middle and the two subdiagrams on the right and the left commute.  The upper parallelogram in the middle commutes due to the pentagon axiom for $\kappa$ in equation \eqref{kappa_assoc} and hence
  the outer paths on the left and on the right are equal. This proves the  commutativity of the diagram \eqref{assoc_strict_2} and completes the proof that the Gray product is strictly associative.

  To conclude that $\uline{\mac G}$ is a  Gray category, we define the tensorator $\uline\sigma_{\uline\alpha,\uline\beta}: (\uline\alpha\,\uline\Box\,\uline K)\uline{\circ}(\uline F\,\uline\Box\,\uline \beta)\Rrightarrow(\uline G\,\uline\Box\,\uline\beta)\uline \circ(\uline\alpha\,\uline\Box\uline H)$ in $\uline{\mac G}$   by 
  $$
  \uline\sigma_{\uline\alpha,\uline\beta}=(\tilde\iota_{\uline G,\uline\beta}\circ \tilde\iota_{\uline\alpha,\uline H})\cdot
  \sigma_{e(\uline\alpha), e(\uline\beta)}\cdot (\tilde\iota_{\uline\alpha, \uline K}^\inv\circ\tilde\iota_{\uline F,\uline\beta}^\inv).$$ It follows from the definition, the properties of the tensorator in $\mac G$ and the identities proved above that $\uline\sigma_{\uline\alpha,\uline\beta}$ satisfies the axioms for the tensorator in Definition \ref{graydef}. 
  This  shows that $\uline{\mac G}$ is a Gray category.

  \bigskip
  2. To construct the strict functor of 2-strict tricategories $\uline\#:\uline{\mac G}\to \uline{\mac G}_{op}$, 
  we set $\uline\#(\mac A)=\mac A$ for objects, $\uline\#(F,z)=(F,-z)$ for basic 1-morphisms $\uline F=(F,z)$ of $\uline{\mac G}$ and extend $\uline\#$ to general 1-morphisms by 
  $$
  \uline\#(\uline{F_n},\ldots,\uline{F_1})=(\uline\#\,\uline{F_1},\ldots, \uline\#\,\uline{F_n}),\quad \uline\#(\emptyset_{\mac A})=\emptyset_{\mac A}.
  $$
  It follows that $\uline \#$ is strictly compatible with the Gray product of 1-morphisms,  preserves the unit-1-morphisms,  satisfies $\uline\#\,\uline\#\uline F=\uline F$ and is  compatible with the evaluation: 
  $e(\uline\#\,\uline F)=\#e(\uline F)$ for all 1-morphisms $\uline F$ in $\uline {\mac G}$.
  
  For  basic 2-morphisms $\uline\alpha=(\alpha,z):\uline F\Rightarrow\uline G$ we set $\uline\# (\alpha, z)=(\alpha,-z)$ and extend $\uline\#$ to general 2-morphisms by
  $$
  \uline\#(\uline{\alpha_n},\ldots,\uline{\alpha_1})=(\uline\#\,\uline{\alpha_1},\ldots, \uline\#\,\uline{\alpha_n}),\quad \uline\#(\emptyset_{\uline F})=\emptyset_{\uline\#\uline F}.
  $$
  This implies that $\uline\#$ is strictly compatible with the horizontal composition of 2-morphisms, preserves the unit 2-morphisms and satisfies  $\uline\#\,\uline\#\uline \alpha=\uline \alpha$ for all 2-morphisms $\uline \alpha$ in $\uline{\mac G}$.  Due to the identity $\#\#F=1$ for all 1-morphisms $F$ in $\mac G$, it also follows that $\uline\#$ is strictly compatible with the Gray product of 1- and 2-morphisms.

  To define the action of  $\uline\#$ on 3-morphisms $\uline\Gamma:\uline\alpha\Rrightarrow\uline\beta$, we note that if $\uline\alpha=(\alpha,z)$ is a basic 2-morphism, 
  the  2-morphisms $\#e(\uline\alpha)$ and $e(\uline\#\,\uline\alpha)$ are related by $\# e(\alpha,1)=e(\uline\#(\alpha,1))$ and $\# e(\alpha,-1)=\#\#e(\uline\#(\alpha,-1))$. 
  For general 2-morphisms $(\uline{\alpha_n},...,\uline{\alpha_1})$, we have
  \begin{align*}&\#e(\uline{\alpha_n},...,\uline{\alpha_1})=\#(e(\uline{\alpha_n})\circ\ldots\circ e(\uline{\alpha_1})),\\
    &e(\uline\#(\uline{\alpha_n},...,\uline{\alpha_1}))=e(\uline\#\,\uline{\alpha_1})\circ\ldots\circ e(\uline\#\,\uline{\alpha_n})
    .\end{align*}
  We  obtain an invertible 3-morphism 
  $\chi_{\uline\alpha}: \#e(\uline\alpha)\Rrightarrow e(\uline\#\,\uline\alpha)$  by setting
  \begin{align*}
    &\chi_{\emptyset_{\uline F}}=\Phi^\inv_{e(\uline F)},\quad \chi_{(\alpha, 1)}=1_{\#\alpha}, \quad \chi_{(\alpha,-1)}=\Theta_\alpha\\
    &\chi_{(\uline\alpha_n,...,\uline\alpha_1)}=(\chi_{\uline\alpha_1}\circ \ldots\circ \chi_{\uline\alpha_n})\cdot \Phi^\inv_{e(\uline\alpha_1),...,e(\uline\alpha_n)},
  \end{align*}
  where $\Phi_{\mu_1,...,\mu_n}:  \#\mu_1\circ\ldots\circ \#\mu_n\Rrightarrow \#(\mu_n\circ\ldots\circ\mu_1)$ 
  denotes the invertible
  3-morphism  determined by  $\Phi_{\mu_1}=1_{\#\mu_1}$ and
  \begin{align}
    &\Phi_{\mu_1,...,\mu_n}=\Phi_{\mu_1,\mu_n\circ ...\circ \mu_2} \cdot (1_{\#\mu_1}\circ \Phi_{\mu_2, \mu_n\circ...\circ \mu_3})\cdots \nonumber\\
    &\cdot (1_{\#\mu_1}\circ ...\circ 1_{\#\mu_{n-3}}\circ \Phi_{\mu_{n-2}, \mu_n\circ\mu_{n-1}})\cdot (1_{\#\mu_1}\circ ...\circ 1_{\#\mu_{n-2}}\circ \Phi_{\mu_{n-1}, \mu_n})\nonumber
  \end{align}
  and the 3-morphisms $\Phi_F: 1_{F^\#}\Rrightarrow \#1_F$ and $\Phi_{\mu,\nu}: \#\mu\circ\#\nu\Rrightarrow \#(\nu\circ\mu)$ are given in the proof of Theorem \ref{graycatduals}.
  Note that it follows from identity \eqref{compat2}  that the bracketing in the definition of $\Phi_{\mu_1,...,\mu_n}$ is irrelevant, and for all composable 2-morphisms $\mu_1$ ,..., $\mu_n$ and all $1\leq k\leq n-1$
  \begin{align*}
    \Phi_{\mu_1,\ldots,\mu_n}=\Phi_{\mu_k\circ \ldots\circ\mu_1, \mu_n\circ \ldots\circ \mu_{k+1}}\cdot (\Phi_{\mu_1, \ldots,  \mu_k}\circ \Phi_{\mu_{k+1}, \ldots,  \mu_{n}}).
  \end{align*}
  From this, it follows  that the 3-morphisms $\chi_{\uline \alpha}: \#e(\uline\alpha)\Rrightarrow e(\uline\#\,\uline\alpha)$ satisfy the relation
  \begin{align}\label{chicomp}
    \chi_{\uline\alpha\circ\uline\beta}=(\chi_{\uline\beta}\circ \chi_{\uline\alpha})\cdot \Phi^\inv_{e(\uline\beta), e(\uline\alpha)}
  \end{align}
  for all composable 2-morphisms $\uline\alpha$, $\uline\beta$.

  For a 3-morphism $\uline\Gamma:\uline\alpha\Rrightarrow\uline\beta$, we define  $\uline\#\,\uline\Gamma:\uline\#\,\uline\alpha\Rrightarrow\uline\#\,\uline\beta$ by
  $$
  e(\uline\#\,\uline\Gamma)=\chi_{\uline\beta} \cdot \#e(\uline\Gamma)\cdot \chi_{\uline\alpha}^\inv.
  $$
  To show that $\uline\#$ defines a strict functor of 2-strict tricategories, we  prove the identities
  $$
  \uline\#(\uline\Psi\,\uline\cdot\,\uline \Phi)=\uline\#\,\uline\Psi\,\uline\cdot\, \uline\#\,\uline \Phi, \; \uline\#(\uline\Psi\,\uline\circ\,\uline \Phi)=\uline\#\,\uline\Phi\,\uline\circ\, \uline\#\,\uline \Psi,\; \uline\#(\uline\Psi\,\uline\Box\,\uline \Phi)=\uline\#\,\uline\Phi\,\uline\Box\, \uline\#\,\uline \Psi$$
  for all 3-morphisms $\uline\Phi,\uline\Psi$ for which these expressions are defined. The first follows directly  from the definition. The identity $\uline\#(\uline\Psi\,\uline\circ\,\uline \Phi)=\uline\#\,\uline\Phi\,\uline\circ\, \uline\#\,\uline \Psi$ follows from  the
  commutative diagram
  \begin{align}\label{hashstrict1}
    &\xymatrix{ e(\uline\#(\uline \alpha\,\uline\circ\, \uline\beta))  \ar[d]_{\chi^\inv_{\uline\alpha\,\uline\circ\,\uline\beta}}  \ar[r]^{\chi_{\uline\beta}^\inv\circ \chi_{\uline\alpha}^\inv} & \#e(\uline\beta)\circ \#e(\uline\alpha) \;\;\ar[r]^{\#e(\uline \Phi)\circ \#e(\uline\Psi)} \ar[ld]^{\Phi_{e(\uline\beta),e(\uline\alpha)}}& \;\;\#e(\uline\beta')\circ \#e(\uline\alpha') \ar[d]^{\chi_{\uline\beta'}\circ \chi_{\uline\alpha'}} \ar[ld]^{\Phi_{e(\uline\beta'), e(\uline\alpha')}}\\
      \#e(\uline\alpha\,\uline \circ\, \uline\beta) \ar[r]_{\#(e(\uline\Psi\uline\circ \uline\Phi))} & \;\;\#e(\uline \alpha'\,\uline\circ\, \uline\beta') \ar[r]_{\chi_{\uline\alpha'\,\uline\circ\,\uline\beta'}} & e(\uline\#(\uline \alpha'\,\uline\circ\, \uline\beta'))
    }
  \end{align}
  for all 3-morphisms $\uline\Phi: \uline\beta\Rrightarrow\uline\beta'$ and $\uline\Psi:\uline\alpha\Rrightarrow \uline\alpha'$. The parallelogram 
  in the middle of this diagram  commutes due to the naturality of  the 3-morphism $\Phi_{\mu,\nu}: \#\mu\circ\#\nu\Rrightarrow \#(\nu\circ\mu)$ and the triangles on the left and right by identity \eqref{chicomp}.
  
  The identity $\uline\#(\uline\Psi\,\uline\Box\,\uline \Phi)=\uline\#\,\uline\Phi\,\uline\Box\, \uline\#\,\uline \Psi$ is equivalent to the commutativity of the diagram
  \begin{align}
    &\label{hashstrict2}
    \xymatrix{ 
      e(\uline\#(\uline\alpha\uline\Box\uline\beta)) 
      \ar@/ _5pc/[dddd]_{\tilde\iota^\inv_{\uline\#\uline\beta, \uline\#\uline F}\circ \tilde\iota^\inv_{\uline\#\uline K,\uline\#\uline\alpha}\;\;} 
      \ar[r]^{e(\uline \#(\uline\Psi\uline\Box\uline\Phi))} 
      & e(\#(\uline\alpha'\uline\Box\uline\beta')) 
      \ar@/ ^5pc/[dddd]^{\;\;\tilde\iota^\inv_{\uline\#\uline\beta', \uline\#\uline F}\circ \tilde\iota^\inv_{\uline\#\uline K,\uline\#\uline\alpha'}}
      \\
      \#e(\uline\alpha\uline\Box\uline\beta)   
      \ar[u]_{\chi_{\uline\alpha\uline\Box\uline\beta}}  
      \ar[d]_{\#(\tilde \iota^\inv_{\uline\alpha,\uline K}\circ \tilde\iota^\inv_{\uline F,\uline\beta})} \ar[r]^{\#e(\uline\Psi\uline\Box\uline\Phi)} 
      & \#e(\uline\alpha'\uline\Box\uline\beta') 
      \ar[d]_{\#(\tilde\iota^\inv_{\uline\alpha',\uline K}\circ\tilde\iota^\inv_{\uline F,\uline\beta'})}\ar[u]^{\chi_{\uline\alpha'\uline\Box\uline\beta'}} 
      \\
      \#(e(\uline\alpha)\Box e(\uline\beta))\;
      \ar[r]^{\#(e(\uline\Psi)\Box e(\uline\Phi))} 
      & \;\#(e(\uline\alpha')\Box e(\uline\beta')) 
      \\
      \#e(\uline\beta)\Box \#e(\uline\alpha) \;
      \ar[u]_{\kappa_{ e(\uline\alpha),e(\uline\beta)}} 
      \ar[r]^{\#e(\uline\Phi)\Box \# e(\uline\Psi)} 
      \ar[d]^{\chi_{\uline\beta}\Box \chi_{\uline\alpha}}
      &  \;\#e(\uline\beta')\Box \#e(\uline\alpha')
      \ar[u]^{\kappa_{ e(\uline\alpha'),e(\uline\beta')}} 
      \ar[d]_{\chi_{\uline\beta'}\Box \chi_{\uline\alpha'}}
      \\
      e(\uline\#\uline\beta)\Box e(\uline\#\uline\alpha) \;\; 
      \ar[r]^{e(\uline\#\uline\Phi)\Box e(\uline\#\uline\Psi)}
      &\;\; e(\uline\#\uline\beta')\Box e(\uline\#\uline\alpha').  
    }
  \end{align}
  In this diagram,  the four rectangles in the middle  commute by definition of the 3-morphisms $\uline\#(\uline\Psi\uline\Box\uline\Phi)$, $\uline\Psi\uline\Box\uline\Phi$ and $\uline\#\Psi$, $\uline\#\Phi$ and due to the naturality of the 3-morphism $\Phi_{\mu,\nu}: \#\mu\circ\#\nu\Rrightarrow\#(\nu\circ\mu)$. It remains to show that the two pentagons on the left and the right commute. 
  As a first step, we reduce the proof of the commutativity of these diagrams  to the cases $\uline\alpha=\emptyset_{\uline F}$ or $\uline\beta=\emptyset_{\uline K}$. For this, we consider the diagram 
  \begin{align}\label{movediag2}
    \xymatrix{
      e(\uline\#(\uline\alpha\uline\Box\uline\beta)) 
      \ar[d]^{\tilde\iota_{\uline\#\uline\beta,\uline \#\uline F}^\inv\circ \tilde\iota_{\uline \#\uline K,\uline\#\uline\alpha}^\inv }
      \ar[rr]^{\chi^\inv_{\uline F\uline\Box\uline\beta}\circ \chi^\inv_{\uline\alpha\uline\Box\uline K}}
      & &  \# e(\uline F\uline\Box\uline\beta)\circ \#e(\uline\alpha\uline\Box\uline K)
      \ar[lddd]^(.7){\#\tilde\iota^\inv_{\uline F,\uline\beta}\circ \#\tilde\iota^\inv_{\uline\alpha,\uline K}}
      \ar[ddd]^{\Phi_{e(\uline F\uline \Box\uline\beta), e(\uline\alpha\uline\Box\uline K)}}
      \\
      e(\uline\#\uline\beta)\Box e(\uline\#\uline\alpha)
      \ar[d]_{\chi^\inv_{\uline\beta}\Box\chi^\inv_{\uline\alpha}}
      & {\begin{array}{l}
          (\#e(\uline\beta)\Box \#1_{e(\uline F)}) 
          \\
          \quad \circ (\#1_{e(\uline K)}\Box \#e(\uline\alpha))\end{array}}
      \ar[dd]^(.3){\substack{  {\kappa_{ 1_{e(\uline F)},e(\uline\beta)}}\\{\circ \kappa_{ e(\uline\alpha),1_{e(\uline K)}}}}}
      \\
      \# e(\uline\beta)\Box \#e(\uline\alpha)
      \ar[ru]_{
        \substack{{(\#e(\uline\beta)\Box \Phi_{e(\uline F)})}\\
          {\circ(\Phi_{e (\uline K)}\Box \#e(\uline\alpha))}}}
      \ar[d]_{\kappa_{ e(\uline\alpha),e(\uline\beta)}}
      & 
      \\
      \#(e(\uline\alpha)\Box e(\uline\beta))\quad
      \ar@/ _3pc/[rr]_{\#(\tilde\iota_{\uline\alpha,\uline K}\circ\tilde\iota_{\uline F,\uline\beta})}
      &  {\begin{array}{l} \#( e(\uline F)\Box e(\uline\beta))\\
          \circ\#(e(\uline\alpha)\Box e(\uline K))\end{array}}
      \ar@/  _1pc/[l]_{\Phi_{e(\uline F)\Box e(\uline\beta), e(\uline\alpha)\Box e(\uline K)}}
      & \#e(\uline\alpha\uline\Box\uline\beta).
    }
  \end{align}
  The path on the outside in this diagram corresponds to the pentagon in diagram \eqref{hashstrict2}.
  The two quadrilaterals in the diagram  commute by the naturality of the 3-morphism $\Phi_{\mu,\nu}:\#\mu\circ\#\nu\Rrightarrow\#(\nu\circ\mu)$ and by identity \eqref{kappanatisom}. The hexagon commutes if and only if the pentagon in diagram \eqref{hashstrict2} commutes for the case where $\uline\alpha=\emptyset_{\uline F}$ or $\uline\beta=\emptyset_{\uline K}$.

  It  is therefore sufficient to prove that the pentagon in the diagram \eqref{hashstrict2} commutes   for   $\uline\alpha=\emptyset_{\uline F}$ or $\uline \beta=\emptyset_{\uline K}$. 
  In the latter, it  reduces to the  diagram
  \begin{align}\label{hashtriv}
    \xymatrix{ 
      e(\uline\#(\uline\alpha\uline\Box \uline K))
      \ar[r]^{\chi^\inv_{\uline\alpha\uline\Box \uline K}}
      \ar[d]_{\tilde\iota^\inv_{\uline\#\uline K,\uline\#\uline\alpha}}
      & \#e(\uline\alpha\uline\Box\uline K)
      \ar[r]^{\#\tilde\iota^\inv_{\uline\alpha,\uline K}}
      & \#(e(\uline\alpha)\Box e(\uline K))
      \ar[d]^{\kappa^\inv_{ e(\uline\alpha),1_{e(\uline K)}}}
      \\
      \#e(\uline K)\Box e(\uline\#\uline\alpha)
      \ar[r]_{\#e(\uline K)\Box \chi^\inv_{\uline\alpha}}
      & \;\;\;\#e(\uline K)\Box \#e(\uline\alpha)\;\;\;
      \ar[r]_{\Phi_{e(\uline K)}\Box \#e(\uline\alpha)}
      & \;\;\;\#1_{e(\uline K)}\Box \#e(\uline\alpha).
    }
  \end{align}
  We start by proving that this diagram commutes for basic 2-morphisms and the empty string of 2-morphisms. For $\uline\alpha=\emptyset_{\uline F}$, the 3-morphism $\tilde\iota_{\uline\alpha,\uline K}$ is trivial, and we have 
  $$\chi_{\uline\alpha}=\Phi_{e(\uline F)}^\inv, \quad \chi_{\uline\alpha\uline\Box\uline K}=\Phi_{e(\uline F\Box\uline K)}^\inv, \quad \kappa_{1_{ e(\uline\alpha),e(\uline K)}}=\Phi_{e(\uline F)\Box e(\uline K)}\cdot (\Phi^\inv_{e(\uline K)}\Box \Phi^\inv_{e(\uline F)}), $$ where the last identity follows directly with  \eqref{kappaunit}. Inserting this into \eqref{hashtriv}, one finds that the diagram commutes.
  For $\uline\alpha=(\alpha,1)$, the 3-morphisms $\chi_{\uline\alpha}$, $\chi_{\uline\alpha\uline\Box\uline K}$, $\tilde\iota_{\uline\alpha,\uline K}$ are trivial and  
  $$\tilde\iota_{\uline\#\uline K,\uline\#,\uline\alpha}=\kappa^\inv_{ e(\uline\alpha),1_{e(\uline K)}}\cdot (\Phi_{e(\uline K)}\Box \#e(\uline\alpha)),$$ which shows that the diagram commutes. 
  For $\uline\alpha=(\alpha,-1)$ 
  diagram \eqref{hashtriv} corresponds to the boundary of the following diagram
  $$
  \xymatrix{
    \#e(\uline K)\Box\alpha 
    & \#\#(\#e(\uline K)\Box \alpha)
    \ar[r]^{\#\kappa^\inv_{ 1_{\#e(\uline K)},\alpha}}
    \ar[l]_{\Theta_{\#e(\uline K)\Box\alpha}}
    & \#(\#\alpha\Box\#1_{\#e(\uline K)})
    \ar[ld]_(.7){\kappa^\inv_{\alpha,\#1_{\#e(\uline K)}}}
    \ar[d]^{\#(\#\alpha\Box \Phi^\inv_{\#e(\uline K)})}
    \\
    \#1_{e(\uline K)}\Box\#\#\alpha 
    \ar[u]^{\Phi^\inv_{e(\uline K)}\Box \Theta_\alpha}
    & \#\#1_{\#e(\uline K)}\Box \#\#\alpha
    \ar[l]_{\#\Phi^\inv_{\#e(\uline K)}\Box \#\#\alpha}
    &\#(\#\alpha\Box 1_{e(\uline K)}),
    \ar@/  ^2pc/[ll]^{\kappa^\inv_{\#\alpha,1_{e(\uline K)}}}
  }
  $$
  where we used identities \eqref{twistun} and \eqref{twistbox} to express $\Theta_{\#e(\uline K)\Box \alpha}$ in terms of $\Theta_\alpha$.  The pentagon in this diagram commutes 
  due to identities \eqref{twistun} and \eqref{twistbox} and the quadrilateral due to the naturality of the 3-morphism $\kappa_{\mu,\nu}: \#\nu\Box\#\mu\Rrightarrow\#(\mu\Box\nu)$. This proves that diagram \eqref{hashtriv}  commutes for basic 2-morphisms and the empty string of 2-morphisms.

  The proof that this identity holds for general 2-morphisms $\uline\alpha=(\uline\alpha_n,...,\uline\alpha_1)$ is by induction over $n$. For $n=1$, $\uline\alpha$ is a basic 1-morphism and this identity was shown above. Suppose that  the commutativity of diagram \eqref{hashtriv} is established for all strings  $\uline\alpha=(\uline\alpha_k,...,\uline\alpha_1)$ of basic 2-morphisms of length  $k\leq n-1$ and let $\uline\gamma=(\uline\gamma_n,...,\uline\gamma_1)$ be a 2-morphism of length $n$. Set $\uline\alpha=\uline\gamma_n$, $\uline\beta=(\uline\gamma_{n-1},...,\uline\gamma_1)$ and consider the
  diagram 
  \begin{align}
    \label{diag1}
    &\xymatrix{
      { \begin{array}{l} \#e(\uline\alpha \uline\Box \uline K) \\ \circ \#e(\uline\beta\uline\Box\uline K) \end{array}}
      \ar[d]_{\chi_{\uline\alpha\uline\Box\uline K}\circ\chi_{\uline\beta\uline\Box\uline K}}
      \ar[r]^{ \Phi_{e(\uline\alpha\uline\Box\uline K),e(\uline\beta\uline\Box\uline K)}}
      & { \begin{array}{l} \#(e(\uline\beta\uline\Box \uline K) \\
          \circ e(\uline\alpha\uline\Box\uline K))\end{array}}
      & {\begin{array}{l} \#( (e(\uline\beta)\Box e(\uline K))\\
          \circ(e(\uline\alpha)\Box e(\uline K)))\end{array}}
      \ar[l]_{\#(\tilde\iota_{\uline\beta,\uline K}\circ \tilde\iota_{\uline\alpha,\uline K} )\;\;}
      \ar[ld]_(.6){\Phi^\inv_{e(\uline\alpha)\Box e(\uline K), e(\uline\beta)\Box e(\uline K)}\qquad}
      \ar[d]^{\kappa^\inv_{e(\uline\alpha\uline\circ\uline\beta),1_{e(\uline K)}}}
      \\
      {\begin{array}{l} e(\uline\#(\uline\alpha\uline\Box\uline K))\\
          \circ e(\uline\#(\uline\beta\uline\Box\uline K))\end{array}}
      \ar[d]_{\tilde \iota^\inv_{\uline\#\uline K,\uline\#\uline\alpha}\circ \tilde\iota^\inv_{\uline\#\uline K,\uline\#\uline\beta}}
      &{\begin{array}{l} \#(e(\uline\alpha)\Box e(\uline K))\\
          \circ \#(e(\uline\beta)\Box e(\uline K))\end{array}}
      \ar[lu]^{\#\tilde\iota_{\uline\alpha,\uline K}\circ \#\tilde\iota_{\uline\beta,\uline K}}
      \ar[d]_{\kappa^\inv_{ e(\uline\alpha),1_{e(\uline K)}}\circ \kappa^\inv_{ e(\uline\beta),1_{e(\uline K)}}}
      &{\begin{array}{l}\#1_{e(\uline K)}\Box\\ \#(e(\uline\alpha)\circ e(\uline\beta))\end{array}}
      \\
      {\begin{array}{l}\#e(\uline K)\Box\\
          (e(\uline\#\uline\alpha)\circ e(\uline\#\uline\beta))\end{array}}
      \ar[rd]_{\#e(\uline K)\Box (\chi_{\uline\alpha}^\inv \circ \chi^\inv_{\uline\beta})\qquad}
      & {\begin{array}{l}(\#1_{e(\uline K)}\Box \#e(\uline\alpha))\\
          \circ(\#1_{e(\uline K)}\Box \#e(\uline\beta))\end{array}}
      & {\quad \begin{array}{l}\#e(\uline K)\Box\\
          \#(e(\uline\beta)\circ e(\uline\alpha))\end{array}}
      \ar[u]^{ \quad \Phi_{e(\uline K)}\Box \#(e(\uline\beta)\circ e(\uline\alpha))}
      \\
      &
      {\begin{array}{l}\#e(\uline K)\Box \\
          (\#e(\uline\alpha)\circ \#e(\uline\beta))\end{array}}
      \ar[u]_(.6){\Phi_{e(\uline K)}\Box ( \#e(\uline\alpha)\circ  \#e(\uline\beta)}
      \ar@/ _1pc/[ru]_{\quad\#e(\uline K)\Box \Phi_{e(\uline\alpha), e(\uline\beta)}}
      & 
    }
  \end{align}
  The outer path in this diagram corresponds to the diagram  \eqref{hashtriv} for $\uline\gamma$. The quadrilateral at the top of the diagram commutes due to the naturality of the 3-morphisms $\Phi_{\mu,\nu}: \#\mu\circ\#\nu\Rrightarrow\#(\nu\circ\mu)$, $\kappa_{\mu,\nu}:\#\nu\Box\#\mu\Rrightarrow\#(\mu\Box\nu)$ and $\Phi_F: 1_{\#F}\Rrightarrow\#1_F$. 
  The hexagon on the left of this diagram commutes because identity \eqref{hashtriv} holds for $\uline\alpha$ and $\uline\beta$.  To show that the hexagon on the right of the diagram commutes, we set $\nu=\tau=1_{K}$ in \eqref{kappanatisom} and use the naturality of the tensorator. This yields  the diagram 
  \begin{align*}
    &\xymatrix{ 
      (\#1_K\Box \#\mu)\circ(\#1_K\Box \#\rho)
      \ar[rr]^{\kappa_{\mu,1_K}\circ \kappa_{\rho,1_K}}
      \ar[d]_{\#1_K\Box \sigma^\inv_{\#1_K,\#\mu}\Box \#\rho}
      &
      & \#(\mu\Box K)\circ \#(\rho\Box K)
      \ar[d]_{\Phi_{\mu\Box K,\rho\Box K}}
      \\
      \#1_K\Box (\#\mu\circ\#\rho) 
      \ar[d]_{\#1_K\Box\Phi_{\mu,\rho}}
      & \#K\Box(\#\mu\Box\#\rho)
      \ar[lu]_{\qquad(\Phi_K\Box \#\mu)\circ (\Phi_K\Box\#\rho)}
      \ar[l]_{\Phi_K\Box (\#\mu\circ\#\rho)}
      &\#((\rho\circ\mu)\Box K)
      \\
      \#1_K\Box\#(\rho\circ\mu).
      \ar[rru]_{\kappa_{ \rho\circ\mu,1_{K}}}
    }
  \end{align*}
  The  triangle in this diagram  commutes due to the naturality of the tensorator and the outer pentagon due to identity  \eqref{kappanatisom}. This implies that the inner hexagon commutes as well and hence the hexagon on the right in diagram \eqref{diag1}.
  This shows that identity \eqref{hashtriv} holds for $\uline\gamma$ and
  concludes the proof  that $\uline\#:\uline{\mac G}\to\uline{\mac G}_{op}$ is a strict functor of 2-strict tricategories for which all coherence data is trivial.

  \medskip
  To prove the identity $\uline\#\uline\#=1$, note that it is obvious for 1- and 2-morphisms. For 3-morphisms $\uline\Psi:\uline\alpha\Rrightarrow\uline\alpha'$, it holds if the following diagram commutes
  \begin{align}\label{hashhash}
    \xymatrix{
      {\begin{array}{l} e(\uline\#\uline\#\uline\alpha)\\
          =e(\uline\alpha)\end{array}} 
      \ar@/ ^3pc/[rrr]^{e(\uline\#\uline\#\uline \Psi)}
      \ar@/ _8pc/[rrr]_{e(\uline\Psi)}
      \ar[r]^{\chi^\inv_{\uline\#\uline\alpha}}
      & \#e(\uline\#\uline\alpha) 
      \ar[d]^{\#\chi^\inv_{\uline\alpha} } 
      \ar[r]^{\#e(\uline\#\uline\Psi)} 
      & \#e(\uline\#\uline\alpha') 
      \ar[r]^{\chi{\uline\#\uline\alpha'}} 
      & {\begin{array}{l}e(\uline\#\uline\#\uline\alpha')\\
          =e(\uline\alpha')\end{array}}
      \\
      & \#\#e(\uline\alpha)
      \ar[lu]_{\Theta_{e(\uline\alpha)}} 
      \ar[r]_{\#\#e(\uline\Psi)}  
      & \#\#e(\uline\alpha')  
      \ar[u]^{\#\chi_{\uline\alpha'} }
      \ar[ru]^{\Theta_{e(\uline\alpha')}}   
    }
  \end{align}
  The three quadrilaterals in  the diagram commute by definition of the 3-morphism $\uline\#\,\uline\#\Psi$, $\uline\#\uline\Psi$ and  by naturality of the 3-morphism $\Theta$. The triangles commute if and only if for all 2-morphisms $\uline\alpha$, we have
  \begin{align}\label{hhidentity}
    \Theta_{e(\uline\alpha)}=\chi_{\uline\#\uline\alpha}\cdot \#\chi_{\uline\alpha}.
  \end{align}
  For $\uline\alpha=\emptyset_{\uline F}$ this follows directly from identity \eqref{twistun}, which implies $$\Theta_{e(\emptyset_{\uline F})}=\Theta_{1_{e(\uline F)}}=\Phi^\inv_{\#e(\uline F)}\cdot \#\Phi^\inv_{e(\uline F)}=\chi_{\uline\#\emptyset_{\uline F}}\cdot \#\chi_{\emptyset_{\uline F}}.$$
  Similarly,  for basic 2-morphisms $\uline\alpha=(\alpha, z)$, we have
  \begin{align*}
    &\chi_{\uline\#(\alpha,1)}\cdot \#\chi_{(\alpha,1)}=\Theta_\alpha\cdot \#1_{\#\alpha}=\Theta_{\alpha}=\Theta_{e(\alpha,1)}\\
    &\chi_{\uline\#(\alpha,-1)}\cdot \#\chi_{(\alpha,-1)}=1_{\#\alpha}\cdot \#\Theta_\alpha=\Theta_{\#\alpha}=\Theta_{e(\alpha, -1)},
  \end{align*}
  where we used the identity $\#\Theta_\alpha=\Theta_{\#\alpha}$ from Lemma \ref{coherence} in the second line. 
  The proof that identity \eqref{hhidentity} holds for strings $\uline\alpha=(\uline\alpha_n,...,\uline\alpha_1)$ of basic 2-morphisms is by induction over the length $n$ of the string. For $n=0,1$ this was shown above. Suppose that the identity \eqref{hhidentity} is established for all 2-morphisms $\uline\alpha=(\uline\alpha_K,...,\uline\alpha_1)$ of length $0\leq k<n$ and let $\uline\gamma=(\uline\gamma_n,...,\uline\gamma_1)$ be a string of basic 2-morphisms of length $n$. Then identity \eqref{hhidentity} holds for $\uline\gamma$ if and only if the following diagram 
  commutes for  $\uline\alpha=\uline\gamma_n$ and $\uline\beta=(\uline\gamma_{n-1},...,\uline\gamma_1)$ 
  \begin{align*}
    \xymatrix{
      \#\#(e(\uline\alpha)\circ e(\uline\beta)) 
      \ar@/ _8pc/[dd]_{\Theta_{e(\uline\alpha\uline\circ \uline\beta)}}
      \ar[r]^{\#\Phi^\inv_{e(\uline\beta), e(\uline\alpha)}}  
      &\#(\#e(\beta)\circ \#e(\alpha))
      \ar[d]^{\#(\chi_{\uline\beta}\circ \chi_{\uline\alpha})}
      \ar[ld]_(.7){\Phi^\inv_{\#e(\uline\alpha),\#e(\uline\beta)}\quad}
      \\
      \#\# e(\uline\alpha)\circ \#\# e(\uline\beta) 
      \ar[d]_{\Theta_{e(\uline\alpha)}\circ\Theta_{e(\uline\beta)}}
      \ar[rd]^(.7){\#\chi_{\uline\alpha}\circ \#\chi_{\uline\beta}}
      & \#(e(\uline\#\uline\beta)\circ e(\uline\#\uline\alpha))
      \ar[d]^{\Phi^\inv_{e(\uline\#\uline\alpha),e(\uline\#\uline\beta)}}
      \\
      e(\uline\alpha)\circ e(\uline\beta) & 
      \quad  \#e(\uline\#\uline\alpha)\circ \#e(\uline\#\uline\beta))
      \ar[l]^{\chi_{\uline\#\uline\alpha}\circ \chi_{\uline\#\uline\beta}\quad}
    }.
  \end{align*}
  The triangle at the bottom of the diagram commutes by induction hypothesis.
  The curved subdiagram at the left commutes due to identity \eqref{twisthorcomp} and the
  subdiagram at the right  due to the naturality of the 3-morphism $\Phi_{\mu,\nu}: \#\mu\circ\#\nu\Rrightarrow\#(\nu\circ\mu)$. This shows that the diagram \eqref{hashhash} commutes for all 3-morphisms $\uline\psi:\uline\alpha\Rrightarrow\uline\alpha'$  in $\uline{\mac G}$ and that  the functor $\uline\#:\uline{\mac G}\to\uline{\mac G}_{op}$ satisfies $\uline\#\uline\#=1$.

  \bigskip
  3. To define the strict functor of 2-strict tricategories $\uline *:\uline{\mac G}\to\uline{\mac G}^{op}$, we set $\uline *$ to be trivial on the  objects and 1-morphisms of $\uline{\mac G}$. For a basic 2-morphism $(\alpha,z)$ we set
  $\uline*(\alpha,z)=(\alpha^*, z)$ and extend $\uline *$ to general 2-morphisms via
  $$
  \uline*(\uline{\alpha_n},...,\uline{\alpha_1})=(\uline*\,\uline{\alpha_1},..., \uline*\,\uline{\alpha_n}),\quad \uline*(\emptyset_{\uline F})=\emptyset_{\uline F}.
  $$
  It follows that $\uline*$ is strictly compatible with the Gray product of 1- and 2-morphisms as well as the horizontal composition, preserves the unit 1- and 2-morphisms and satisfies $\uline*\,\uline*(\uline \alpha)=\uline\alpha$ for all 2-morphisms $\uline\alpha$.   We also have the identities
  \begin{align*}
    &*e(\emptyset_{\uline F})=e(\uline * \emptyset_{\uline F})=1_{e(\uline F)},  &  &*e(\alpha,1)=e(\uline *(\alpha,1))=\alpha^*,\\
    &*e(\alpha,-1)=*\#\alpha, & &e(\uline *(\alpha, -1))=\#*\alpha.\end{align*}
  To obtain a 3-morphism $\uline *\uline\Psi: \uline*\uline\alpha'\Rrightarrow \uline *\uline\alpha$ for each 3-morphism  $\uline\Psi: \uline\alpha\Rrightarrow \uline\alpha'$ we consider
  the 3-morphism  $\xi_{\uline\alpha}: *e(\uline\alpha)\Rrightarrow e(\uline*\,\uline\alpha)$ given by 
  \begin{align*}
    &\xi_{\emptyset_{\uline F}}=1_{1_{e(\uline F)}}, \quad \xi_{(\alpha,1)}=1_{*\alpha},\quad \xi_{(\alpha,-1)}=\Delta_{*\alpha}^\inv\\
    &\xi_{(\uline{\alpha_n},...,\uline{\alpha_1})}=\xi_{\uline{\alpha_1}}\circ ... \circ \xi_{\uline{\alpha_n}}
  \end{align*} 
  and set  
  $$e(\uline *\uline\Psi)= \xi_{\uline\alpha}\cdot *e(\uline\Psi)\cdot \xi_{\uline\alpha'}^\inv.$$
  The strict compatibility of $\uline *$ with the vertical composition of 3-morphisms is a direct consequence of the definition. The strict compatibility of $\uline *$ with the horizontal composition  is equivalent to the commutativity of the diagram
  \begin{align*}
    &\xymatrix{ e(\uline{*}(\uline \alpha\,\uline\circ\, \uline\beta))   & \ar[l]_{\xi_{\uline\beta}\circ \xi_{\uline\alpha}}  \ar[ld]^{1_{*e(\uline\alpha\,\uline\circ\,\uline\beta)}} {*}e(\uline\beta)\circ {*}e(\uline\alpha) \;\;& \ar[l]_{{*}e(\uline \Phi)\circ {*}e(\uline\Psi)}   \ar[ld]^{1_{*e(\uline\alpha'\,\uline\circ\,\uline\beta')}}\;\;{*}e(\uline\beta')\circ {*}e(\uline\alpha')\\
      {*}e(\uline\alpha\,\uline \circ\, \uline\beta) \ar[u]^{\xi_{\uline\alpha\,\uline\circ\,\uline\beta}} &  \ar[l]^{{*}(e(\uline\Psi \uline \circ \uline\Phi))} \;\;{*}e(\uline \alpha'\,\uline\circ\, \uline\beta')  &   \ar[l]^{\xi^\inv_{\uline\alpha'\,\uline\circ\,\uline\beta'}} e(\uline{*}(\uline \alpha'\,\uline\circ\, \uline\beta') ) \ar[u]_{\xi^\inv_{\uline\beta'}\circ \xi^\inv_{\uline\alpha'}}.
    }.\end{align*}
  By definition, the paths on the boundary correspond to the 3-morphisms $e(\uline*(\uline\Psi\uline\circ\uline\Phi))$ and $e(\uline*\uline\Phi\uline\circ \uline *\uline\Psi)$, and
  the parallelogram in the middle commutes due to the identity $e(\uline\Psi\uline\circ \uline\Phi)=e(\uline\Psi)\circ e(\uline\Phi)$. As we have $\xi_{\uline\beta}\circ \xi_{\uline\alpha}=\xi_{\uline\alpha\uline\circ\uline\beta}$ for all 2-morphisms $\uline\alpha,\uline\beta$ by definition, the diagram commutes and we obtain  $\uline*(\uline\Psi\uline\circ\uline\Phi)=e(\uline*\uline\Phi\uline\circ \uline *\uline\Psi)$ for all  composable 3-morphisms $\uline\Psi,\uline\Phi$. The identity $\uline *\uline *\uline\Psi=1$ then follows from the identity 
  $\Delta^*_{*\alpha}=\Delta_\alpha$ for all 2-morphisms $\alpha$ in a spatial Gray category $\mac G$ and the compatibility of $\uline *$ with the horizontal composition.

  The strict compatibility of the functor $\uline *$ with the Gray product  corresponds to the commutativity of the diagram 
  \begin{align}\label{starbox}
    &\xymatrix{
      e(\uline *(\uline\alpha\uline \Box\uline\beta)) 
      \ar[r]^{\xi_{\uline\alpha\uline\Box\uline\beta}^\inv} 
      \ar[rdd]^(.3){\tilde\iota^\inv_{\uline F,\uline*\uline\beta}\circ \tilde\iota^\inv_{\uline K,\uline*\uline \alpha}} 
      & {*e}(\uline\alpha\uline\Box\uline\beta)  
      \ar[d]^{\!\!\!{*}(\tilde \iota_{\uline \alpha, \uline K}\circ \tilde\iota_{\uline F, \uline\beta})} 
      & {*e}(\uline\alpha'\uline\Box\uline\beta') 
      \ar[d]_{{*}(\tilde \iota_{\uline \alpha', \uline K}\circ \tilde\iota_{\uline F, \uline\beta'})} 
      \ar[l]_{*e(\uline\Psi \uline\Box \uline\Phi)} 
      &e(\uline *(\uline\alpha' \uline \Box\uline\beta')) 
      \ar@/ _3pc/[lll]_{e(\uline *(\uline\Psi\uline\Box\uline\Phi))}
      \ar[l]_{\xi_{\uline\alpha'\uline\Box\uline\beta'}^\inv} 
      \ar[ldd]_(.3){\tilde\iota^\inv_{\uline F,\uline*\uline\beta'}\circ \tilde\iota^\inv_{\uline K,\uline*\uline\alpha'}} 
      \ar@/ ^13pc/[lll]^{e((\uline *\uline\Phi)\uline\Box(\uline *\uline\Psi)) } 
      \\
      & {*e}(\uline\beta) {\widehat\Box} {*e}(\uline\alpha)\quad &  
      \quad {*e}(\uline\beta) {\widehat\Box} {*e}(\uline\alpha)
      \ar[l]_{{*e}(\uline\Phi){\widehat\Box}*e(\uline\Psi)}   
      \\
      &  e(\uline *\uline\beta){\widehat\Box} e(\uline *\uline\alpha)   
      \ar[u]_{\xi_{\uline\beta}^\inv{\widehat\Box}\xi_{\uline\alpha}^\inv}
      &  e(\uline *\uline\beta'){\widehat\Box} e(\uline *\uline\alpha')  
      \ar[l]_{e(\uline*\uline\Phi){\widehat\Box}e(\uline*\uline\Psi)}
      \ar[u]^{\xi_{\uline\beta'}^\inv{\widehat\Box}\xi_{\uline\alpha'}^\inv}
    }
  \end{align}
  for  3-morphisms $\uline\Phi: \uline\beta'\Rrightarrow\uline\beta$, $\uline\Psi:\uline\alpha'\Rrightarrow\uline\alpha$ and 2-morphisms $\uline\alpha,\uline\alpha': \uline F\Rightarrow\uline G$, $\uline\beta,\uline\beta': \uline H\Rightarrow\uline K$. In this diagram the expression ${\widehat\Box}$
  denotes the opposite Gray product of 2-morphisms
  $
  \beta{\widehat\Box}\alpha=( G\Box\beta)\circ (\alpha\Box H)$ from  \eqref{eq:product2} and from  Corollary \ref{lemma:cubical-to-opcubical-tricat}.  
  The  two rectangles in the middle  of the diagram  and the curved quadrilaterals at the top and bottom of the diagram  commute by definition of the 3-morphisms $\uline *(\uline\Psi\uline\Box\uline\Phi)$, $\uline\Psi\uline\Box\uline\Phi$, $\uline *\uline\Phi\uline \Box\uline *\uline \Psi$.
  To show that  the two  quadrilaterals at the left and right of the diagram commute, we note that it is sufficient to prove this for the case where either $\uline\alpha=\emptyset_{\uline F}$ or $\uline\beta=\emptyset_{\uline K}$.
  In the latter, the diagram reduces to
  \begin{align}\label{redbox}
    \xymatrix{ 
      e(\uline *(\uline\alpha\uline\Box\uline K))
      \ar[d]_{\tilde\iota^\inv_{\uline*\uline\alpha,\uline K}}
      &{*}e(\uline\alpha\uline\Box\uline K)
      \ar[l]_{\xi_{\uline\alpha\uline\Box\uline K}}
      \ar[d]_{*\tilde\iota_{\uline\alpha,\uline K}}
      \\
      e(\uline*\uline\alpha)\Box e(\uline K)
      &{*}e(\uline\alpha)\Box e(\uline K)
      \ar[l]^{\xi_{\uline\alpha}\Box \uline K},
    }
  \end{align}
  which clearly commutes if $\uline\alpha=\emptyset_{\uline F}$ or $\uline\alpha=(\alpha,1)$. For $\uline\alpha=(\alpha,-1)$, we consider the following diagram whose boundary corresponds to the diagram \eqref{redbox}
  \begin{align}\label{boxstar}
    \xymatrix{ 
      \#*(\#K\Box \alpha) 
      \ar[r]^{\Delta^\inv_{*(\#K\Box\alpha))}}
      \ar[d]_(.65){\#*\Theta_{\#K\Box\alpha}}
      \ar@/ _6pc/[dd]^(.85){\kappa^\inv_{1_{\#K},*\alpha}}
      & {*}(\#\alpha\Box\#1_{\#K})
      \ar[ld]_(.6){*\Gamma_{\#(\#K\Box\alpha)}}
      \ar[r]^{*\kappa_{1_{\#K},\alpha}}
      & {*}(\#\alpha\Box\#1_{\#K})
      \ar[ld]_(.7){*\Gamma_{\#\alpha\Box \#1_{\#K}}}
      \ar@/ ^5pc/[dd]_(.2){*(\#\alpha\Box \Phi_{\#K})}
      \\
      \#*\#\#(\#K\Box \alpha)\quad
      \ar[r]_{\#*\#\kappa_{ \#1_{\#K},\alpha}}
      &\quad \#*\#(\#\alpha\Box\#1_{\#K})\quad
      \ar[r]_{\quad \#*\#(\#\alpha\Box\Phi_{\#K})}
      & \quad \#*\#(\#\alpha\Box K)
      \\
      \#*\alpha\Box \#1_{\#K}
      \ar[rd]_{\#*\alpha\Box\Phi^\inv_{\#K}}
      \ar[r]^{\#*\Theta_\alpha\Box \Phi_{\#K}^\inv}
      & \#*\#\#\alpha\Box K
      & {*}\#\alpha\Box K
      \ar[l]_{*\Gamma_{\#\alpha}\Box K}
      \ar[ld]^{\Delta^\inv_{*\alpha}\Box K}
      \ar[u]^{*\Gamma_{\#\alpha\Box K}}
      \\
      & \#*\alpha\Box K
    }
  \end{align}
  The   triangle and the bottom quadrilateral in this diagram commute due to identity  \eqref{Deltagamma}. 
  The parallelogram at the top of the diagram commutes due to the naturality of the 3-morphisms
  $\kappa_{\mu,\nu}: \#\nu\Box\#\mu\Rrightarrow\#(\mu\Box\nu)$ and the quadrilateral on the right of the diagram due to the naturality of the 3-morphism $\Gamma_\mu: *\#{*}\#\mu\Rrightarrow\mu$.
  The heptagon in this diagram can be subdivided  as 
  \begin{align*}
    \xymatrix{
      & \;\;\#*\#\#(\#K\Box \alpha)\;\;
      \ar[rd]^{\#*\#\kappa_{ \#1_{\#K},\alpha}}
      \\
      \#*(\#K\Box \alpha)
      \ar[d]^{\kappa^\inv_{ 1_{\#K},\*\alpha}}
      \ar[ru]^{{\#*\Theta_{\#K\Box\alpha}}}
      \ar[r]^{\#*(\Theta_{1_{\#K}}\Box \Theta_\alpha)}
      & \#*(\#\# 1_{\#K}\Box \#\#\alpha)
      \ar[d]^{\#*(\Phi_{\#K}\Box\#\#\alpha)}
      &\#*\#(\#\alpha\Box \#1_{\#K})
      \ar[l]_{\#*\kappa_{ \#\alpha,\#1_{\#K}}}
      \ar[d]_{\#*\#(\#\alpha\Box\Phi_{\#K})}
      \\
      \#*\alpha\Box \#1_{\#K}
      \ar[ddd]^{\#*\Theta_\alpha\Box \Phi^\inv_{\#K}}
      & \#*(\#1_K\Box \#\#\alpha)
      \ar[d]^{\#\sigma^\inv_{*\#1_K, *\#\#\alpha}}
      & \#*\#(\#\alpha\Box K)
      \ar[l]_{\#*\kappa_{\#\alpha,1_{K}}}
      \\
      & \#(*\#1_K\Box *\#\#\alpha)
      \ar[d]^{\kappa^\inv_{ *\# 1_K,*\#\#\alpha}}
      &
      \\
      & \#*\#\#\alpha\Box \#*\#1_K
      &
      \\
      \#*\#\#\alpha\Box K\;\;
      \ar[ru]^(.4){\#*\#\#\alpha\Box *\Gamma_{1_K}\;\;}
      &\;\;{*(*\#{*}\#\#\alpha\Box *\#{*}\#1_K)}
      \ar[u]_{\sigma^\inv_{\#*\#\#\alpha, \#*\#1_K}}
      &{*\#\alpha\Box K}
      \ar[l]_{\qquad\qquad *(\Gamma_{\#\alpha}\Box \Gamma_{1_K})}
      \ar@/ ^3pc/[ll]^{*\Gamma_{\#\alpha}\Box K}
      \ar[uuu]^{*\Gamma_{\#\alpha\Box K}}.
    }
  \end{align*}
  The upper quadrilateral in this diagram  commutes due to relation \eqref{twistbox} and the rectangle below it by naturality of the 3-morphism $\kappa_{\mu,\nu}: \#\nu\Box\#\mu\Rrightarrow\#(\mu\Box\nu)$. The quadrilateral at the bottom of  the diagram commutes due to the naturality of the tensorator and the hexagon on the right  due to 
  identity \eqref{gamgraycomp}.  The heptagon on the left  commutes by naturality of the tensorator and of the 3-morphisms $\kappa_{\mu,\nu}: \#\nu\Box\#\mu\Rrightarrow\#(\mu\Box\nu)$, $\Phi_F: 1_{\#F}\Rrightarrow\#1_F$, $\Gamma_\mu: *\#{*}\#\mu\Rrightarrow\mu$ and 
  $\Theta_\mu:\#\#\Rrightarrow\mu$ together with identity \eqref{gamuncomp}.
  Hence, the diagram commutes, which implies that diagram \eqref{boxstar} commutes. This in turn proves the commutativity of diagram \eqref{redbox} for 2-morphisms $\uline\alpha=(\alpha, -1)$.

  For general 2-morphisms $\uline\alpha=(\uline\alpha_n,...,\uline\alpha_1)$ the commutativity of diagram \eqref{redbox}  follows directly from the identities \eqref{iotacomp} and 
  $\xi_{\uline\alpha\uline\circ\uline\beta}=\xi_{\uline\beta}\circ\xi_{\uline\alpha}$. This  proves that
  $\uline *$ defines a strict functor of 2-strict tricategories $\uline *:\uline{\mac G}\to\uline{\mac G}^{op}$ with trivial  coherence data and  $\uline*\uline*=1$.

  \bigskip
  4. It remains to prove the identity $\uline*\uline\#\uline *\uline \#=1$.  It is obvious that this identity holds for 1- and 2-morphisms.  To prove that it holds  for  3-morphisms $\uline\Psi:\uline\alpha\Rrightarrow\uline\beta$,
  we consider
  the  diagram
  \begin{align*}
    \xymatrix{ 
      e(\uline{*}\,\uline\#\,\uline {*}\,\uline\#\, \uline\alpha) 
      \ar@/ ^3pc/[rrr]^{e(\uline *\uline\#\uline *\uline \#\uline\Psi)}
      \ar@/ _17pc/[rrr]_{e(\uline\psi)}
      \ar[r]^{\xi^\inv_{\uline\#\,\uline {*}\,\uline\#\, \uline\alpha}}
      & {*}e(\uline\#\,\uline {*}\,\uline\#\, \uline\alpha) 
      \ar[d]^{*\chi_{\uline {*}\,\uline\#\, \uline\alpha}}  
      \ar[r]^{*e(\uline\#\,\uline {*}\,\uline\#\, \uline\Psi)} 
      & {*}e(\uline\#\,\uline {*}\,\uline\#\, \uline\beta)  
      & e(\uline{*}\,\uline\#\,\uline {*}\,\uline\#\, \uline\beta)
      \ar[l]_{\xi_{\uline\#\,\uline {*}\,\uline\#\, \uline\beta}}  
      \\
      & {*}\#e(\uline {*}\,\uline\#\, \uline\alpha)  
      \ar[d]^{*\# \xi_{\uline\#\, \uline\alpha}}  
      \ar[r]^{*\#e(\uline {*}\,\uline\#\, \uline\Psi)} 
      & {*}\#e(\uline {*}\,\uline\#\, \uline\beta) 
      \ar[u]^{\chi^\inv_{\uline {*}\,\uline\#\, \uline\beta}}
      \\
      & {*}\#{*}e(\uline\#\, \uline\alpha) 
      \ar[d]^{*\#* \chi^\inv_{\uline\alpha}}  
      \ar[r]^{*\#*e(\uline\#\, \uline\Psi)} 
      & {*}\#{*}e(\uline\#\, \uline\beta) 
      \ar[u]^{*\#\xi^\inv_{\uline\#\, \uline\beta}}
      \\
      & {*}\#{*}\#e(\uline\alpha) 
      \ar[luuu]_{\Gamma_{e(\uline\alpha)}} 
      \ar[r]^{*\#{*}\# e(\Psi)} 
      & {*}\#{*}\#e(\uline\beta) 
      \ar[u]^{*\#*\chi_{ \uline\beta}} 
      \ar[ruuu]^{\Gamma_{e(\uline\beta)}}
    }
  \end{align*}
  The  three rectangles  and the curved quadrilateral at the top of this diagram commute by definition of $\uline\#\,\uline\Psi$, $\uline*\uline\Psi$,  and the curved quadrilateral at the bottom commutes due to the naturality of $\Gamma$. It is therefore sufficient to show that  the curved subdiagrams at the left and the right  commute, which amounts to the relation
  $$
  \Gamma_{e(\uline\alpha)}\cdot *\#*\chi_{\uline\alpha}^\inv\cdot *\#\xi_{\uline\#\uline\alpha}\cdot *\chi_{\uline*\uline\#\uline\alpha}\cdot \xi^\inv_{\uline\#\uline *\uline\# \uline\alpha}=1_{e(\uline\alpha)}
  $$
  for all 2-morphisms $\uline\alpha$ in $\uline{\mac G}$. For $\uline\alpha=\emptyset_{\uline F}$, the 3-morphisms $\xi_{\uline\#\uline\alpha}$, $\xi_{\uline\#\uline*\uline\#\uline\alpha}$ are trivial, and this relation reduces to $$\Gamma_{1_{e(\uline F)}}\cdot *\# *\Phi_{e(\uline F)}\cdot *\Phi^\inv_{\#e(\uline F)}=1_{1_{e(\uline F)}},$$ 
  which holds by \eqref{gamuncomp}. For $\uline\alpha=(\alpha, 1)$, $\xi_{\uline\#\uline*\uline\#\uline\alpha}$ and $\chi_{\uline\alpha}$ are trivial, and from equation \eqref{Deltagamma} one obtains
  $$
  \Gamma_{\alpha}\cdot *\#\Delta_{*\alpha}^\inv\cdot *\Theta_{\alpha*}=1.
  $$
  For $\uline\alpha=(\alpha,-1)$, equation \eqref{Deltagamma} together with the identity $\#\Theta_{\alpha}=\Theta_{\#\alpha}$ in Lemma \ref{coherence} and the naturality of $\Delta_{\alpha}$
  implies
  $$
  \Gamma_{\#\alpha}\cdot *\#*\Theta_{\alpha}^\inv\cdot \Delta_{\alpha}=1_{\#\alpha}.
  $$

  To prove the identity for general 2-morphisms $\uline\alpha=(\uline\alpha_n,...,\uline\alpha_1)$, it is sufficient to show that the following diagram commutes
  $$
  \xymatrix{
    {*}\#e(\uline*\uline\#\uline\alpha_n)\circ\ldots \circ *\#e(\uline*\uline\#\uline\alpha_1)\quad 
    \ar[d]_{*\#\xi_{\uline\#\uline\alpha_n}\circ\ldots\circ  *\#\xi_{\uline\#\uline\alpha_1}}
    & \quad {*}\#(e(\uline*\uline\# \uline\alpha_n)\circ\ldots\circ e(\uline*\uline\#\uline\alpha_1))
    \ar[l]_{{*}\Phi_{e(\uline*\uline\#\uline\alpha_n),...,e(\uline*\uline\#\uline\alpha_1)}}
    \ar[d]_{*\#(\xi_{\uline*\uline\#\uline\alpha_n}\circ\ldots\circ \xi_{\uline*\uline\#\uline\alpha_1})}
    \\
    {*\#*}e(\uline\#\uline\alpha_n)\circ\ldots\circ *\#*e(\uline\#\uline\alpha_1) \quad 
    \ar[d]_{*\#*\chi^\inv_{\uline\alpha_n}\circ\ldots\circ *\#*\chi^\inv_{\uline\alpha_1}}
    & \quad  {*}\#*(e(\uline\# \uline\alpha_1)\circ\ldots\circ e(\uline\#\uline\alpha_n))
    \ar[l]_{*\Phi_{*e(\uline\#\uline\alpha_n),...,*e(\uline\#\uline\alpha_1)}}
    \ar[d]_{*\#*(\chi^\inv_{\uline\alpha_1}\circ\ldots\circ \chi^\inv_{\uline\alpha_n})}
    \\
    {*\#{*}\#}e(\uline\alpha_n)\circ\ldots\circ {*\#{*}\#}e(\uline\alpha_1) \quad 
    \ar[d]_{\Gamma_{e(\uline\alpha_n)}\circ\ldots\circ \Gamma_{e(\uline\alpha_1)}}
    &\quad  {*}\#*(\#e(\uline\alpha_1)\circ\ldots\circ \#e(\uline\alpha_n))
    \ar[l]_{*\Phi_{*\#e(\uline\alpha_n),...,*\#e(\uline\alpha_1)}}
    \ar[d]_{*\#*\Phi_{\uline\alpha_1,...,\uline\alpha_n}}
    \\
    e(\uline\alpha_n)\circ\ldots\circ e(\uline\alpha_1)\quad 
    & \quad  {*}\#*\#(e(\uline\alpha_n)\circ\ldots\circ e(\uline\alpha_1))
    \ar[l]_{\Gamma_{e(\uline\alpha)}}.
  }
  $$
  The upper two rectangles commute due to the naturality of the 3-morphism $\Phi_{\mu,\nu}:\#\mu\circ\#\nu\Rrightarrow\#(\nu\circ\mu)$. The rectangle at the bottom commutes due to identity 
  \eqref{gamhorcompeq}. This shows that the diagram commutes and the functors of 2-strict tricategories $\uline\#: \uline{\mac G}\to\uline{\mac G}_{op}$, $\uline *:\uline{\mac G}\to\uline{\mac G}^{op}$ satisfy $\uline *\uline\#\uline *\uline \#=1$.

  \bigskip
  5. To show that the Gray category $\uline{\mac G}$ is equivalent to ${\mac G}$, we note that the evaluation defines a \Glax functor of 2-strict tricategories $e:\uline{\mac G}\to\mac G$. As the evaluation is strictly compatible with the horizontal and the vertical composition and with all unit morphisms, the only coherence data of this functor is given by the 3-morphisms
  $\tilde\iota_{\uline\alpha,\uline K}: e(\uline \alpha)\Box e(\uline K)\Rrightarrow e(\uline\alpha\uline\Box \uline K)$, $\tilde\iota_{\uline F,\uline\beta}: e(\uline F)\Box e(\uline\beta)\Rrightarrow e(\uline F\uline \Box\uline\beta)$. As in the proof of Theorem \ref{graycatduals}, it is therefore sufficient to show that for all composable 3-morphisms $\uline\alpha:\uline F\Rightarrow\uline G$, $\uline\beta:\uline H\Rightarrow\uline K$
  the 3-morphisms $\tilde\iota_{\uline\alpha,\uline\beta}=\tilde\iota_{\uline\alpha, \uline K}\circ \tilde\iota_{\uline F,\uline\beta}$ are natural in both arguments and satisfy conditions analogous to 
  \eqref{kappanatisom}, \eqref{kappaunit} and \eqref{kappa_assoc}
  as well as
  $
  \tilde\iota_{{1_\mathcal C}, \uline \alpha}=\tilde\iota_{\uline\alpha, {1_\mathcal D}}=1_{e(\uline\alpha)}
  .$
  The naturality and the compatibility with the unit-morphisms are a direct consequence of the definitions. Condition \eqref{kappanatisom} follows from the commutative diagram  \eqref{graycirc}
  and Condition \eqref{kappa_assoc} from the commutative diagrams \eqref{assoc_strict_1} and \eqref{assoc_strict_2}.
  This shows that the evaluation defines a \Glax functor of 2-strict tricategories  $e:\uline{\mac G}\to\mac G$.

  We construct an embedding functor $f: \mac G\to\uline{\mac G}$ that will be a \Glax functor of 2-strict tricategories. For this, we set  $f(\mac A)=\mac A$ for all objects, $f(F)=(F,1)$ for all 1-morphisms , $f(\alpha)=(\alpha,1)$ for all 2-morphisms of $\mac G$ and $f(\Gamma)=\Gamma$ for all 3-morphisms of $\mac G$. This defines for all objects $\mac A,\mac B$ a weak 2-functor $f_{\mac A,\mac B}: \mac G(\mac A,\mac B)\to\uline{\mac G}(\mac A,\mac B)$ with the coherence data given by the invertible 3-morphisms  $1_{\mu\circ \nu}: f(\mu)\circ f(\nu)\Rrightarrow f(\mu\circ \nu)$ and  $1_{1_F}: \emptyset_{f( F)}\Rrightarrow f(1_F)$ for all 1-morphisms $F$ and composable 2-morphisms $\mu,\nu$ in $\mathcal G$. 
  
  The 2-morphisms $\iota_{\mathcal{C}}: f(1_{\mathcal{C}}) \Rightarrow 1_{f(\mathcal{C})}$  from Definition \ref{grayfunc} and their inverses are given by 
  $\iota_{\mathcal{C}}= \iota_{\mathcal{C}}^{-1}= (1_{1_{\mathcal{C}}},1)$. 
  The invertible pseudo-natural transformation $\kappa_{\mac A,\mac B,\mac C}: \uline\Box ( f_{\mac B,\mac C}\times f_{\mac A,\mac B})\to f_{\mac A,\mac C}\Box$ is determined by the 2-morphisms $(1_{F\Box G}, 1):f(G)\uline \Box f(F)\Rightarrow f(G\Box F)$ and the invertible 3-morphisms 
  $1_{\mu\Box\nu}: (1_{G\Box K},1)\uline \circ (f(\mu)\uline \Box f(\nu))  \Rrightarrow f(\mu\Box\nu)\uline \circ (1_{F\Box H},1)$ for all pairs of composable 1-morphisms $G,F$ and $H,K$ and 2-morphisms $\mu: F\Rightarrow G$, $\nu:H\Rightarrow K$ in $\mac G$. It is easy to show that the coherence conditions in Definitions \ref{grayfunc} and  \ref{laks2func} are satisfied and, consequently,  
  $f$ defines a \Glax functor of 2-strict tricategories $f: \mac G\to\uline{\mac G}$.

  It follows directly that  $ef=1_{\mac G}$. The \Glax functor of 2-strict tricategories $f e:\uline{\mac G}\to\uline{\mac G}$ is given by $f e(\mac A)=\mac A$, $f e (\uline F)=(e(\uline F), 1)$, $f e(\uline\alpha)=(e(\uline \alpha), 1)$ and  $f e(\uline\Gamma)=\uline\Gamma$.
  A natural isomorphism of \Glax functors of 2-strict tricategories $\eta: fe\to 1$ is given by the trivial 1-morphism $\emptyset_{\mac A}:\mac A\to\mac A$ for each object $\mac A$ of $\uline{\mac G}$ together with the  
  invertible pseudo-natural transformation of weak 2-category functors $1_{\mac A,\mac B}\to (fe)_{\mac A,\mac B}$ that is determined by the 2-morphisms  $\eta_{\uline F}=(1_{e(\uline F)}, 1): fe(\uline F)\Rightarrow  \uline F$ for each 1-morphism $\uline F:\mac A\to\mac B$ and the invertible 3-morphism  $\eta_{\uline \alpha}=1_{e(\uline\alpha)}:  (1_{e(\uline G)},1)  \uline \circ fe(\uline\alpha)\Rrightarrow \uline \alpha  \uline \circ (1_{e(\uline F)},1)$ for each 2-morphism $\uline\alpha:\uline F\Rightarrow\uline G$.  A direct calculation shows that the consistency conditions in Definitions \ref{graynat} and \ref{nat2transformations}  are satisfied and that this defines a natural isomorphism $fe\to 1$ of  \Glax functors of 2-strict tricategories. It also follows directly that $e\eta=1_e: e\to e$ and $\eta fe=1_{fe}: fe\to fe$. 

  The
  invertible pseudo-natural transformation $\eta f: f\to f$ is determined by the 2-morphisms $(1_F,1): f(F)\Rightarrow f(F)$ for each 1-morphism $F$ and the 3-morphisms $1_{\alpha}: ((\alpha,1), (1_{F},1))\Rrightarrow ((1_G,1), (\alpha,1))$. A modification $\Psi:\eta f\Rightarrow 1_f$ is therefore given by the trivial 2-morphism $\emptyset_{1_{\mac A}}$ for each object $\mac A$ of $\mac G$ and the invertible 3-morphisms $1_{1_F}: (\eta f)_F\Rrightarrow \emptyset_{f(F)}$. This   implies $e\Psi=1_{1_{\uline{\mac G}}}:
  e\eta f=1_ef=1_{ef}=1_{1_{\mac G}}\Rrightarrow e1_f=1_{ef}=1_{1_{\mac G}}$
  and 
  concludes the proof that the Gray categories $\mac G$ and $\uline{\mac G}$ are equivalent.

  \bigskip
  6. By definition, the \Glax functor of 2-strict tricategories $f: \mac G\to\uline{\mac G}$ satisfies $\uline*f=f*$. As the functors are the identity on the objects, a
  natural isomorphism $\tilde\chi: \uline\# f\to f\#$ is determined by an invertible pseudo-natural transformation of weak 2-functors $\tilde \chi: (\uline\# f)_{\mac A,\mac B}\to (f\#)_{\mac A,\mac B}$  for each pair of objects $\mac A,\mac B$. This natural isomorphism is given by the 2-morphisms $\tilde\chi_{F}=(1_{\#F},1): f\#(F)\Rightarrow \uline\# f(F)$ for each 1-morphism $F$ in $\mac G$ and the invertible 3-morphisms 
  $\tilde\chi_{\mu}=1_{\#\mu}:  (1_{\#G}, 1)\uline\circ_{op} \uline\#f(\mu)   \Rrightarrow  f\#(\mu)\uline\circ_{op} (1_{\#F}, 1)$ for each 2-morphism $\mu:F\Rightarrow G$.
  It follows directly that all coherence conditions in Definitions  \ref{graynat} and \ref{nat2transformations}
  are satisfied. 

  The natural isomorphisms $\chi: \#e\to e\uline\#$ and $\xi: *e\to e\uline *$ are 
  obtained from the coherence data of $\mac G$. As the functors of 2-strict tricategories $\#e:\uline{\mac G}\to \mac G_{op}$ and $e\uline\#:\uline{\mac G}\to \mac G_{op}$, as well as  $*e:\uline{\mac G}\to \mac G^{op}$ and $e\uline*:\uline{\mac G}\to \mac G^{op}$ agree on the objects and 1-morphisms of $\uline{\mac G}$, such natural isomorphisms are specified uniquely by 
  natural isomorphisms between the functors $(\#e)_{\uline F,\uline G}, (e\uline\#)_{\uline F,\uline G}: \uline{\mac G}(\uline F,\uline G)\to \mac G_{op}(\#e(\uline F), \#e(\uline G))$ and between the functors $(*e)_{\uline F,\uline G}, (e\uline *)_{\uline F,\uline G}: \uline{\mac G}(\uline F,\uline G)\to \mac G^{op}(*e(\uline F), *e(\uline G))$. They are determined by the   invertible 3-morphisms
  $\chi_{\uline\alpha}: \#e(\uline\alpha)\Rrightarrow e(\uline\#\,\uline\alpha)$ and $\xi_{\uline\alpha}: *e(\uline\alpha)\Rrightarrow e(\uline*\,\uline\alpha)$ for each 2-morphism $\uline \alpha:\uline F\Rightarrow\uline G$. That they satisfy the consistency conditions  in Definitions  \ref{graynat} and \ref{nat2transformations} was shown, respectively, in the second and third part of the proof. 
\end{proof}

Theorem \ref{thm:strictification} explicitly constructs a  Gray category $\uline{\mac G}$ and strictifications $\uline*:\uline{\mac G}\to\uline{\mac G}^{op}$,  $\uline \#: \uline {\mac G}\to\uline{\mac G}_{op}$ of the functors of 2-strict tricategories $*: \mac G\to\mac G^{op}$, $\#:\mac G\to\mac G_{op}$. 
This construction has the benefit that it is conceptually clear and  concrete and allows one to  verify the properties of the strictified functors by direct computations. 
It  remains to show that the Gray category $\uline{\mac G}$  with the strict functors of 2-strict tricategories  $\uline*:\uline{\mac G}\to\uline{\mac G}^{op}$,  $\uline \#: \uline {\mac G}\to\uline{\mac G}_{op}$ is again a Gray category with duals in the sense of Definition \ref{dualgray} and to clarify  which additional relations hold in the strictified Gray category.

\begin{theorem}\label{lem:graystrict} For every spatial Gray category  $\mac G$,  the associated Gray category $\uline{\mac G}$ from Theorem \ref{thm:strictification} is a Gray category with strict duals in the sense of Definition \ref{graystrictdual}.
\end{theorem}
\begin{proof} 
  $\quad$ \newline
  1.  For each pair of objects $\mac C$, $\mac D$ of $\uline{\mac G}$, the functor $\uline *:\uline{\mac G}\to\uline{\mac G}^{op}$ defines a strict 2-functor $\uline *:\uline{\mac G}(\mac C,\mac D)\to \uline{\mac G}(\mac C,\mac D)^{op}$ that is trivial on the objects of $\uline{\mac G}(\mac C,\mac D)$ and satisfies  $\uline *\uline *=1$. To show that this gives $\uline{\mac G}(\mac C,\mac D)$ the structure of a planar 2-category, it is sufficient to construct 
  for each
  2-morphism $\uline\mu$ a 3-morphism
  $\uline\epsilon_{\uline\mu}: \emptyset_{\uline G}\Rrightarrow \uline\mu \uline\circ \uline * \uline\mu$
  that satisfies the conditions  in Definition \ref{dualgray} and in \eqref{pivcond}. This 3-morphism is defined by
  $$
  e(\uline\epsilon_{\uline\mu})=(1_{e(\uline\mu)}\circ \xi_{\uline\mu})\cdot \epsilon_{e(\uline\mu)}, $$
  where $\epsilon_{e(\uline\mu)}$ denotes the corresponding 3-morphism in $\mac G$. 
  The identity
  $\uline H\uline\Box \uline\epsilon_{\uline\mu}\uline\Box\uline K=\uline\epsilon_{\uline H\uline\Box\uline \mu\uline\Box\uline K}$ from Definition \ref{dualgray} follows from the commutative diagram
  \eqref{strpr1}   and the analogous diagram with the 1-morphism on the left.

  \begin{align}\label{strpr1}
    \xymatrix{ 1_{e(\uline F\uline\Box\uline K)}
      \ar@/ _12pc/[rr]^{e(\uline\epsilon_{\uline\mu}\uline\Box\uline K)}
      \ar@/ _1.9pc/[rdd]^(.7){e(\uline\epsilon_{\uline\mu})\Box e(\uline K)}
      \ar[r]^{\epsilon_{e(\uline\mu\uline\Box\uline K)}\qquad}
      \ar@/ ^3pc/[rr]^{e(\uline\epsilon_{\uline\mu\uline\Box\uline K})}
      \ar[rd]^(.6){\epsilon_{e(\uline\mu)}\Box e(\uline K)}
      & e(\uline\mu\uline\Box\uline K)\circ * e(\uline\mu\uline\Box\uline K)
      \ar[r]^{e(\uline\mu\uline\Box\uline K)\circ \xi_{\uline\mu\uline\Box\uline K}}
      &e( (\uline \mu\uline\circ \uline *\uline\mu)\uline\Box\uline K)
      \\
      & (e(\uline\mu) \circ *e(\uline\mu))\Box e(\uline K)
      \ar[d]^{(e(\uline\mu)\circ \xi_{\uline\mu})\Box e(\uline K)}
      \ar[u]_{\tilde\iota_{\uline\mu. \uline K}\circ *\tilde\iota^\inv_{\uline\mu,\uline K}}
      & 
      \\
      &e(\uline\mu\uline\circ \uline *\uline\mu)\Box e(\uline K)
      \ar@/ _2pc/[ruu]^(.7){\tilde\iota_{\uline\mu\uline\circ\uline *\uline\mu, \uline K}}
    }
  \end{align}

  That the 3-morphism $\uline\epsilon_{\uline\mu}: \emptyset_{\uline G}\Rrightarrow \uline\mu\,\uline\circ\,\uline*\,\uline \mu$ satisfies  the conditions \eqref{pivcond} is a consequence 
  of  the following three  commutative diagrams
  and the analogue of the second diagram for the composite $(\uline *\epsilon_{\uline\mu}\uline\circ \uline\mu)\uline \cdot(\uline\mu\uline\circ\uline \epsilon_{\uline *\uline \mu})$
  \begin{align*}
    &\xymatrix{
      1_G 
      \ar@/ _4pc/[rrdd]_{e(\epsilon_{\uline\mu\uline\circ\uline\nu})}
      \ar[r]^{e(\uline\epsilon_{\uline\mu})\qquad} 
      \ar[rd]_(.6){\epsilon_{e(\uline\mu)}}
      & e(\uline\mu\uline\circ\uline*\uline \mu) 
      \ar[r]^{{\uline\mu}\uline\circ \uline\epsilon_{\uline\nu}\uline\circ {\uline*\uline \mu}\qquad}
      \ar[rd]_(.4){e(\uline\mu)\circ \epsilon_{e(\uline\nu)}\circ e(\uline *\uline\mu)}
      & e(\uline\mu\uline\circ\uline\nu\uline\circ\uline*\uline\nu\uline\circ \uline *\uline\mu) 
      \\
      & e(\uline\mu)\circ *e(\uline\mu) 
      \ar[u]^(.6){e(\uline\mu)\circ \xi_{\uline\mu}}
      \ar[rd]_(.4){e(\uline\mu)\circ \epsilon_{e(\uline\nu)}\circ *e(\uline\mu)}
      & e(\uline\mu \uline\circ\uline\nu)\circ* e(\uline\nu)\circ e(\uline *\uline\mu)
      \ar[u]_{e(\uline\mu\uline\circ\uline\nu)\circ \xi_{\uline\nu}\circ e(\uline *\uline\mu)}
      \\
      & & e(\uline\mu\uline\circ\uline\nu)\circ *e(\uline\nu)\circ *e(\uline\mu)
      \ar[u]_{e(\uline\mu\uline\circ\uline\nu)\circ *e(\uline\nu)\circ \xi_{\uline\mu}}
    }
  \end{align*}

  \begin{align*}
     &\xymatrix{
      e(\uline\mu)
      \ar[d]_{\epsilon_{e(\uline\mu)}\circ e(\uline\mu)}
      \ar[r]^{e(\uline \epsilon_{\uline\mu}\uline \circ {\uline\mu})}
      & e(\uline\mu\uline\circ\uline*\uline\mu\uline\circ\uline\mu)
      \ar[rd]^{e(\uline\mu\uline\circ \uline *\uline\epsilon_{\uline *\uline\mu})}
      \\
      e(\uline\mu)\circ *e(\uline\mu)\circ e(\uline\mu)
      \ar@/ _7pc/[rr]_{e(\uline\mu)\circ \epsilon_{*e(\uline\mu)}}
      \ar[ru]^{e(\uline\mu)\circ \xi_{\uline\mu}\circ e(\uline\mu)}
      & e(\uline\mu)\circ *e(\uline\mu)\circ *e(\uline *\uline\mu)
      \ar[u]^{e(\uline\mu)\circ \xi_{\uline\mu\uline\circ \uline *\uline\mu}}
      \ar[r]^{\qquad e(\uline\mu)\circ * e(\epsilon_{\uline*\uline\mu})}
      \ar[d]_(.4){e(\uline\mu)\circ \xi_{\uline\mu}\circ *e(\uline *\uline\mu)}
      &\quad e(\uline\mu)
      \\
      & e(\uline\mu)\circ e(\uline *\mu)\circ*e(\uline *\mu)
      \ar[ru]^{e(\uline\mu)\circ *\epsilon_{e(\uline*\uline\mu)}}
      \ar[lu]^{e(\uline\mu)\circ e(\uline *\uline\mu)\circ \xi_{\uline\mu}}
    }
  \end{align*}
  \begin{align*}
     &\xymatrix{ 1_{e(\uline G)} 
      \ar[d]_{e(\uline\epsilon_{\uline\nu})}
      \ar[rd]_{\epsilon_{e(\uline\nu)}}
      \ar[rr]^{e(\uline\epsilon_{\uline\mu})}
      \ar[rrd]^{\epsilon_{e(\uline\mu)}}
      & & e(\uline\mu)\circ e(\uline *\uline \mu)
      \ar@/ ^5pc/[dd]^{e(\uline\Psi)\circ e(\uline *\uline\mu)}
      \\
      e(\uline\nu)\circ e(\uline *\uline \nu)
      \ar@/ _4pc/[rrd]_{e(\uline\nu)\circ e(\uline *\uline\Psi)}
      & e(\uline\nu)\circ *e(\uline\nu)
      \ar[l]^{e(\uline\nu)\circ \xi_{\uline\nu}}
      \ar[d]_{e(\uline\mu)\circ *e(\uline\Psi)}
      & e(\uline\mu)\circ *e(\uline\mu)
      \ar[u]_{e(\uline\mu)\circ \xi_{\uline\mu}}
      \ar[ld]^(.4){e(\uline\Psi)\circ *e(\uline\mu)}
      \\
      &
      e(\uline\nu)\circ *e(\uline\mu)
      \ar[r]_{e(\uline\nu)\circ \xi_{\uline\mu}}
      & e(\uline\nu)\circ e(\uline*\uline\mu).
    }
  \end{align*}
  This shows that for all objects $\mac C$, $\mac D$ the 2-category  $\uline{\mac G}(\mac C,\mac D)$
  is planar and that the first condition in Definition \ref{dualgray} is satisfied. 

  \bigskip
  2. The functor $\uline\#:\uline{\mac G}\to\uline{\mac G}_{op}$ defines the dual of each 1-morphism $\uline F:\mac C\to\mac D$  and by definition satisfies $\uline\#\,\uline\#\,\uline F=\uline F$, $\uline\#(\uline F\,\uline\Box\,\uline G)=\uline\#\uline G\,\uline\Box\,\uline\#\uline F$, $\uline\#\emptyset_{\mac C}=\emptyset_{\mac C}$. It remains to construct 
  the fold  2-morphisms $\uline\eta_{\uline F}: \emptyset_{\mac D}\Rightarrow \uline F\uline\Box \uline \# \uline F$ and the triangulator  3-morphisms $\uline T_{\uline F}: (\uline *\uline\eta_{\uline F}\uline \Box\uline F)\uline\circ (\uline F\uline \Box\uline\eta_{\uline F})\Rrightarrow \emptyset_{\uline F}$
  and to show that they satisfy the conditions in Definition \ref{dualgray}. We define
  $$
  \uline\eta_{\emptyset_{\mac C}}=\emptyset_{\emptyset_{\mac C}},\qquad \uline\eta_{\uline G}=(\eta_{e(\uline G)}, 1): \emptyset_{\mac D}\Rightarrow \uline F\uline\Box \uline \#\uline F,\qquad e(\uline T_{\uline F})=T_{e(\uline F)}
  $$
  for all 1-morphisms $\uline F:\mac C\to\mac D$ and all non-empty 1-morphisms $\uline G:\mac C\to\mac D$.
  Conditions (2) (b) in Definition \ref{dualgray} then hold by definition.
  As  the 3-morphisms $\tilde\iota_{\uline \eta_{\uline F}, \uline H}$, $\tilde\iota_{\uline K,\uline\eta_F}$
  and $\xi_{\uline\eta_{\uline F}}$ and the analogous 3-morphisms for $\uline *\uline\eta$  are trivial, the remaining  identities in  (2) (c), (d) then follow directly from the corresponding properties of the fold 2-morphisms and the triangulator in $\mac G$. 
  This shows that $(\uline{\mac G}, \uline *,\uline\#)$ is a Gray category with duals in the sense of Definition \ref{dualgray}.

  \bigskip
  3. To prove the  identity $\uline\#\uline\epsilon_{\uline\mu}=\uline\epsilon_{\uline *\uline\#\uline\mu}$, consider the diagram
  \begin{align}\label{strpr2}
    \xymatrix{
      1_{\#e(\uline G)} 
      \ar@/  ^4pc/[rrr]^{e(\uline \epsilon_{\uline *\uline\#\uline\mu})}
      \ar[r]^{\epsilon_{e(\uline *\uline\#\uline\mu)}\qquad}
      \ar[rd]^{\epsilon_{*e(\uline\#\uline\mu)}}
      \ar[rdd]_(.75){\qquad\qquad\epsilon_{*\# e(\uline\mu)}}
      \ar@/ _18pc/[rrr]^{e(\uline\#\uline\epsilon_{\uline\mu})}
      & e(\uline *\uline \#\uline\mu )\circ *e(\uline *\uline\#\uline\mu)
      \ar[rr]^{e(\uline *\uline \#\uline\mu)\circ \xi_{\uline *\uline\#\uline\mu}}
      & & e(\uline *\uline\#\uline\mu\uline\circ \uline\#\uline \mu)
      \\
      & {*}e(\uline\#\uline\mu)\circ e(\uline\#\uline\mu)
      \ar[u]_{\xi_{\uline\#\uline\mu}\circ *\xi^\inv_{\uline\#\uline\mu}}
      & \#*e(\uline\mu)\circ \#e(\uline\mu)
      \ar@/ _3pc/[dd]_(.7){\Phi_{*e(\uline\mu), e(\uline\mu)}}
      \ar[ld]^(.6){\Delta_{*e(\uline\mu)}\circ \#e(\uline\mu)}
      \\
      & {*}\# e(\uline\mu)\circ \#e(\uline\mu)
      \ar[u]_{*\chi_{\uline\mu}^\inv\circ \chi_{\uline\mu}}
      & \#(e(\uline\mu\uline\circ\uline *\uline\mu))
      \ar[ruu]_{\chi_{\uline\mu\uline\circ\uline *\uline\mu}}
      \\
      & \#1_{e(\uline G)}
      \ar@/ ^1pc/[luuu]_(.15){\chi_{\emptyset_{\uline G}}\qquad\qquad}
      \ar[r]_{\#\epsilon_{e(\uline\mu)}}
      &  \#(e(\uline\mu)\circ *e(\uline\mu)).
      \ar[u]_{\#(e(\uline\mu)\circ \xi_{\uline\mu})}
    }
  \end{align}
  The triangles at the top and the pentagon at the bottom of the diagram  commute by definition of $\uline\epsilon_{\uline *\uline\#\uline\mu}$ and $\uline\#\uline\epsilon_{\uline\mu}$. The two triangles on the left commute by \eqref{pivcond}. The upper polygon in the middle commutes since it can be decomposed into diagrams whose commutativity was established in the proof of Theorem \ref{thm:strictification}. The lower polygon in the middle  which involves the 3-morphisms $\#\epsilon_{e(\uline\mu)}$ and $\epsilon_{*\#e(\uline\mu)}$ commutes by definition of the 3-morphism $\Delta_{*e(\uline\mu)}$.  This shows that the diagram \eqref{strpr2} commutes and proves the identity 
  $\uline\#\uline\epsilon_{\uline\mu}=\uline\epsilon_{\uline *\uline\#\uline\mu}$.
\end{proof}

\medskip
The conditions in  Definition \ref{graystrictdual} on a Gray category with strict duals have a clear geometrical interpretation in terms of Gray category diagrams. In the diagrams, the functor of 2-strict tricategories $*$ corresponds to a $180$ degree rotation around the 
$w$-axis, and the condition $**=1$ ensures that the evaluation is invariant under a $360$ degree rotation of the diagram. Similarly, the functor $\#$ in a Gray category with strict duals corresponds to a $180$ degree rotation around the $y$-axis and the condition $\#\#=1$ ensures that the evaluation is invariant under  a $360$ degree rotation. The condition $*\#{*}\#=1$ corresponds to the fact that the $180$ degree rotations around the $w$- and $y$-axis commute. Together with the strictness of the functors $*$ and $\#$, these conditions on $*$ and $\#$  ensure that the functors of 2-strict tricategories $*$ and $\#$ correspond to the symmetries of a cube.
In contrast to the original Gray category $\mac G$, where these symmetries were realised up to higher morphisms, in the strictified Gray category $\uline{\mac G}$ these symmetries are realised exactly.

The  condition in Definition \ref{graystrictdual} that the 3-morphisms $\#\epsilon_\mu, \epsilon_{*\#\mu}: \#1_{G}\to \#\mu\circ *\#\mu$  agree ensures that the labelling of the minima and maxima of lines in the associated Gray category diagrams does not become ambiguous.  

Note also that the strictification theorem implies a coherence result for the 3-morphisms $\Phi_F$, $\Phi_{\mu,\nu}$ and $\kappa_{\mu,\nu}$ from Theorem \ref{graycatduals} and the 3-morphisms $\Theta_\mu$, $\Gamma_\mu$ that characterise the natural isomorphisms $\Theta:\#\#\to 1$ and $\Gamma: *\#{*}\#\to 1$ in Theorem \ref{graydualnat}.  As shown in the proof of Theorem  \ref{thm:strictification}, these are precisely the 3-morphisms associated with the evaluation functor 
$e: \uline{\mac G}\to\mac G$. The strictification theorem implies that any two 3-morphisms $\Psi,\Omega:\mu\Rrightarrow\nu$ which are constructed from these 3-morphisms, their inverses and their $*$- and $\#$-duals via the Gray product, the horizontal and vertical composition and the tensorator
are equal. 

With respect to the discussion in Section \ref{subsec:geom-interpret} this suggests that  in a spatial Gray category, it  should be possible to omit the labelling by 3-morphisms $\Omega_\mu$ at the points where lines in the diagrams cross the fold lines, as these labellings are canonical. Similarly,  the evaluation of two diagrams that can be transformed into each other by sliding lines over folds and cusps as in Figure \ref{foldcrossgamma2}, \ref{foldcrosstheta} should be related by a unique 3-morphism that is constructed from the 
3-morphisms $\Gamma_\mu$ and $\Theta_\mu$, their $*$- and $\#$-duals and their inverses via the  Gray product, the horizontal and  vertical composition and the tensorator. However, a detailed exploration of this idea is beyond the scope of this paper.

\section{Diagrams for Gray categories with duals}
\label{sec:diagrams_duals}

In this section, it is shown that, under suitable additional assumptions, the evaluation of a Gray category diagram
is invariant under homomorphisms of Gray category diagrams. The first assumption is that the Gray category is spatial in the sense of Definition 4.8. The second assumption is that the 2-skeleton of the diagram is a two-manifold.

\subsection{Non-progressive Gray category diagrams}
\label{non-progressive}

In this subsection more general diagrams than the progressive ones are defined.  This is familiar from two-dimensional diagrams, where the maxima and minima can be interpreted as points on a single line that changes direction. Exactly the same thing happens for lines in three-dimensional diagrams. 

The analogous situation for surfaces is that they can have singularities of the projection.  So the diagrams in Figure \ref{baez} can be interpreted as containing a single surface with folds and cusps. A singular point of the projection of a surface is called a fold if it is locally isomorphic to a subdiagram of any of the diagrams in Figure \ref{saddle}. A singular point is called a cusp if it is locally isomorphic to Figure \ref{baez} g), \ref{BL_consisteny} a), or their rotations by $\pi$ around the $w$-axis. Note that the two folds meeting a cusp are either both input lines (lines meeting the top face of a small cube around the vertex) or both output lines (lines meeting the bottom face of a small cube around the vertex).

A diagram with folds and cusps can be subdivided to give a progressive diagram by introducing additional lines along the folds and additional vertices at the cusps, at maxima and minima of folds and at the intersection of folds and lines. This remark is not entirely trivial because it is necessary to check that the additional lines at folds meet the boundary correctly according to Definition  \ref{3ddiag} \eqref{3d1}. This follows from condition \eqref{3d2} of the same definition, since a fold line meeting a side face would imply that the side face is not a progressive diagram.

According to Whitney's classification \cite{W}, the generic singularities of the projection of a smooth surface are the smooth folds and cusps. This motivates the analogous definition of generic singularity in the PL case.

\begin{definition}
A three-dimen\-sio\-nal diagram is called {\bf generic} if the only singularities of the projection $p_2$ on surfaces are folds and cusps, and the only singularities of $p_1\circ p_2$ on lines and folds are maxima and minima. In addition, the subdivision obtained by additional lines at the folds and additional vertices at the cusps and the intersection of folds and lines is required to be a generic progressive diagram according to Definition \ref{3ddiag}.  
\end{definition}

The definition and evaluation of a diagram for a Gray category with duals follows the same pattern as for two-dimensional diagrams in Section \ref{2catdualssection}. Again it will be assumed without further mention that diagrams are generic where this is appropriate.

 \begin{definition}\label{ev_diag}
   Let $\mathcal G$ be a Gray category with duals. 
   A {\bf diagram for $\mathcal G$}  is a three-dimensional diagram $D$ together with a labelling of its minimal progressive subdivision  $S$ with elements of $\mathcal G$ that makes $S$ into a Gray category diagram.
   This labelling is required to be such that the fold lines are labelled appropriately with either $\eta_F$ or $\eta^*_F$, the cusps with $T_F$, $T_F^{-1}$, $*T_F$, or $*T_F^{-1}$, 
   the minima and maxima for lines by $\epsilon_\mu$ or  $\epsilon_\mu^*$, and the minima and maxima for folds by 
   $\epsilon_{\eta_F}$,  $\epsilon_{\eta_F}^*$, $\epsilon^*_{\eta_F}$ or  $\epsilon^*_{\eta_F^*}$, as shown in Figure \ref{saddle}.  The {\bf evaluation of $D$} is defined as the evaluation of $S$.
 \end{definition}
 According to this definition, any new vertices introduced at the intersection of folds and lines are labelled with a 3-morphism but these are not constrained.  It follows from the discussion in Section \ref{subsec:geom-interpret}
 that, in one case of the intersection of a fold and a line,  a canonical label is given by the 3-morphism $\Omega_\mu\colon \eta_G^*\circ (\mu\Box G^\#)\Rrightarrow \eta_{F}^*\circ (F\Box\#\mu)$  defined in \eqref{omega_mu}.  As explained there, the interaction of $\Omega_\mu$ with the unit 2-morphisms, the horizontal composition and the Gray product is  determined  by the coherence data of the functor of 2-strict tricategories $\#$ (see Figure \ref{foldcrossconsist3}) and the evaluations of diagrams related by sliding lines over a cusp are related by the  natural isomorphisms $\Theta\colon\#\#\to 1$, $\Gamma\colon *\#{*}\#\to 1$ (see Figures \ref{foldcrossgamma2} and \ref{foldcrosstheta}). Together with the strictification result from Section \ref{sec:strictification}, which implies that this data is coherent, this suggests that the labels $\Omega_\mu$ at the folds are canonical and can be omitted, since any two diagrams constructed from such labellings should be related by a unique 3-morphism. 
 However, this aspect is not  analysed systematically in this paper, and in the following we restrict attention to diagrams where there are no such vertices. 

\begin{definition}\label{def:standard}
  A three-dimensional diagram is called {\bf standard} if there are no folds meeting vertices or lines.
\end{definition}

Standard diagrams are familiar from knot theory. A ribbon knot is given by an embedding  $e\colon S^1\times[-1,1]\to\mathbb R^3$.  A ribbon knot with no folds is called blackboard-framed in knot theory. 
Ribbon knots were generalised to ribbon graphs by Reshetikhin and Turaev \cite{RT}. A ribbon graph consists of a graph, called the core, and a 
compact surface with boundary which contains the core. The ribbon graph is considered to be a thickening of the core. 

The ribbon graphs can be realised as diagrams in the sense of this paper in the following way.
Given a diagram $D$ for a Gray category with duals, as in  Definition \ref{3ddiag},  we denote the $k$-skeleton of $D$ (not of its minimal progressive subdivision) by $X^{k}$.

\begin{definition}
  A {\bf ribbon diagram}  is a three-dimen\-sio\-nal  diagram $D$ with an embedded graph $\gamma\subset X^1$,  called the {\bf core}, such that 
  \begin{enumerate}
  \item The core $\gamma$ is the union of all of the vertices $X^0$ and a subset of the set of lines.
  \item There  exists a  two-dimensional PL manifold\footnote{Note that, to avoid confusion, the definition uses the terminology `two-dimensional PL manifold' because the word `surface' is reserved for the 2-dimensional strata of a diagram.} $\Sigma\subset [0,1]^3$ such that $X^2$ is a regular neighbourhood of $\gamma$ in $\Sigma$.
  \item $X^1\subset \gamma \cup\partial X^2$.
  \end{enumerate}
\end{definition}

\begin{figure}
  \centering
  \includegraphics[scale=0.4]{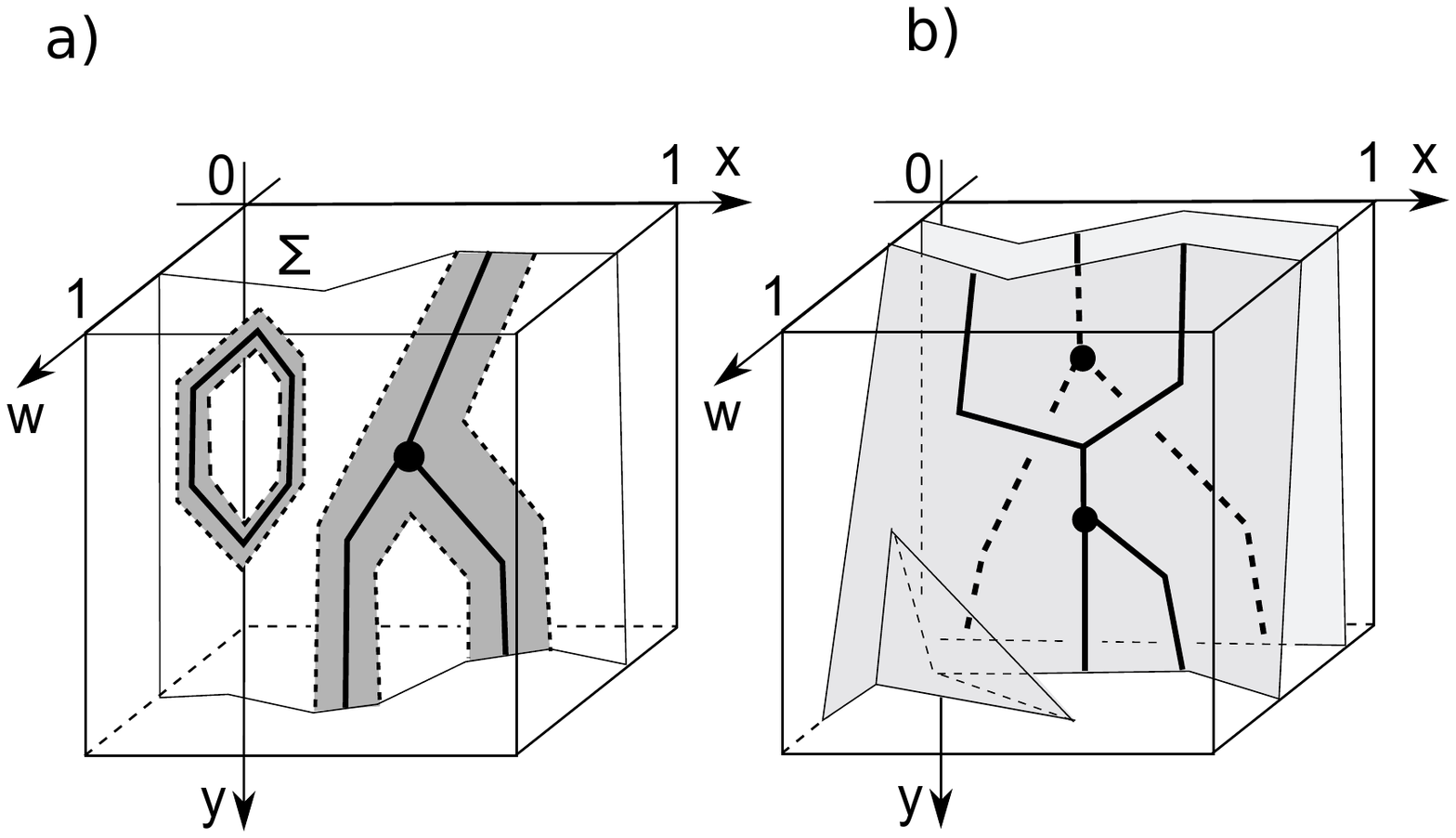}
  \caption{a) A ribbon diagram. The dashed lines denote lines in the subset $l\subset X^1$, the thick solid lines and vertices  the core $\gamma$ of the ribbon. The 2-skeleton $X^2$ is drawn in grey, the auxiliary surface $\Sigma$ in white. \newline
    b)  A  surface diagram.}
  \label{surface_diagram}
\end{figure}

A standard result on regular neighbourhoods is that $X^2$ is a compact two-manifold with boundary \cite{RS}. The lines of $D$ that are not in the core form a subset 
$l\subset  X^1$. 
The surfaces of $D$ are the components of $X^2\setminus\{l\cup\gamma\}$. Note that $\Sigma$ is not part of the structure of the ribbon diagram. It is just required that a suitable PL manifold $\Sigma$ exists.
An example of a ribbon diagram is given  in Figure \ref{surface_diagram} a). 

As shown in Corollary \ref{cor:spatial},  if $\mac G$ is spatial, then for all objects $\mac C$ of $\mac G$  the 2-category $\mac G(1_{\mac C}, 1_{\mac C})$ is a ribbon category. Conversely,  a  
ribbon category can be viewed as a spatial Gray category, with only one object $\mac C$ and one 1-morphism $1_{\mac C}$.  This is the appropriate category data for labelling a ribbon diagram with no folds.  

It is easy to see that the evaluation of a ribbon diagram labelled with such data
coincides with the Reshetikhin-Turaev evaluation of the associated ribbon. For this, one labels 
the regions of the diagram  with the object  $\mac C$, its  surfaces with the 1-morphism $1_{\mac C}$ and  assigns the trivial 2-morphism $1_{1_{\mac C}}$ to the lines in $l$. 
The vertices and lines in $\gamma$ are labelled with data from the ribbon category $\mac G(1_{\mac C}, 1_{\mac C})$. The evaluation of such a ribbon diagram with no folds 
according to Definition \ref{ev_diag} then  
coincides with the Reshetikhin-Turaev evaluation of the associated ribbon labelled with data from the ribbon category $\mac G(1_{\mac C}, 1_{\mac C})$.

The ribbon diagrams can be modified to provide another interesting class of examples of Gray category diagrams.

\begin{definition} \label{surf_diagram}A {\bf surface diagram} is a three-dimensional diagram such that $X^2$ is a two-dimensional PL manifold whose boundary is contained in the boundary of the cube, $\partial X^2\subset\partial [0,1]^3$. 
\end{definition}

An example of a surface diagram is given in Figure \ref{surface_diagram} b).
Surface diagrams and ribbon diagrams are closely related.
For every surface diagram $D$,  taking a regular neighbourhood of $X^1\subset X^2$ yields a ribbon diagram.
 This ribbon graph is called a ribbon neighbourhood of $X^1$. 

Given a ribbon diagram with ribbon $X^2$,  one can construct  a surface diagram by embedding $\partial(X^2\times[0,1])$ in such a way that   $X^2\times\{0\}$ coincides with the ribbon. 
This corresponds to doubling the ribbon, placing one copy of the ribbon in front of the other and gluing the two copies of the ribbon together at their boundaries. The resulting surface is the boundary of a tubular neighbourhood of the core. 
For standard diagrams, labelling with data from a ribbon category will result in the same evaluation.

\subsection{Invariance}\label{sec:invariance}

For a standard diagram one can choose a sufficiently small ribbon neighbourhood of $X^1$ so that it does not contain any folds. Then the projection of this ribbon to the projection plane is a regular mapping. The projection plane is assumed to have a canonical orientation, and so the projection map induces an orientation of the ribbon.

With these definitions, it is possible to consider the behaviour of the evaluation under mappings of the diagrams. The mappings of interest are the following.

\begin{definition} A homomorphism of standard surface diagrams $f\colon D\to D'$ is called an {\bf oriented homomorphism} if $f$ is an orientation-preserving map of a ribbon neighbourhood of  $X^1\subset D$ to a ribbon neighbourhood ${X'}^1\subset D'$.
\end{definition}

For diagrams that are labelled with a Gray category with duals, a mapping of diagrams determines the relation between the labels. The discussion is parallel to the two-dimensional case in Section \ref{2catdualssection}. By subdividing the diagrams,  one can restrict attention to isomorphisms of progressive diagrams. Surfaces are oriented by the projection $p_2$, and lines are oriented by the projection $p_1\circ p_2$, by comparing with standard orientations of the coordinates. 

For general diagrams, the mapping of vertices has a complicated structure. However, for standard surface diagrams the situation simplifies. Let $q$ be the projection map $p_2$ restricted to $X^2\subset D$ and $q'$ the corresponding map for  ${X'}^2\subset D'$.  The mapping of the projection plane $q'\circ f\circ q^{-1}$ is a local isomorphism near each vertex, and for an oriented homomorphism it is also orientation-preserving.

\begin{definition}\label{3dlabels}
  Let $f\colon D\to D'$ be an oriented isomorphism of progressive surface diagrams. The labels of $D'$ are called {\bf induced} from the labels of $D$ by $f$ if
  \begin{enumerate}
  \item  the labels on a region of $D$ and its image in $D'$ are equal
  \item  the labels on a line and its image are equal if the orientation of the line is preserved; they are related by $*$ if the orientation is reversed
  \item  the labels on a surface and its image are equal if the orientation of the surface is preserved; they are related by $\#$ if the orientation of the surface is reversed
  \item the label on a vertex of $D'$ is induced from the label on the corresponding vertex of $D$ by the mapping of the projection plane in a neighbourhood around the vertex, using Definition \ref{rotatedvertex}.
  \end{enumerate}
\end{definition}

The main invariance result follows.  In its most general form, it relies
 on the conjecture that the moves on a projection of a PL two-manifold are the analogues of the moves in the smooth case. To our knowledge, this problem has not been investigated in the literature.

  \begin{conjecture} \label{conjecture} An isotopy of a surface projection with singularities in a compact subset of the surface can be adjusted so that the moves are the Reidemeister II and III moves for the folds, sliding of folds over or under cusps, or the moves in which cusps appear or disappear in pairs as given in Figures \ref{whitney} and \ref{whitney_proj} or similar figures obtained by reversing any of the axes.
  \end{conjecture} 

We prove the result for isomorphisms of surface diagrams that can be transformed to the identity by an isotopy that has the properties listed in Conjecture \ref{conjecture}.

\begin{theorem} \label{invariance_surface} Let $D$ and $D'$ be standard surface diagrams that are labelled with a spatial Gray category. Let $f\colon D\to D'$ be an oriented isomorphism of standard surface diagrams that is related to the identity by an isotopy that projects to a composite of isotopies of 2d diagrams and the following moves  for surface projections with singularities
\begin{enumerate}
\item the Reidemeister II and III moves for the ribbon neighbourhood of $X^1\subset X^2$ in Figure \ref{rtmoves} a) and \ref{rtmoves} c), 
\item Kauffman's double twist cancellation move for the ribbon neighbourhood of $X^1\subset X^2$ Figure \ref{rtmoves} b),
\item The move sliding a line over a vertex for the ribbon neighbourhood of $X^1\subset X^2$  from Figure \ref{rtmoves} d),
\item the Reidemeister II and III moves for the folds, 
\item sliding of folds over or under cusps, 
\item  the moves in which cusps appear or disappear in pairs as given in Figures \ref{whitney} and \ref{whitney_proj}, 
\item  similar figures obtained by reversing any of the axes.
\end{enumerate}
If  the labels of $D'$ are  induced from $D$ by $f$, then
 the evaluations of $D$ and $D'$ are equal.   
\end{theorem}

\begin{proof} 
     Invariance under (1) the Reidemeister II and III moves for the ribbon neighbourhood of $X_1\subset X_2$ follows from Theorem \ref{th_progressive_invariant}.
  Invariance under (2) the Kauffman double-twist cancellation is the equation
  \begin{equation} \left(*\Gamma_{*\nu}\right)^{-1}\cdot \#\left(\Delta^*_{\nu^*}\cdot\Delta^{-1}_\nu\right)\cdot *\Gamma_{*\nu}=1_\nu\end{equation}
  for all 2-morphisms $\nu$, in which $\#$ is the  functor from  \eqref{hashdualformula}. This equation, whose left-hand side is shown in Figure \ref{kauffman} and \ref{kauffman_slice}, follows from the spatial condition. Invariance under (3) the sliding move is a consequence of the naturality of the tensorator. It can be seen directly  that  each of these moves maps standard diagrams to standard diagrams.  The invariance under (4) the Reidemeister II and III moves for folds  and  (5)  the sliding moves is proved in the same way as above for the ribbon moves.

The evaluation is invariant under the moves (6) in which cusps appear or disappear in pairs according to the equations depicted in Figure \ref{whitney}. For the diagrams obtained by reversing axes, the corresponding equations are obtained by applying to both sides the $*$ operation or the inverse of vertical composition. 

Invariance under isotopies of 2d diagrams follows, because the projected diagrams are labelled by the planar 2-categories  $\mac G(\mac C,\mac D)$.  Hence, their  evaluation  is invariant if the labels for this two-dimen\-sio\-nal diagram are induced according to Definition \ref{rotatedvertex}. For the regions, lines and vertices projected from regions, lines and vertices of $D$ this follows from Definition \ref{3dlabels}. For folds and cusps, this follows from the fact that $*$ is a rotation in the projection plane, and for crossings from Lemma \ref{sigma_ids}.

  \begin{figure}
    \centering
    \includegraphics[scale=0.4]{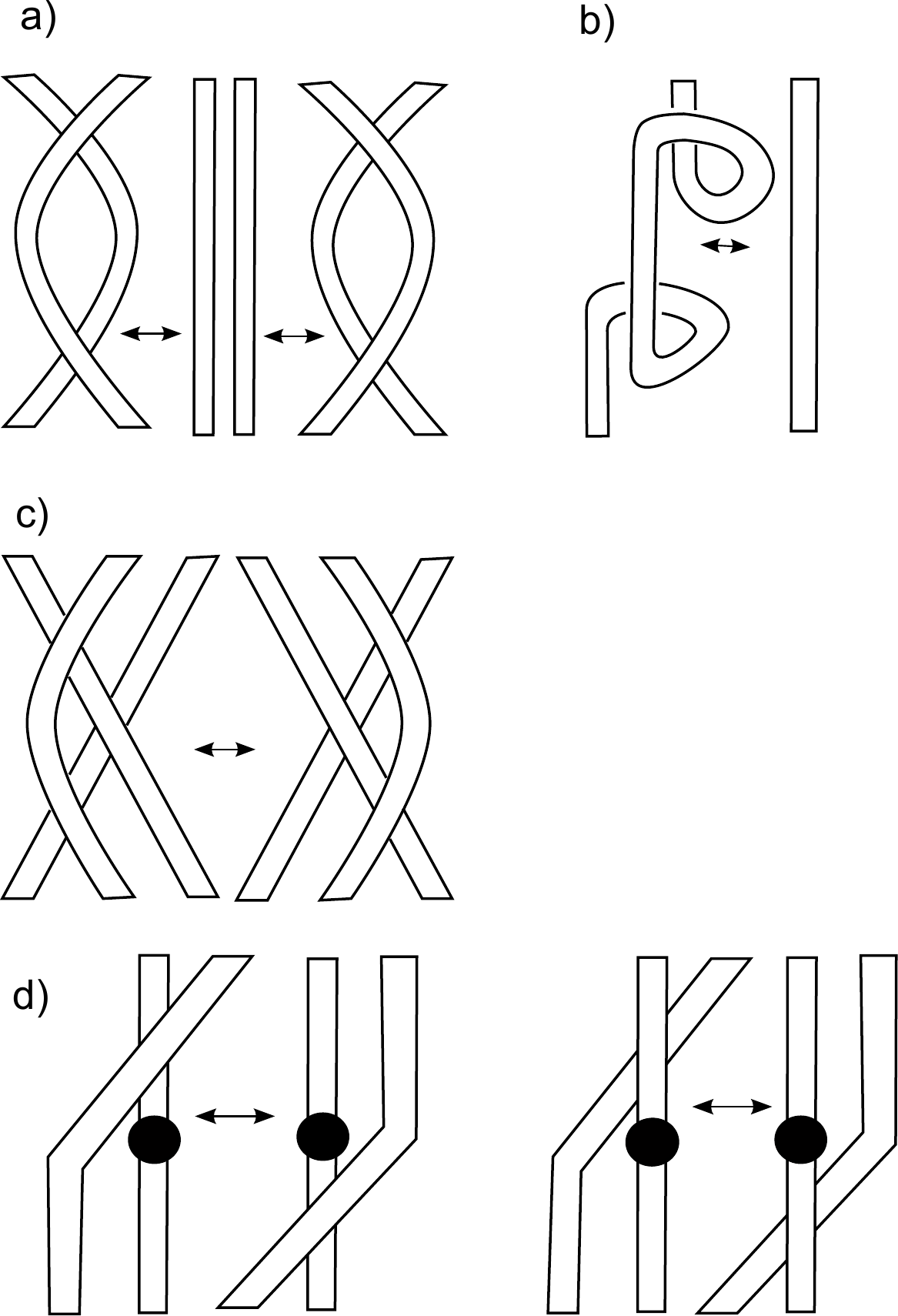}
    \caption{ Moves relating different  projections of \newline ribbon diagrams:\newline
      a) Reidemeister II move. \newline
      b) Double twist cancellation move.\newline
      c) Reidemeister III move.\newline
      d) The two additional moves from \cite {RT}.
    }
    \label{rtmoves}
  \end{figure}

  \begin{figure}
    \centering
    \includegraphics[scale=0.4]{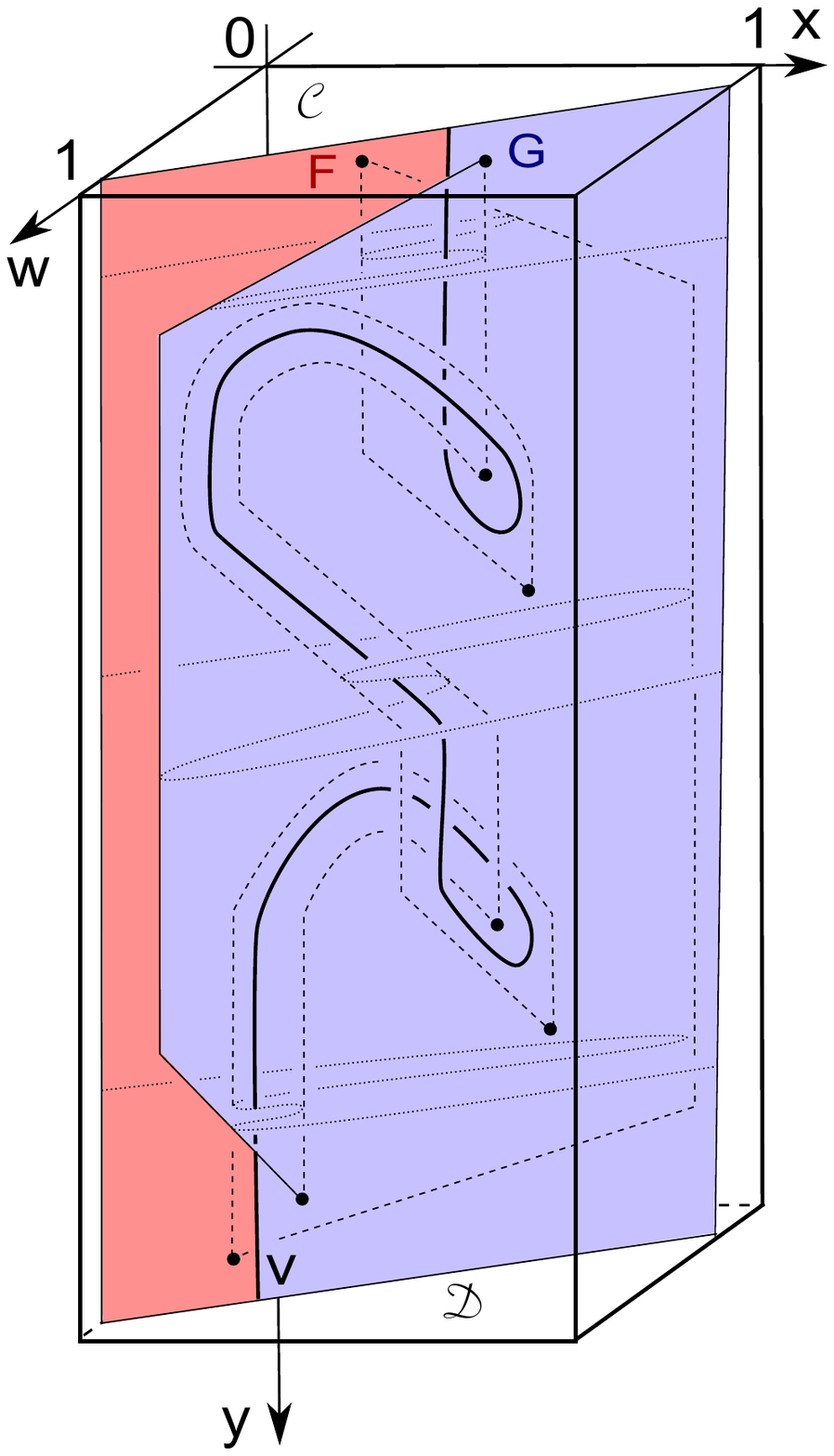}
    \caption{The Kauffman double-twist cancellation.}
    \label{kauffman}
  \end{figure}

  \begin{figure}
    \centering
    \includegraphics[scale=0.7]{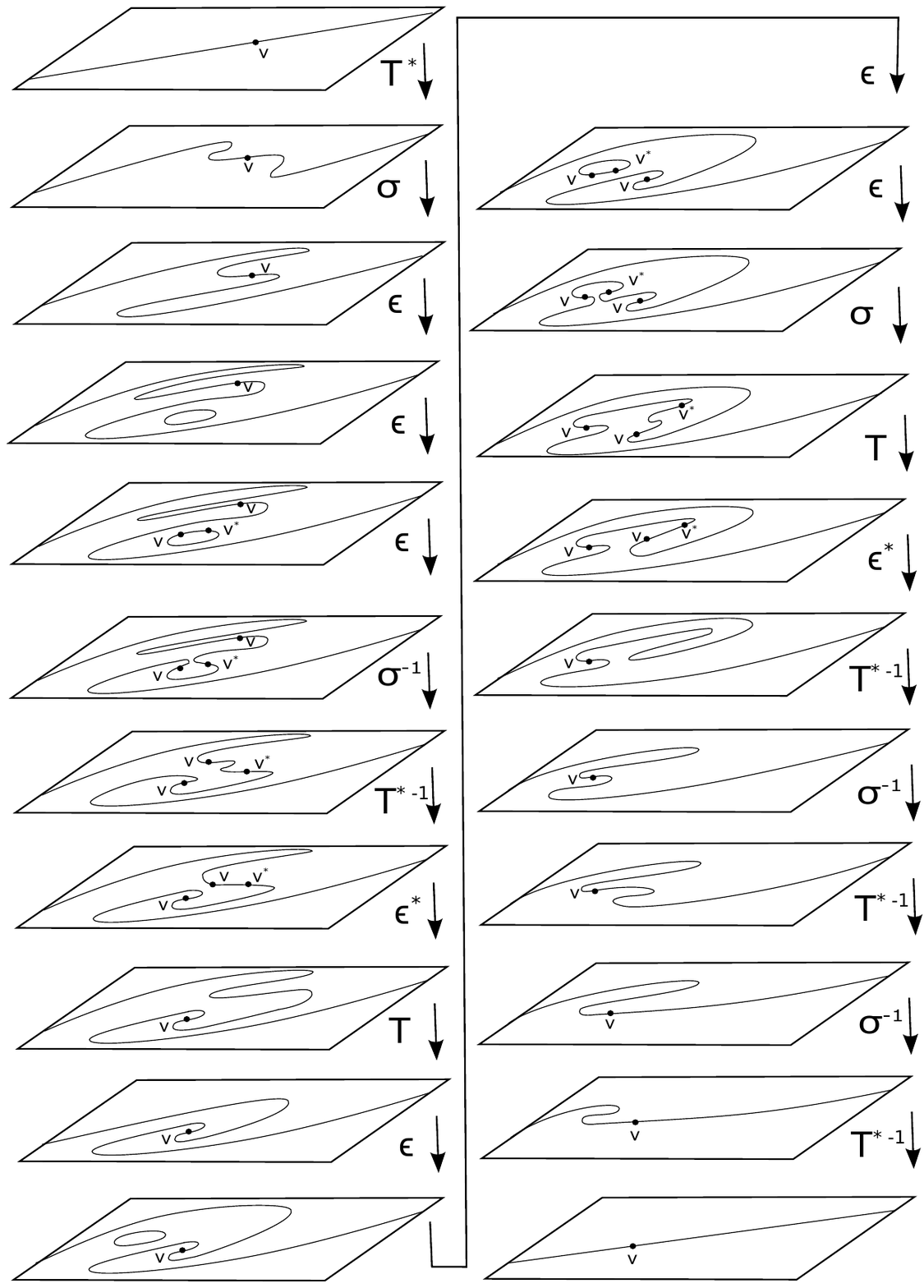}
    \caption{The Kauffman double-twist cancellation in movie representation, obtained by taking constant height slices in Figure \ref{kauffman}.}
    \label{kauffman_slice}
  \end{figure}

\end{proof}

\begin{corollary}\label{cor:fullinv} Let $D$ and $D'$ be standard surface diagrams that are labelled with a spatial Gray category. Let $f\colon D\to D'$ be an oriented isomorphism of standard surface diagrams and the labels of $D'$ induced from $D$ by $f$.
 If Conjecture \ref{conjecture} is true, then the evaluations of $D$ and $D'$ are equal. 
\end{corollary}

\begin{proof}
  For each homomorphism $f$, the `Alexander trick' guarantees that there is an isotopy of the diagram from the identity to $f$ \cite[Proposition 3.22]{RS}. As in the proof of Theorem \ref{th_progressive_invariant}, this isotopy can be chosen in such a way that its  effect on the two-dimensional diagram obtained by projection with $p_2$ is an isotopy of the two-dimensional projection plane punctuated by a finite sequence of moves that generalise the Reidemeister moves. 

  First, the isotopy is factored into a product of moves for the ribbon neighbourhood of $X^1\subset X^2$. These moves will be called ribbon moves. At this stage, the action of the isotopy on the surfaces is not considered, except for the requirement that for each move the projection of the surface singularities avoids the move in the projection plane. If this condition is not satisfied, then the isotopy can be adjusted by a small perturbation so that it is satisfied.

  Reidemeister moves for blackboard-framed links projected to the plane were given by Kauffman \cite{K-RI, K} and, more explicitly, by Freyd and Yetter \cite{FY}. These moves are the (1) Reidemeister II and III moves (Figure \ref{rtmoves} a) and \ref{rtmoves} c)), plus (2) Kauffman's double-twist cancellation move (Figure \ref{rtmoves} b). The additional moves for ribbon graphs were given by Reshetikhin and Turaev \cite{RT,TBK} and consist of (3) sliding a line over or under a vertex (Figure \ref{rtmoves} d).
  
   It remains to consider the effect of the isotopies between the ribbon moves.  
  Between these moves, the projection of the ribbon neighbourhood changes by an isotopy of the projection plane, and this also preserves the fact that the diagram is standard. However, the singularities of the projection of the surfaces change during these isotopies.

If Conjecture \ref{conjecture} holds, these changes are given by  
 (4) the Reidemeister II and III moves for the folds, (5) sliding of folds over or under cusps, or (6) the moves in which cusps appear or disappear in pairs as given in Figures \ref{whitney} and \ref{whitney_proj} or (7) similar figures obtained by reversing any of the axes.
  
   Between surface moves, the diagram can be subdivided to a diagram $E$ by introducing an additional vertex for each cusp and an additional line for each fold, and labelling with the canonical 2- and 3-morphisms for folds and cusps. Note that $E$ need not be progressive, but the projection of each surface in $E$ is regular. Between the surface moves, the subdivision to $E$ is preserved by the isotopy.  

  The topology of the two-dimensional diagram that is given by the projection of $E$ is determined by the 1-skeleton of $E$. The associated moves induced by the isotopy are the moves for a projection of a graph. These moves are (4) the Reidemeister II and III moves and (5) the Reshetikhin-Turaev sliding move. Note that the Kauffman double-twist cancellation move does not occur. This is because this move is always accompanied by  pairs of cusps appearing or disappearing, and all such cusp cancellations have already been accounted for in the surface moves or ribbon moves. 
 Finally, between these moves, the diagram changes by isotopies of the projection plane.

  By Theorem \ref{invariance_surface} the evaluation of a diagram is invariant under such two-dimensional isotopies as well as the moves (1) to (7) considered above, and this proves the claim.

\end{proof}

Note that in the proof the spatial condition is only required 
when regarding the 1-skeleton $X^{1}$. Thus we immediately obtain 
\begin{corollary}
  \label{corollary:Invariant-gen-Gray-cat}
Let $\mac G$ be a (not necessarily spatial) Gray category with duals and $D$ a surface diagram with $X^{1} = \emptyset$. Then if Conjecture \ref{conjecture} is true, the evaluation of $D$ is invariant under oriented surface isomorphisms. 
\end{corollary}

This result can also be understood from the fact that the invariant uses only a subcategory of $\mathcal G$ and this subcategory is in fact spatial, as summarised in the following Remark.

\begin{remark}
  Given a Gray category $\mac G$ with duals, it is a lengthy but straightforward computation to check that the following construction gives   a subcategory $\widetilde{ \mac G}$, which is a spatial Gray category with duals. 
The objects and 1-morphisms of $\widetilde{\mac G}$ are those of $\mac G$, the 2-morphisms  are 
 all morphisms  generated 
from the identities and the fold morphisms $\eta_{F}$ and  $\eta_{F}^{*}$, using the compositions $\Box$ and $\circ$. 
As 3-morphisms we take all composites of identities, $\epsilon$ and $\epsilon^*$ for all 2-morphisms,  tensorators and their inverses and triangulators with their duals and inverses. 
\end{remark}


\section*{Acknowledgements} The work of Catherine Meusburger and Gregor Schaumann was supported by the Emmy-Noether Fellowship ME 3425/1-1 of the German Research Foundation. The work of John Barrett was supported by the STFC Particle Theory grant  ST/L000393/1, and the Simons Foundation as visiting scientist at the GGI. All three authors would like to thank the Galileo Galilei Institute for Theoretical Physics for support during the workshop Emergent Geometries from Strings and Quantum Fields in June 2023.

The results of this paper were presented at `New Perspectives in Topological Field Theories' in Hamburg 27 - 31 August 2012, with 
support from the European Science Foundation Research Networking Programme ITGP (Interactions of Low-Dimentional Topology and Geometry with Mathematical Physics).

\appendix

\section{Functors of strict tricategories, natural transformations and modifications}
\label{2functors}

In this appendix, we define functors of strict tricategories  and their  natural transformations and modifications by  specialising the associated definitions for tricategories in \cite{GPS, gurski}.
For completeness, we also record the standard definitions for  functors of (strict) 2-categories \cite{KS,BasicBicat}.

\begin{definition}
\label{laks2func}
  A {\bf lax 2-functor} $F\colon\mathcal C\rightarrow \mathcal D$ between 2-categories $\mathcal C,\mathcal D$ is given by the following data
  \begin{itemize}
  \item A function $F_0\colon \text{Ob}(\mathcal C)\rightarrow \text{Ob}(\mathcal D)$.
  \item For all objects $G,H$ of $\mathcal C$, a functor $F_{G,H}\colon \mathcal C_{G,H}\rightarrow \mathcal D_{F_0(G), F_0(H)}$.
  \item For all objects $G,H,K$ of $\mathcal C$ a natural transformation $\Phi_{GHK}\colon\circ\;(F_{H,K}\times F_{G,H})\to F_{G,K}\;\circ$. These determine, for all 1-morphisms $\nu\colon G\rightarrow H$, $\mu\colon H\rightarrow K$, a 2-morphism $\Phi_{\mu,\nu}\colon F_{H,K}(\mu)\circ F_{G,H}(\nu)\rightarrow F_{G,K}(\mu\circ\nu)$. 
  \item For all objects $G$, a 2-morphism $\Phi_G\colon 1_{F_0(G)}\to F_{G,G}(1_G)$.
  \end{itemize}
  The function $F_0$, the functors $F_{G,H}$ and the 2-morphisms $\Phi_{\mu,\nu}$ and $\Phi_G$ are required to satisfy the following consistency conditions
  \begin{enumerate}
  \item For all  1-morphisms $\nu\colon G\rightarrow H$
    $$\Phi_{1_H, \nu}\cdot(\Phi_H\circ 1_{F_{G,H}(\nu)})=\Phi_{\nu, 1_G}\cdot(1_{F_{G,H}(\nu)}\circ\Phi_G)=1_{F_{G,H}(\nu)}.$$
  \item  For all 1-morphisms $\nu\colon G\rightarrow H$, $\mu\colon H\rightarrow K$, $\rho\colon K\rightarrow L$, the following diagram commutes
    \begin{align*}
      \xymatrix{
        F_{K,L}(\rho)\circ F_{H,K}(\mu)\circ F_{G,H}(\nu)\ar[d]^{1\circ \Phi_{\mu,\nu}} \ar[r]^{\qquad\Phi_{\rho,\mu}\circ 1} & F_{H,L}(\rho\circ\mu)\circ F_{G,H}(\nu)\ar[d]^{\Phi_{\rho\circ \mu,\nu}}\\
        F_{K,L}(\rho)\circ F_{G,K}(\mu\circ \nu) \ar[r]^{\Phi_{\rho,\mu\circ \nu}} & F_{G,L}(\rho\circ\mu\circ \nu).
      }
    \end{align*}
  \end{enumerate}

  A {\bf weak 2-functor} (also called a strong 2-functor, pseudo-functor, or homomorphism) is a lax 2-functor in which all 2-morphisms $\Phi_{\mu,\nu}$ and $\Phi_G$ are invertible. A weak 2-functor is said to have {\bf strict units} if the 2-morphisms  $\Phi_G$ are all identities, and it is
  called {\bf strict} if the 2-morphisms $\Phi_{\mu,\nu}$ and $\Phi_G$ are all identities. In this case, one has
  $$F_{G,K}(\mu\circ\nu)=F_{H,K}(\mu)\circ F_{G,H}(\nu),\qquad 1_{F_0(G)}=F_{G,G}(1_G).$$

  There is an analogous  definition with the arrows labelled by $\Phi_{\mu,\nu}$ and $\Phi_G$ reversed. In this case, the functor is called an  {\bf op-lax 2-functor}.
\end{definition}

In the following, we will also require the notion of cubical and opcubical functors between certain 2-categories. Our definition is a special case of the 
definition of cubical and opcubical functors from \cite{GPS, gurski}.
\begin{definition}
  \label{definition:cubical-opcubical-functors}
  Let $\mathcal C, \mathcal D$ and $\mathcal E$ be 2-categories. A functor $F\colon \mathcal C \times \mathcal D \rightarrow \mathcal E$ with coherence  isomorphisms 
  \begin{equation*}
    \Phi_{\mu,\nu}\colon F( (\mu_{1} , \mu_{2})) \circ F( (\nu_{1}, \nu_{2})) \rightarrow F((\mu_{1} \circ \nu_{1}), (\mu_{2} \circ \nu_{2})),
  \end{equation*}
  for $\circ$-composable 1-morphisms $\mu=(\mu_1,\mu_2)$ and $\nu=(\nu_1,\nu_2)$ in $\mac C\times\mac D$ 
  is called {\bf (op)cubical}, if  the 2-morphism $ \Phi_{\mu,\nu}$ is the identity in case $\mu_{1}$ or $\nu_{2}$ ($\mu_{2}$ or $\nu_{1}$) is an identity 1-morphism. 
\end{definition}

The following notion of natural transformation of 2-functors  adopts  the convention of \cite{GPS,gurski} and is sometimes also referred to as `op-lax 2-transformation'.

\begin{definition}
\label{nat2transformations}

  A {\bf natural transformation} $\rho\colon F\to G$ between lax 2-functors \newline
  $F=(F_0,F_{A,B}, \Phi_{\mu,\nu}, \Phi_A)\colon\mac C\to\mac D, G=(G_0,G_{A,B}, \Psi_{\mu,\nu}, \Psi_{A})\colon\mathcal C\to\mathcal D$ \newline  is given by the following data:
  \begin{itemize}
  \item For all objects $A$ of $\mathcal C$, a 1-morphism $\rho_{A}\colon F_0(A)\to G_0(A)$.
  \item For all objects $A,B$ of $\mathcal C$ a natural transformation $$\rho_{A,B}\colon (\rho_B\circ -) F_{A,B} \rightarrow  (-\circ \rho_A) G_{A,B} ,$$ where  $-\circ \rho_A\colon\mathcal D_{G_0(A), G_0(B)}\to\mathcal D_{F_0(A), G_0(B)}$ and  $\rho_B\circ -\colon \mathcal D_{F_0(A),F_0(B)}$ $\to \mathcal D_{F_0(A), G_0(B)}$ denote the functors given by pre- and post-composition with $\rho_A$ and $\rho_B$. These natural transformations determine for all 1-morphisms $\mu\colon A\to B$ a 2-morphism
    $\rho_\mu\colon \rho_B\circ F_{A,B}(\mu) \to  G_{A,B}(\mu)\circ \rho_A$.
  \end{itemize}
  The 1-morphisms $\rho_A$ and 2-morphisms $\rho_{\mu}$ are required to satisfy  the following consistency conditions:
  \begin{enumerate}
  \item For all 1-morphisms $\nu\colon A\to B$ and $\mu\colon B\to C$ the following diagram commutes
    \begin{align}
      \xymatrix{\rho_C\circ F_{B,C}(\mu)\circ F_{A,B}(\nu)  \ar[d]^{1\circ \Phi_{\mu,\nu}}  \ar[r]^{\rho_\mu \circ 1} & G_{B,C}(\mu)\circ \rho_B\circ F_{A,B}(\nu)
        \ar[d]^{1 \circ \rho_\nu} \\
        \rho_C\circ F_{A,C}(\mu\circ\nu) \ar[d]^{\rho_{\mu\circ\nu}} &  G_{B,C}(\mu)\circ G_{A,B}(\nu)\circ \rho_A  \ar[ld]^{\Psi_{\mu,\nu}\circ 1 }\\
        G_{A,C}(\mu\circ\nu)\circ \rho_A    .
      }\nonumber
    \end{align}
  \item For all objects $A$ of $\mathcal C$ the following diagram commutes
    \begin{align}
      \xymatrix{
        1_{G_0(A)}\circ \rho_A=  \rho_A=\rho_A\circ 1_{F_0(A)} \ar[d]^{ 1\circ \Phi_A} \ar[rd]^{\Psi_A\circ 1}\\
        \rho_A\circ F_{A,A}(1_A)  \ar[r]^{\rho_{1_A}} &  G_{A,A}(1_A)\circ \rho_A .
      }\nonumber
    \end{align}
  \end{enumerate}
  A {\bf pseudo-natural transformation} $\rho\colon F\to G$ of lax 2-functors $F,G\colon\mathcal C\to\mathcal D$ is a natural transformation of lax 2-functors in which all  
  2-morphisms $\rho_{\mu}\colon \rho_A\circ F_{A,B}(\mu) \to G_{A,B}(\mu)\circ \rho_A$ are isomorphisms. A pseudo-natural transformation is called {\bf invertible} if all the 1-morphisms $\rho_{A}$ are invertible.
  A {\bf natural isomorphism} is a pseudo-natural transformation in which for every object $A$, $F_0(A)= G_0(A)$ and the 1-morphism $\rho_A$ is the identity.
\end{definition}

It is easy to see that an invertible pseudo-natural transformation has indeed a unique inverse pseudo-natural transformation.

\begin{definition}
  \label{modi_gen}
  Let $\rho=(\rho_{A}, \rho_{A,B})\colon F\to G$ and  
$\tau=(\tau_A,\tau_{A,B})\colon F\to G$ be natural transformations between lax 2-functors $F=(F_0, F_{A,B}, \Phi_{\mu,\nu}, \Phi_A)$ and $G=(G_0,G_{A,B}, \Psi_{\mu,\nu}, \Psi_{A})\colon\mathcal C\to\mathcal D$. A {\bf modification} $\Psi\colon \rho\Rightarrow\tau$ is a collection of 2-morphisms $\Psi_A\colon \rho_A\Rightarrow\tau_A$ for every object $A$ of $\mac G$ such that for all 1-morphisms $\mu\colon A\to B$ 
  $$
  \tau_\mu\cdot (\Psi_A\circ 1_{F_{A,B}(\mu)})= (1_{G_{A,B}(\mu)}\circ \Psi_B) \cdot \rho_\mu.
  $$
  A modification is called {\bf invertible} if all 2-morphisms $\Psi_A$ are invertible.
\end{definition}

In terms of these definitions for 2-categories,  the concepts of functors of strict tricategories, natural transformations and modifications can be formulated. 

  The definition of a strict tricategory is given in \cite{GPS}. This can be summarised informally as follows:

\begin{definition}
 A {\bf strict  tricategory} is a tricategory $(\mac G, \Box,\circ,\cdot)$ in which the composition $\Box$ is strictly associative and unital.
 \end{definition}

Note that in general there are some differences between the definition of tricategory in \cite{GPS} and `algebraic tricategory' in \cite{gurski}.  However, for the strict  tricategories considered in this paper  these definitions coincide. 

It is important to note that a strict tricategory is \emph{not} a 3-category. In particular,  $\mac G(C,D)$ is a bicategory for all objects $C,D$ and $\Box$ is a set of weak 2-functors. In this paper, the only cases of interest are where these bicategories are 2-categories, and there is the additional condition that the tricategories are either cubical or opcubical.

\begin{definition}
\label{strictopcubical}
  A {\bf 2-strict  tricategory} is a strict tricategory that satisfies the following additional condition 
  \begin{enumerate}
  \item For all objects $\mac C,\mac D$ the bicategory $\mac G(\mac C,\mac D)$ is a strict 2-category.
     \end{enumerate}
     A {\bf strict (op)cubical tricategory}
     is a 2-strict tricategory that also satisfies the conditions
      \begin{enumerate}
       \item[(2)] 
    $
    1_{1_{\mac C}}\circ 1_{1_\mac C}=1_{1_{\mac C}}.
    $
  \item[(3)] Each functor $\Box\colon \mac G(\mac D,\mac E)\times\mac G(\mac C,\mac D)\to\mac G(\mac C,\mac E)$ is (op)cubical, i.e., the invertible coherence 3-morphisms
    $$
    \Box_{\mu,\nu}\colon (\mu_1\Box\mu_2)\circ (\nu_1\Box \nu_2)\to (\mu_1\circ \nu_1)\Box (\mu_2\circ\nu_2)
    $$
    for $\mu=(\mu_1,\mu_2), \nu=(\nu_1,\nu_2)\in\mac G(\mac D,\mac E)\times \mac G(\mac C,\mac D)$
    are identity 3-morphisms  if $\mu_1$ or $\nu_2$  (in the cubical case), or $\mu_2$ or $\nu_1$ (in the opcubical case) is an identity 2-morphism.
  \end{enumerate}
\end{definition}

In the following we call a 1-morphism $F \colon \mathcal{C} \rightarrow \mathcal{D}$ in a strict tricategory  {\bf invertible}, if there exists a 1-morphism $G\colon \mathcal{D} \rightarrow \mathcal{C}$ with 
$F \Box G=1_{\mathcal{D}}$ and $G \Box F= 1_{\mathcal{C}}$. Similarly, a 2-morphism $\mu\colon F \rightarrow G$ in a strict tricategory  is called {\bf invertible}, if there exists a 2-morphism $\nu\colon G \rightarrow F$ with $\mu \circ \nu = 1_{G}$ and $\nu \circ \mu= 1_{F}$.
Note that the inverse 1-morphism $G$ and 2-morphism  $\nu$ are  determined uniquely.

A functor of strict tricategories  is a  functor of tricategories  that is compatible with the extra requirements that hold in case of a strict  tricategory. As we consider only tricategory functors between strict (op)cubical  tricategories, we will not give the most general definition of a tricategory functor, but refer the reader to  
\cite{GPS, gurski}  for details. Note that the definitions in \cite{GPS, gurski} differ slightly and that the definition of \cite{gurski} is stronger than 
that of \cite{GPS}, since the 
coherence data of a functor $F$  consists of adjoint equivalences instead of just equivalences in certain bicategories.  In particular, it contains an adjoint 
equivalence $\kappa\colon \widetilde\Box\circ (F \times F) \rightarrow F \circ \Box$. A functor $F\colon \mathcal{G} \rightarrow \tilde{\mathcal{G}}$  of 2-strict tricategories according to our definition is a 
functor of tricategories according to \cite{gurski}, where all coherence data consists of identities and $\kappa$ is a natural isomorphism, which is automatically an adjoint equivalence.

\begin{definition}
\label{grayfunc}
A {\bf \Glax functor}  $F\colon\mathcal G\rightarrow \tilde {\mathcal G}$ {\bf between 2-strict tricategories}
$\mathcal G, \tilde {\mathcal G}$ consists of \begin{itemize}
  \item a function $F_0\colon\text{Ob}(\mathcal G)\rightarrow \text{Ob} (\tilde {\mathcal G})$,
  \item weak 2-functors $F_{\mathcal C,\mathcal D}\colon \mathcal G(\mathcal C,\mathcal D)\rightarrow \tilde {\mathcal G}(F_0(\mathcal C), F_0(\mathcal D))$ for all objects $\mathcal C,\mathcal D$ of $\mathcal G$,
  \item    an  invertible pseudo-natural transformation of weak 2-functors $$\kappa_{\mathcal C,\mathcal D, \mathcal E}\colon \widetilde\Box (F_{\mathcal D,\mathcal E}\times F_{\mathcal C,\mathcal D})\rightarrow F_{\mathcal C,\mathcal E}\,\Box$$ for all objects $\mathcal C,\mathcal D,\mathcal E$ of $\mathcal G$,
  \item  an invertible 2-morphism $\iota_{\mathcal{C}}\colon F_{\mathcal{C},\mathcal{C}}(1_{\mathcal{C}}) \rightarrow 1_{F_{0}(\mathcal{C})}$ for all objects $\mac C$ of $\mac G$, 
  \end{itemize}
  such that the following consistency conditions are satisfied
  \begin{enumerate} 
  \item For all objects $\mathcal{B,C,D,E}$ of $\mathcal G$
    $$\left(\kappa_\mathcal{BCE}(\Box\times I)\right)\circ\left( \widetilde\Box(\kappa_\mathcal{CDE}\times I)\right)=\left(\kappa_\mathcal{BDE}(I\times\Box)\right)\circ\left(\widetilde\Box(I\times\kappa_\mathcal{BCD})\right).$$ In this formula $I$ is the identity functor, the unnamed product is the Gray product in 2Cat, and $\circ$ is the horizontal composition of pseudo-natural transformations. 
  \item For all objects $\mathcal{C,D}$ of $\mathcal G$, $\kappa_\mathcal{CCD}(I\times F_{\mathcal{C},\mathcal{C}}(I_{\mathcal{C}}))= \widetilde \Box (1_{F_{\mathcal C,\mathcal D}} \times \iota_{\mathcal{C}})$, where
    $I_{\mathcal C}$ is the strict functor from the trivial 2-category that has as image the object $1_{\mathcal C}$ of  $\mathcal G(\mathcal C ,\mathcal C)$ and
    $\iota_{\mathcal{C}}$ is considered as a natural transformation of 2-functors $\iota_{\mathcal{C}}\colon F_{\mathcal{C},\mathcal{C}}(I_{\mathcal{C}}) \rightarrow I_{F_{0}(\mathcal{C})}$.
  \item For all objects $\mathcal{C,D}$ of $\mathcal G$, $\kappa_\mathcal{CDD}(I_{\mathcal D}\times I)=\widetilde \Box (\iota_{\mathcal{D}} \times 1_{F_{\mathcal C,\mathcal D}})$. 
  \end{enumerate}
  The \Glax functor is called a {\bf  functor} (or  {\bf weak functor with strict units}) {\bf of 2-strict} tricategories,   if additionally 
  \begin{enumerate}
    \addtocounter{enumi}{4}
  \item $\kappa$ is a natural isomorphism,
  \item for  all objects $\mathcal C$ of $\mathcal G$, $F_{\mathcal C,\mathcal C}(1_{\mathcal C})=1_{F_0(\mathcal C)}$ and $\iota_{\mathcal{C}}=\iota_{\mathcal{C}}^{-1}=1_{1_{\mathcal{C}}}$.
  \end{enumerate}
  We  call the map $F_{0}$ and the maps of $F_{\mathcal{C},\mathcal{D}}$ the {\bf mappings} of $F$, while all the other data is called {\bf coherence morphisms} or {\bf coherence data} of $F$.
\end{definition}

Unpacking the definition of a functor of 2-strict tricategories  leads to the following explicit description of  the coherence data for $F$. 
The weak 2-functors $F_{\mac C,\mac D}$  have as coherence data a collection of invertible 3-morphisms 
$$\Phi_{\mu,\nu}\colon F_{\mac C,\mac D}(\mu)\mathbin{\widetilde\circ} F_{\mac C,\mac D}(\nu)\to F_{\mac C,\mac D}(\mu\circ\nu)$$
for all $\circ$-composable 2-morphisms $\mu$, $\nu$ in $\mathcal G(\mac C,\mac D)$ and for each 1-morphism $G$ in $\mac G(\mac C,\mac D)$ an invertible  3-morphism
$$\Phi_G\colon 1_{F_{\mac C,\mac D}(G)}\to F_{\mac C,\mac D}(1_G),$$
which satisfy the axioms in Definition \ref{laks2func}.
The natural isomorphisms of weak 2-functors $\kappa_{\mathcal C,\mathcal D, \mathcal E}$ 
are characterised by invertible 3-morphisms
$$\kappa_{\mu,\nu}\colon  F_{\mac D,\mac E}(\mu)\mathbin{\widetilde\Box} F_{\mac C,\mac D}(\nu) \to F_{\mac C,\mac E}(\mu\Box\nu)$$
for  all  $\Box$-composable 2-morphisms $\mu\in \mac G(\mac D,\mac E)$, $\nu\in \mac G(\mac C,\mac D)$.
The  conditions in Definition  \ref{nat2transformations}  take the following form:

\medskip
$\bullet$ For all 1-morphisms $G_1\in\mac G(\mac D,\mac E)$, $G_2\in\mac G(\mac C,\mac D)$ one has
\begin{equation} F_{\mac D,\mac E}(G_1)\widetilde\Box F_{\mac C,\mac D}(G_2)=F_{\mac C,\mac E}(G_1\Box G_2).\end{equation}

$\bullet$ For all  2-morphisms $\mu,\rho\in\mac G(\mac D,\mac E)$ and $\nu,\tau\in\mac G(\mac C,\mac D)$ such that $\mu,\rho$ and $\nu,\tau$ are $\circ$-composable,  the following diagram commutes
\begin{align}\label{tensoratorcompatibility}
  &\xymatrix{ (F_{\mac D,\mac E}(\mu)\widetilde\Box F_{\mac C,\mac D}(\nu))\mathbin{\widetilde\circ} (F_{\mac D,\mac E}(\rho)\widetilde\Box F_{\mac C,\mac D}(\tau))
    \ar[d]_*{ \widetilde\Box_{F_{\mac D,\mac E}(\mu),F_{\mac C,\mac D}(\nu),F_{\mac D,\mac E}(\rho),F_{\mac C,\mac D}(\tau)}}
    \ar[r]<4pt>^-*{\quad\kappa_{\mu,\nu}\widetilde\circ \kappa_{\rho,\tau}}
    &  F_{\mac C,\mac E}(\mu\Box\nu)\mathbin{\widetilde\circ} F_{\mac C,\mac E}(\rho\Box\tau) 
    \ar[d]^*{\Phi_{\mu\Box\nu,\rho\Box\tau}}
    \\
    (F_{\mac D,\mac E}(\mu)\widetilde\circ F_{\mac D,\mac E}(\rho))\mathbin{\widetilde\Box}(F_{\mac C,\mac D}(\nu)\widetilde\circ F_{\mac C,\mac D}(\tau)) 
    \ar[d]_*{\Phi_{\mu,\rho}\widetilde\Box\Phi_{\nu,\tau}}
    &  F_{\mac C,\mac E}((\mu\Box\nu)\circ(\rho\Box\tau)) 
    \ar[d]^*{F_{\mac C,\mac E}(\Box_{\mu,\nu,\rho,\tau}) }
    \\
    F_{\mac D,\mac E}(\mu\circ\rho)\mathbin{\widetilde\Box} F_{\mac C,\mac D}(\nu\circ\tau) 
    \ar[r]_*{\kappa_{\mu\circ\rho,\nu\circ\tau}}
    &    F_{\mac C,\mac E}((\mu\circ\rho)\mathbin\Box(\nu\circ\tau)).
   }
\end{align}

$\bullet$ All 1-morphisms $G\in\mac  G(\mac D,\mac E)$, $H\in\mac G(\mac C,\mac D)$ satisfy 
\begin{align}
  &\xymatrix{
    1_{F_{\mac C,\mac E}(G\Box H)}=1_{F_{\mac D,\mac E}(G)}\widetilde\Box 1_{F_{\mac C,\mac D}(H)} 
    \ar[d]_ {\Phi_G\widetilde\Box \Phi_H}
    \ar[rd]^ {\Phi_{G\Box H}}
    \\
    F_{\mac D,\mac E}(1_G)\widetilde\Box F_{\mac C,\mac D}(1_H) \ar[r]_*{\kappa_{1_G,1_H}\qquad} 
    &   F_{\mac C,\mac E}(1_{G\Box H})=F_{\mac C,\mac E}(1_G\Box 1_H).
  }\label{unitcompatibility}
\end{align}
Condition (1)  in Definition \ref{grayfunc} states  that
the  diagram 
\begin{align}\label{Boxcompatibility}
  & \xymatrix{F_{\mac D,\mac E}(\nu) \widetilde\Box F_{\mac C,\mac D}(\mu) \widetilde\Box F_{\mac B,\mac C}(\rho)   \ar[r]^{1 \widetilde \Box \kappa_{\mu,\rho} } \ar[d]^{\kappa_{\nu,\mu} \widetilde\Box 1}  &  F_{\mac D,\mac E}(\nu) \widetilde\Box F_{\mac B,\mac D}(\mu \Box \rho)  \ar[d]^{ \kappa_{\nu,\mu \Box \rho}} \\
    F_{\mac C,\mac E}(\nu \Box \mu ) \widetilde\Box F_{\mac B,\mac C}(\rho)  \ar[r]^{\kappa_{\nu \Box \mu,\rho}} &  F_{\mac B,\mac E}(\nu \Box \mu \Box \rho) 
  } 
\end{align}
commutes for all 2-morphisms $\rho\in\mac G(\mac B,\mac C)$, $\mu\in\mac G(\mac C,\mac D)$, $\nu\in\mac G(\mac D,\mac E)$, and conditions (2), (3)  in Definition \ref{grayfunc} read
\begin{align*}
  & \kappa_{\mu,1_{\mathcal{C}}}=1_{F_{\mac C,\mac D}(\mu)} \colon F_{\mac C,\mac D}(\mu) = F_{\mac C,\mac D}(\mu) \widetilde \Box F_{\mac C,\mac C}(1_{\mathcal{C}}) \rightarrow  F_{\mac C,\mac D}(\mu \Box 1_{\mathcal{C}})=F_{\mac C,\mac D}(\mu), \\
  &\kappa_{1_{\mathcal{D}}, \mu}=1_{F_{\mac C,\mac D}(\mu)}  \colon  F_{\mac C,\mac D}(\mu)= F_{\mac D,\mac D}(1_{\mathcal{D}} ) \widetilde \Box F_{\mac C,\mac D}(\mu) \rightarrow  F_{\mac C,\mac D}(1_{\mathcal{D}}\Box \mu) =F_{\mac C,\mac D}(\mu).
\end{align*}

The notion of a functor of 2-strict tricategories in Definition \ref{grayfunc} thus corresponds to a trihomomorphism in \cite[Def 3.3.1]{gurski} for which
the adjoint equivalence $\chi$ in  \cite[Def 3.3.1]{gurski}  is a natural isomorphism given by the invertible 3-morphisms $\kappa_{\mu,\nu}$ and  for which  the adjoint equivalence $\mathbf \iota$ and the invertible modifications $\omega,\gamma,\delta$ in \cite[Def 3.3.1]{gurski} are trivial. 

There is an obvious composition of functors of 2-strict  tricategories, that is a special case of the general composition of functors between tricategories considered in \cite{gurski}.  Although in general the composition of functors between tricategories is not strictly associative, this is the case for functors between 2-strict tricategories. 
\begin{lemma}
  \label{composition-functors-asso}
  Let  $F\colon\mathcal{G} \rightarrow \mathcal{H}$, $G\colon\mathcal{H} \rightarrow \mathcal{K}$, $H\colon \mathcal{K} \rightarrow \mathcal{L}$
  be weak  functors of 2-strict tricategories. Then $H  (G  F)= (H G) F $.
\end{lemma}
\begin{proof}
  In \cite[Prop. 4.2.3]{gurski}, explicit expressions for a natural transformation  $\alpha\colon H  (G F) \rightarrow (H  G)   F $ are given.  It is easy to see that for 2-strict tricategories, the data from 
  which $\alpha$ is constructed consists entirely of identity mappings and morphisms. 
\end{proof}

In the following, we also require the notion of a strict functor of 2-strict tricategories. The standard definition of a strict functor of 2-strict tricategories is that of a functor of 2-strict tricategories 
for which all 2-functors $F_{\mathcal C,\mathcal D}$ in Definition \ref{grayfunc} are strict and all natural isomorphisms $\kappa_{\mathcal C,\mathcal D, \mathcal E}$ in Definition  \ref{grayfunc} are identities. 
However, as we only consider functors between strict cubical or opcubical tricategories, we change this definition slightly to adapt it to our setting. 

For this, note that a functor $F\colon \mac G\to\tilde{\mac G}$ between an opcubical  strict tricategory $\mac G$ and a cubical strict tricategory $\tilde {\mac G}$ can never be strict in the usual sense
unless the coherence morphisms $\Box_{\mu,\nu}$ from Definition \ref{strictopcubical} are trivial.  The strictness of $F$ and the fact that $\mac G$ is opcubical  imply
\begin{align*}
  F_{\mac C,\mac E}(\mu\Box\nu)=&F_{\mac C,\mac E}((\mu\Box 1_{H_2})\circ (1_{G_1}\Box\nu))\\
  =&(F_{\mac D,\mac E}(\mu)\widetilde {\Box} 1_{F_{\mac C,\mac D}(H_2)})\widetilde{\circ} (1_{F_{\mac D,\mac E}(G_1)}\widetilde {\Box} F_{\mac C,\mac D}(\nu))
\end{align*}
for all composable 2-morphisms $\mu\colon G_1\Rightarrow G_2\in\mac G(\mac D,\mac E)$, $\nu\colon H_1\Rightarrow H_2\in\mac G(\mac C,\mac D)$.  As $\tilde {\mac G}$ is cubical, one has
$$
F_{\mac D,\mac E}(\mu)\widetilde \Box F_{\mac C,\mac D}(\nu)=(1_{F_{\mac D,\mac E}(G_2)}\widetilde \Box F_{\mac C,\mac D}(\nu))\widetilde\circ (F_{\mac D,\mac E}(\mu)\widetilde \Box 1_{F_{\mac C,\mac D}(H_1)}).
$$
The two expressions cannot agree for all composable 2-morphisms unless the coherence morphisms $\widetilde {\Box}_{\mu,\nu}$ from Definition \ref{strictopcubical}  are trivial. 
For this reason, we modify the notion of strictness for the case of functors between opcubical and cubical strict tricategories and call such a functor strict if and only if its composition with the functor $\Sigma$ from Corollary \ref{lemma:cubical-to-opcubical-tricat} is strict in the usual sense. This amounts  to the requirement that the coherence morphism $\kappa_{\mac C,\mac D,\mac E}\colon \widetilde\Box (F_{\mac D,\mac E}\times F_{\mac C,\mac D})\to F_{\mac C,\mac E}\Box$ from Definition \ref{grayfunc} is given by the coherence morphisms $\tilde \Box^\inv_{\mu,\nu}$ from Definition \ref{strictopcubical}.

\begin{definition}
 \label{strictdef}
  Let  $\mac G$ and $\tilde {\mac G}$ be 2-strict tricategories  that are either cubical or opcubical and $F\colon\mac G\to\tilde {\mac G}$ a functor of 2-strict tricategories.  Then the functor $F$ is called {\bf strict} if
  for all objects $\mac C,\mac D$ of $\mac G$ the 2-functors  $F_{\mathcal C,\mathcal D}$ are strict and the natural isomorphisms $\kappa_{\mac C,\mac D,\mac E}$ from Definition \ref{grayfunc} are
  \begin{itemize}
  \item  the  identity morphisms in case $\mac G$ and $\tilde{\mac G}$ are both cubical or both opcubical,
  \item given by the 3-morphisms
    $$
    \kappa_{\mu,\nu}=\widetilde{\Box}^\inv_{({F_{\mac D,\mac E}(1_{G_2})}, F_{\mac C,\mac D}(\nu)), (F_{\mac D,\mac E}(\mu), {F_{\mac C,\mac D}(1_{H_1})})}
    $$
    for 2-morphisms  $\mu\colon G_1\Rightarrow G_2\in\mac G(\mac D,\mac E)$, $\nu\colon H_1\Rightarrow H_2\in\mac G(\mac C,\mac D)$
    in case $\mac G$ is cubical and $\tilde{\mac G}$ is opcubical,

  \item  given by  the 3-morphisms
    $$
    \kappa_{\mu,\nu}=\widetilde{\Box}^\inv_{({F_{\mac D,\mac E}(\mu)}, F_{\mac C,\mac D}(1_{H_2})), (F_{\mac D,\mac E}(1_{G_1}), {F_{\mac C,\mac D}(\nu)})}
    $$
    for 2-morphisms  $\mu\colon G_1\Rightarrow G_2\in\mac G(\mac D,\mac E)$, $\nu\colon H_1\Rightarrow H_2\in\mac G(\mac C,\mac D)$ in case $\mac G$ is opcubical and $\tilde {\mac G}$ is cubical.
  \end{itemize}
\end{definition}

\begin{definition}
\label{graynat}
  A {\bf natural transformation}  $\omega\colon F\to G$ between \Glax functors $F=(F_0, F_{\mathcal C,\mathcal D}, \kappa^F_{\mathcal B,\mathcal C,\mathcal D}, \iota^{F}_{\mac C})$, 
  $G=(G_0,G_{\mathcal C,\mathcal D}, \kappa^G_{\mathcal B,\mathcal C,\mathcal D}, \iota^{G}_{\mac C})\colon\mathcal G\to\tilde {\mathcal G}$ of 2-strict  tricategories consists of the following data: 
  \begin{itemize}
  \item For all objects $\mathcal C$ of $\mathcal G$ a 1-morphism $\omega_{\mathcal C}\colon F_0(\mathcal C)\to G_0(\mathcal C)$
  \item For all pairs of objects $\mathcal C,\mathcal D$ of $\mathcal G$ a natural transformation of weak 2-functors $\omega_{\mathcal C,\mathcal D}\colon  (\omega_\mathcal D\Box -)F_{\mathcal C,\mathcal D}  \to
    (-\Box\omega_\mathcal C)G_{\mathcal C,\mathcal D}$
  \end{itemize}
  such that for all triples $\mathcal B$, $\mathcal C$, $\mathcal D$ the following diagrams commute
  \begin{align}
    \xymatrix{
      {\begin{array}{r}
          (\omega_{\mathcal D}\Box -)\Box (F_{\mathcal B,\mathcal C}\times F_{\mathcal C,\mathcal D})\\
          =\Box (1\times (\omega_\mathcal D\Box -))(F_{\mathcal B,\mathcal C}\times F_{\mathcal C,\mathcal D})\end{array}}\quad \ar[r]^{\qquad\qquad(\omega_{\mathcal D}\Box -)\kappa^F_{\mathcal B, \mathcal C,\mathcal D}} 
      \ar[d]^{\Box(1\times \omega_{\mathcal C,\mathcal D})}
      & \quad(\omega_\mathcal D\Box -)F_{\mathcal B,\mathcal D}\Box  \ar[d]_{\omega_{\mathcal B,\mathcal D}\Box}\\
      {\begin{array}{r}\Box (1\times (-\Box \omega_\mathcal C))(F_{\mathcal B,\mathcal C}\times G_{\mathcal C,\mathcal D})\\ = \Box((\omega_\mathcal C\Box -)\times 1)(F_{\mathcal B,\mathcal C}\times G_{\mathcal C,\mathcal D})\end{array}} \ar[d]^{\Box (\omega_{\mathcal B,\mathcal C}\times 1)}  & (-\Box\omega_\mathcal B)G_{\mathcal B,\mathcal D}\Box\\
      {\begin{array}{r}\Box ((-\Box \omega_{\mathcal B})\times 1)(G_{\mathcal B,\mathcal C}\times G_{\mathcal C,\mathcal D})\\= (-\Box \omega_\mathcal B)\Box (G_{\mathcal B,\mathcal C}\times G_{\mathcal C,\mathcal D}),\end{array}}\quad \ar[ru]_{\qquad\quad(-\Box \omega_\mathcal B)\kappa^G_{\mathcal B,\mathcal C,\mathcal D}} } \nonumber
  \end{align}
  \begin{align}
    \xymatrix{
      \omega_{\mathcal{C}} \Box F_{\mathcal{C},\mathcal{C}}(1_{\mathcal{C}}) \ar[r]^{\omega_{1_{\mathcal{C}}}} \ar[d]_{1\Box (\iota_{F})_{\mathcal{C}}} &  G_{\mathcal{C},\mathcal{C}}(1_{\mathcal{C}}) \Box \omega_{\mathcal{C}} \ar[dl]^{(\iota_{G})_{\mathcal{C}} \Box 1} \\
      \omega_{\mathcal{C}} 
    }\qquad 
    \xymatrix{
      &  \omega_{\mathcal{C}} \ar[dl]_{1\Box (\iota_{F}^{-1})_{\mathcal{C}}} \ar[d]^{(\iota_{G}^{-1})_{\mathcal{C}} \Box 1}\\
      \omega_{\mathcal{C}} \Box F_{\mathcal{C},\mathcal{C}}(1_{\mathcal{C}}) \ar[r]^{\omega_{1_{\mathcal{C}}}} & G_{\mathcal{C},\mathcal{C}}(1_{\mathcal{C}}) \Box \omega_{\mathcal{C}}. 
    }
    \nonumber
  \end{align}    
  Here,  $-\Box\omega_\mathcal C\colon\tilde{\mathcal G}( G_0(\mathcal C), \tilde{\mathcal D})\to \tilde{\mathcal G}(F_0(\mathcal C),\tilde{\mathcal D})$ and $\omega_\mathcal D\Box-\colon\tilde{\mathcal G}(\tilde{\mathcal C}, F_0(\mathcal D))\to \tilde{\mathcal G}(\tilde{\mathcal C}, G_0(\mathcal D))$ denote the weak 2-functors defined by pre- and postcomposition with $\omega_C$ with respect to the Gray product.
  The natural transformations $\omega_{\mac C,\mac D}$ determine for all 1-morphisms $H\colon\mac C\to\mac D$ in $\mac G$ a 2-morphism $\omega_H\colon\omega_{\mac D}\Box F_{\mac C,\mac D}(H) \Rightarrow  G_{\mac C,\mac D}(H)\Box \omega_{\mac C}$.

  A natural transformation $\omega$ is called a {\bf natural isomorphism}, if the 1-morphisms $\omega_{\mathcal{C}}$ are invertible and all natural transformations $\omega_{\mathcal{C},\mathcal{D}}$ are invertible pseudo-natural transformations.
\end{definition}

\begin{definition}

  Let $F=(F_0,F_{\mac C,\mac D},  \kappa^F_{\mac C,\mac D,\mac E},\iota^F_{\mac C})\colon\mac G\to\tilde{\mac G}$ and $G=(G_0,G_{\mac C,\mac D}, \kappa^G_{\mac C,\mac D,\mac E},\iota^G_{\mac C})\colon\mathcal G\to\tilde{\mathcal G}$ be functors of 2-strict  tricategories and  $\omega=(\omega_{\mac C},\omega_{\mac C,\mac D}), \eta=(\eta_{\mac C}, \eta_{\mac C,\mac D})\colon F\Rightarrow G$  natural transformations. A {\bf modification}
  $\Psi\colon \omega\Rrightarrow\eta$ consists of the following data:
  \begin{itemize}
  \item For every object $\mac C$ of $\mac G$ a 2-morphism $\Psi_{\mac C}\colon\omega_{\mac C}\Rightarrow \eta_{\mac C}$ in $\tilde{\mac G}$.
  \item For every pair of objects $\mac C,\mac D$ of $\mac G$, an invertible 3-morphism 
    $\Psi_{\mac C,\mac D}\colon ( 1 \Box \Psi_{\mac C})\circ \omega_{\mac C,\mac D}\Rrightarrow \eta_{\mac C,\mac D}\circ( \Psi_{\mac D} \Box 1)$.
  \end{itemize}
  These determine for all 1-morphisms $H\colon\mac C\to\mac D$ a 3-morphism 
  $$\Psi_H\colon (1_{F_{\mac C,\mac D}(H)} \Box \Psi_{\mac C})\circ \omega_{H}\Rrightarrow \eta_H\circ (\Psi_{\mac D} \Box 1_{G_{\mac C,\mac D}(H)}).$$
  A modification $\Psi$ is called invertible if all 2-morphisms $\Psi_{\mathcal{C}}$ are invertible.
\end{definition}

\begin{definition}
  \label{Equivalence-Gray-Cat}
  Two 2-strict  tricategories $\mathcal{G}$, $\widetilde{\mathcal{G}}$  are called {\bf equivalent}, if there exist \Glax functors of  2-strict tricategories $F\colon \mathcal{G}\rightarrow \widetilde{\mathcal{G}}$ and 
  $G\colon \widetilde{\mathcal{G}} \rightarrow \mathcal{G}$ together with invertible pseudo-natural transformations $\eta\colon F G \rightarrow 1_{\widetilde{\mathcal{G}}}$ and $\varphi\colon GF \rightarrow 1_{\mathcal{G}}$. 
\end{definition}


\end{document}